\documentclass[11pt]{article}
\usepackage[applemac]{inputenc}
\usepackage{latexsym,amsmath,amscd,amssymb,framed,graphics,color,mathrsfs,stmaryrd}
\usepackage{enumerate}
\usepackage{geometry}
\usepackage[colorlinks]{hyperref}
\usepackage[T1]{fontenc}
\usepackage{graphicx}
\definecolor{sand}{rgb}{0.76, 0.7, 0.5}
\definecolor{taupegray}{rgb}{0.55, 0.52, 0.54}
\usepackage[colorlinks]{hyperref}
\hypersetup{%
pdfstartview={FitH},
urlcolor=sand,
citecolor=blue,
linkcolor=taupegray,
}

\usepackage{array}
\usepackage{float}
\usepackage{multicol}
\usepackage{aeguill}
\usepackage[dvipsnames]{xcolor}
\usepackage{fancyhdr}
\usepackage[arrow, matrix, curve]{xy}
\usepackage{url}
\usepackage{colordvi}
\usepackage{framed}
\usepackage{color}
\usepackage{algorithmic}
\usepackage{algorithm}

\usepackage{manfnt}
\usepackage[bbgreekl]{mathbbol}
\usepackage{tikz}

\newcommand{\mathsym}[1]{{}}

\newtheorem{theorem}{Theorem}[section]
\newtheorem{definition}[theorem]{Definition}

\newtheorem{remark}[theorem]{Remark}

\begin{document}
\title{Multisymplectic variational integrators \\ for barotropic and incompressible fluid models with constraints}

\author{Fran\c{c}ois Demoures$^{1}$ and Fran\c{c}ois Gay-Balmaz$^{2}$}
\addtocounter{footnote}{1} \footnotetext{EPFL, Doc \& Postdoc Alumni. Av Druey 1, Lausanne 1018, Switzerland \\\texttt{francois.demoures@alumni.epfl.ch}}
\addtocounter{footnote}{2}
\footnotetext{CNRS \& \'Ecole 
Normale Sup\'erieure, Laboratoire de M\'et\'eorologie Dynamique, Paris, France.
\texttt{francois.gay-balmaz@lmd.ens.fr}}

\maketitle

\begin{abstract}
We present a structure preserving discretization of the fundamental spacetime geometric structures of fluid mechanics in the Lagrangian description in 2D and 3D. Based on this, multisymplectic variational integrators are developed for barotropic and incompressible fluid models, which satisfy a discrete version of Noether theorem. We show how the geometric integrator can handle regular fluid motion in vacuum with free boundaries and constraints such as the impact against an obstacle of a fluid flowing on a surface. Our approach is applicable to a wide range of models including the Boussinesq and shallow water models, by appropriate choice of the Lagrangian.
\end{abstract} 


 \section{Introduction}

This paper presents a multisymplectic variational integrator for barotropic fluids and incompressible fluids with free boundaries in the Lagrangian description. The integrator is derived from a spacetime discretization of the Hamilton principle of fluid dynamics and is based on a discrete version of the multisymplectic geometric formulation of continuum mechanics. As a consequence of its variational nature, the resulting scheme preserves exactly the momenta associated to symmetries, it is symplectic in time, and energy is well conserved. In addition to its conservative poperties, the variational scheme can be naturally extended to handle constraints, such as the impact against an obstacle of fluid flowing on a surface, by augmenting the discrete Lagrangian with penalty terms.

\medskip

Multisymplectic geometry is the natural geometric setting for classical field theories and is the appropriate spacetime extension of the symplectic formulation of classical mechanics. Important properties of Lagrangian and Hamiltonian systems in classical mechanics, such as the symplecticity of the flow and the preservation of the momentum maps associated to symmetries, have corresponding statements for field theories that are intrinsically formulated via multisymplectic geometry. These are the \textit{multisymplectic form formula} and the \textit{covariant Noether theorem} for the solution of Euler-Lagrange field equations. Of particular importance in these formulations are the \textit{Cartan forms} associated to the Lagrangian density of the theory.

\medskip 

Multisymplectic variational integrators were developed in \cite{MaPaSh1998} via a spacetime discretization of the Hamilton principle of field theories, which results in numerical schemes that satisfy a discrete version of the multisymplectic form formula and a discrete covariant Noether theorem. The discrete framework also allows the definitions of the concepts of \textit{discrete Cartan forms} and \textit{discrete covariant momentum maps}. Other approaches to multisymplectic integrators
have been also developed, in, for example, \cite{BrRe2001}. 
We refer to \cite{LeMaOrWe2003,DeGBRa2014,DeGBKoRa2014,DeGBRa2016,DGBDRA2017} for the development of multisymplectic variational integrators for several mechanical systems of interest in engineering. Examples include the simulation of the dynamics of rotor blades via asynchronous variational integrators where it is necessary to compute accurate solutions for long periods of time, the dynamics of geometrically exact (Cosserat) beams, or the simulation of elastodynamic frictionless impact problems.

\medskip

In this paper, we develop this method towards its application to compressible and incompressible fluid dynamics by using, at the continuous level, the multisymplectic variational formulation of continuum mechanics as described in \cite{MaPeShWe2001,FeMaWe2003}. The main ingredients in our discrete approach are the concepts of discrete deformation gradient and discrete Jacobian, defined both in the 2D and 3D cases. They enter in a fundamental way in the definition of the spacetime discretized Lagrangian and they allow to exactly impose discrete incompressibility via an augmented Lagrangian approach.
Besides it conservative properties, thanks to its variational nature, our scheme can be naturally extended to include constraints. This is illustrated with fluid flowing or impacting on a surface.

\medskip 

The variational discretization in this paper is carried out in the Lagrangian frame and for fluid dynamics interpreted as a special class of field theory on spacetime. Geometric variational discretizations for fluids have also been developed in the Eulerian description and for fluid dynamics interpreted as an infinite dimensional dynamical system on diffeomorphism groups, as opposed to the spacetime covariant description carried out here. This variational approach is based on a discretization of groups of diffeomorphisms, see \cite{Pa2009,PMTKMD2011,BaGB2017,GaGB2020} for both incompressible and compressible models.

\medskip 

This paper is a first step towards the development of dynamic mesh update from a structure preserving point of view, inspired by arbitrary Lagrangian-Eulerian methods. Several approaches have been proposed in the literature, such as \cite{HuLiZi1981,DoGiHa1982,Do1983,SoFoDhOu1991,MaHu1997,FaRaWaBe2010}.

\medskip

The organization of the paper is as follows. Section \S\ref{baro_field_th} first briefly reviews the variational formulation of barotropic and incompressible fluid models in the Lagrangian description in a classical way. We mention in particular the case of isentropic perfect gas, the shallow water and Boussinesq equations, and the ideal fluid. This variational setting is then recasted in the multisymplectic variational formalism, which is fundamental for the discretization carried out later. The multisymplectic form formula and the covariant Noether theorems are recalled.
The two dimensional discrete fluid models are formulated in Section \S\ref{2D_wave_prop_fluid}.  In \S\ref{2D_mult_discret}, the discrete configuration bundle and jet bundle are recalled, and the discrete deformation gradient, the discrete Jacobian as well as the discrete Lagrangian for barotropic models are defined. The discrete Euler-Lagrange equations are obtained from the discrete version of the Hamilton principle. The algorithmically conserved quantities (discrete multisymplectic form formula and Noether theorem) are written. In \S\ref{2D_penalty_method} discrete incompressibility is treated via a Lagrange multiplier constraint and via a penalty term. Numerical results are presented in \S\ref{2D_examples} to demonstrate the basic properties of the method and to validate it, for both compressible and incompressible fluids, with free boundary or flowing on a surface and impacting against an obstacle.
Section \S\ref{3D_wave_prop_fluid} develops the three dimensional discrete multisymplectic formulation for barotropic and incompressible ideal fluids with the same class of examples than in the two dimensional situation. 
The paper concludes with the Appendix \ref{appendix} where several expressions needed to implement the integrators are given.

\section{Barotropic and incompressible fluids}\label{baro_field_th}

In this section we briefly review the variational formulation of barotropic and incompressible fluid models in the Lagrangian (or material) description in Cartesian coordinates. This formulation is then recasted in a multisymplectic variational setting, which allows to formulate intrinsically the Hamilton principle, the multisymplectic property of the solutions, and the covariant Noether theorem with the help of Cartan forms. This gives the geometric framework to be discretized in a structure preserving way later.

\medskip

Assume that the reference configuration of the fluid is a compact domain $ \mathcal{B}  \subset  \mathbb{R} ^n $ with piecewise smooth boundary, and the fluid moves in the ambient space $\mathcal{M} = \mathbb{R} ^n$. We denote by $ \varphi : \mathbb{R} \times \mathcal{B} \rightarrow \mathcal{M} $ the fluid configuration map, which indicates the location $ m=\varphi (t,X)$ at time $t$ of the fluid particle with label $X \in \mathcal{B}$. The deformation gradient is denoted $\mathbf{F}(t,X)$, given in coordinates by $ \mathbf{F} ^a{}_i= \varphi ^a{}_{,i}$, with $X^i, i=1,...,n$ the Cartesian coordinates on $ \mathcal{B} $ and $m^a, a=1,...,n$ the Cartesian coordinates on $\mathcal{M} $. We assume that the fluid configuration is regular enough so that all the computations below are valid.

\subsection{Barotropic fluids}

\subsubsection{Definition} 
A fluid is \textit{barotropic} if it is compressible and the surfaces of constant pressure $p$ and constant density $\rho$ coincide, i.e., we have a relation
\begin{equation} \label{P(rho)}
p = p(\rho).
\end{equation}
The internal energy $W$ of barotropic fluids in the material description depends on the deformation gradient $\mathbf{F}$ only through the Jacobian $J$ of $ \varphi $, given in Cartesian coordinates by
\[
J(t,m)= \mathrm{det}(\mathbf{F}(t,X)),
\]
hence in the material description we have $W=W( \rho  _0, J)$, with $\rho_0(X)$ the mass density of the fluid  in the reference configuration. 
The pressure in the material description is
\begin{equation}\label{baro_pressure} 
P_W( \rho _0, J)= - \rho  _0 \frac{\partial W}{\partial J} ( \rho  _0, J).
\end{equation} 

The \textit{continuity equation for mass} can be written as
\begin{equation}\label{cont_mass_eq}
\rho_0(X)= \rho(t, \varphi(t,X)) J(t,X),
\end{equation}
with $ \rho  (t,m)$ the Eulerian mass density. The internal energy $w( \rho  )$ in the Eulerian description satisfies the relation
\[
W( \rho  _0, J)= w \left( \frac{\rho  _0}{J}\right)
\]
and one notes that $P_W= p_w \circ \varphi $, with $p_w= \rho  ^2 \frac{\partial w}{\partial \rho  }$ the Eulerian pressure $p$ in \eqref{P(rho)}.

\subsubsection{Hamilton's principle for barotropic fluids}

The Lagrangian of the barotropic fluid evaluated on a fluid configuration map $ \varphi (t,X)$ has the standard form
\begin{equation}\label{Lagrangian_barotropic} 
L( \varphi , \dot \varphi , \nabla \varphi )= \frac{1}{2} \rho  _0 |\dot \varphi | ^2 - \rho  _0 W( \rho  _0,J) - \rho  _0 \Pi ( \varphi ),
\end{equation} 
with $ \Pi $ a potential energy, such as the gravitational potential $ \Pi ( \varphi )= \mathbf{g} \cdot \varphi $.

\medskip 

Hamilton's principle
\[
\delta \int_0^T\!\int_ \mathcal{B} L( \varphi , \dot \varphi , \nabla \varphi ) {\rm d}t\,{\rm d}X=0
\]
for variations of $\varphi $ vanishing at $t=0,T$ yields the Euler-Lagrange equations
\[
\frac{\partial }{\partial t} \frac{\partial L}{\partial \dot \varphi } + \frac{\partial }{\partial x ^i } \frac{\partial L}{\partial \varphi _{, i}}= \frac{\partial L}{\partial \varphi },
\]
together with the natural boundary conditions
\[
\frac{\partial L}{\partial \varphi ^a _{,i}} n_i \delta \varphi ^a =0 \quad\text{on}\quad \partial \mathcal{B} ,
\]
for allowed variations $ \delta \varphi $. Here $n$ denotes the outward pointing unit normal vector field to $ \partial \mathcal{B} $.

\medskip

From the Lagrangian of the barotropic fluid \eqref{Lagrangian_barotropic} and the material pressure $P_W$ defined in \eqref{baro_pressure} we get  the barotropic fluid equations in the Lagrangian description as
\begin{equation}\label{CEL_barotropic} 
\rho  _0 \ddot \varphi + \frac{\partial }{\partial x^i }  \left( P_W J \mathbf{F} ^{-1} \right) ^i = - \rho  _0 \frac{\partial \Pi }{\partial \varphi } 
\end{equation} 
together with the natural boundary conditions
\begin{equation}\label{zero_pressure_BC} 
P_W J \, n _i (\mathbf{F} ^{-1} ) ^i _a \delta \varphi ^a =0 \quad\text{on}\quad \partial \mathcal{B} ,
\end{equation}
for all allowed variations $ \delta \varphi $. For instance for a free boundary problem, the variations $ \delta \varphi $ are arbitrary on $ \partial \mathcal{B} $, hence the boundary condition \eqref{zero_pressure_BC} yields the zero pressure condition
\begin{equation}\label{BC_pW} 
P_W|_{ \partial \mathcal{B} }=0.
\end{equation} 
Boundary conditions with surface tension can be deduced from the Hamilton principle by adding an area term in the Lagrangian, see \cite{GBMaRa2012}.

\medskip 

Using the relations $ \dot \varphi = u \circ \varphi $, $P_W=p_w \circ \varphi $, and $ \rho  _0 = (\rho  \circ \varphi )J$, between Lagrangian and Eulerian quantities, one deduces from \eqref{CEL_barotropic} the familiar Eulerian form of barotropic fluids as
\[
\rho  (\partial _t u + u \cdot \nabla u)= - \nabla p_w - \rho  \nabla \Pi , \qquad \partial _t \rho  + \operatorname{div}( \rho  u)=0. 
\]

\subsubsection{Example: isentropic perfect gas and rotating shallow water}

Let us consider the following general barotropic expression for the internal energy and pressure
\begin{equation}\label{Barotropic_general} 
w( \rho  )= \frac{A}{ \gamma -1}\rho   ^{\gamma -1} + B \rho  ^{-1}  , \qquad p_w( \rho  )= A \rho  ^ \gamma - B,
\end{equation}
for constants $A$, $B$, and adiabatic coefficient $ \gamma $, see \cite{CoFr1948}. The material internal energy to be used in the Lagrangian \eqref{Lagrangian_barotropic} is 
\begin{equation}\label{energy_water} 
W( \rho  _0, J)= \frac{A}{ \gamma -1} \left(\frac{J}{ \rho  _0}  \right)    ^{1- \gamma } + B  \left(\frac{J}{ \rho  _0}  \right).
\end{equation} 
For an \textit{isentropic perfect gas} we have $B=0$.

\medskip 

In our tests, we shall use the expression \eqref{Barotropic_general} for the treatment of an isentropic perfect gas, where the value of the constant $B\neq 0$ does not affect the dynamics, while it allows to naturally impose from \eqref{BC_pW} the boundary condition
\[
P|_{ \partial \mathcal{B} }=B,
\]
with $P= A \left( \frac{\rho  _0}{J} \right) ^ \gamma $ the pressure of the isentropic perfect gas. This is crucial for the discretization, since it allows to find the appropriate discretization of the boundary condition directly from the boundary terms of the discrete variational principle.

\medskip 

The \textit{rotating shallow water model} can also be recasted in the formulation above, in which case the variable $ \rho  _0$ is interpreted as the water depth in the reference configuration. The Lagrangian is
\[
L( \varphi ,\dot \varphi , \nabla \varphi ) = \frac{1}{2} \rho  _0| \dot \varphi | ^2 + \rho  _0 \dot \varphi \cdot R( \varphi ) - \rho  _0W( \rho  _0,J) , 
\]
where $R$ is the vector potential of the angular velocity of the Earth and $W$ is chosen as
\[
W( \rho  _0, J)=g \frac{1}{2}\frac{ \rho  _0}{J}.
\]

\subsection{Incompressible fluid models}

\subsubsection{Hamilton principle with incompressibility constraint}\label{const_Incompressible}

Incompressible models are obtained by inserting the constraint $J=1$ in the Hamilton principle as 
\begin{equation}\label{HP_constraint} 
\delta \int_0^T\!\int_ \mathcal{B} \big( L( \varphi , \dot \varphi , \nabla \varphi )  + \lambda (J -1) \big) {\rm d}t\,{\rm d}X=0,
\end{equation} 
where $ \lambda (t,X)$ is the Lagrange multiplier. With the Lagrangian \eqref{Lagrangian_barotropic}, this results in the system
\begin{equation}\label{CEL_barotropic_incomp} 
\rho  _0 \ddot \varphi + \frac{\partial }{\partial x^i }  \left( (P_W + \lambda ) J \mathbf{F} ^{-1} \right) ^i = - \rho  _0 \frac{\partial \Pi }{\partial \varphi } , \qquad J=1.
\end{equation} 
With the relations $ \dot \varphi = u \circ \varphi $, $P_W=p_w \circ \varphi $, and $ \rho  _0 = (\rho  \circ \varphi )J$, we get from \eqref{CEL_barotropic_incomp} the familiar Eulerian formulation
\[
\rho  ( \partial _t u  + u \cdot \nabla u)= - \nabla (p_w+ p_ \lambda ) - \rho  \nabla \Pi , \qquad \operatorname{div} u=0, \qquad \partial _t \rho  + u \cdot \nabla \rho  =0.
\] 
In this case $p_w+p_ \lambda $ is determined from the incompressibility constraint via a Poisson equation.

\subsubsection{Example: Boussinesq model, nonhomogeneous Euler equations, and ideal fluid}

The \textit{Boussinesq model} is obtained from the Hamilton principle with incompressibility constraint \eqref{HP_constraint} by interpreting $ \rho  _0$ as the buoyancy in the reference configuration and taking the Lagrangian
\begin{equation}\label{Lagrangian_Boussinesq} 
L( \varphi , \dot \varphi , \nabla \varphi )= \frac{1}{2} |\dot \varphi | ^2 - \rho  _0 \varphi \cdot \mathbf{g} 
\end{equation}
with gravitational acceleration vector $\mathbf{g}$. For the \textit{nonhomogeneous Euler fluid}, the Lagrangian is the kinetic energy
\begin{equation}\label{Lagrangian_nonhomogeneous} 
L( \varphi , \dot \varphi , \nabla \varphi )= \frac{1}{2} \rho  _0  |\dot \varphi | ^2 ,
\end{equation} 
for some non-constant density $ \rho  _0( X)$. For the \textit{ideal fluid}, one takes
\begin{equation}\label{Lagrangian_ideal} 
L( \varphi , \dot \varphi , \nabla \varphi )= \frac{1}{2} |\dot \varphi | ^2 
\end{equation} 
in \eqref{HP_constraint}, which gives
\begin{equation}\label{CEL_ideal} 
\ddot \varphi + \frac{\partial }{\partial x^i } \left(   \lambda  J \mathbf{F} ^{-1} \right) ^i = 0 , \qquad J=1,
\end{equation} 
and hence $ \partial _t u + u \cdot \nabla u= - \nabla p_\lambda $, $ \operatorname{div}u=0$ is obtained in the Eulerian formulation.

\subsection{Multisymplectic variational continuum mechanics}\label{MVCM}

In this paragraph, we briefly review the geometric variational framework of classical field theory, as it applies to continuum mechanics, following \cite{MaPeShWe2001}. This setting will be discretized in a structure preserving way which allows the identification of the notion of discrete multisymplecticity, discrete momentum map, and discrete Noether theorems.

\subsubsection{Configuration bundle, jet bundle, and Lagrangian density}

The geometric formulation of classical field theories starts with the identification of the configuration bundle of the theory, denoted $\pi _{ \mathcal{Y} , \mathcal{X} }: \mathcal{Y}  \rightarrow \mathcal{X} $, such that the fields $ \varphi $ of the theory are sections of this fiber bundle, i.e., they are smooth maps $ \varphi :\mathcal{X}  \rightarrow \mathcal{Y} $ such that $ \pi _{\mathcal{Y} , \mathcal{X} } \circ \varphi =\operatorname{id}_ \mathcal{X} $, where $ \operatorname{id}_ \mathcal{X} $ denotes the identity map on $ \mathcal{X} $. 
We assume $ \operatorname{dim} \mathcal{X} =n+1$ and denote by $x ^\mu $, $ \mu =0,1,2,...,n$, the coordinates on $ \mathcal{X} $. The fiber coordinates on $\mathcal{Y} $ are $y^a$, $a=1,...,N$, hence coordinates on the manifold $ \mathcal{Y} $ are $(x^\mu, y^a)$, $\mu=0,...,n$, $a=1,...,N$.
While the configuration bundle for continuum mechanics is a trivial bundle, it is advantageous to use the general setting of fiber bundles since it allows to efficently particularise to continuum mechanics the intrinsic geometric formulation and structures of field theories.

\medskip 

The \textit{first jet bundle} of the configuration bundle $\pi _{\mathcal{X} ,\mathcal{Y} }: \mathcal{Y}  \rightarrow \mathcal{X} $ is the field theoretic analogue of the tangent bundle of classical mechanics, i.e., its fiber at $y$  contains the first derivatives $ \varphi^a{} _{,\mu}(x)$ of a field $ \varphi $ at $x$  with $ \varphi (x)=y$. It is defined as the fiber bundle $\pi _{\mathcal{Y} ,J ^1 \mathcal{Y} }: J ^1 \mathcal{Y}  \rightarrow \mathcal{Y} $ over $\mathcal{Y} $, whose fiber at $y \in \mathcal{Y} $ consists of linear maps $ \gamma :T_x\mathcal{X}  \rightarrow T_y\mathcal{Y} $ satisfying $T \pi _{\mathcal{Y} ,\mathcal{X} } \circ \gamma = \operatorname{id}_{T_x\mathcal{X} }$, where $x= \pi _{\mathcal{X} ,\mathcal{Y} }(y)$.   The induced coordinates on the fiber of $J ^1 \mathcal{Y} \rightarrow \mathcal{Y} $ are denoted $v^a{}_\mu$. We note that $ J ^1 \mathcal{Y} $ can also be regarded as the total space of a bundle over $\mathcal{X} $, namely $ \pi _{\mathcal{X} , J ^1 \mathcal{Y} }:= \pi _{\mathcal{X} ,\mathcal{Y} } \circ \pi _{\mathcal{Y} ,J ^1 \mathcal{Y} }: J ^1 \mathcal{Y}  \rightarrow \mathcal{X} $. Natural coordinates on the manifold $ J^1\mathcal{Y} $ are hence $(x^\mu, y^a, v^a{}_\mu)$, $\mu=0,...,n$, $a=1,...,N$.

\medskip

The derivative of a field $ \varphi $ can be regarded as a section of $\pi _{\mathcal{X} ,J ^1 \mathcal{Y} }:J ^1 \mathcal{Y}  \rightarrow \mathcal{X} $, by writing $x \in \mathcal{X}  \mapsto j ^1 \varphi (x):=T_x \varphi \in J ^1 _{ \varphi (x)}\mathcal{Y} $, with $T_x \varphi :T_x\mathcal{X}  \rightarrow T_{ \varphi (x)}\mathcal{Y} $ the tangent map (or first derivative) of $ \varphi $. The section $ j^1 \varphi $ is called the \textit{first jet extension} of $\varphi $ and is the intrinsic object corresponding to the value of a field and of its first derivatives, at the points in $\mathcal{X} $.  In the natural coordinates $(x^\mu, y^a, v^a{}_\mu)$ of $J^1\mathcal{Y} $, the first jet extension reads $ j^1 \varphi :x ^\mu \mapsto ( x ^\mu , \varphi ^a (x), \varphi ^a {}_{, \mu }(x))$.

\medskip

A \textit{Lagrangian density} is a smooth bundle map $ \mathcal{L} : J ^1 \mathcal{Y}  \rightarrow \Lambda ^{n+1}\mathcal{X} $ over $\mathcal{X} $, where $ \Lambda ^{n+1}\mathcal{X}  \rightarrow \mathcal{X} $ is the vector bundle of $(n+1)$-form on $\mathcal{X} $. In coordinates we write $ \mathcal{L} ( j ^1 \varphi (x) )= L( x ^\mu , \varphi^a , \varphi ^a {}_{, \mu } )d^{n+1}x$. The associated \textit{action functional} is
\begin{equation}\label{action_funct}
S( \varphi ):=\int_\mathcal{X}  \mathcal{L} ( j ^1 \varphi (x)).
\end{equation}

\subsubsection{The case of continuum mechanics}

For continuum mechanics, the configuration bundle is the trivial fiber bundle
\[
\mathcal{Y} = \mathcal{M}  \times \mathcal{X}  \rightarrow \mathcal{X}, \quad\text{with}\quad \mathcal{X} = \mathbb{R} \times \mathcal{B},
\]
where $ \mathcal{B} $ is the reference configuration of the continuum and $\mathcal{M} $ is the ambient space, see the beginning of \S\ref{baro_field_th}. We have the equalities $x=(t,X)$ and $y=(x,m)=(t,X,m)$ between the variables of the general theory and those of continuum mechanics.

\medskip 

A section of this bundle is a map $ \varphi : \mathcal{X} \rightarrow \mathcal{X} \times \mathcal{M}$, whose first component is $ \operatorname{id}_ \mathcal{X} $. It is canonically identified with a map $\varphi : \mathcal{X} =\mathbb{R} \times \mathcal{B} \rightarrow  \mathcal{M}$ referred to as the fluid configuration map above.

\medskip 

The first jet bundle is canonically identified with the vector bundle $L(T\mathcal{X}  , T\mathcal{M} ) \rightarrow \mathcal{Y} = \mathcal{X}  \times \mathcal{M} $, whose fiber at $y=(x,m)$ is the vector space $L(T_x\mathcal{X} , T_m\mathcal{M} )$ of linear maps from $T_x \mathcal{X} $ to $T_m \mathcal{M}$. The first jet extension is $j^ 1\varphi (t,X)= (\varphi (t,X), \dot \varphi (t,X), \nabla \varphi (t,X))$ and the Lagrangian density reads
\[
\mathcal{L} (\varphi , \dot \varphi , \nabla \varphi )= L (\varphi , \dot \varphi , \nabla \varphi) dt \wedge d^nX
\]
with $L$ the Lagrangian of barotropic fluids given in \eqref{Lagrangian_barotropic}.

\subsubsection{Multisymplectic form and Cartan forms}\label{MF_CF}

Without entering into the details, we recall that the dual jet bundle $J^1\mathcal{Y} ^\star \rightarrow \mathcal{Y} $, defined as the bundle of affine maps $J^1\mathcal{Y}  \rightarrow \Lambda ^{n+1}\mathcal{X} $, is endowed with a \textit{canonical $(n+1)$ form} $ \Theta _{\rm can}$ and a \textit{canonical multisymplectic $(n+2)$-form} $ \Omega  _{\rm can}= - \mathbf{d}\Theta _{\rm can}$. These are the field theoretic analogue to the canonical one-form and canonical symplectic form on the phase space (cotangent bundle of the configuration manifold) in classical mechanics. By pulling back these canonical forms with the Legendre transform $\mathbb{F} \mathcal{L} :J^1Y \rightarrow J^1Y^\star$ of a given Lagrangian density $ \mathcal{L} :J^1\mathcal{Y}  \rightarrow \Lambda ^{n+1}\mathcal{X} $, one gets the \textit{Cartan forms} $ \Theta _ \mathcal{L} $ and $ \Omega _ \mathcal{L} $ on $J^1\mathcal{Y} $, see \cite{GiMmsy}. These forms appear naturally in the Hamilton principle, in the multisymplectic form formula, and in the Noether theorem, as will shall explain below. All these three notions have discrete analogues, that we shall deeply use in \S\ref{2D_wave_prop_fluid} and \S\ref{3D_wave_prop_fluid}.

\medskip 

The Cartan forms arise in the Hamilton principle as follows. Using the relation $ \mathcal{L}(j^1 \varphi )= (j^1 \varphi ) ^* \Theta _ \mathcal{L} $, \cite{GiMmsy}, the derivative of the action functional \eqref{action_funct} takes the intrinsic form
\begin{equation}\label{HP_intrinsic} 
\begin{aligned}
\mathbf{d} S( \varphi ) \cdot V( \varphi ) &= \left. \frac{d}{d\varepsilon}\right|_{\varepsilon=0} \int_ \mathcal{X}  \mathcal{L} ( j^1( \phi _ \varepsilon \circ \varphi )) \\
&= - \int_\mathcal{X}  (j^1 \varphi ) ^* \mathbf{i} _{j^1V}  \Omega_ \mathcal{L} + \int_{ \partial \mathcal{X} } (j^1 \varphi ) ^* \mathbf{i} _{j^1V} \Theta  _ \mathcal{L},
\end{aligned} 
\end{equation} 
where $ \phi _ \varepsilon $ is the flow of a vertical vector field $V$ on $\mathcal{Y} $, i.e., $T \pi _{\mathcal{X} ,\mathcal{Y}  } \circ V=0$, and $j^1V$ denotes the first jet extension of $V$ to $J^1\mathcal{Y} $ defined as $j^1V= \left. \frac{d}{d\varepsilon}\right|_{\varepsilon=0} j^1 \phi ^ \varepsilon $.

\subsubsection{Multisymplectic form formula and Noether theorem}\label{MFF_NT}

The multisymplectic form formula is a property of the solution of the Euler-Lagrange field equations that extends the symplectic property of the solution of the Euler-Lagrange equations of classical mechanics. It is obtained from the identity \eqref{HP_intrinsic}, by evaluating the action functional at a solution of the Euler-Lagrange equations and taking its derivative along variations of solutions, see \cite{MaPaSh1998}.
Let $ \varphi $ be a solution of the Euler-Lagrange field equations and $V$, $W$ solutions of the first variation of the Euler-Lagrange equations at $ \varphi $. Then $V$, $W$, $\varphi $ satisfy the \textit{multisymplectic form formula}:
\begin{equation}\label{MFF}
\int_ { \partial U} (j^1\varphi ) ^*  \mathbf{i} _{j^1V} \mathbf{i} _{j^1W} \Omega _ \mathcal{L} =0,
\end{equation} 
for all open subset $U \subset \mathcal{X} $ with with piecewise smooth boundary.

\medskip

We now recall the general statement of the \textit{covariant Noether theorem}. Let a Lie group $G$ act on $\mathcal{Y} $ an assume that the action covers a diffeomorphism of $\mathcal{X} $. Assume that the Lagrangian density $ \mathcal{L} $ is $G$-equivariant with respect to this action, see later in \S\ref{Sym_barotropic} for a concrete example. Then, considering only variations along the Lie group action, and restricting the action functional to an arbitrary open subset $U \subset \mathcal{X} $ with piecewise smooth boundary, formula \eqref{HP_intrinsic} shows that a solution of the Euler-Lagrange field equations satisfy the \textit{covariant Noether theorem}
\begin{equation}\label{NT} 
\int_{ \partial U} ( j^1 \varphi ) ^* J^ \mathcal{L} ( \xi )=0, \qquad \text{for all $ \xi \in \mathfrak{g} $},
\end{equation}
where $ J^ \mathcal{L} ( \xi )= \mathbf{i} _{j^1 \xi _\mathcal{Y} } \Theta _ \mathcal{L}:J^1\mathcal{Y}  \rightarrow \mathfrak{g} ^* \otimes \Lambda ^nJ^1\mathcal{Y} $ is the \textit{covariant momentum map} associated to $ \mathcal{L} $ and $ \xi _\mathcal{Y} $ is the infinitesimal generator of the Lie group action associated to the Lie algebra element $ \xi \in \mathfrak{g} $.

\section{2D Discrete barotropic and incompressible fluid models} \label{2D_wave_prop_fluid}

In this section we propose a multisymplectic variational discretization of fluid mechanics, by focusing on compressible barotropic models and incompressible models. We consider free boundary fluids, as well as fluid impacting on a surface.
A main step in our construction is the definition of discrete deformation gradient and discrete Jacobian.

\subsection{Multisymplectic discretizations} \label{2D_mult_discret}

We consider the geometric setting of continuum mechanics with the configuration bundle $\mathcal{Y} = \mathcal{X}  \times \mathcal{M} \rightarrow \mathcal{X} = \mathbb{R} \times \mathcal{B} $.
We assume that $ \mathcal{B} $ is a rectangle in $ \mathbb{R} ^2 $ and take $ \mathcal{M}= \mathbb{R} ^2 $.

\subsubsection{Discrete configuration bundle}

\begin{figure}[H] {
\begin{displaymath}
\begin{xy}
\xymatrix{   &    \mathcal{Y} _d=  \mathcal{X} _d \times \mathcal{M}     \ar@<2pt>[dd]^{\pi_d}  \\ 
&\\
\mathcal{U} _d  \ar[ruu]^-{\phi_d}  \ar[r]_-{\phi_{X_d}}  &   \phi _{\mathcal{X} _d}(\mathcal{U} _d) = \mathcal{X} _d\subset \mathcal{X} \ar@<2pt>[uu]^{\varphi_d }  }
\end{xy}
\end{displaymath} }
\caption{\footnotesize Discrete configuration and discrete configuration bundle}\label{D_conf_bundle} 
\end{figure}
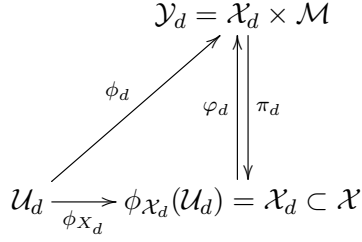

The general discrete setting is the following. One first considers a \textit{discrete parameter space} $\mathcal{U} _d$ and a \textit{discrete base-space configuration}, which is a one-to-one map
\[
\phi _{\mathcal{X} _d}:  \mathcal{U}_d \rightarrow \phi _{ \mathcal{X} _d}( \mathcal{U} _d)= \mathcal{X} _d \subset  \mathcal{X} 
\]
whose image is the discrete spacetime $ \mathcal{X} _d$. The \textit{discrete configuration bundle} is defined as $\pi _d: \mathcal{Y} _d= \mathcal{X} _d \times \mathcal{M}  \rightarrow \mathcal{X} _d$. The \textit{discrete fields} are the sections of the discrete configuration bundle, identified with maps $ \varphi _d: \mathcal{X} _d\rightarrow  \mathcal{M}$. In order to describe both the discrete spacetime as well as the discrete field, one introduces the \textit{discrete configuration} $\phi _d: \mathcal{U} _d \rightarrow  \mathcal{Y} $, from which the discrete base-space configuration and the discrete physical deformation are obtained as $ \phi_{ \mathcal{X} _d}= \pi _d \circ \phi _d$ and $ \varphi _d= \phi _d \circ \phi _{ \mathcal{X} _d} ^{-1} $, see Fig.\,\ref{D_conf_bundle}. This setting is particularly well adapted to situations where the discrete spacetime is also variable, see \cite{LeMaOrWe2003,DeGBRa2016}.

\begin{figure}[H] \centering 
\includegraphics[width=2.4 in]{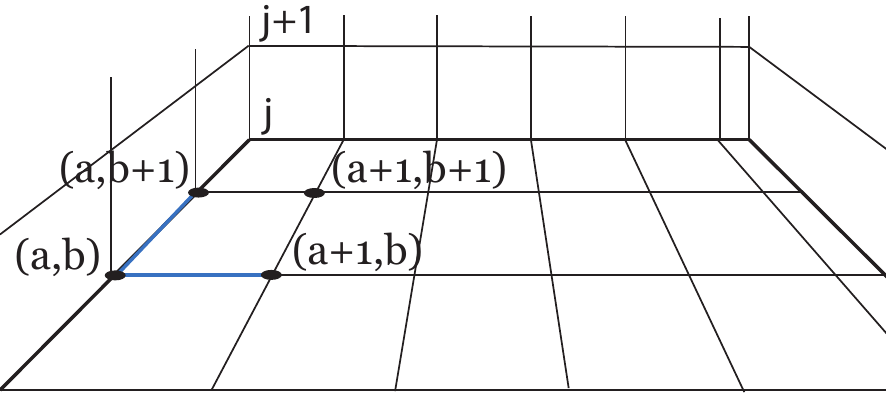}  \vspace{-3pt}  
\caption{\footnotesize Discrete spacetime domain $\mathcal{U} _d$.} \label{pyramid} 
\end{figure}

We consider the discrete parameter space defined by $ \mathcal{U} _d:=\{0,...,j, ...,N\}\times \mathbb{B}_d$, where $\{0,...,j, ...,N\}$ encodes an increasing sequence of time and $\mathbb{B}_d$ parameterizes the nodes and simplexes of the discretization of $ \mathcal{B} $. In this paper we restrict to the case $\mathbb{B}_d= \{0, \ldots, A\}\times \{0, \ldots, B\}$, where $A$ and $B$ are the number of spatial grid points. Therefore, $ \mathcal{U} _{d}=\{0, \ldots, N\} \times \{0, \ldots, A\}\times \{0, \ldots, B\}$ with elements denoted $(j,a,b) \in  \mathcal{U} _d$. The discrete parameter space determines a set of parallelepipeds, denoted $\mbox{\mancube}_{a,b}^j$, and defined by the following eight pairs of indices (see Fig.\,\ref{pyramid})
\begin{equation}\label{cube}
\begin{aligned}
\mbox{\mancube}_{a,b}^j &= \big\{ (j,a,b), (j+1,a,b),(j,a+1,b),(j,a, b+1), (j,a+1,b+1), \\
&\hspace{3cm}  (j+1,a+1, b),(j+1,a,b+1),  (j+1,a+1,b+1) \big\},
\end{aligned} 
\end{equation}
$j= 0, ...,N-1$, $a=0, ..., A-1$, $b=0, ..., B-1$. The set of all such parallelepipeds is denoted $ \mathcal{U} _d^{\,\mbox{\mancube}}$. 

\subsubsection{Discrete Jacobian}\label{DJ2D}

\begin{figure}[H] \centering 
\includegraphics[width=4.8 in]{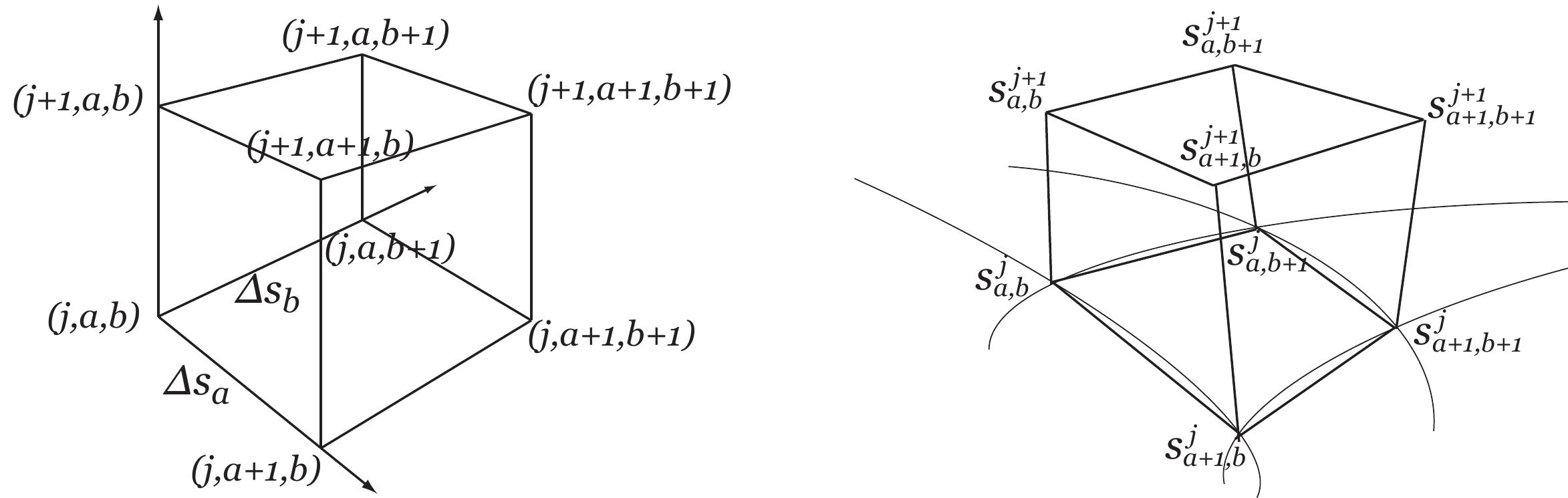} \vspace{-3pt}  
\caption{\footnotesize \textit{On the left}: discrete coordinate system. \textit{On the right}: Nodes of the mesh with Euclidean coordinates. } \label{2D_basis_vectors}
\end{figure}

\begin{figure}[H] \centering 
 \includegraphics[width=4 in]{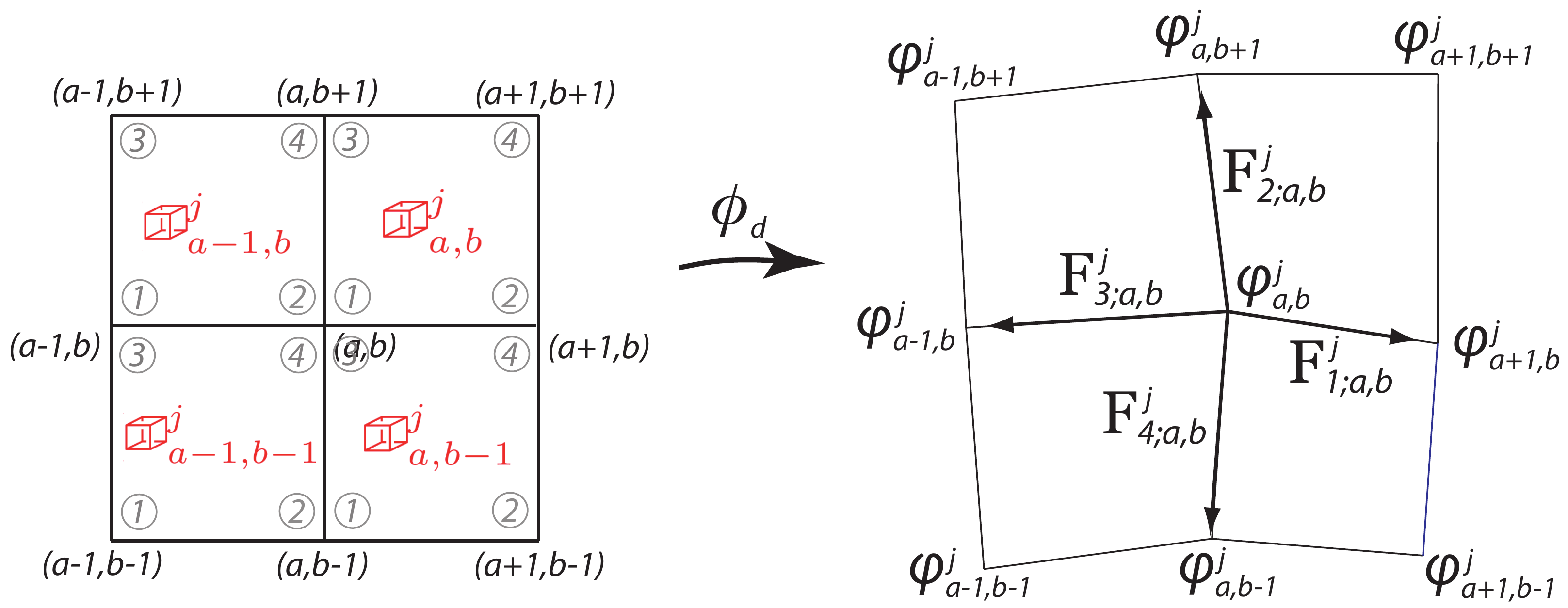} \vspace{-3pt}  
 \caption{\footnotesize Discrete field $\phi_d=\varphi_d\circ  \phi _{ \mathcal{X} _d}$ evaluated on $\mbox{\mancube}_{a,b}^j $, $\mbox{\mancube}_{a,b}^j $, $\mbox{\mancube}_{a,b}^j $, $\mbox{\mancube}_{a,b}^j $ at time $t^j$.} \label{2D_tessellation} 
 \end{figure}

As recalled above, in the continuous setting, the material internal energy function $W( \rho  _0,J)$ of the barotropic fluid depends on the deformation gradient only through its Jacobian. To define the discrete deformation gradient and the discrete Jacobian, we assume that the discrete base space configuration is of the form
\begin{equation}\label{choice_of_DBSC} 
\phi _{\mathcal{X} _d}(j,a,b) = s^j_{a,b}= (t^j, z_a^j, z_b ^j) \in \mathbb{R} \times \mathcal{B},
\end{equation} 
see Fig.\,\ref{2D_basis_vectors}. The discrete field $\varphi_d$ evaluated at $s_{a,b}^j$ is denoted $\varphi_{a,b}^j:= \varphi_d(s_{a,b}^j)$, see Fig.\,\ref{2D_tessellation}.

Given a discrete base space configuration $ \phi _{ \mathcal{X} _d}$ and a discrete field $ \varphi _d$, we define the following four vectors $\mathbf{F}_{\ell;a,b}^j \in \mathbb{R} ^2 $, $\ell=1,2,3,4$ at each node $(j,a,b) \in  \mathcal{U} _d$, see Fig.\,\ref{2D_tessellation} on the right:
\begin{equation}\label{basis_vectors}
\mathbf{F}_{1;a,b}^j= \frac{\varphi_{a+1,b}^j - \varphi_{a,b}^j}{ | s^j_{a+1,b} - s^j_{a,b}| } \quad \text{and} \quad  \mathbf{F}_{2;a,b}^j = \frac{\varphi_{a,b+1}^j- \varphi_{a,b}^j}{ |s^j_{a,b+1} -s^j_{a,b}| }
\end{equation}
\[
\mathbf{F}_{3;a,b}^j= \frac{\varphi_{a-1,b}^j - \varphi_{a,b}^j}{|s^j_{a,b} - s^j_{a-1,b}| } = -  \mathbf{F}_{1;a-1,b}^j \quad \text{and} \quad \mathbf{F}_{4;a,b}^j = \frac{\varphi_{a,b-1}^j- \varphi_{a,b}^j}{ |s^j_{a,b} - s^j_{a,b-1}|} = - \mathbf{F}_{2;a,b-1}^j.
\]
Based on these definitions, the discrete gradient is constructed as follows.

\begin{definition}\label{gradient_definition}
The  \textit{discrete gradient deformations} of a discrete field $ \varphi _d$ at the parallelepiped $\,\mbox{\mancube}_{a,b}^j$ are the four $2 \times 2$ matrices $\mathbf{F}^\ell(\mbox{\mancube}_{a,b}^j)$, $\ell=1,2,3,4$, defined at the four nodes at time $t^j$ of $\;\mbox{\mancube}_{a,b}^j$, as follows:
\begin{equation}\label{gradient_def}
\begin{aligned}
\mathbf{F}_1(\mbox{\mancube}_{a,b}^j) &= \left[\mathbf{F}_{1;a,b}^j \; \; \mathbf{F}_{2;a,b}^j \right], \qquad & \mathbf{F}_2(\mbox{\mancube}_{a,b}^j) &= \left[\mathbf{F}_{2;a+1,b}^j  \; \; \mathbf{F}_{3;a+1,b}^j \right] ,
 \\
\mathbf{F}_3(\mbox{\mancube}_{a,b}^j)&= \left[\mathbf{F}_{4;a,b+1}^j \; \; \mathbf{F}_{1;a,b+1}^j \right], \qquad & \mathbf{F}_4(\mbox{\mancube}_{a,b}^j) &= \left[\mathbf{F}_{3;a+1,b+1}^j  \; \; \mathbf{F}_{4;a+1,b+1}^j \right].
\end{aligned}
\end{equation}
The ordering $\ell=1$ to $\ell=4$ is respectively associated to the nodes $(j,a,b)$, $(j,a,b+1)$, $(j,a+1,b)$, $(j,a+1,b+1)$, see Fig.\,\ref{2D_tessellation} on the left.
\end{definition}

\medskip

It is assumed that the discrete field $ \varphi _d $ is such that the determinant of the discrete gradient deformations are positive.

\begin{definition}\label{2D_discrete_J}
The discrete Jacobians of a discrete field $ \varphi _d $ at the parallelepiped $\,\mbox{\mancube}_{a,b}^j$ are the four numbers $J_\ell(\mbox{\mancube}_{a,b}^j)$, $\ell=1,2,3,4$, defined at the four nodes at time $t^j$ of $\;\mbox{\mancube}_{a,b}^j$ as follows:
\begin{equation}\label{2D_disc_jacob}
 \begin{aligned}
J_1(\mbox{\mancube}_{a,b}^j) &=  | \mathbf{F}_{1;a,b}^j \times \mathbf{F}_{2;a,b}^j |=\mathrm{det} \big( \mathbf{F}_1(\mbox{\mancube}_{a,b}^j)\big),\\
J_2(\mbox{\mancube}_{a,b}^j) &=  | \mathbf{F}_{2;a+1,b}^j \times \mathbf{F}_{3;a+1,b}^j | = \mathrm{det} \big( \mathbf{F}_2(\mbox{\mancube}_{a,b}^j)\big) \\
J_3(\mbox{\mancube}_{a,b}^j) &=  | \mathbf{F}_{4;a,b+1}^j \times \mathbf{F}_{1;a,b+1}^j | =\mathrm{det} \big( \mathbf{F}_3(\mbox{\mancube}_{a,b}^j)\big)\\
J_4(\mbox{\mancube}_{a,b}^j)& =  | \mathbf{F}_{3;a+1,b+1}^j \times \mathbf{F}_{4;a+1,b+1}^j | =\mathrm{det} \big( \mathbf{F}_4(\mbox{\mancube}_{a,b}^j)\big).
\end{aligned}
\end{equation} 
\end{definition}

As a consequence, from relations \eqref{2D_disc_jacob}, the variation of the discrete Jacobian is given by
\begin{equation*}
\delta J_\ell = \frac{\partial \,\mathrm{det}(\mathbf{F}_\ell)}{\partial\, \mathbf{F}_\ell}: \delta \mathbf{F}_\ell =J_\ell(\mathbf{F}_\ell)^{-\mathsf{T}}: \delta \mathbf{F}_\ell,
\end{equation*}
at each $\,\mbox{\mancube}_{a,b}^j$, which is used in the derivation of the discrete Euler-Lagrange equations.

\subsubsection{Discrete Lagrangian}

Recall that the set of all parallelepipeds in the discrete parameter space is denoted $ \mathcal{U} _d^{\,\mbox{\mancube}}$. We write
\[
\mathcal{X} _d^{\,\mbox{\mancube}}:=\phi _{ \mathcal{X} _d}\big( \mathcal{U} _d^{\,\mbox{\mancube}}\big)
\]
the set of all parallelepipeds in $ \mathcal{X} _d$. The \textit{discrete version of the first jet bundle} is given by
\begin{equation}\label{discrete_J1} 
J^1 \mathcal{Y} _d:= \mathcal{X} _d^{\,\mbox{\mancube}} \times \underbrace{ \mathcal{M} \times ... \times \mathcal{M}}_{ \text{$8$ times}} \rightarrow \mathcal{X} _d^{\,\mbox{\mancube}}.
\end{equation} 
Given a discrete field $\varphi _d $, its \textit{first jet extension} is the section of \eqref{discrete_J1} defined by
\begin{equation}\label{discrete_jet_extension} 
j^1 \varphi _d( \mbox{\mancube}_{a,b}^j) = \big( \varphi _{a,b}^j, \varphi _{a,b}^{j+1},  \varphi _{a+1,b}^j, \varphi _{a+1,b}^{j+1},  \varphi _{a,b+1}^j ,  \varphi _{a,b+1}^{j+1}, \varphi _{a+1,b+1}^j ,  \varphi _{a+1,b+1}^{j+1} \big),
\end{equation} 
which associates to each parallelepiped, the values of the field at its nodes. A \textit{discrete Lagrangian} is a map
\[
\mathcal{L}  _d : J ^1 \mathcal{Y} _d \rightarrow \mathbb{R},
\]
see \cite{MaPaSh1998}. The discrete Lagrangian evaluated on a discrete field is denoted as
\[
\mathcal{L} _d\big( j ^1 \varphi _d( \mbox{\mancube})\big).
\]

\medskip

We now consider the case of the barotropic fluid. We assume for simplicity that the mass density $\rho_0$ of the fluid in the reference configuration is a constant number. The case of a Lagrangian density with a nonconstant mass density $ \rho  _0$ is important for applications to stratified flows and can be easily treated by our approach. We consider a class of discrete Lagrangians associated to \eqref{Lagrangian_barotropic} of the form
\begin{equation}\label{Discrete_Lagrangian_2D_fluid_metric}
\mathcal{L}_d\big(j^1 \varphi _d(\mbox{\mancube})\big) = \text{vol}\big( \mbox{\mancube}\big)  \Big( \rho  _0 K_d\big(j^1 \varphi _d(\mbox{\mancube})\big) - \rho  _0 W_d\big( \rho  _0, j^1 \varphi _d(\mbox{\mancube})\big) - \rho  _0 \Pi _d  \big(j^1 \varphi _d(\mbox{\mancube})\big) \Big),
\end{equation}
where $\text{vol}\big( \mbox{\mancube}\big)$ is the volume of the parallelepiped $\mbox{\mancube} \in \mathcal{X} _d^{\,\mbox{\mancube}}$. Examples of $K_d, W_d, \Pi _d$ are given as follows.

\begin{itemize}
\item[--] The discrete kinetic energy $K_d: J^1 \mathcal{Y} _d \rightarrow \mathbb{R}$ is defined as
\begin{equation}\label{kinetic_energy}
K_d\big(j^1 \varphi _d(\mbox{\mancube}_{a,b}^j)\big)  : = \frac{1}{4} \sum_{\alpha=a}^{a+1} \sum_{\beta=b}^{b+1} \frac{1}{2} \big| v_{\alpha,\beta}^j\big|^2 ,
\end{equation}
with $v_{ \alpha , \beta }^j= ( \varphi_{ \alpha , \beta }^{j+1}- \varphi_{ \alpha , \beta }^j)/\Delta t^j$ the discrete velocity.

\item[--] The discrete internal energy $W_d:J^1 \mathcal{Y} _d \rightarrow \mathbb{R}$ is defined as
\begin{equation}\label{2D_internal_energy}
W_d\big( \rho  _0,j^1 \varphi _d(\mbox{\mancube}_{a,b}^j)\big) := \frac{1}{4}  \sum_{\ell=1}^4 W\big( \rho _0 , J_\ell(\mbox{\mancube}_{a,b}^j)\big),
\end{equation} 
where $W$ is the material internal energy of the continuous model and $J_\ell(\mbox{\mancube}_{a,b}^j)$ are the discrete Jacobians associated to $\mbox{\mancube}_{a,b}^j$ at time $t^j$.

\item[--] The discrete potential energy $\Pi_d:J^ 1 \mathcal{Y} _d\rightarrow \mathbb{R}$ is given by
\begin{equation}\label{general_potential}
\Pi_d\big(j^1 \varphi _d(\mbox{\mancube}_{a,b}^j)\big):= \frac{1}{4}  \sum_{\alpha=a}^{a+1} \sum_{\beta=b}^{b+1}\Pi ( \varphi ^j_{\alpha,\beta}),
\end{equation}
where $ \Pi  $ is the potential energy of the continuous model. We shall focus on the gravitation potential $ \Pi  ( \varphi ) = \mathbf{g} \cdot \varphi $, with gravitational acceleration vector $\mathbf{g}$, in which case
\begin{equation}\label{grav_potential}
\Pi_d\big(j^1 \varphi _d(\mbox{\mancube}_{a,b}^j)\big) =  \frac{1}{4}  \sum_{\alpha=a}^{a+1} \sum_{\beta=b}^{b+1}  \mathbf{g} \cdot  \varphi_{\alpha,\beta}^j .
\end{equation}
\end{itemize}
We will also consider a mid-point rule discretization later in \S\ref{conv}.

\subsubsection{Discrete variations and discrete Euler-Lagrange equations}\label{DVDEL_2D}

To simplify the exposition, we assume that the discrete base space configuration is fixed and given by $\phi_{ \mathcal{X} _d} (j,a,b) = (j \Delta  t, a \Delta  s_1, b \Delta s_2)$, for given $ \Delta t$, $ \Delta s_1$, $ \Delta s_2$, that is, we assume that the mesh is fixed\footnote{Mesh deformations can be also considered in this setting and will be explored in a future work.} and matches with the standard basis axis of the Euclidean space (reference frame). In this case, we have $\text{vol}\big( \mbox{\mancube}\big) = \Delta t \Delta s_1 \Delta s_2$ in the discrete Lagrangian \eqref{Discrete_Lagrangian_2D_fluid_metric} and the mass of each $2D$ cell in $ \phi _{ \mathcal{X} _d}( \mathbb{B}_d)$ is  $M= \rho_0 \Delta s_1 \Delta s_2$.

\medskip

The discrete action functional associated to $ \mathcal{L} _d $ is obtained as
\begin{equation}\label{Disc_act_sum}
S_d(\varphi _d)= \sum_{ \mbox{\mancube} \,\in \mathcal{X} _d^{\,\mbox{\mancube}}} \mathcal{L} _d\big(j^1 \varphi_d (\mbox{\mancube}) \big)  = \sum_{j=0}^{N-1} \sum_{a=0}^{A-1} \sum_{b=0}^{B-1} \mathcal{L}_d\big(j^1 \varphi_d (\mbox{\mancube}_{a,b}^j) \big).
\end{equation}
In order to apply the discrete Hamilton principle, we compute the variation $\delta S_d(\varphi _d) $ of the action sum and we get

\vspace{-0.3cm}{\small\begin{equation}
\begin{aligned}
&  \sum_{j=0}^{N-1} \sum_{a=0}^{A-1} \sum_{b=0}^{B-1} \Big[  \frac{M}{4}\left( v_{a,b}^{j} \cdot \delta \varphi_{a,b}^{j+1}  + v_{a+1,b}^j \cdot \delta \varphi_{a+1,b}^{j+1} + v_{a,b+1}^j \cdot \delta \varphi_{a,b+1}^{j+1} + v_{a+1,b+1}^j \cdot \delta \varphi_{a+1,b+1}^{j+1}   \right) \nonumber
\\
& \hspace{3cm}  + A_{a,b}^j \cdot \delta \varphi_{a,b}^j + B_{a,b}^j \cdot \delta \varphi_{a+1,b}^j  + C_{a,b}^j \cdot \delta \varphi_{a,b+1}^j + D_{a,b}^j \cdot \delta \varphi_{a+1,b+1}^j  \Big], \label{2D_action_sum} 
\end{aligned}
\end{equation}}
where we have used the following expressions of the partial derivative of $ \mathcal{L} _d$:
\begin{equation}\label{PD_2}
\begin{aligned} 
D_2 \mathcal{L} _{a,b}^j&= \frac{M}{4} v _{a,b}^j & \qquad  D_6 \mathcal{L} _{a,b}^j&= \frac{M}{4} v _{a,b+1}^j\\
D_4 \mathcal{L} _{a,b}^j&= \frac{M}{4} v _{a+1,b}^j &  \qquad D_8 \mathcal{L} _{a,b}^j&= \frac{M}{4} v _{a+1,b+1}^j
\end{aligned}
\end{equation} 
and we have introduced the following notations for the other partial derivatives
\begin{equation}\label{PD_1}
\begin{aligned} 
D_1 \mathcal{L} _{a,b}^j&= A_{a,b}^j & \qquad   D_5 \mathcal{L} _{a,b}^j&=C _{a,b}^j\\
D_3 \mathcal{L} _{a,b}^j&=B_{a,b}^j & \qquad  D_7 \mathcal{L} _{a,b}^j&= D _{a,b}^j,
\end{aligned}
\end{equation} 
whose expressions are given in Appendix \ref{water_DCEL} for an arbitrary internal energy function $W$.
Note that $ D_k \mathcal{L} _{a,b}^j$ is the partial derivative of $ \mathcal{L} _d$, at $\mbox{\mancube}_{a,b}^j$, with respect to the $k^{th}$ variable, in the order listed in \eqref{discrete_jet_extension}.

\medskip 

Rearranging the expression \eqref{2D_action_sum} we get the discrete Euler-Lagrange equations
\begin{equation} \label{2D_DCEL}
M v_{a,b}^{j-1} + A_{a,b}^j + B_{a-1,b}^j + C_{a,b-1}^j + D_{a-1,b-1}^j =0,
\end{equation}
which correspond to variations $ \delta \varphi _{a,b}^j$ at the interior of the domain. Variations at the spatial boundary gives the boundary conditions
\begin{equation}\label{2D_boundary_cond_space} {\small
\left\{
\begin{aligned}
 & \frac{M}{2} v_{0,b}^{j-1} + A_{0,b}^j + C_{0,b-1}^j =0, &\qquad  &\frac{M}{2} v_{a,0}^{j-1} + A_{a,0}^j + B_{a-1,0}^j =0,  
 \\
 &  \frac{M}{2} v_{A,b}^{j-1} + B_{A-1,b}^j + D_{A-1,b-1}^j =0, &\qquad  &\frac{M}{2} v_{a,B}^{j-1} + C_{a,B-1}^j + D_{a-1, B-1}^j =0, 
 \\
& \frac{M}{4} v_{0,0}^{j-1} +  A_{0,0}^j=0, &\qquad &\frac{M}{4} v_{A,0}^{j-1} +  B_{A-1,0}^j =0, 
\\
& \frac{M}{4} v_{0,B}^{j-1} +  C_{0,B-1}^j =0, &\qquad &\frac{M}{4} v_{A,B}^{j-1} + D_{A-1,B-1}^j =0,
 \end{aligned} \right. }
\end{equation}
while variations at the temporal boundary gives
\begin{equation}\label{2D_boundary_cond_time} {\small
\left\{
\begin{aligned}
 & A_{a,b}^0 + B_{a-1,b}^0 + C_{a,b-1}^0 + D_{a-1,b-1}^0 =0, &\qquad  &M v_{a,b}^{N-1} =0, 
 \\
& A_{0,b}^0 + C_{0,b-1}^0 =0, &\qquad &A_{a,0}^0 +B_{a-1,0}^0 =0, 
\\
&  B_{A-1,b}^0 + D_{A-1,b-1}^0 =0, &\qquad  &C_{a,B-1}^0+ D_{a-1,B-1}^0 =0.
 \end{aligned} \right. }
\end{equation}

We assume that the variations of the discrete field at the spatial boundary are arbitrary so we get the boundary conditions \eqref{2D_boundary_cond_space}. This corresponds to the discrete version of the boundary condition \eqref{zero_pressure_BC}. We assume that the variations at the temporal extremity vanish, hence \eqref{2D_boundary_cond_time} is not imposed.

\subsubsection{Discrete Cartan forms}

In a similar way with the continuous case recalled in \S\ref{MVCM}, the discrete multisymplectic form formula and the discrete Noether theorem are efficiently derived and written by using discrete analogues to the Cartan forms $ \Theta _ \mathcal{L} $ and $ \Omega _ \mathcal{L} $ on the first jet bundle $J^1 \mathcal{Y} $, see \S\ref{MF_CF}, and by using differential exterior calculus.  The discrete Cartan forms of multisymplectic variational integrators are the natural spacetime generalizations of the discrete Cartan forms appearing in time variational integrators, \cite{MaWe2001}.

\medskip 

Given a discrete Lagrangian $ \mathcal{L} _d: J^1 \mathcal{Y} _d \rightarrow \mathbb{R} $, the \textit{discrete Cartan one-forms} are defined on the discrete first jet bundle \eqref{discrete_J1}  as
\begin{equation}\label{DCF} 
\Theta ^\mathtt{p}_{ \mathcal{L}_d}= D_\mathtt{p} \mathcal{L} _d \,{\rm d} \varphi _d^{(\mathtt{p})}, \quad \mathtt{p}=1,...,8,
\end{equation}
see \cite{MaPaSh1998,LeMaOrWe2003,DeGBRa2016}. In \eqref{DCF} we have used the notation
\begin{equation}\label{D_conf_2D}
\varphi_{d}^{(\mathtt{p})} \in \big\{\varphi_{a,b}^j, \varphi_{a,b}^{j+1},\varphi_{a+1,b}^j, \varphi_{a+1,b}^{j+1} , \varphi_{a,b+1}^j,  \varphi_{a,b+1}^{j+1}, \varphi_{a+1,b+1}^j, \varphi_{a+1,b+1}^{j+1}\big\}.
\end{equation}

\medskip 

For the discrete Lagrangian \eqref{Discrete_Lagrangian_2D_fluid_metric} of the barotropic fluid, using \eqref{PD_1} and \eqref{PD_2} we get the following expressions of the discrete Cartan one-forms evaluated on the first jet extension $j^1 \varphi _d( \mbox{\mancube}_{a,b}^j) \in J^1 \mathcal{Y} _d$ of discrete field $ \varphi _d $:
\begin{equation}\label{DCartanForm_extended}
 \begin{aligned}
\Theta_{\mathcal{L}_d }^{1}& = A_{a,b}^j \, {\rm d} \varphi_{a,b}^j, &\qquad &  \Theta_{\mathcal{L}_d }^{2}  =  \frac{M}{4}  v_{a,b}^{j} \, 
{\rm d} \varphi_{a,b}^{j+1},
\\
\Theta_{\mathcal{L}_d }^{3} & = B_{a,b}^j \, {\rm d} \varphi_{a+1,b}^j,&\qquad   &  \Theta_{\mathcal{L}_d }^{4}  =  \frac{M}{4}  v_{a+1,b}^{j} \, 
{\rm d} \varphi_{a+1,b}^{j+1},
\\
\Theta_{\mathcal{L}_d }^{5} & = C_{a,b}^j \, {\rm d} \varphi_{a,b+1}^j, &\qquad  &  \Theta_{\mathcal{L}_d }^{6}  = \frac{M}{4}  v_{a,b+1}^{j+1} \, 
{\rm d} \varphi_{a,b+1}^{j+1},
\\
\Theta_{\mathcal{L}_d }^{7}& = D_{a,b}^j \, {\rm d} \varphi_{a+1,b+1}^j, &\qquad  &  \Theta_{\mathcal{L}_d }^{8} = \frac{M}{4}  v_{a+1,b+1}^{j} \, 
{\rm d} \varphi_{a+1,b+1}^{j+1}.
\end{aligned}
\end{equation}

\medskip 

In order to present the multisymplectic form formula and the discrete covariant Noether theorem, we shall rewrite the differential of the discrete action functional \eqref{Disc_act_sum} in an intrinsic form using the discrete Cartan one-forms. Given a vector field $V_d$ tangent to the discrete configuration $\varphi_d$, we consider its first jet extension $j^1V_d$ which attributes to the set of nodes in $\mbox{\mancube}$ the set of values of $V_d$ on these nodes. With this definition, for a given $ \mbox{\mancube} \in \mathcal{X} ^{ \,\mbox{\mancube}}_d$, we can write the partial derivatives of $ \mathcal{L} _d$ in terms of the discrete Cartan forms as
\begin{equation} \label{pull_back_contrac}
D_\mathtt{p} \mathcal{L}_d\big(j^1 \varphi _d(  \mbox{\mancube})\big) \cdot V_d^{(\mathtt{p})} = \Theta _{ \mathcal{L} _d}^\mathtt{p} \big(j^1 \varphi _d(  \mbox{\mancube})\big) \cdot j^1V_d = \left[ (j^1 \varphi _d)^\ast \big(\mathbf{i}_{j^1 V_d} 
\Theta_{\mathcal{L}_d}^{\mathtt{p}}\big)\right] (\mbox{\mancube}),
\end{equation}
for $\mathtt{p}=1,...,8$. In the last term we have used standard notations from differential calculus: the notation $(j ^1 \varphi ) ^*$ for the pull-back by $j^1 \varphi _d $ of $k$-forms from $ J^1 \mathcal{Y} _d$ to $ \mathcal{X} _d^{\mbox{\mancube}}$ and the notation $ \mathbf{i} _{ j^1V_d}$ for the insertion of a vector in a $k$-form. With these notations, the total derivative of the discrete action functional \eqref{Disc_act_sum} is
\begin{equation} \label{disc_Ham_princ_Cartan}
\mathbf{d} S_d( \varphi _d  ) \cdot V _d = \sum_{\mbox{\mancube} \, \in \mathcal{U} _d^{ \, \mbox{\mancube}}} \;\sum_{\mathtt{p} \in \,\mbox{\mancube}} \;
 \left[(j^1\varphi_d )^* \big(\mathbf{i}_{j^1V_d} 
\Theta_{\mathcal{L}_d}^{\mathtt{p}}\big) \right](\mbox{\mancube}).
\end{equation}
Here $\mathtt{p} \in \mbox{\mancube}$ denotes a node $\mathtt{p}$ of the parallelepiped $\mbox{\mancube}$.
Such a formula is true on any subdomain $ \mathcal{U} '_d \subset  \mathcal{U} _d$, by considering the restricted action $S'_d= S_d|_{ \mathcal{U} '_d}$.

\subsubsection{Discrete multisymplectic form formula and discrete Noether theorem}

When restricted to a solution $ \varphi _d $ of \eqref{2D_DCEL}, the total derivative of $S_d$ reads 
\begin{equation} \label{reduced_disc_Ham_princ}
\mathbf{d} S_d(\varphi_d ) \cdot V _d =\sum_{\mbox{\mancube} \, \in \,\mathcal{U} _d^{ \, \mbox{\mancube}}} \;  \sum_{\mathtt{p};\, \mbox{\mancube}^{(\mathtt{p})} \in \,\partial \mathcal{U} _d} 
 \left[(j^1\varphi_d )^* \big(\mathbf{i}_{j^1V_d} 
\Theta_{\mathcal{L}_d}^{\mathtt{p}}\big) \right](\mbox{\mancube}),
\end{equation}
similarly on any subdomains $\mathcal{U} '_d \subset \mathcal{U} _d$. Note the difference with formula \eqref{disc_Ham_princ_Cartan}. From \eqref{reduced_disc_Ham_princ} two important results are obtained:
\begin{enumerate}
\item  The \textit{discrete multisymplectic form formula}.

It is obtained by taking the exterior derivative of \eqref{reduced_disc_Ham_princ}, evaluating it on the first variations $V_d$, $W_d$ of a solution $ \varphi _d $, and using the rules of exterior differential calculus, which gives
\begin{equation} \label{D_multi_form_formula}
\mathbf{d} \mathbf{d} S_d(\varphi_d ) ( V _d, W_d)= \sum_{\mbox{\mancube} \, \in \,{ \mathcal{U} '}_d^{ \, \mbox{\mancube}}}\; \sum_{\mathtt{p};\, \mbox{\mancube}^{(\mathtt{p})} \in \,
\partial  {\mathcal{U} '}_d}   \left[(j^1\varphi _d  )^* \big(\mathbf{i}_{j^1V_d}\mathbf{i}_{j^1W_d} 
\Omega_{\mathcal{L}_d}^{\mathtt{p}}\big) \right](\mbox{\mancube} )=0,
\end{equation}
for any subdomains $ \mathcal{U} '_d \subset \mathcal{U} _d$.
Here $\Omega_{\mathcal{L}_d}^{\mathtt{p}} =  - \mathbf{d} \Theta_{\mathcal{L}_d}^{\mathtt{p}}$, $\mathtt{p}=1,...,8$ are the \textit{discrete Cartan 2-forms} on $J^1 \mathcal{Y} _d$, see \cite{MaPaSh1998}. This is the discrete version of the multisymplectic form formula \eqref{MFF}. 
It extends to spacetime discretization, the symplectic property of variational integrators, \cite{MaWe2001}.
This formula encodes a discrete version of the reciprocity theorem of continuum mechanics, as well as discrete time symplecticity of the solution flow, see \cite{LeMaOrWe2003}.

\item The \textit{discrete covariant Noether theorem}. 

Consider an action $\Phi: G \times \mathcal{M}  \rightarrow \mathcal{M} $ of a Lie group $G$ on $ \mathcal{M}$. For $ \xi \in \mathfrak{g} $, the Lie algebra of $G$, we denote by $ \xi _ \mathcal{M} $ the infinitesimal generator of the action, i.e. the vector field on $ \mathcal{M} $ defined by
\[
\xi _ \mathcal{M} (m):= \left. \frac{d}{d\varepsilon}\right|_{\varepsilon=0} \Phi _{\exp( \varepsilon \xi )}(m),
\]
for every $m \in \mathcal{M} $. Assume that $ \mathcal{L} _d$ is $G$-invariant with respect this action. As a consequence, the discrete action is also $G$-invariant and we get
\begin{equation}\label{equivariance}
\mathbf{d} S'_d( \varphi _d ) \cdot \xi _ \mathcal{M} ( \varphi _d)=0 \quad \text{for all $\xi \in \mathfrak{g}$}.
\end{equation}
From \eqref{reduced_disc_Ham_princ}, it follows
\begin{equation} \label{global_disc_Noether_th}
\mathbf{d} S'_d(\varphi_d ) \cdot \xi _ \mathcal{M} ( \varphi _d)  = \sum_{\, \mbox{\mancube} \, \in {\mathcal{U} '}_d^{\,\mbox{\mancube}}}\;\sum_{\mathtt{p};\, \mbox{\mancube}^{(\mathtt{p})} \in 
\partial \mathcal{U} '_d}   \left[(j^1\varphi_d )^* \big\langle J ^{\mathtt{p}} _{ \mathcal{L}_d  }, \xi \big\rangle \right](\mbox{\mancube}) =0,
\end{equation}
for every $ \xi \in \mathfrak{g} $, where the \textit{discrete covariant momentum maps} are defined by
\begin{equation}\label{disc_mom_map}
J ^{\mathtt{p}}_{ \mathcal{L} _d } : J ^1 \mathcal{Y}  _d \rightarrow \mathfrak{g}  ^\ast , \quad \langle J ^{\mathtt{p}} _{ \mathcal{L}_d  }, \xi \rangle :=\mathbf{i} _{ \xi _{ J ^1 \mathcal{Y}  _d }} \Theta ^{\mathtt{p}} _{ \mathcal{L}_d },\quad \xi  \in \mathfrak{g},\quad \mathtt{p}=1,\ldots,8.
\end{equation}
In \eqref{disc_mom_map} $ \xi _{J^1 \mathcal{Y} _d}$ is the infinitesimal generator of the action of $G$ induced on $ J^ 1\mathcal{Y} _d$ by the action $ \Phi $ on $\mathcal{M} $. It is given at each $ j^1 \varphi _d (\mbox{\mancube}_{a,b}^j) \in J^ 1 \mathcal{Y} _d $ by 
\[
\begin{aligned}
\small
\xi _{J^1Y_d}\big( j^1 \varphi _d (\mbox{\mancube}_{a,b}^j)\big)= & \left(\mbox{\mancube}_{a,b}^j,
\xi_ \mathcal{M}  (\varphi_{a,b}^j), \xi_ \mathcal{M}  (\varphi_{a, b}^{j+1}),\xi_ \mathcal{M}  (\varphi_{a+1,b}^j), 
\xi _\mathcal{M}  ( \varphi _{a+1,b}^{j+1}) ,  \right.
\\
& \qquad \left. \xi_ \mathcal{M}  (\varphi_{a,b+1}^j), \xi_ \mathcal{M}   (\varphi_{a, b+1}^{j+1}) , \xi_ \mathcal{M}  (\varphi_{a+1,b+1}^j), \xi_\mathcal{M}   (\varphi_{a+1, b+1}^{j+1})   \right).
\end{aligned}
\]
From \eqref{global_disc_Noether_th}, we thus obtain the \textit{discrete covariant Noether theorem}
\begin{equation}\label{cons_law}
\sum_{\mbox{\mancube} \, \in  {\mathcal{U}'} _d^{\, \mbox{\mancube}}} \; \sum_{\mathtt{p};\, \mbox{\mancube}^{(\mathtt{p})} \in 
\partial {\mathcal{U} '}_d} J ^{\mathtt{p}} _{ \mathcal{L}_d  }(\mbox{\mancube}) =0,
\end{equation}
for every subdomain  $ \mathcal{U}' _d \subset \mathcal{U}  _d $ and for $ \varphi _d$ a solution of the discrete Euler-Lagrange equations. This is the discrete version of the covariant Noether theorem \eqref{NT}.
\end{enumerate}

We refer to \cite{LeMaOrWe2003}, \cite{DeGBRa2014}, \cite{DeGBRa2016} for more explanations concerning discrete conservation laws for multisymplectic variational discretizations.

\subsubsection{Symmetries for barotropic fluids}\label{Sym_barotropic}

In absence of the gravitation potential, the discrete Lagrangian \eqref{Discrete_Lagrangian_2D_fluid_metric} is invariant under rotation and translation, i.e., the action of the special Euclidean group $SE(2)$. This follows from inspection of the expressions \eqref{kinetic_energy} and \eqref{2D_internal_energy}, and the expression of the discrete Jacobian.

\medskip

From this invariance, the discrete covariant Noether theorem \eqref{cons_law} is satisfied with the discrete covariant momentum maps $J^\mathtt{p}_{ \mathcal{L} _d}: J^1 \mathcal{Y} _d \rightarrow \mathfrak{se}(2) ^* $ given by
\begin{equation}\label{cov_momap} 
J^\mathtt{p}_{ \mathcal{L} _d}\big( j^1 \varphi _d (\mbox{\mancube}_{a,b}^j)\big) = \big( \varphi ^{(\mathtt{p})} \times D_\mathtt{p} \mathcal{L} _{a,b}^j,D_\mathtt{p} \mathcal{L} _{a,b}^j\big), \qquad \mathtt{p}=1,...,8.
\end{equation}
A consequence of this discrete covariant Noether theorem is the conservation of the \textit{classical} discrete momentum map given in terms of the discrete \textit{covariant} momentum map as
\begin{equation}\label{total_momentum_map} 
\begin{aligned} 
\mathbf{J} _d^j=\mathbf{J} _d ( \boldsymbol{\varphi }^j, \boldsymbol{\varphi }^{j+1}) &= \sum_{a=0}^{A-1}\sum_{b=0}^{B-1} \left( J_{ \mathcal{L} _d}^2 + J_{ \mathcal{L} _d}^4+J_{ \mathcal{L} _d}^6+J_{ \mathcal{L} _d}^8 \right) \\
&= - \sum_{a=0}^{A-1}\sum_{b=0}^{B-1} \left( J_{ \mathcal{L} _d}^1 + J_{ \mathcal{L} _d}^3+J_{ \mathcal{L} _d}^5+J_{ \mathcal{L} _d}^7 \right),
\end{aligned}
\end{equation} 
i.e., $ \mathbf{J} _d^{j+1}=\mathbf{J} _d^j$.
On the left hand side $ \boldsymbol{\varphi }^j=\{ \varphi ^j_{a,b}\mid 0\leq a \leq A-1, \, 0\leq b\leq B-1\}$ is the collection of all positions at time $t^j$. On the right hand sides each of the discrete momentum maps $J^\mathtt{p}_{ \mathcal{L} _d}$ are evaluated on $j^1 \varphi _d (\mbox{\mancube}_{a,b}^j)$. 
We refer to \cite{DeGBRa2014} for details regarding the link between discrete \textit{classical} and discrete \textit{covariant} momentum maps underlying formulas like \eqref{total_momentum_map}.
Boundary conditions play an important role in this correspondence.

\medskip 

From \eqref{cov_momap}, \eqref{total_momentum_map}, and Appendix \ref{water_DCEL}, we get the expression
\begin{equation}\label{3D_momentum_map_fluid} \small
\mathbf{J}_d^j =
\begin{bmatrix}
\vspace{0.2cm}\displaystyle\sum_{a=0}^{A-1} \sum_{b=0}^{B-1}\mathbf{J}_r\big(j^1 \varphi _d (\mbox{\mancube}_{a,b}^j)\big)  
\\ 
\displaystyle\sum_{a=0}^{A-1} \sum_{b=0}^{B-1} \mathbf{J}_l\big(j^1 \varphi _d (\mbox{\mancube}_{a,b}^j)\big)
\end{bmatrix} 
\quad  \text{with} \quad 
\begin{aligned}
\mathbf{J}_r\big(j^1 \varphi _d (\mbox{\mancube}_{a,b}^j)\big) & = \sum_{\alpha =a}^{a+1} \sum_{\beta =b}^{b+1} \varphi_{\alpha,\beta}^j \times \frac{M}{4} v_{\alpha,\beta}^j  \in \mathbb{R},\\
\mathbf{J}_l\big(j^1 \varphi _d (\mbox{\mancube}_{a,b}^j)\big)& = \sum_{\alpha =a}^{a+1} \sum_{\beta =b}^{b+1} \; \frac{M}{4} v_{\alpha,\beta}^j \in \mathbb{R}^2.
\end{aligned}
\end{equation}

\subsection{Incompressible models and penalty method} \label{2D_penalty_method}

In this section we adapt the multisymplectic variational integrator obtained above to the case of incompressible models.

\paragraph{Equality constraint.} As recalled in \S\ref{const_Incompressible}, incompressible models can be obtained from a Lagrange multiplier approach which imposes the equality constraint $J=1$. In the discrete case, one similarly adds to the discrete action \eqref{Disc_act_sum} the corresponding Lagrange multiplier term to get
\begin{equation}\label{const_func_f1}
\widehat{S}_d( \varphi _d, \lambda _d)=  \sum_{ \mbox{\mancube} \,\in \mathcal{X} _d^{\,\mbox{\mancube}}} \Big[ \mathcal{L} _d\big(j^1 \varphi_d (\mbox{\mancube}) \big) + \sum_{\ell=1}^4  \lambda_d ^\ell( \mbox{\mancube}) \big(J_\ell( \mbox{\mancube})  - 1\big) \Big] ,
\end{equation}
which imposes the equality constraint $J_\ell( \mbox{\mancube})  =1$ for the discrete Jacobian, for all parallelepiped $\mbox{\mancube}$ and all $\ell=1,2,3,4$. The critical point condition associated to \eqref{const_func_f1} reads
\begin{equation} \label{equality_constraint}
\nabla_{\varphi_d} \widehat{S}_d( \overline{\varphi}_d, \overline{\lambda} _d)=0 \quad \text{and} \quad  J_d( \mbox{\mancube})  - 1 =0 \;\; \text{on all} \;\; \mbox{\mancube}\in \mathcal{X} _d^{\,\mbox{\mancube}}.
\end{equation}

It is well-known, \cite[p.187]{Rockafellar1993}, that if $\overline{\varphi}_d$ is a local optimal solution of the function $S_d(\varphi _d)$ in \eqref{Disc_act_sum} restricted to the discrete incompressibility equality constraint $\mathcal{C} _d=\{ \varphi _d \mid J_\ell( \mbox{\mancube})  = 1,\;\forall\;\mbox{\mancube}\in \mathcal{X} _d^{\,\mbox{\mancube}}\}$, then there must be a Lagrangian multiplier $\overline{\lambda}_d$ such that \eqref{equality_constraint} holds. However, solving the equations in \eqref{equality_constraint} may not be practical, because inequality constraints also naturally appear, as we will see in the following examples. Let us thus consider a constraint set $\mathcal{C}_d$ associated to the equality constraint $J_\ell=1$ and to inequality constraints $g_i \leq 0$ for $i=1,...m$, i.e.,
\begin{equation}\label{feasible_set}
\mathcal{C}_d = \{ \varphi _d \mid J_\ell(\mbox{\mancube})  =1 \;\; \text{and} \; \; g_i (\mbox{\mancube})\leq 0  \;\; \text{for} \;\; i=1,...m \}.
\end{equation}
Under appropriate conditions, see \cite{RoWe1998}, generalizations of the Lagrangian multiplier rule allow to find the critical points of the action $S_d$ defined in \eqref{Disc_act_sum} under constraints of the form \eqref{feasible_set}, see \cite{DeGBRa2016} for an application to variational integrators. In particular, we have the following necessary condition for $\overline{\varphi}_d \in \mathcal{C}_d$ to be locally optimal: $-\nabla S_d(\overline{\varphi}_d) \in N_{\mathcal{C}_d}(\overline{\varphi}_d)$, where $N_{\mathcal{C}_d}(\overline{\varphi}_d)$ is the normal cone to $\mathcal{C}_d$ at $\overline{\varphi}_d$, which can be viewed as a special case of the calculus of subgradients, i.e., $N_{\mathcal{C}_d}(\overline{\varphi}_d)=\partial I_{\mathcal{C}_d}(\overline{\varphi}_d)$, where $I_{\mathcal{C}_d}$ is the indicator function of $ \mathcal{C}_d$.

When $\mathcal{C}_d$ and $S_d$ are convex, the locally optimal condition is sufficient for $\overline{\varphi}_d$ to be globally optimal. If the ambient space is $\mathcal{M}=\mathbb{R}^n$, this relation reduces to 
\begin{equation} \label{locally_optimal_R^n}
\nabla S _d (\overline{\varphi}_d) + \overline{\lambda}_0^\ell \nabla J_\ell(\overline{\varphi}_d)+\overline{\lambda}_1 \nabla g_1(\overline{\varphi}_d) +.... + \overline{\lambda}_m \nabla g_m(\overline{\varphi}_d) =0,
\end{equation}
for $ \overline{\lambda}_0^\ell \in \mathbb{R} $ and where $\overline{\lambda}_i \geq 0$, for $i=1,...,m$, is non-vanishing only when  $ g_i ( \overline{\varphi }_d)=0$.
When $S_d$ is not convex it is difficult, on a practical viewpoint, to find the global optimal among the set of local optimals, see e.g., \cite{DeGBRa2016}.

\medskip 

Also, given the examples that we will study, where $\mathcal{C}_d$ is convex, instead of solving our problem with the Lagrangian multiplier approach we will introduce quadratic penalty functions $ r \alpha(\varphi_d)$ associated with \eqref{feasible_set}, with penalty parameter $r$, which may be considered as an approximation of the indicator function $I_{\mathcal{C}_d}$, see \cite[p.280]{Mo1973}.
Moreover, if we suppose that for each $r$ there exists a solution $\varphi_r\in \mathcal{M}$ to the problem to minimize $S _d (\varphi_d) + r \,\alpha(\varphi_d)$ with $\varphi_d \in \mathcal{M}$, and that the sequence $\{\varphi_r\}$ is contained in a compact subset of $\mathcal{M}$, we know (see \cite[p.477]{BaShSh2006}) that the limit $\overline{\varphi}_d$ of any convergent subsequence of $\{\varphi_r\}$ when $r \rightarrow \infty$ is an optimal solution to the original problem. If $\mathcal{C}_d$ is nonconvex, a large enough penalty parameter $r$ must be used to get sufficiently close to an optimal solution. In this case computational difficulties could appear to solve the penalty problem and an augmented Lagrangian penalty function can be considered, which enjoys several advantegous properties, see, e.g., \cite{Rockafellar1993,RoWe1998,BaShSh2006}.

\paragraph{Penalty method.} Given the discrete action $S_d$ defined in \eqref{Disc_act_sum}, the Hamilton principle subject to constraints is approximated by a penalty scheme where one seeks the critical points of the action 
\begin{equation}\label{2D_disc_penalty_action} {\small
\begin{aligned}
\widetilde{S}_d( \varphi _d)  = S_d( \varphi _d )- \!\!\!\!\!\!  \sum_{ \mbox{\mancube} \,\in \mathcal{X} _d^{\,\mbox{\mancube}}} \!\! \operatorname{vol}( \mbox{\mancube}) \left(  \Phi_{d0}\big(j^1 \varphi _d (\mbox{\mancube}) \big) +    \Phi_{d1}\big(j^1 \varphi _d (\mbox{\mancube}) \big) + ... + \Phi_{dm}\big(j^1 \varphi _d (\mbox{\mancube}) \big)\right),
\end{aligned}}
\end{equation} 
with quadratic penalty term associated to the incompressibility (equality) constraint
\begin{equation}\label{2D_disc_penalty_function}
\Phi_{d0}\big(j^1 \varphi _d (\mbox{\mancube})\big) := \frac{1}{4}  \sum_{\ell=1}^4 \frac{r}{2} \big(J_\ell(\mbox{\mancube})- 1\big)^2,
\end{equation} 
where $r$ is the penalty parameter, and with quadratic penalty terms $\Phi_{di}\big(j^1 \varphi _d (\mbox{\mancube}) \big)$, $ i=2,...m$, associated to inequality constraints.

\subsubsection{Implicit version} \label{conv}

To test the convergence of our multisymplectic integrator, we will use an implicit integrator obtained from the Lagrangian \eqref{Discrete_Lagrangian_2D_fluid_metric} discretized through the mid-point rule, that is, the discrete internal energy $W_d\big( \rho  _0,j^1 \varphi _d(\mbox{\mancube}_{a,b}^j)\big)$ is evaluated at time $(t^j + t^{j+1})/2$. In this case the four vectors \eqref{basis_vectors} at each node $(j,a,b) \in  \mathcal{U} _d$ are now defined by

{\footnotesize
\begin{equation}\label{basis_vectors_midpoint}
\begin{aligned}
& \mathbf{F}_{1;a,b}^{j+1/2}= \frac{(\varphi_{a+1,b}^j +\varphi_{a+1,b}^{j+1}) - (\varphi_{a,b}^j + \varphi_{a,b}^{j+1}) }{ 2 \Delta s_1 }, \quad   \mathbf{F}_{2;a,b}^{j+1/2} = \frac{(\varphi_{a,b+1}^j +\varphi_{a,b+1}^{j+1})- (\varphi_{a,b}^j + \varphi_{a,b}^{j+1})}{ 2\Delta s_2 }
\\
&\mathbf{F}_{3;a,b}^{j+1/2}= \frac{(\varphi_{a-1,b}^j + \varphi_{a-1,b}^{j+1}) - (\varphi_{a,b}^j + \varphi_{a,b}^{j+1})}{2\Delta s_1 }, \quad \mathbf{F}_{4;a,b}^{j+1/2} = \frac{(\varphi_{a,b-1}^j + \varphi_{a,b-1}^{j+1})- (\varphi_{a,b}^j + \varphi_{a,b}^{j+1})}{2\Delta s_2} .
\end{aligned}
\end{equation}}

\noindent The variation $\delta S_d(\varphi_d)$ of the action sum is found as

\vspace{-0.3cm}{\small\begin{equation}
\begin{aligned}
&  \sum_{j=0}^{N-1} \sum_{a=0}^{A-1} \sum_{b=0}^{B-1} \left[  \left( \frac{M}{4} v_{a,b}^{j} + \mathbb{A}_{a,b}^j\right) \cdot \delta \varphi_{a,b}^{j+1}  + \left( \frac{M}{4} v_{a+1,b}^j + \mathbb{B}_{a,b}^j \right) \cdot \delta \varphi_{a+1,b}^{j+1} \right.  \nonumber
\\
& \hspace{2.2 cm} + \left( \frac{M}{4} v_{a,b+1}^j + \mathbb{C}_{a,b}^j\right) \cdot \delta \varphi_{a,b+1}^{j+1} +  \left( \frac{M}{4}v_{a+1,b+1}^j + \mathbb{D}_{a,b}^j \right) \cdot \delta \varphi_{a+1,b+1}^{j+1}    \nonumber
\\
& \hspace{2.2 cm}  + \left( -\frac{M}{4} v_{a,b}^{j} + \mathbb{A}_{a,b}^j\right) \cdot \delta \varphi_{a,b}^{j}  + \left( -\frac{M}{4} v_{a+1,b}^j + \mathbb{B}_{a,b}^j \right) \cdot \delta \varphi_{a+1,b}^{j}   \nonumber
\\
& \hspace{2.2 cm} \left. + \left( - \frac{M}{4} v_{a,b+1}^j + \mathbb{C}_{a,b}^j\right) \cdot \delta \varphi_{a,b+1}^{j} +  \left(- \frac{M}{4}v_{a+1,b+1}^j + \mathbb{D}_{a,b}^j \right) \cdot \delta \varphi_{a+1,b+1}^{j} \right], \label{2D_action_sum_midpoint} 
\end{aligned}
\end{equation}}

\vspace{-0.2cm}\noindent with coefficients $\mathbb{A}_{a,b}^j$, $\mathbb{B}_{a,b}^j$, $\mathbb{C}_{a,b}^j$, $\mathbb{D}_{a,b}^j$ given in Appendix \ref{water_DCEL} for an arbitrary internal energy function $W$.
It yields the implicit discrete Euler-Lagrange equations 
{\footnotesize
\begin{equation} \label{2D_DCEL_midpoint}
M v_{a,b}^{j-1} -M v_{a,b}^{j} + \mathbb{A}_{a,b}^j + \mathbb{A}_{a,b}^{j-1} + \mathbb{B}_{a-1,b}^j + \mathbb{B}_{a-1,b}^{j-1} + \mathbb{C}_{a,b-1}^j + \mathbb{C}_{a,b-1}^{j-1} + \mathbb{D}_{a-1,b-1}^j + \mathbb{D}_{a-1,b-1}^{j-1} =0,
\end{equation}}

\vspace{-0.4cm}\noindent with the spatial boundary conditions
{\footnotesize
\begin{equation} \label{2D_spat_bound_midpoint}
\begin{aligned}
&\frac{M}{2} v_{0,b}^{j-1} -\frac{M}{2} v_{0,b}^{j} + \mathbb{A}_{0,b}^j + \mathbb{A}_{0,b}^{j-1} + \mathbb{C}_{0,b-1}^j + \mathbb{C}_{0,b-1}^{j-1}  =0,
\\
& \frac{M}{2} v_{a,0}^{j-1} -\frac{M}{2} v_{a,0}^{j} + \mathbb{A}_{a,0}^j + \mathbb{A}_{a,0}^{j-1} + \mathbb{B}_{a-1,0}^j + \mathbb{B}_{a-1,0}^{j-1}  =0,
\\
& \frac{M}{2} v_{A,b}^{j-1} - \frac{M}{2} v_{A,b}^{j}  + \mathbb{B}_{A-1,b}^j + \mathbb{B}_{A-1,b}^{j-1} +  \mathbb{D}_{A-1,b-1}^j + \mathbb{D}_{A-1,b-1}^{j-1} =0,
\\
& \frac{M}{2} v_{a,B}^{j-1} - \frac{M}{2} v_{a,B}^{j}  + \mathbb{C}_{a,B-1}^j + \mathbb{C}_{a,B-1}^{j-1} + \mathbb{D}_{a-1,B-1}^j + \mathbb{D}_{a-1,B-1}^{j-1} =0,
\\
&\frac{M}{4} v_{0,0}^{j-1} -\frac{M}{4} v_{0,0}^{j} + \mathbb{A}_{0,0}^j + \mathbb{A}_{0,0}^{j-1}   =0, \qquad  \frac{M}{4} v_{A,0}^{j-1} -\frac{M}{4} v_{A,0}^{j}  + \mathbb{B}_{A-1,0}^j + \mathbb{B}_{A-1,0}^{j-1}  =0,
\\
&\frac{M}{4} v_{0,B}^{j-1} -\frac{M}{4} v_{0,B}^{j} + \mathbb{C}_{0,B-1}^j + \mathbb{C}_{0,B-1}^{j-1}  =0, \qquad \frac{M}{4} v_{A,B}^{j-1} - \frac{M}{4} v_{A,B}^{j}  + \mathbb{D}_{A-1,B-1}^j + \mathbb{D}_{A-1,B-1}^{j-1} =0.
\end{aligned}
\end{equation}}

\vspace{-0.3cm}\noindent The corresponding temporal boundary conditions can be computed similarly.

\subsection{Numerical simulation}\label{2D_examples}

We evaluate the properties of the proposed multisymplectic integrator for barotropic and incompressible ideal fluid models with the case of a free boundary fluid and with the case of a fluid flowing on a surface and impacting an obstacle. In the two cases we present an explicit integrator, while we consider an implicit integrator (mid-point rule) for the convergence tests.

\subsubsection{Example 1: fluid motion in vacuum with free boundaries} \label{Example1_2D}

Consider a barotropic fluid with properties $\rho_0= 997 \, \mathrm{kg/m}^2$, $\gamma =6 $, $A = \tilde{A} \rho_0^{-\gamma}$ with $\tilde A = 3.041\times 10^4$ Pa, and $B = 3.0397\times 10^4$ Pa. The size of the discrete reference configuration at time $t^0$ is $1\mathrm{m} \times 1\mathrm{m}$, with space-steps $\Delta s_1= \Delta s_2= 0.0714$m. We consider both the compressible barotropic fluid ($r=0$) and the incompressible case with penalty parameters $r= 10^6$ and $r= 10^7$. The time-steps are $\Delta t=10^{-3}$ when $r\in \{0,10^6\}$, and $\Delta t=5 \times 10^{-4}$ when $r= 10^7$. Initial perturbations (tiny compression) are applied at time $t^1$ on nodes $(4,0)$ and $(5,1)$. 
Note that in the incompressible case, using the penalty approach allows to treat this slight compression as initial condition.
\medskip

Regarding the incompressible models, in the continuous setting the internal energy $W$ plays no role since its effect is absorbed into the gradient of the pressure. In the discrete case, when using the penalty method for incompressible fluids it is advantageous to include the internal energy of the isentropic perfect fluid, as the case $W=0$ needs to deal with a much higher penalty term.

\medskip

\medskip 


\begin{figure}[H] \centering 
\includegraphics[width=1.75 in]{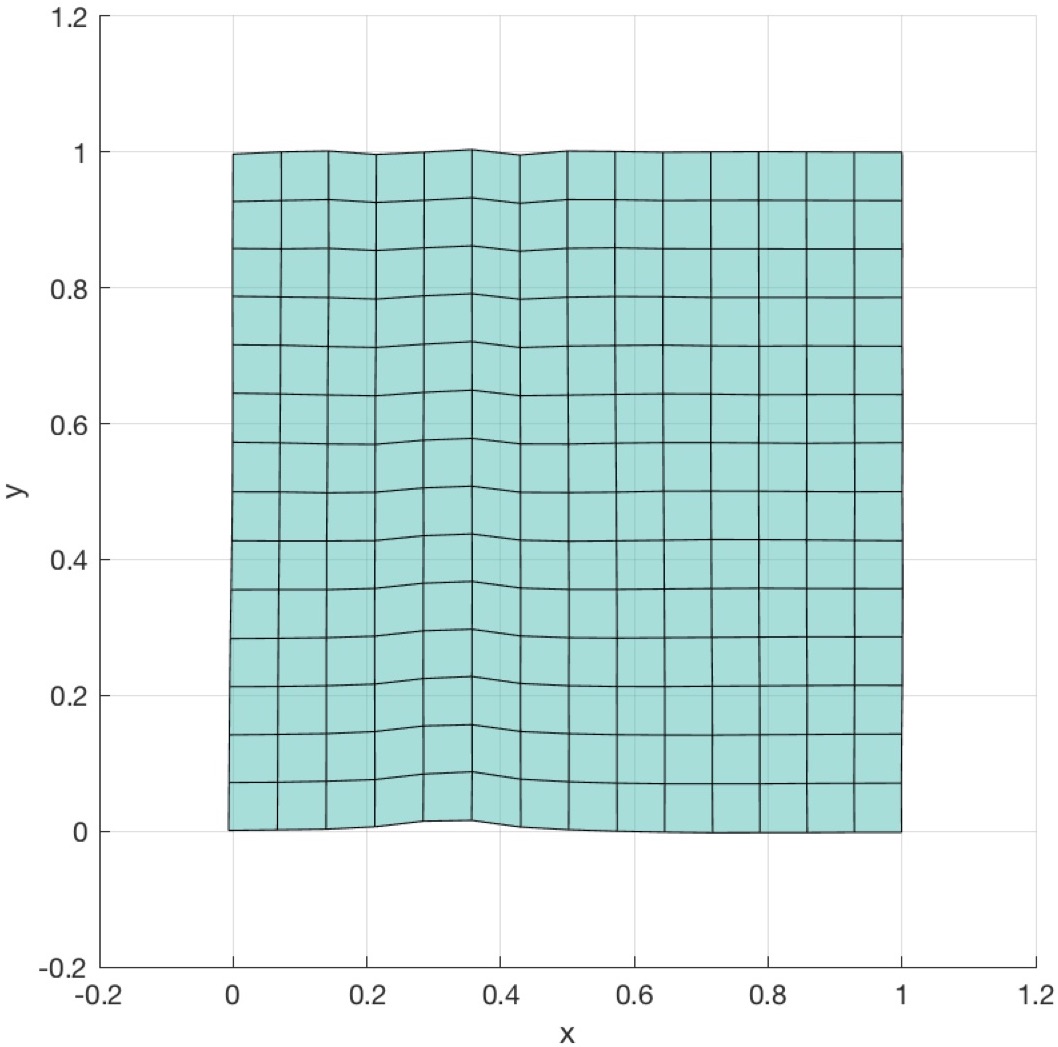} \vspace{-0pt} \;
\includegraphics[width=1.8 in]{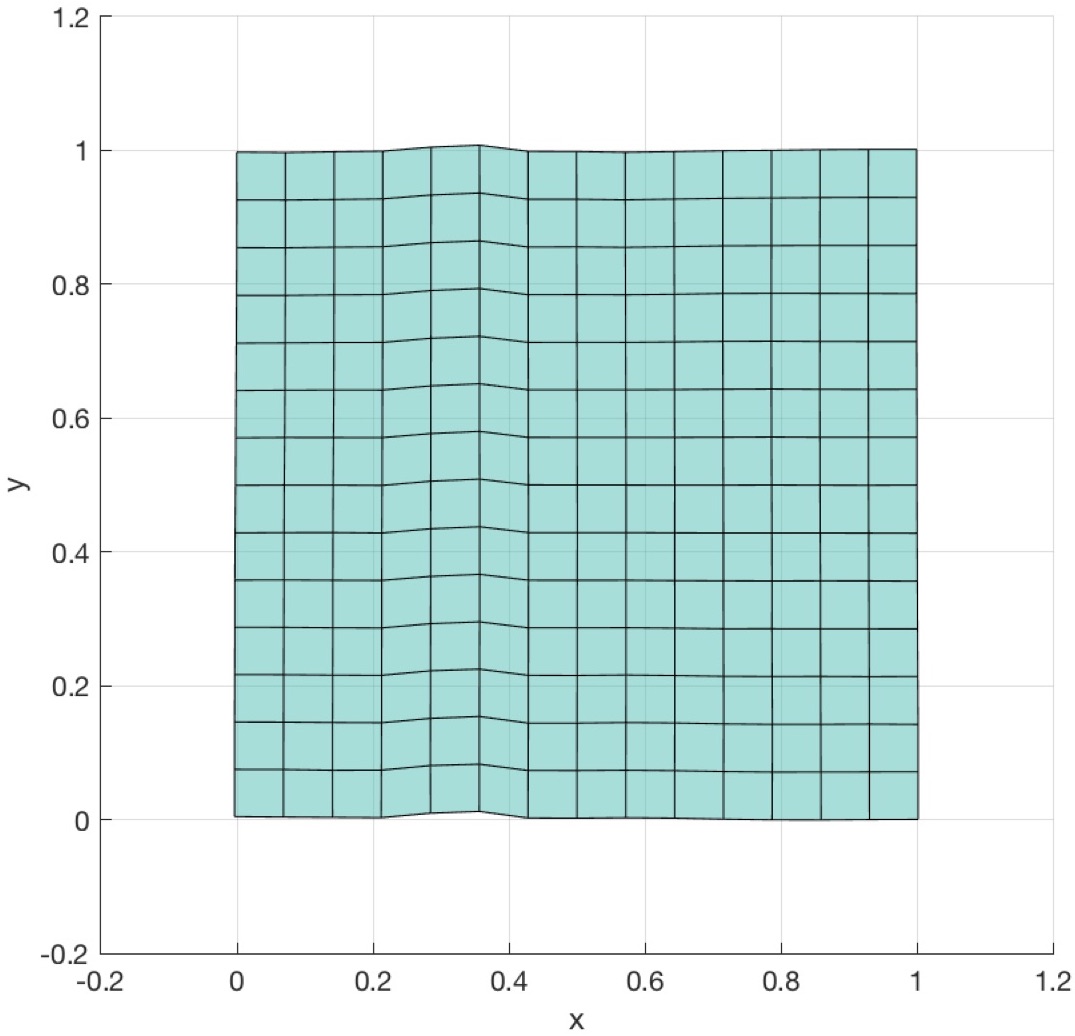} \vspace{-3pt} \;
\includegraphics[width=1.77 in]{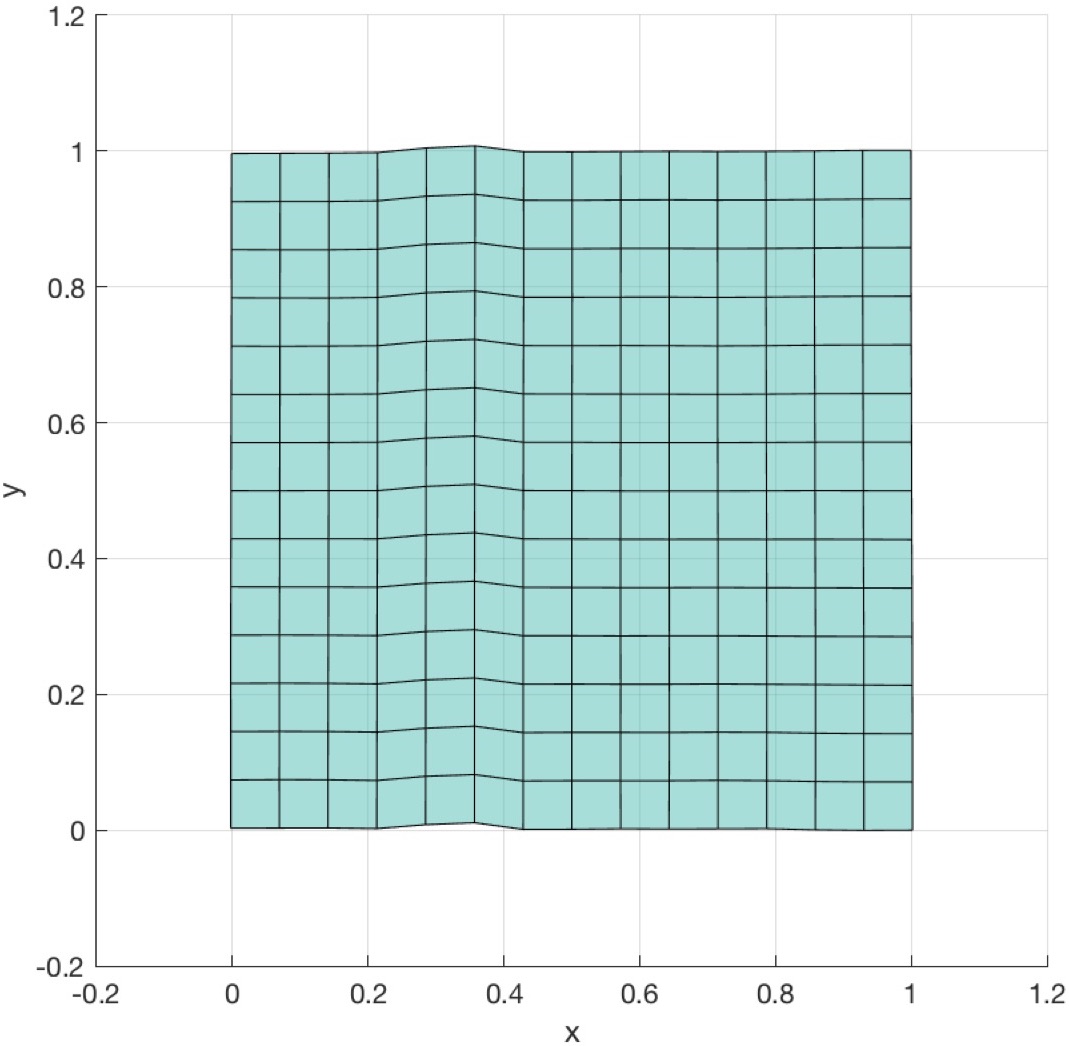} \vspace{-3pt} \\
\includegraphics[width=1.3 in]{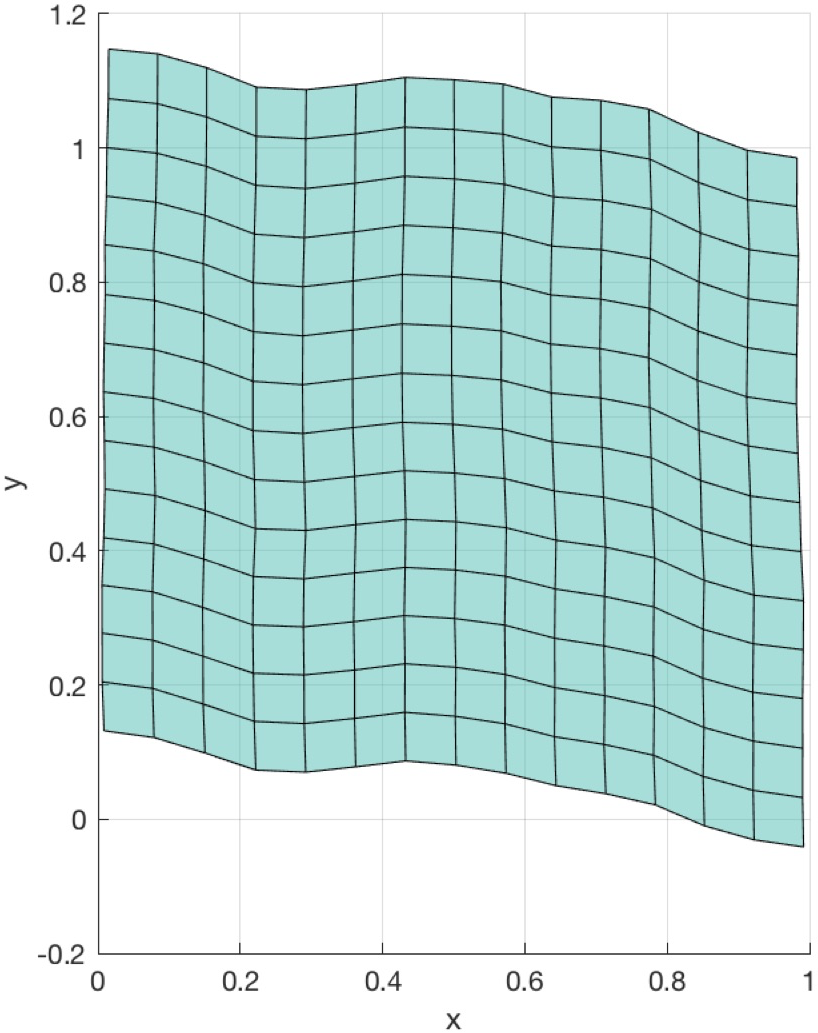} \vspace{-3pt} \qquad  \qquad
\includegraphics[width=1.3 in]{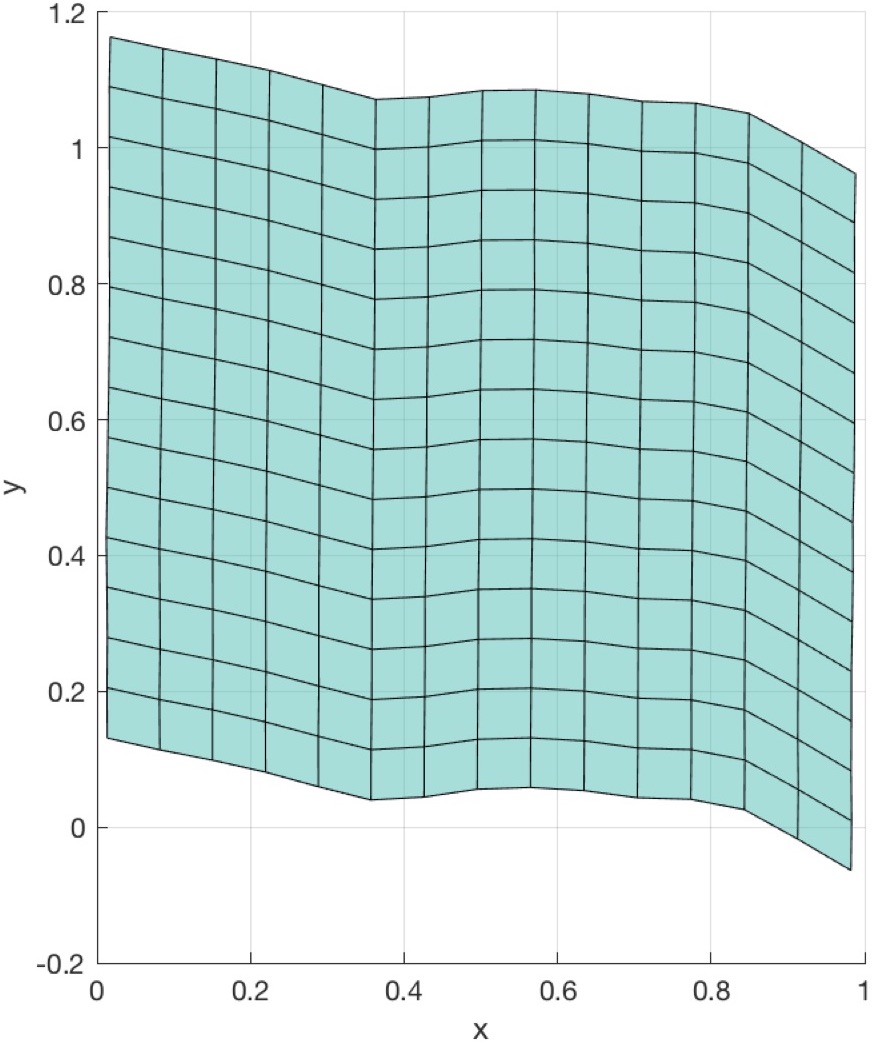} \vspace{-3pt} \qquad \qquad
\includegraphics[width=1.25 in]{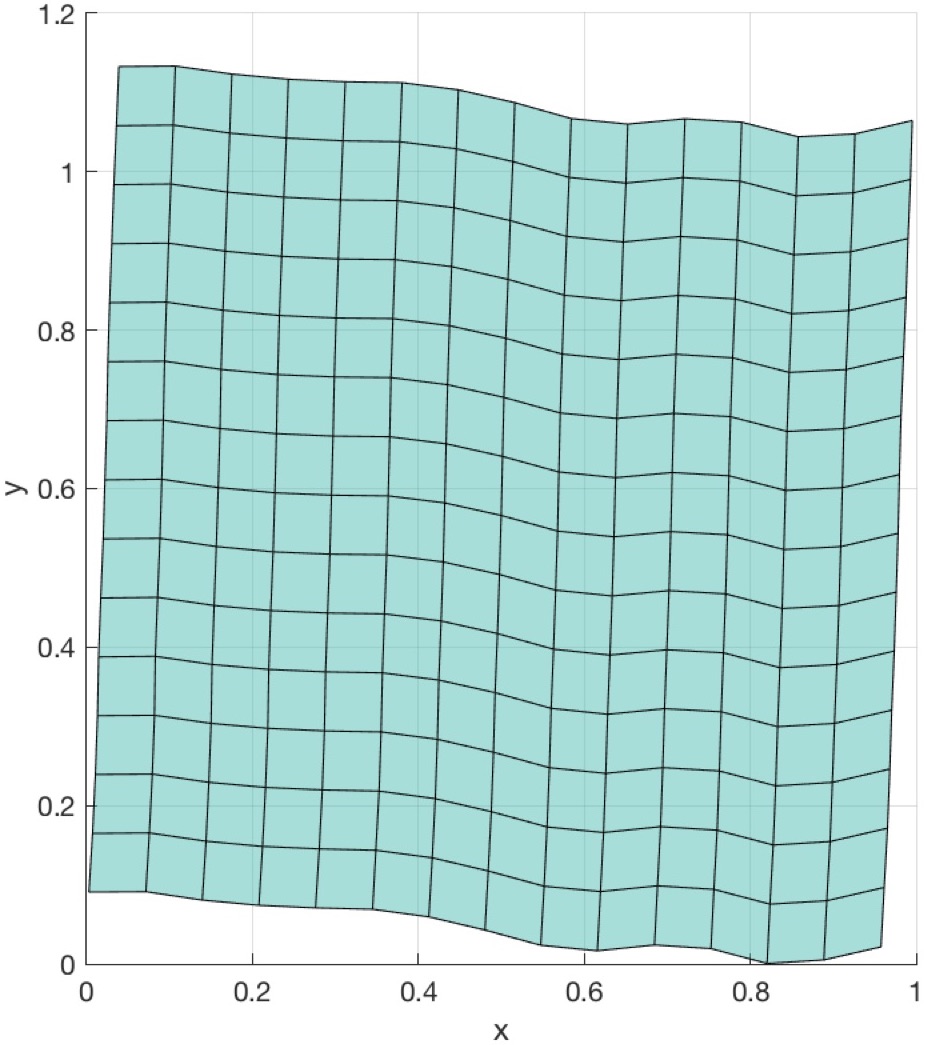} \vspace{-3pt}
\caption{\footnotesize\textit{Left to right}: barotropic and incompressible ideal fluid models ($r= 10^6$ and $r= 10^7$). \textit{Top to bottom}: after $0.1$s and $6$s.} \label{water_perturb} 
\end{figure}
We observe that  the compressible model exhibits enhanced deformation. The different behavior of the compressible and incompressible model will be even more noticeable in the following test, see Fig.\,\ref{water_contact_incompressible} with respect to Fig.\,\ref{water_contact}.

\medskip
 
\begin{figure}[H] \centering 
   \includegraphics[width=1.8 in]{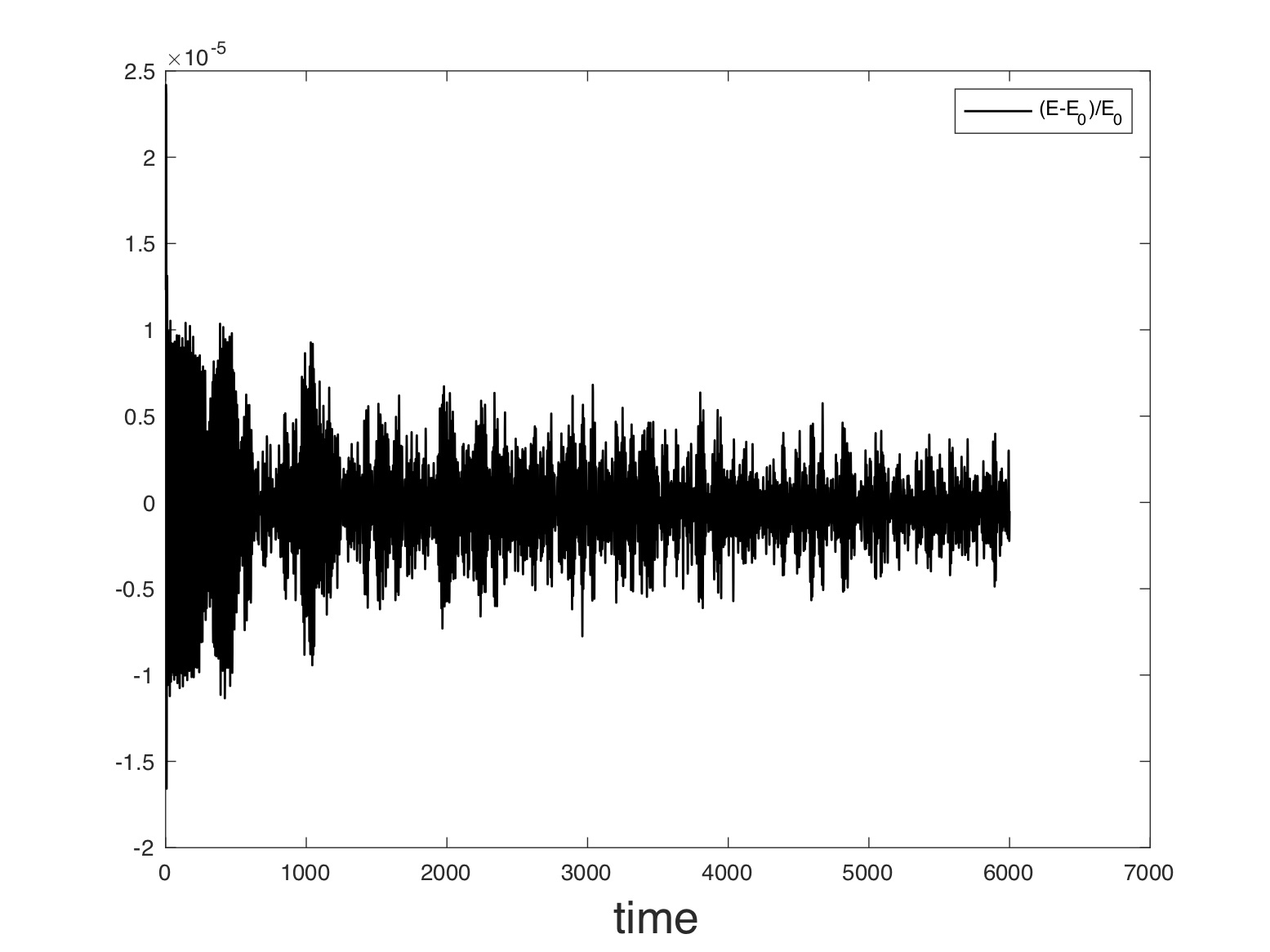} \vspace{-3pt} \; 
  \includegraphics[width=1.8 in]{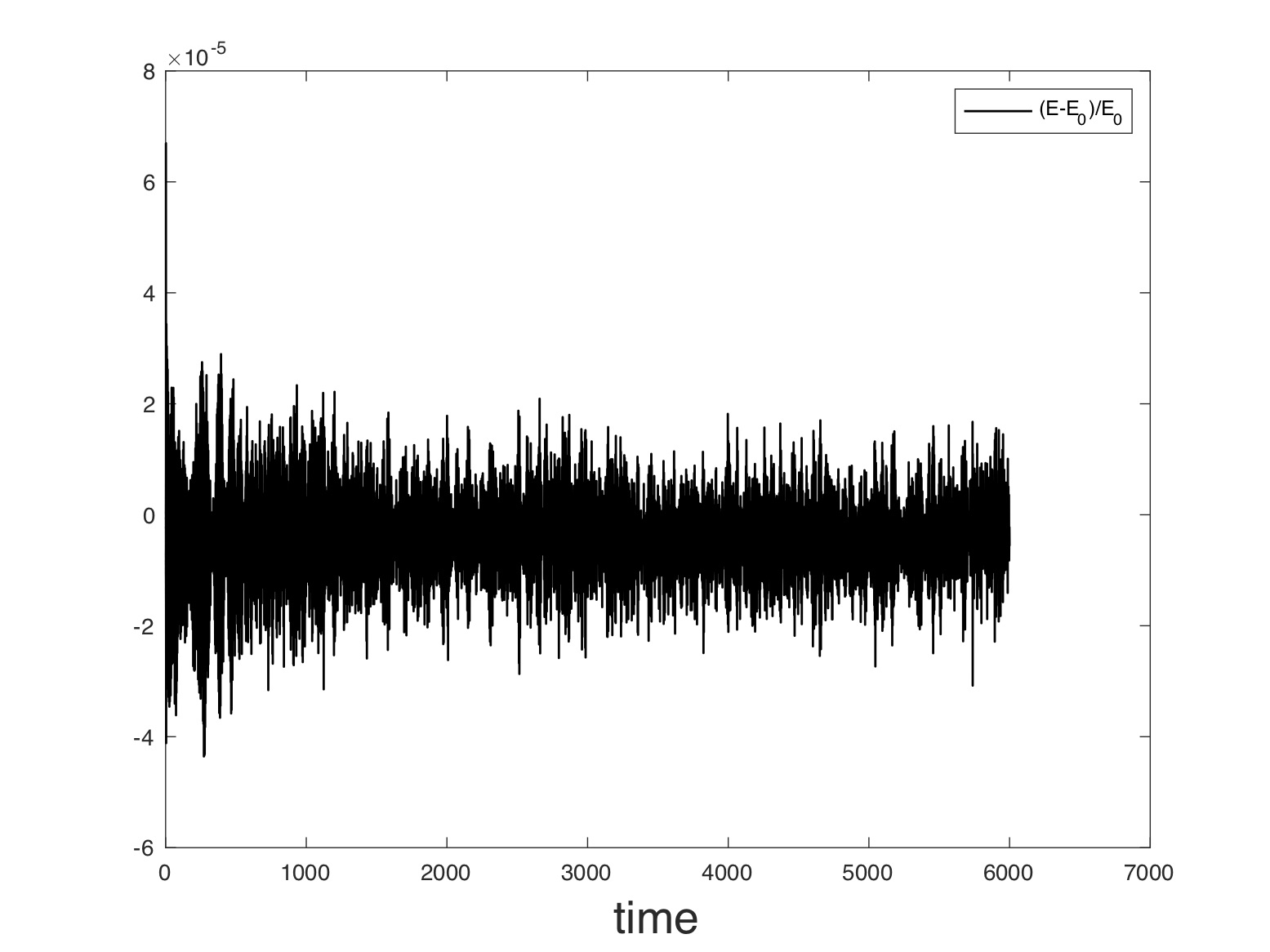} \vspace{-3pt}  \;
  \includegraphics[width=1.8 in]{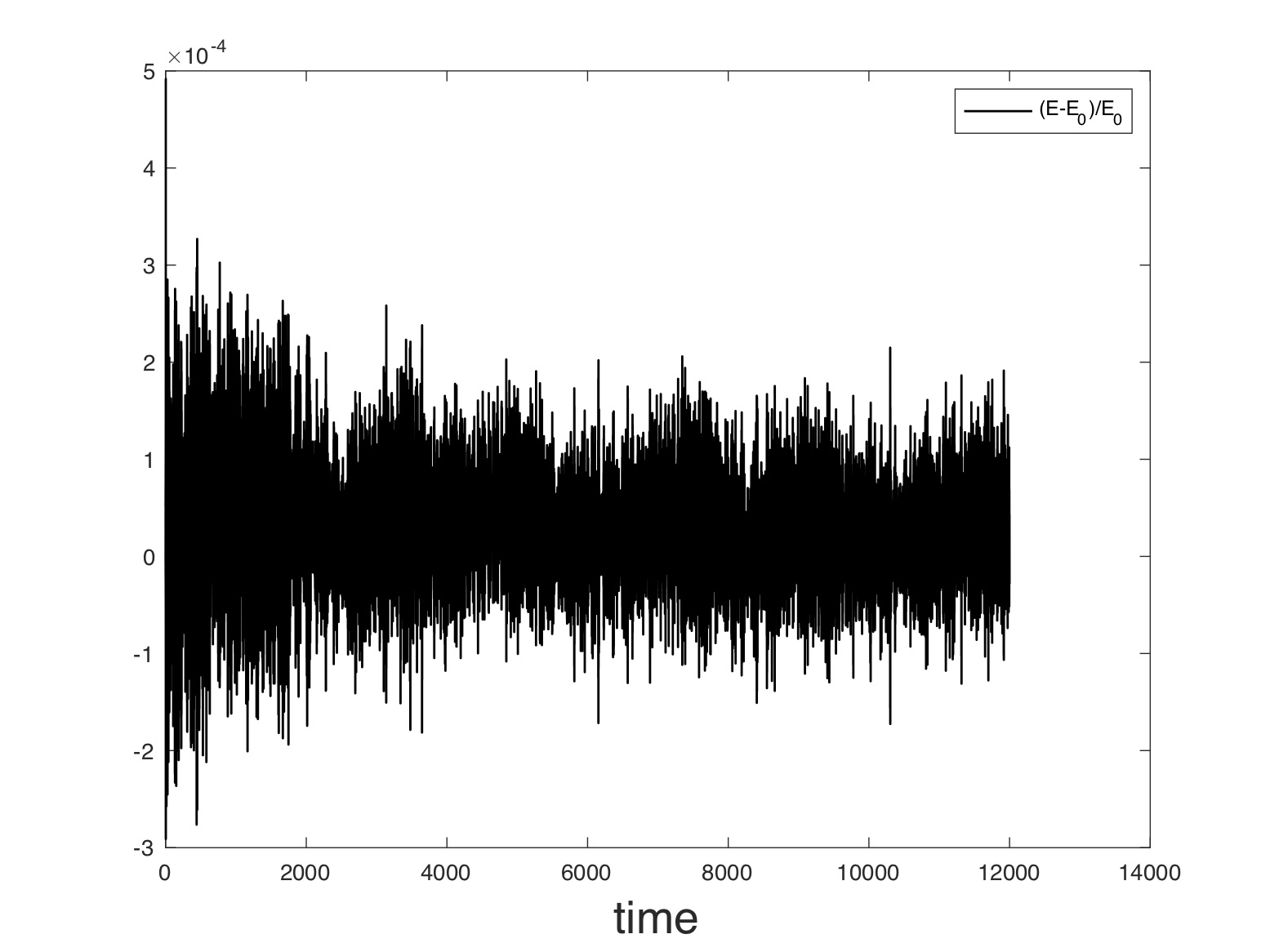} \vspace{-3pt}  
  \\  \includegraphics[width=1.8 in]{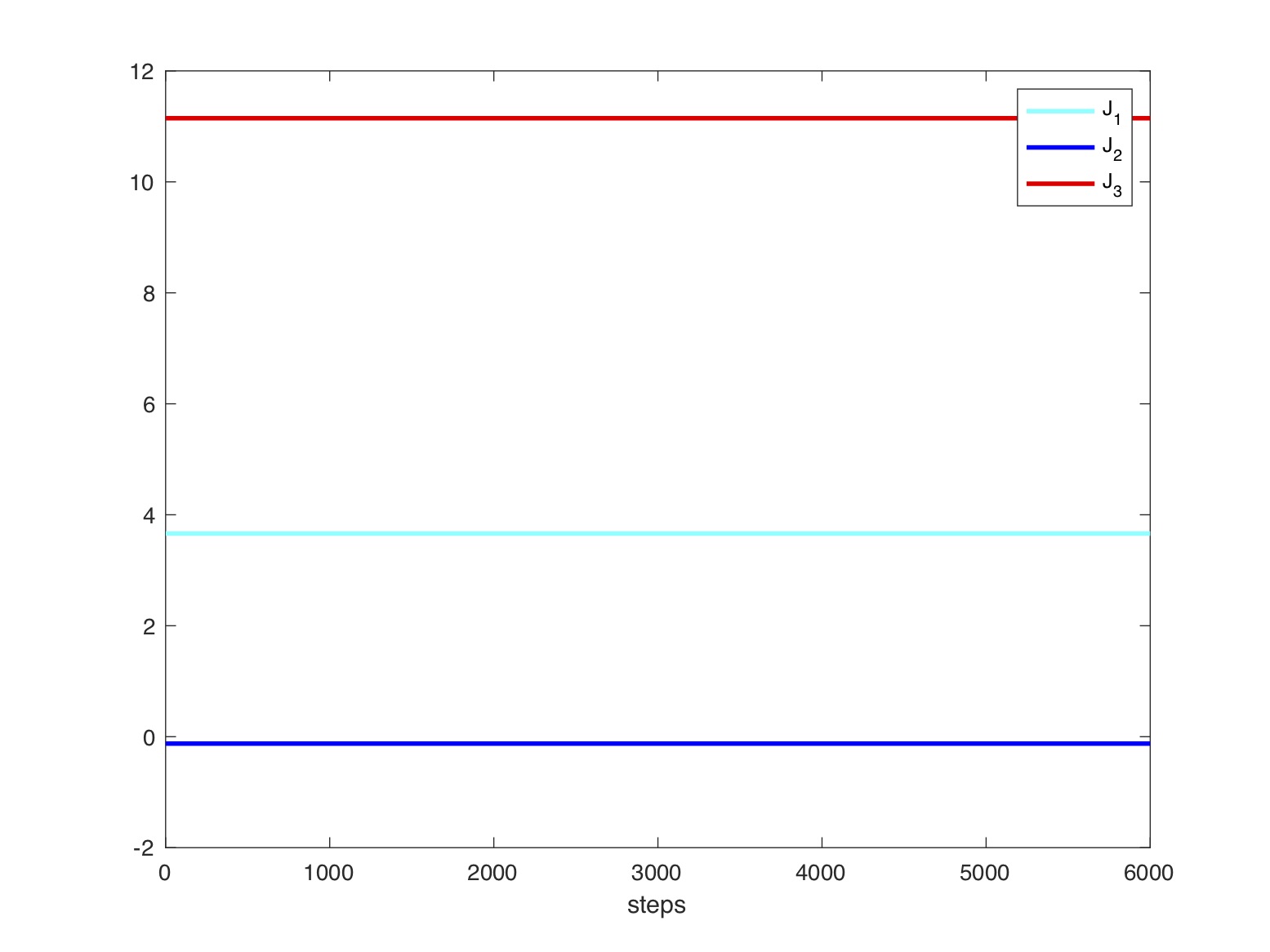} \vspace{-3pt} \;  \includegraphics[width=1.8 in]{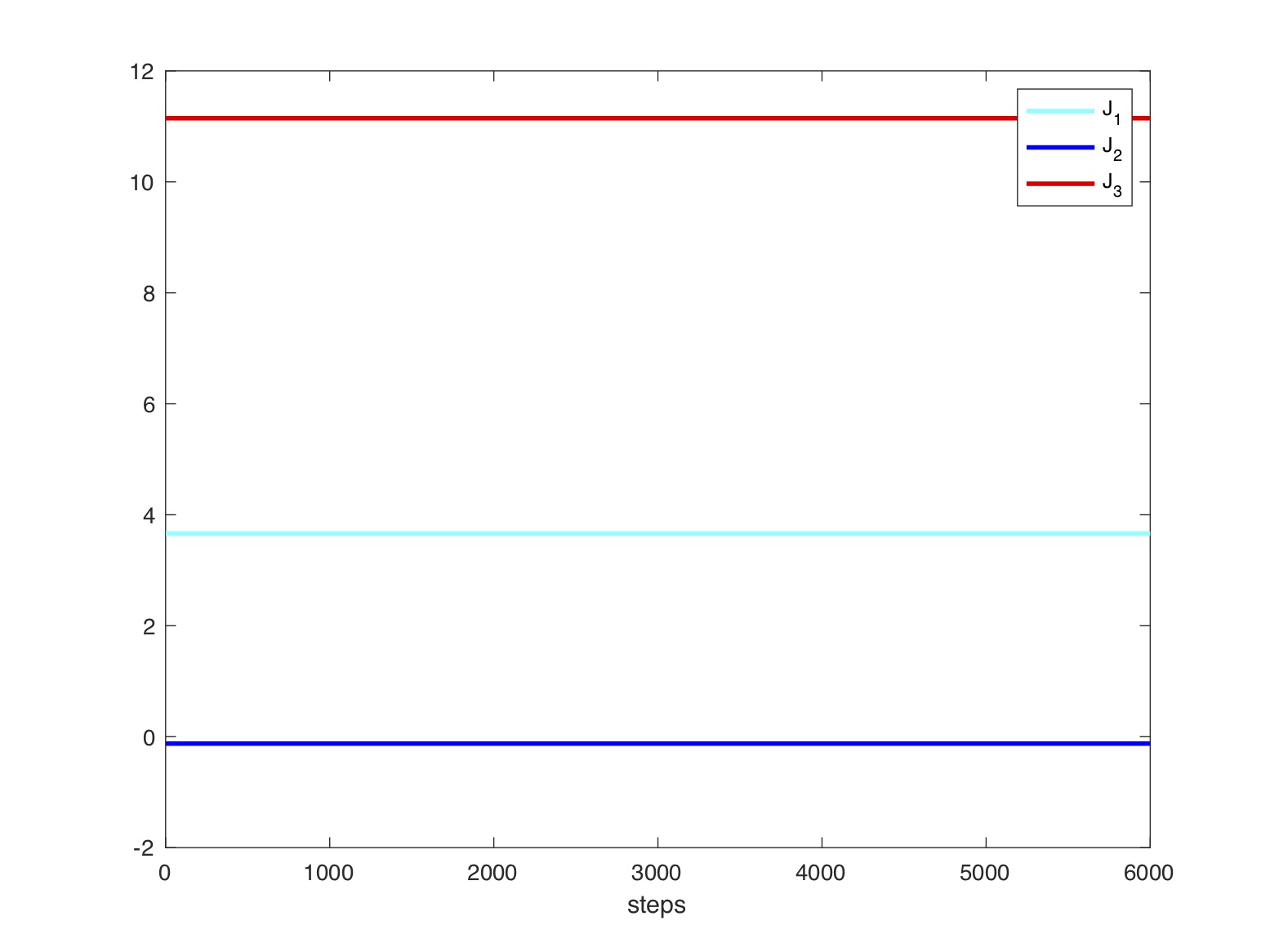} \vspace{-3pt} \;
  \includegraphics[width=1.8 in]{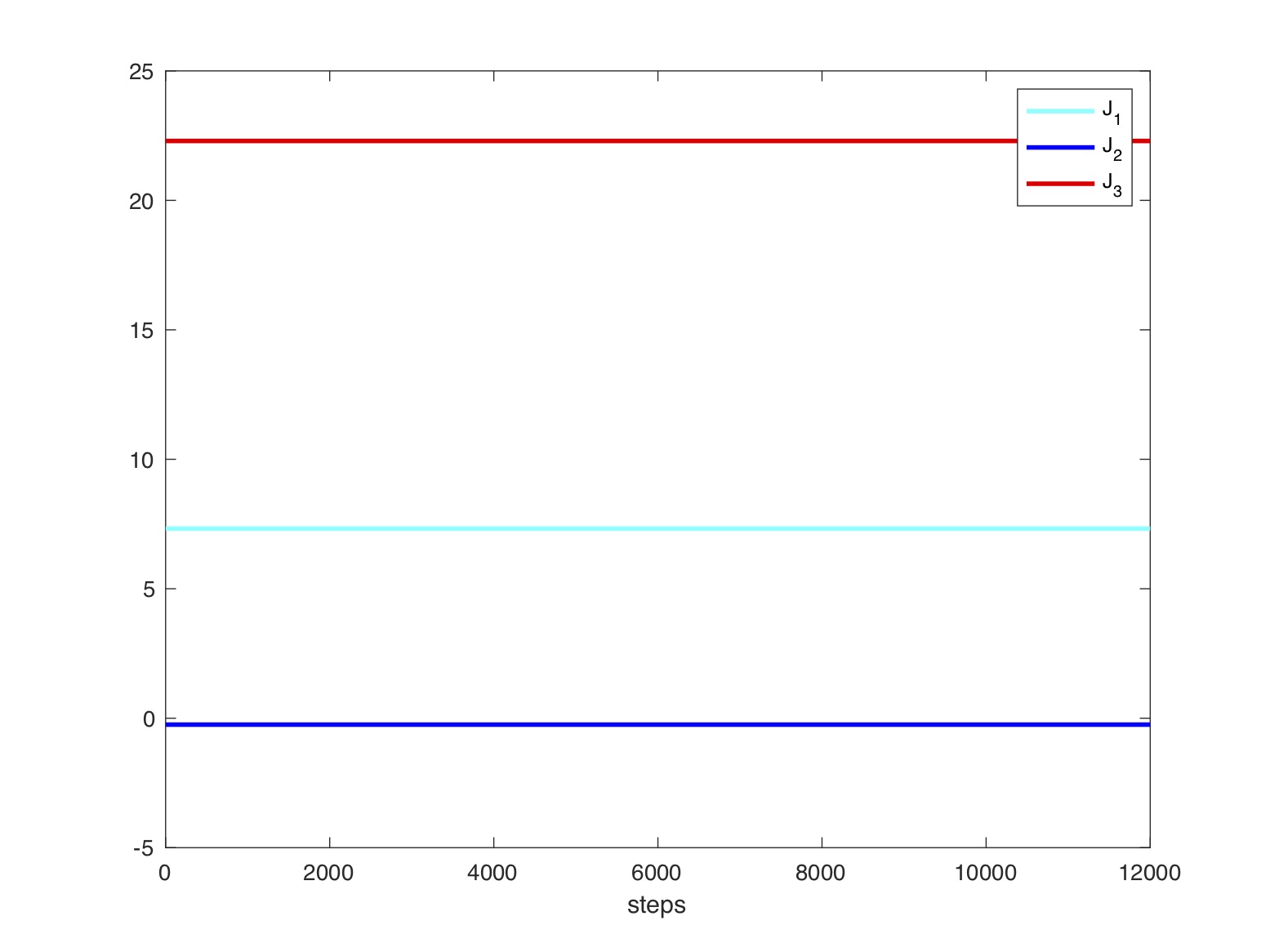} \vspace{-3pt}
  \caption{\footnotesize \textit{Left to right}: barotropic and incompressible ideal fluid models ($r= 10^6$ and $r= 10^7$). \textit{Top to bottom}: Relative energy and momentum map evolution during $6$s.}\label{water_perturb_sym_energ} 
  \end{figure}

The discrete Lagrangian is invariant under rotation and translation, hence from the discrete Noether theorem the angular and linear momentum map \eqref{3D_momentum_map_fluid} are preserved. Energy and momentum preservation is illustrated in Fig.\,\ref{water_perturb_sym_energ}.
  
\subsubsection{Example 2: impact against an obstacle of a fluid flowing on a surface} \label{2D_water_impact}

\paragraph{Inequality constraint.} Let us consider a fluid subject to gravity and flowing without friction on a surface until it comes into contact with an obstacle. 
The gravitational potential is described by \eqref{grav_potential}
with $\mathbf{g}=  \mathrm{g}  \mathbf{E}_2$. We consider both a barotropic fluid and an incompressible ideal fluid.

\medskip 

We impose the following constraints on the configuration:
\begin{itemize}
\item[--] The fluid is bounded below by a rigid surface, defined by the inequality constraint $\Psi_1(\varphi_{a,b}^j) \leq 0$ verified for all $\varphi_{a,b}^j$. 
\item[--] There is a second inequality constraint $\Psi_2(\varphi_{a,b}^j) \leq 0$, verified for all $\varphi_{a,b}^j$, which forces the fluid to stay outside of the obstacle.  
\end{itemize}

For barotropic fluid and incompressible ideal fluid the problems to solve are respectively described as follows
\begin{itemize}
\item[--] $(\mathcal{P}_1)$ \quad \textit{Find the critical points of the action $S_d$ defined in \eqref{Disc_act_sum} subject to the inequality constraints $\Psi_\alpha(\varphi_{a,b}^j) \leq 0$, $\alpha =1,2$, for all nodes.} 
\item[--] $(\mathcal{P}_2)$ \quad \textit{Find the critical points of the action $S_d$ defined in \eqref{Disc_act_sum} subject to the equality constraint $J_\ell( \mbox{\mancube})  =1$, see \eqref{const_func_f1}, and the inequality constraints $\Psi_\alpha(\varphi_{a,b}^j) \leq 0$, $\alpha =1,2$, for all nodes.}
\end{itemize}

We solve the previous problems via the penalty method. For problem $\mathcal{P}_1$, we must find the critical points of the action 
\begin{align}
\widetilde{S}_d(\varphi_d)  & = S_d(\varphi_d) - \sum_{j=0}^{N-1} \sum_{a=0}^{A-1} \sum_{b=0}^{B-1}\Delta t \Delta s_1\Delta s_2  \left( \Phi_{d1}(\varphi_{a,b}^j) + \Phi_{d2}(\varphi_{a,b}^j)\right) , \label{Discrete_penalty_contact_2D}
\\
\Phi_{d\alpha}(\varphi_{a,b}^j) & = \frac{1}{2} K_\alpha |\Psi_\alpha(\varphi_{a,b}^j)|^2  \quad  \text{with} \quad \left\{\begin{array}{lc} K_\alpha \in \;] 0, \infty [ & \text{if \; $\Psi_\alpha(\varphi_{a,b}^j) \geq 0$} \\
K_\alpha=0 & \text{if \; $\Psi_\alpha(\varphi_{a,b}^j) < 0$} \end{array}\right. .  \label{2D_inequality_penalty_term}
\end{align} 
For problem $\mathcal{P}_2$ we add the penalty function \eqref{2D_disc_penalty_function}, associated to the equality constraint $J^\ell( \mbox{\mancube})=1$, into the discrete action \eqref{Discrete_penalty_contact_2D}.
 
\paragraph{Test.} Consider a barotropic fluid model with properties $\rho_0= 997 \, \mathrm{kg/m}^2$, $\gamma =6 $, $A = \tilde{A} \rho_0^{-\gamma}$ with $\tilde A = 3.041\times 10^4$ Pa, and $B = 3.0397\times 10^4$ Pa. The size of the discrete reference configuration at time $t^0$ is $2 \, \mathrm{m} \times 0.4 \, \mathrm{m}$, with time-step $\Delta t=10^{-4}$ and  space-steps $\Delta s_1=0.0625$m, $\Delta s_2= 0.033$m. The values of the impenetrability penalty parameters are chosen as $K_1=4.8\times10^{10}$, $K_2=4.8 \times10^{6}$. For the incompressible case we consider the penalty parameter $r=5\times10^{8}$.

\medskip 

\begin{figure}[H] \centering 
\hspace{1.6cm} \includegraphics[width=4.1 in]{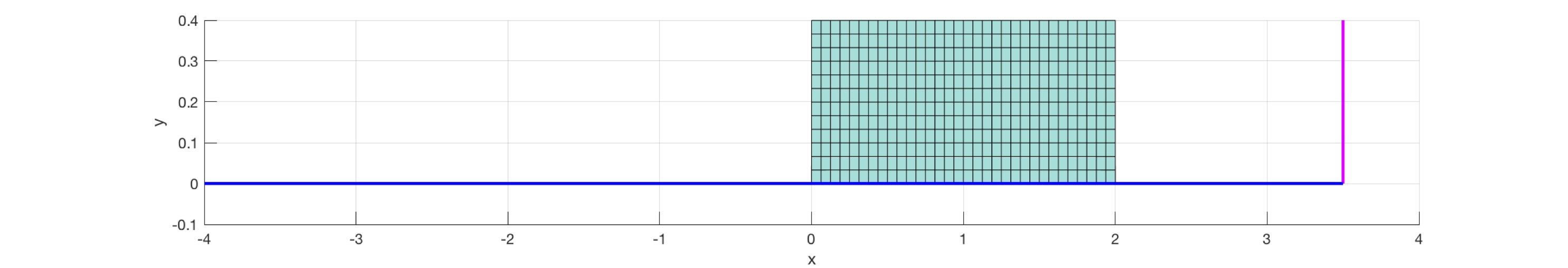} \vspace{-3pt}  
     \\
\hspace{1.6cm}  \includegraphics[width=4.1 in]{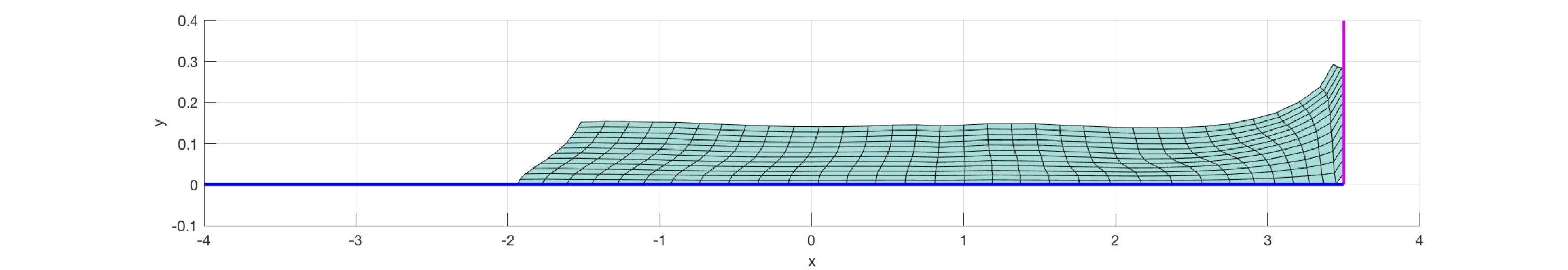} \vspace{-3pt} 
    \\
     
\hspace{1.6cm} \includegraphics[width=4.1 in]{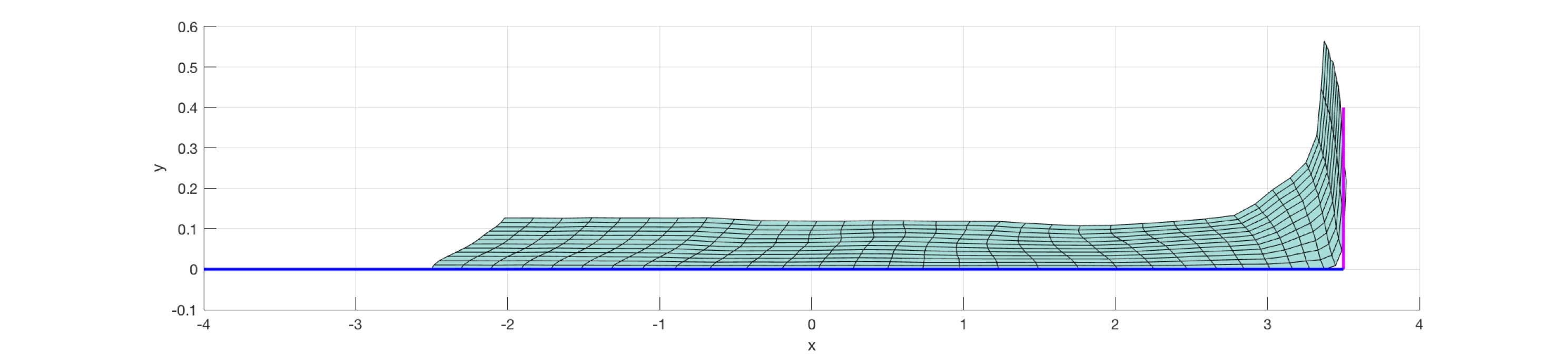} \vspace{-3pt}
     \\
 \hspace{1.6cm}  \includegraphics[width=4.1 in]{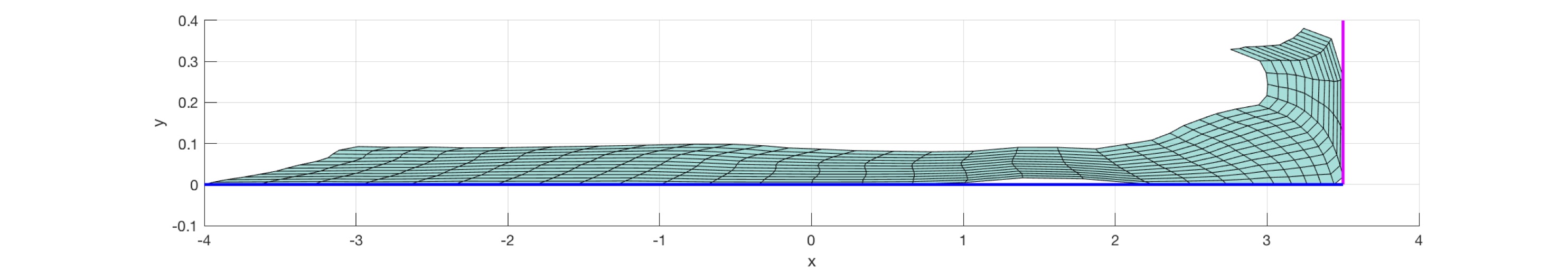} \vspace{-3pt}
    \\
     \includegraphics[width=5 in]{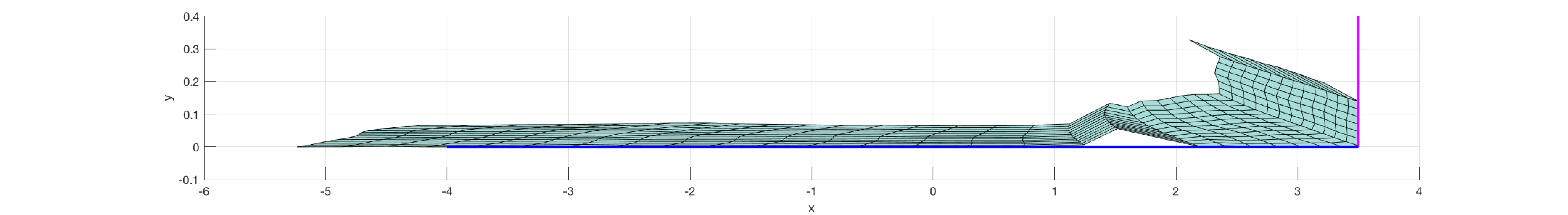} \vspace{-3pt}
   \caption{\footnotesize Barotropic fluid with contact. \textit{Top to bottom}: after $0.01$s, $1$s, $1.2$s, $1.6$s, $2$s.} \label{water_contact} 
 \end{figure}

\begin{figure}[H] \centering 
\hspace{1.6cm} \includegraphics[width=4.1 in]{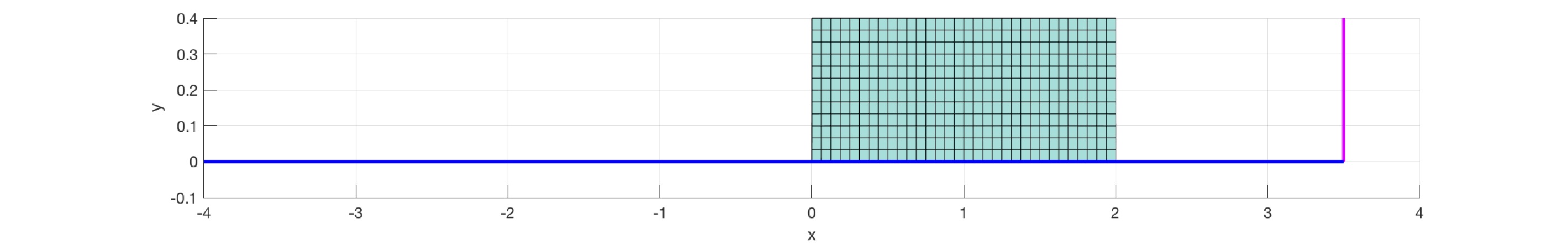} \vspace{-3pt}  
    \\
 \hspace{1.6cm} \includegraphics[width=4.1 in]{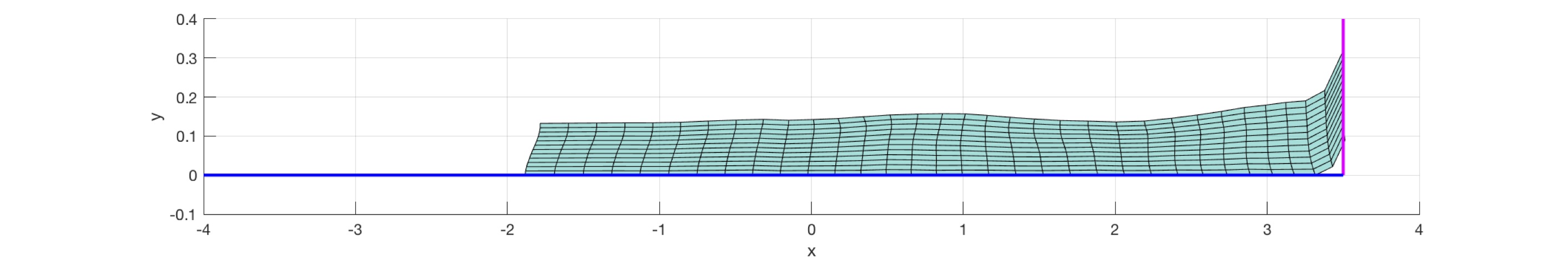} \vspace{-3pt} 
    \\    
\hspace{1.6cm} \includegraphics[width=4.1 in]{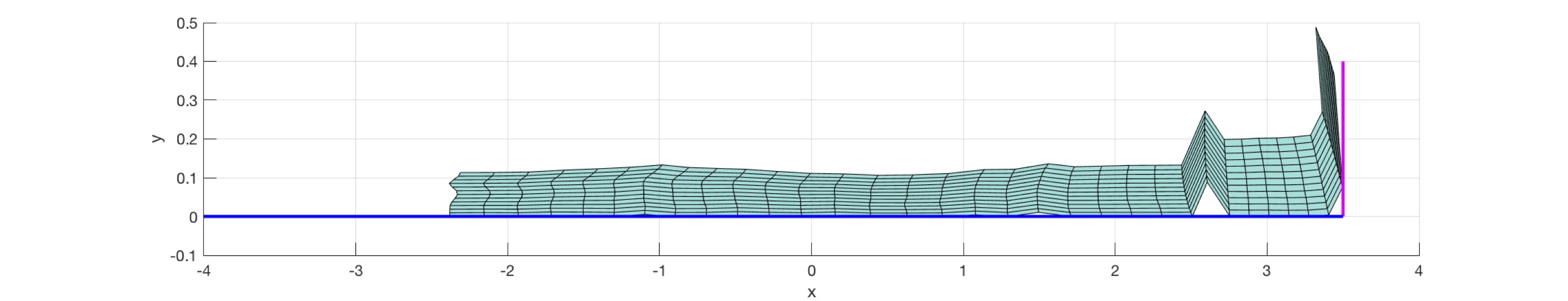} \vspace{-3pt}
     \\
\hspace{1.6cm}  \includegraphics[width=4.1 in]{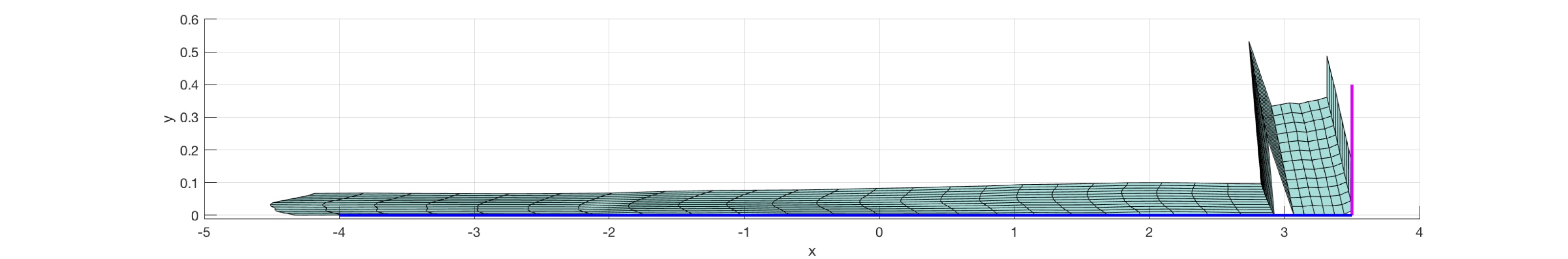} \vspace{-3pt}
    \\
     \includegraphics[width=5 in]{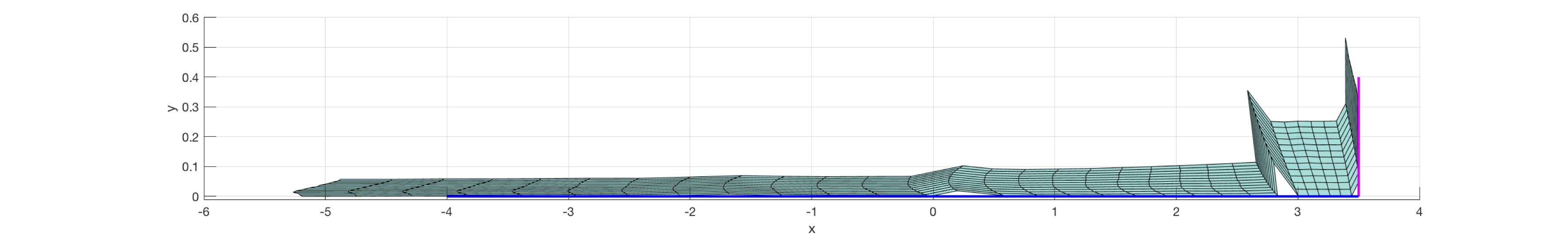} \vspace{-3pt}
   \caption{\footnotesize Incompressible ideal fluid with contact. \textit{Top to bottom}: after $0.01$s, $1$s, $1.2$s, $1.8$s, $2$s.}\label{water_contact_incompressible} 
 \end{figure}

The initial motion of the fluid is only due to the gravity. There are no other perturbations so that there is no expansion or compression imposed in the initial conditions. The evolution in the barotropic and incompressible cases are illustrated in Fig.\,\ref{water_contact} and Fig.\,\ref{water_contact_incompressible},  with the incompressibility conditions imposed by the penalty term.
 
 \medskip


\begin{figure}[H] \centering 
\includegraphics[width=1.9 in]{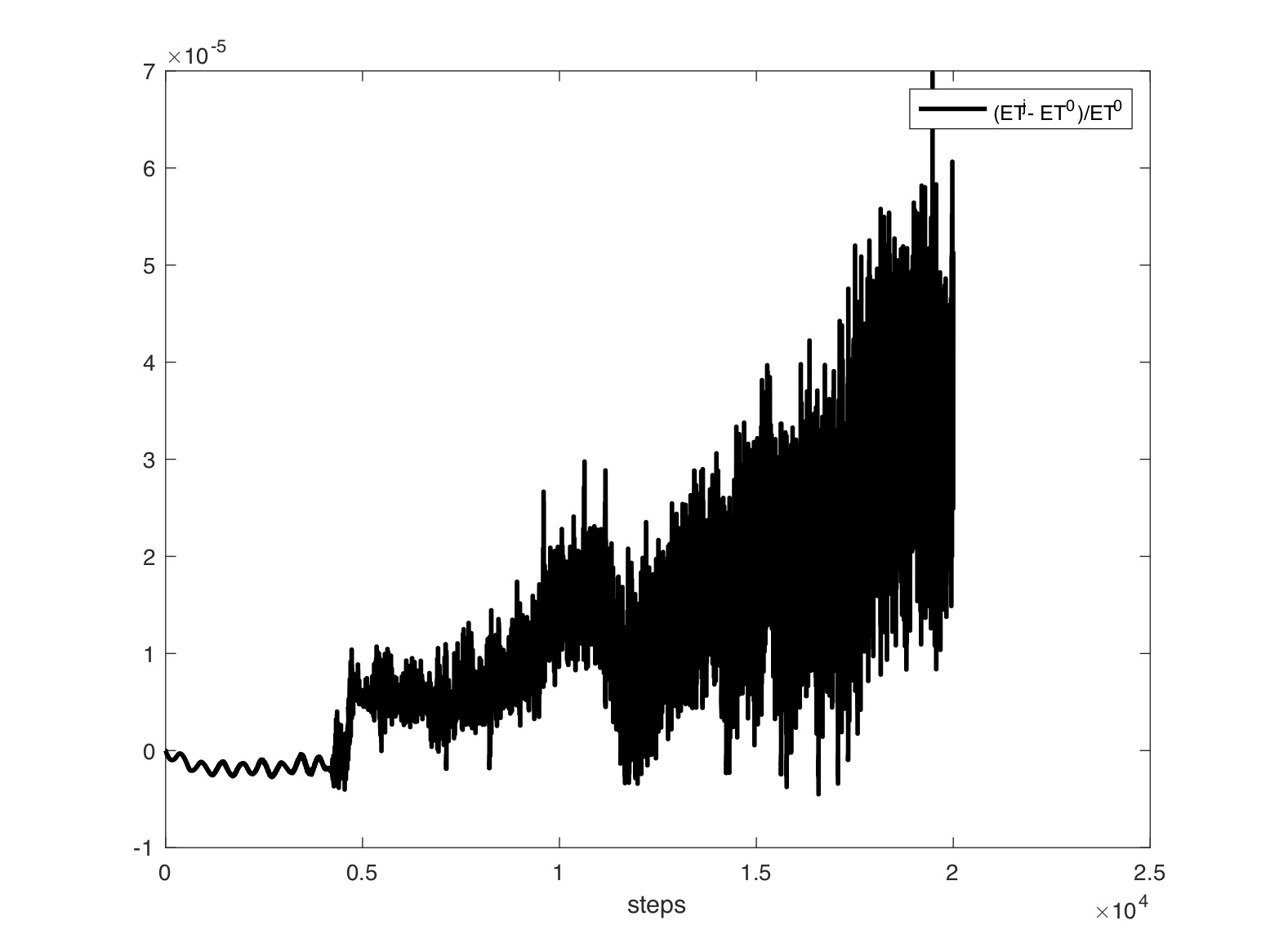} \vspace{-3pt} \qquad \includegraphics[width=1.9 in]{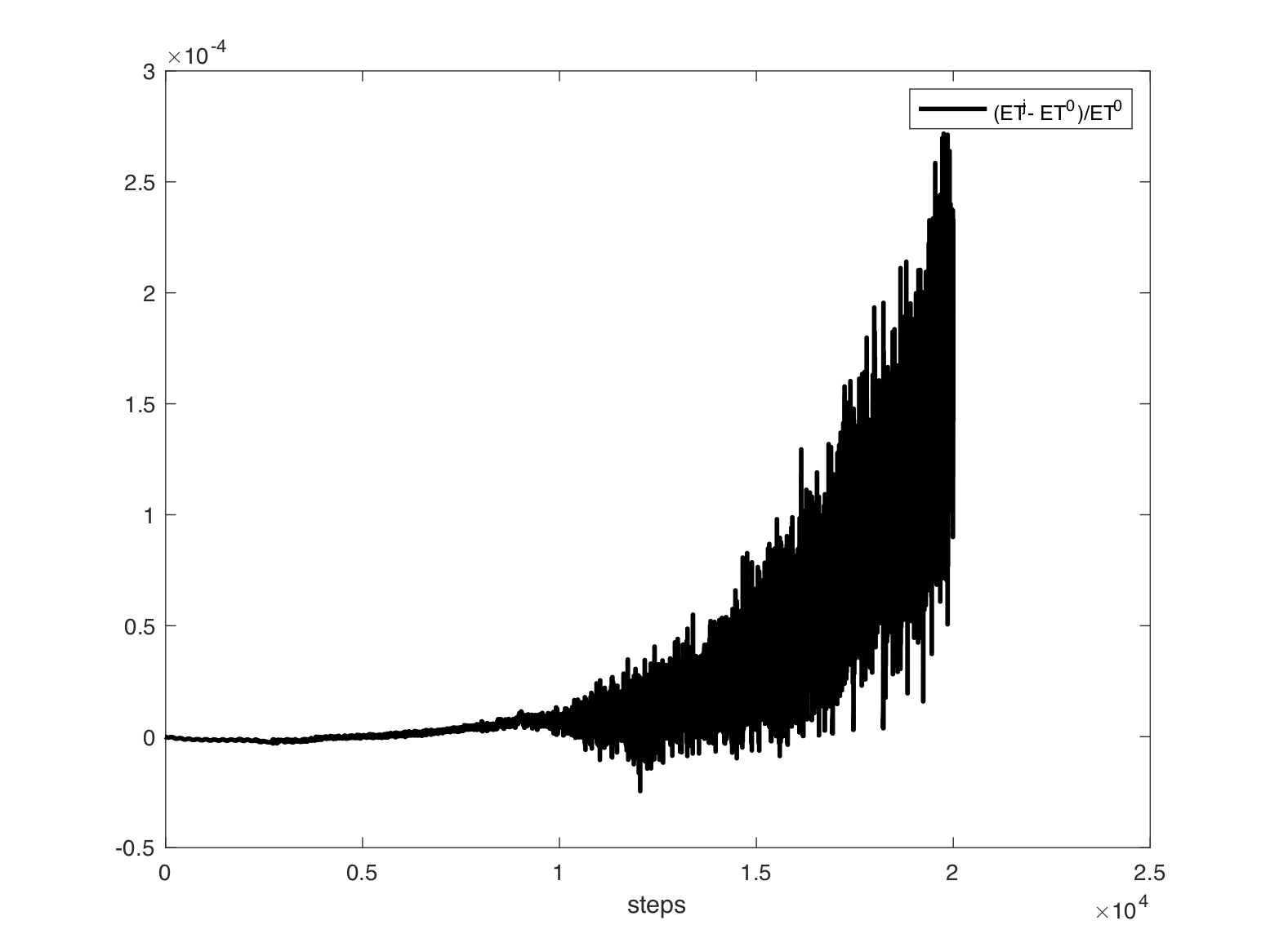} \vspace{-3pt}\\
\includegraphics[width=1.9 in]{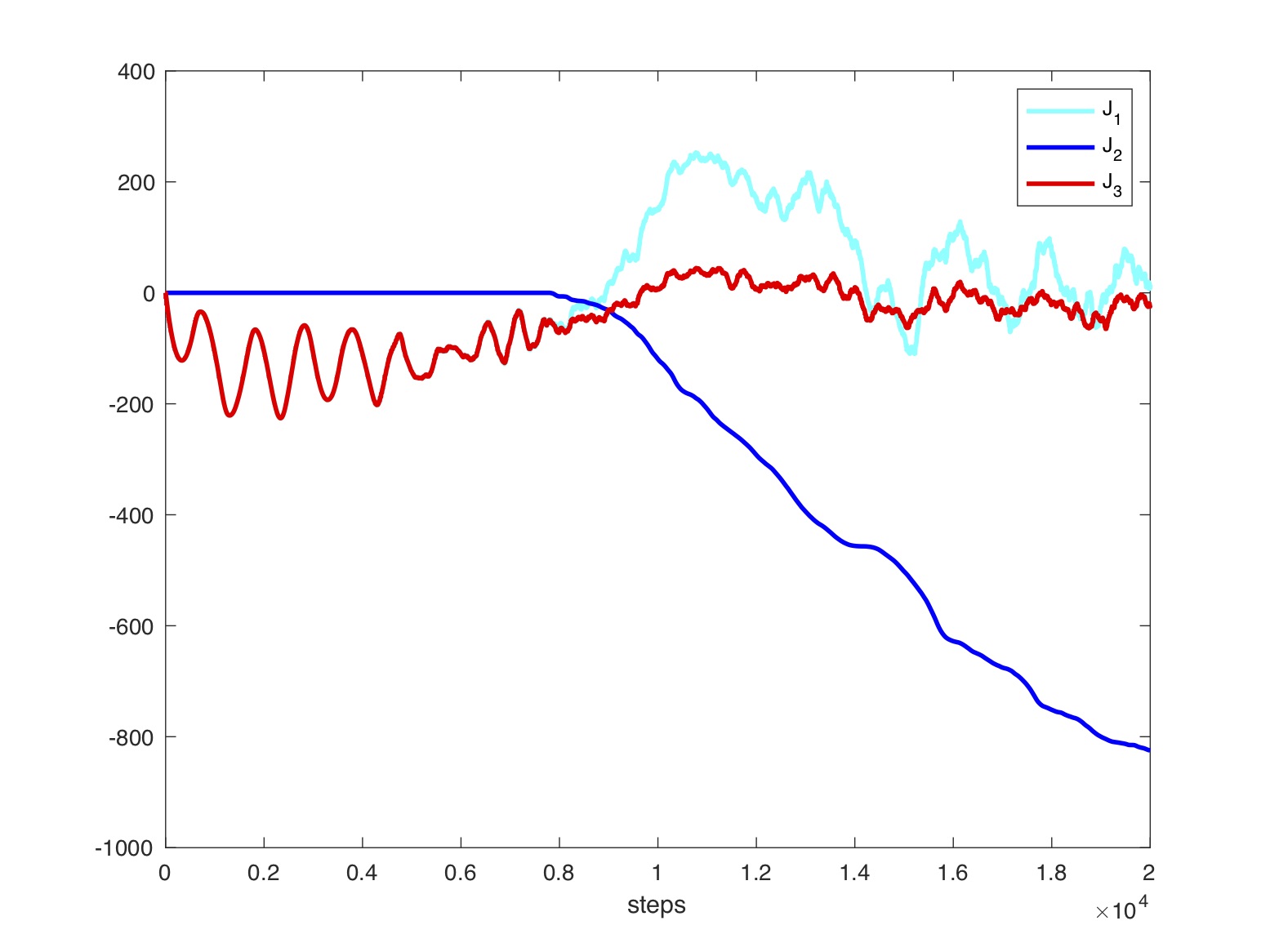} \vspace{-3pt} \qquad \includegraphics[width=1.9 in]{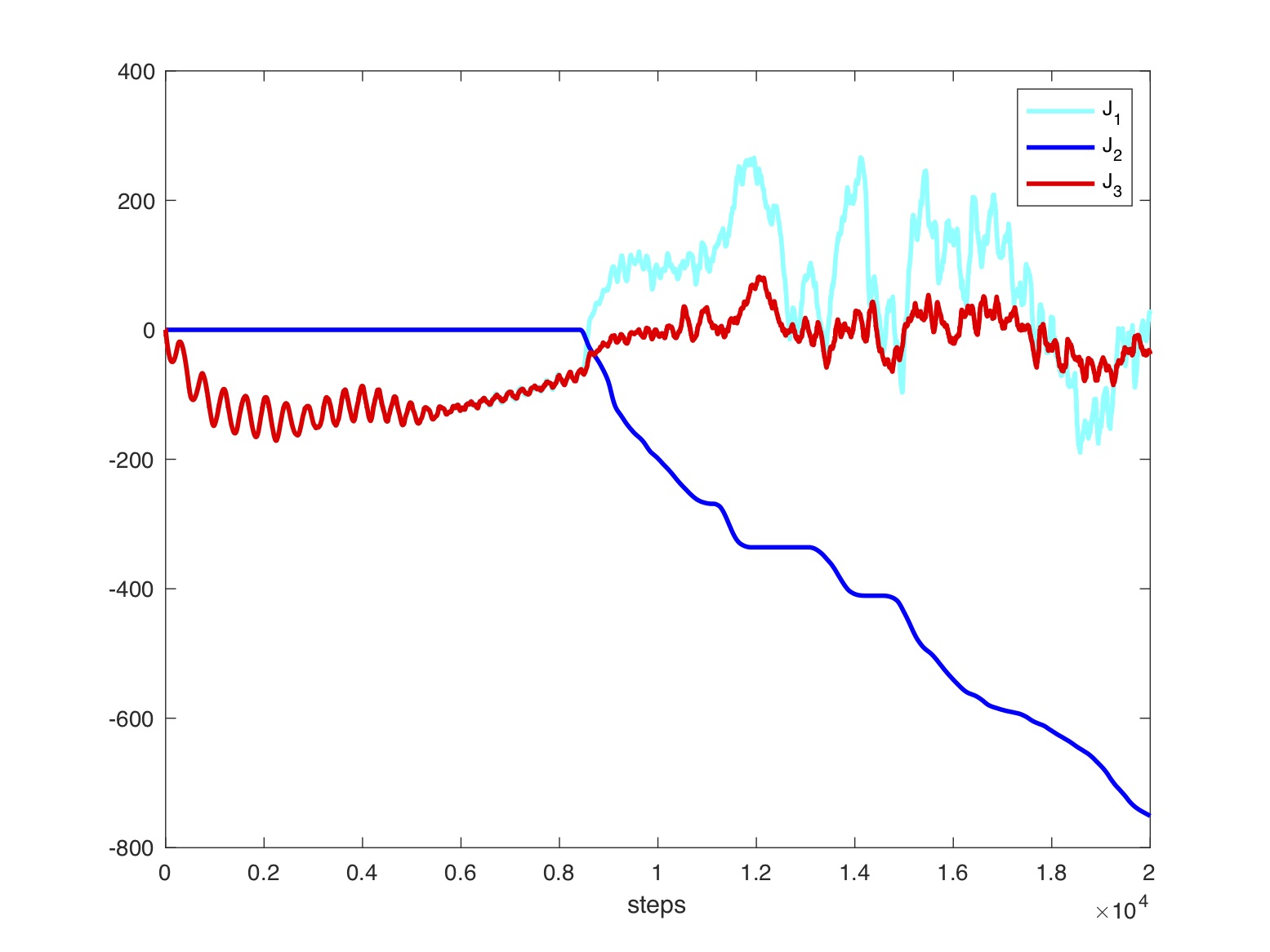} \vspace{-3pt}
\caption{\footnotesize   \textit{Left to right} Barotropic fluid and incompressible ideal fluid with contact. \textit{Top to bottom}: Relative energy and momentum map evolution during $2$s.} \label{water_contact_energy_J}
\end{figure}

The momentum map evolution is given in Fig.\,\ref{water_contact_energy_J}, where we note that only the component of the momentum map associated with vertical translation is preserved before the impact because of the presence of the gravity term. The energy perturbation increases after the contact, while complex phenomena similar to those encountered in breaking waves appear, like plunging waves, giving rise to a turbulent motion.

\subsubsection{Convergence tests} \label{2D_convergence}

Consider a barotropic fluid model with properties $\rho_0= 997 \, \mathrm{kg/m}^2$, $\gamma =6 $, $A = \tilde{A} \rho_0^{-\gamma}$ with $\tilde A = 3.041\times 10^4$ Pa, and $B = 3.0397\times 10^4$ Pa. The size of the discrete reference configuration at time $t^0$ is $0.4 \, \mathrm{m} \times 0.4 \, \mathrm{m}$. We consider the \textit{implicit} integrator to study the convergence with respect to $\Delta t$ and $\Delta s_1=\Delta s_2$.

\paragraph{Barotropic fluid motion in vacuum with free boundaries.} Given a fixed mesh, with values $\Delta s_1=\Delta s_2= 0.057$m, we impose an initial speed $V_{a,0}^0= (0,\; 0.163 \times a)^T$ on the boundary $b=0$ for all $a\in╩\{0,...,A\}$, and vary the time-steps as $\Delta t \in \{ 5 \times 10^{-3}, \, 2.5 \times 10^{-3}, \, 1.25 \times 10^{-3}, \, 6.25 \times 10^{-4} \}$.
We compute the $L^2$-errors in the position $\varphi_d$ at time $t^N=0.25$s, by comparing $\varphi_d$ with an ``exact solution'' obtained with the time-step $\Delta t_{\rm ref}=3.125 \times 10^{-4}s$. That is, for each value of $\Delta t$ we calculate
\begin{equation}\label{L_2norm}
\| \varphi_d - \varphi_{\rm ref} \|_{L^2} = \left( \sum_a \sum_b \| \varphi_{a,b}^N - \varphi_{{\rm ref};a,b}^N \|^2 \right)^{1/2}.
\end{equation}
This yields the following convergence with respect to $\Delta t$ 
\begin{figure}[H] \centering 
\begin{tabular}{| c | c | c | c | c |}
\hline
$\Delta t$ & $5 \times 10^{-3}$ & $ 2.5 \times 10^{-3}$ & $1.25 \times 10^{-3}$ & $6.25 \times 10^{-4}$ \\
\hline
$\| \varphi_d - \varphi_{\rm ref} \|_{L^2}$ & $1.5\times 10^{-2}$  & $ 7 \times 10^{-3}$  &  $3.6\times 10^{-3}$ &  $1\times 10^{-3}$ \\
\hline
$ \text{rate} $  &    & 1.106 &  0.964 &  1.815 \\
\hline
\end{tabular}
\end{figure}

Given a fixed time-step $\Delta t=2\times 10^{-3}$ we impose an initial speed\footnote{Note that, we need to take care of the initial sum of momentum $\sum_a m_a V_{a,0}^0$ which must be of the same value regardless of the number of nodes in the mesh.} $V_{a,0}^0= (0,\; 0.163 \times a)^T$, on the boundary $b=0$ for all $a\in╩\{0,...,A\}$, and vary the space-steps as $\Delta s_1=\Delta s_2$ $\in \{0.4, \, 0.2,\, 0.1,\, 0.05 \}$. The ``exact solution'' is chosen with $\Delta s_{1;\rm ref}=\Delta s_{2;\rm ref} = 0.025$m. We compute the $L^2$-errors in the position $\varphi_d$ at time $t^N=0.1$s. We get the following convergence with respect to $\Delta s_1=\Delta s_2$
\begin{figure}[H] \centering 
\begin{tabular}{| c | c | c | c | c | c |c|} 
\hline 
$\Delta s_1=\Delta s_2 $ &  $0.4$  &$0.2$  & $0.1$ & $0.05$   \\
\hline
$\| \varphi_d - \varphi_{\rm ref} \|_{L^2}$ & $3\times 10^{-2}$  & $2.16\times 10^{-2}$   &  $1.53\times 10^{-2} $ &  $7.1\times 10^{-3} $   \\
\hline
$ \text{rate} $  &    & 0.475 &  0.493 &  1.12 \\
\hline 
\end{tabular}
\end{figure}

\begin{figure}[H] \centering 
\includegraphics[width=1.45 in]{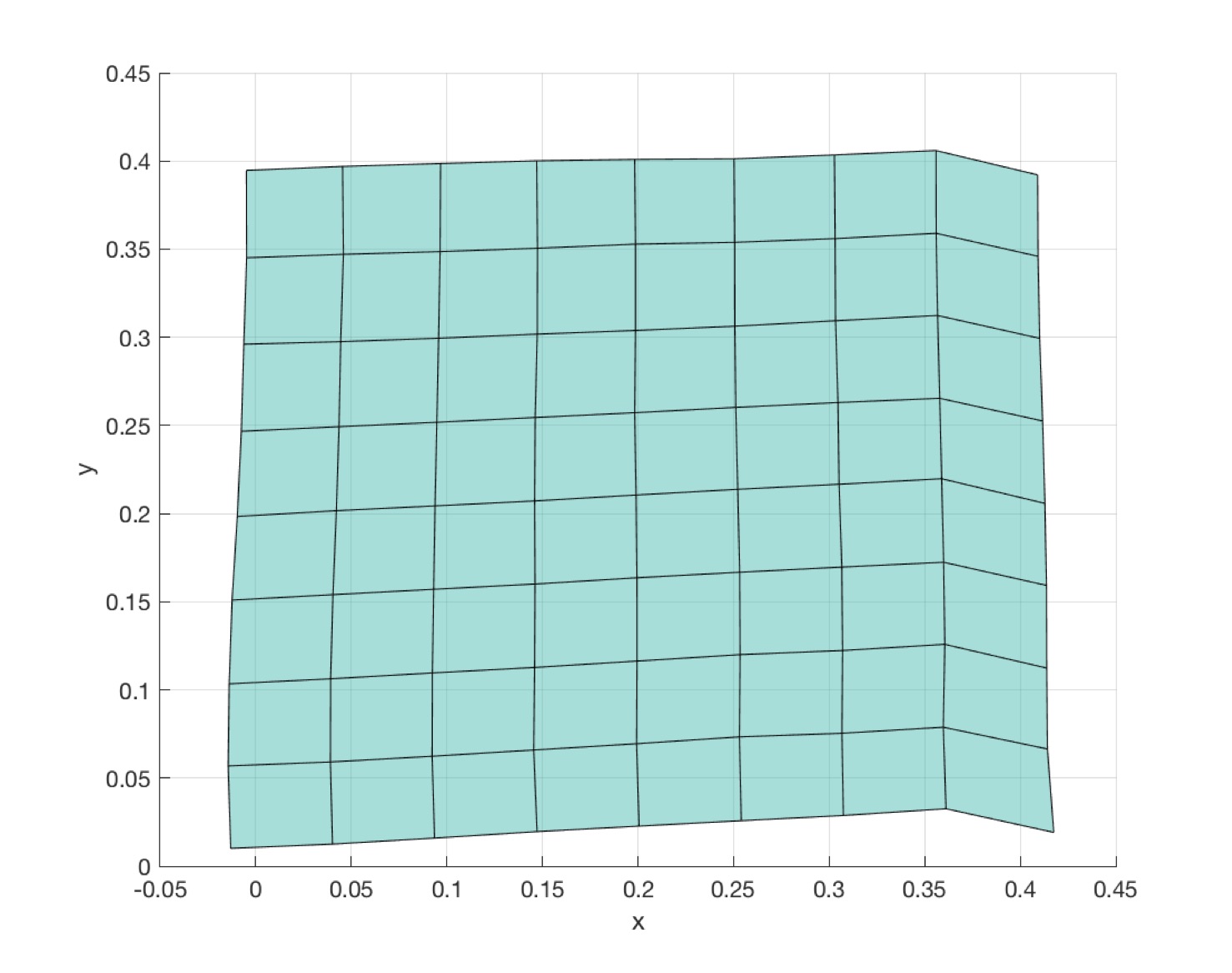} \vspace{-1pt}  \qquad
 \includegraphics[width=1.35 in]{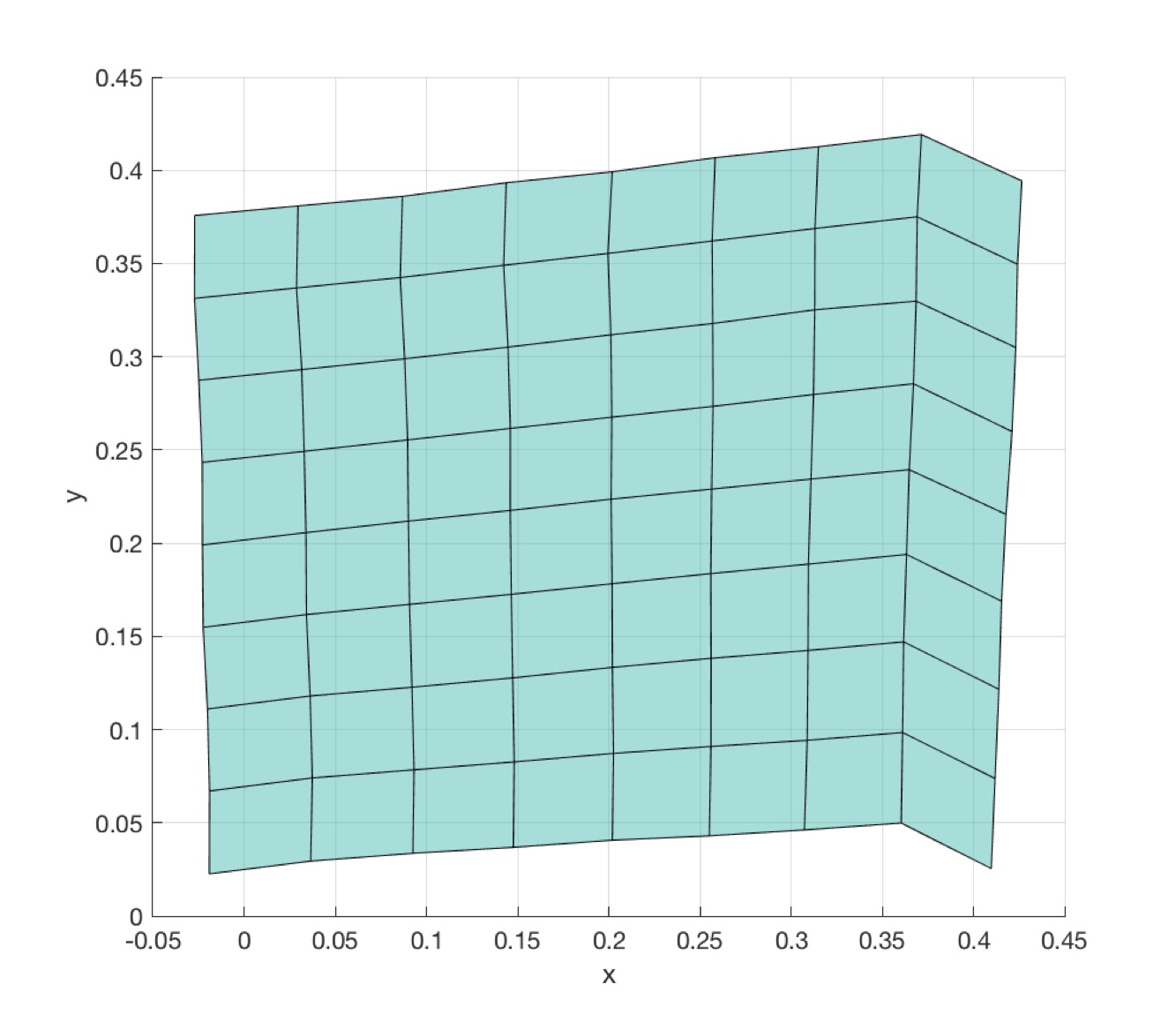}  \vspace{-1pt}  \qquad
 \includegraphics[width=1.38 in]{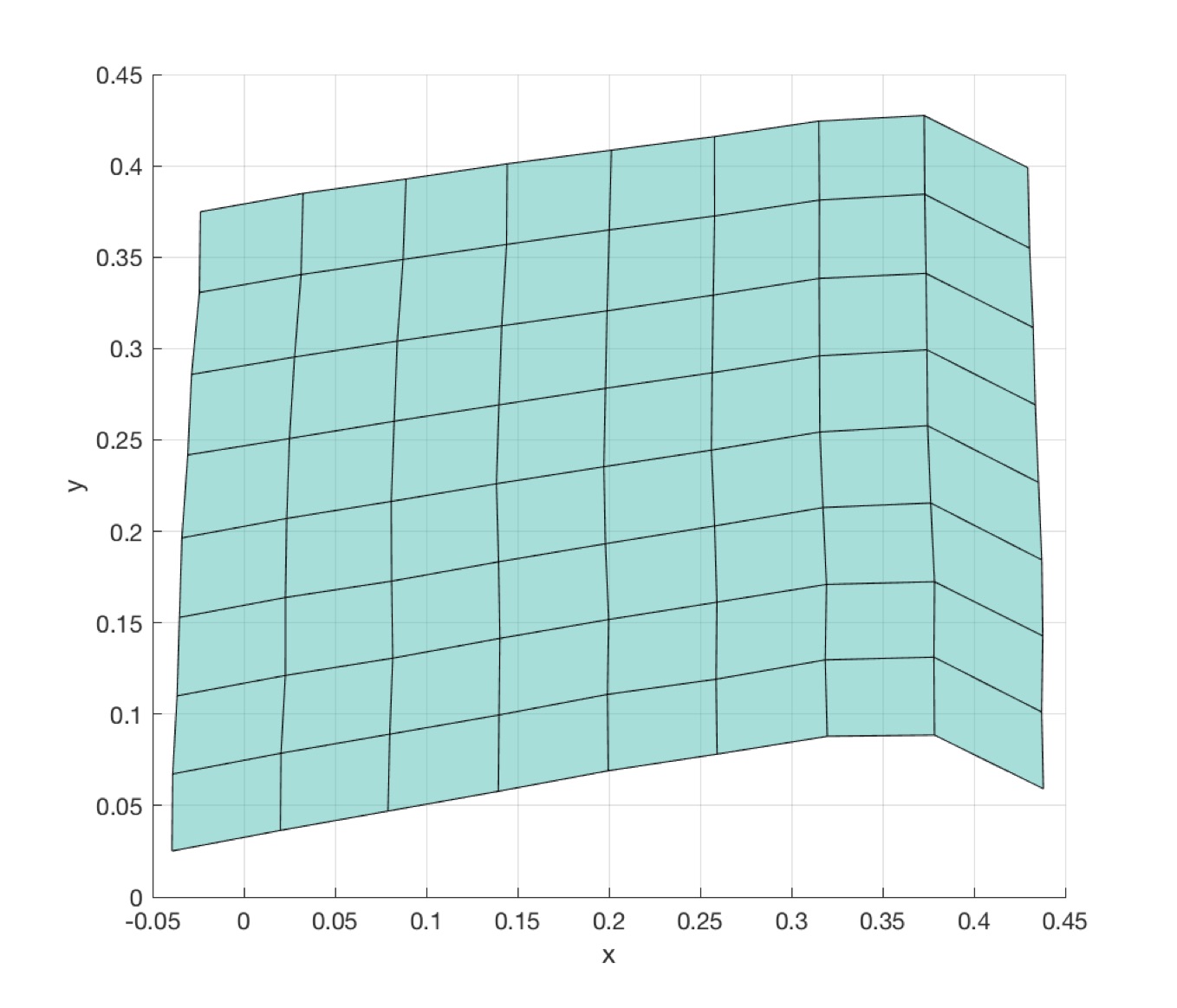} \vspace{-1pt} 
\\
  \includegraphics[width=1.88 in]{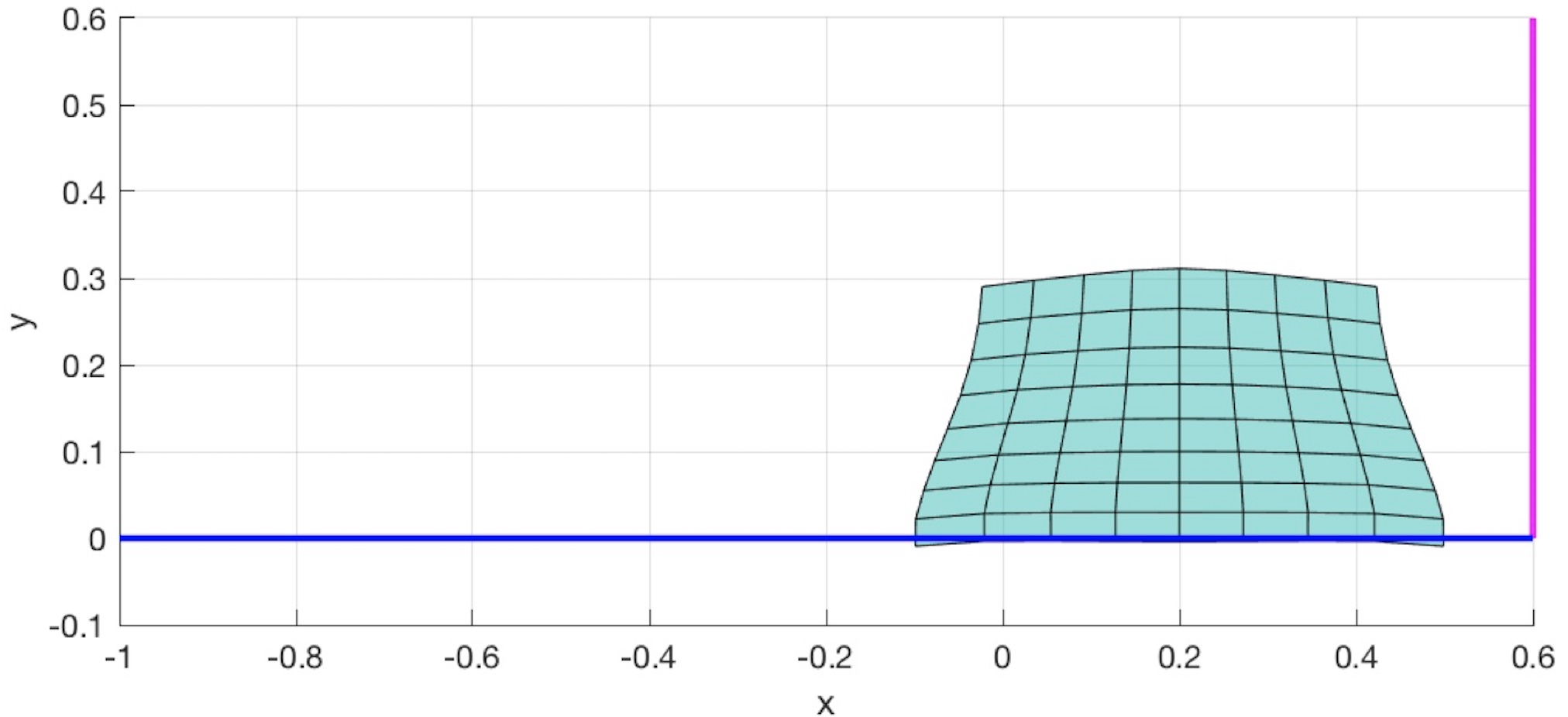} \vspace{-3pt} 
   \includegraphics[width=1.88 in]{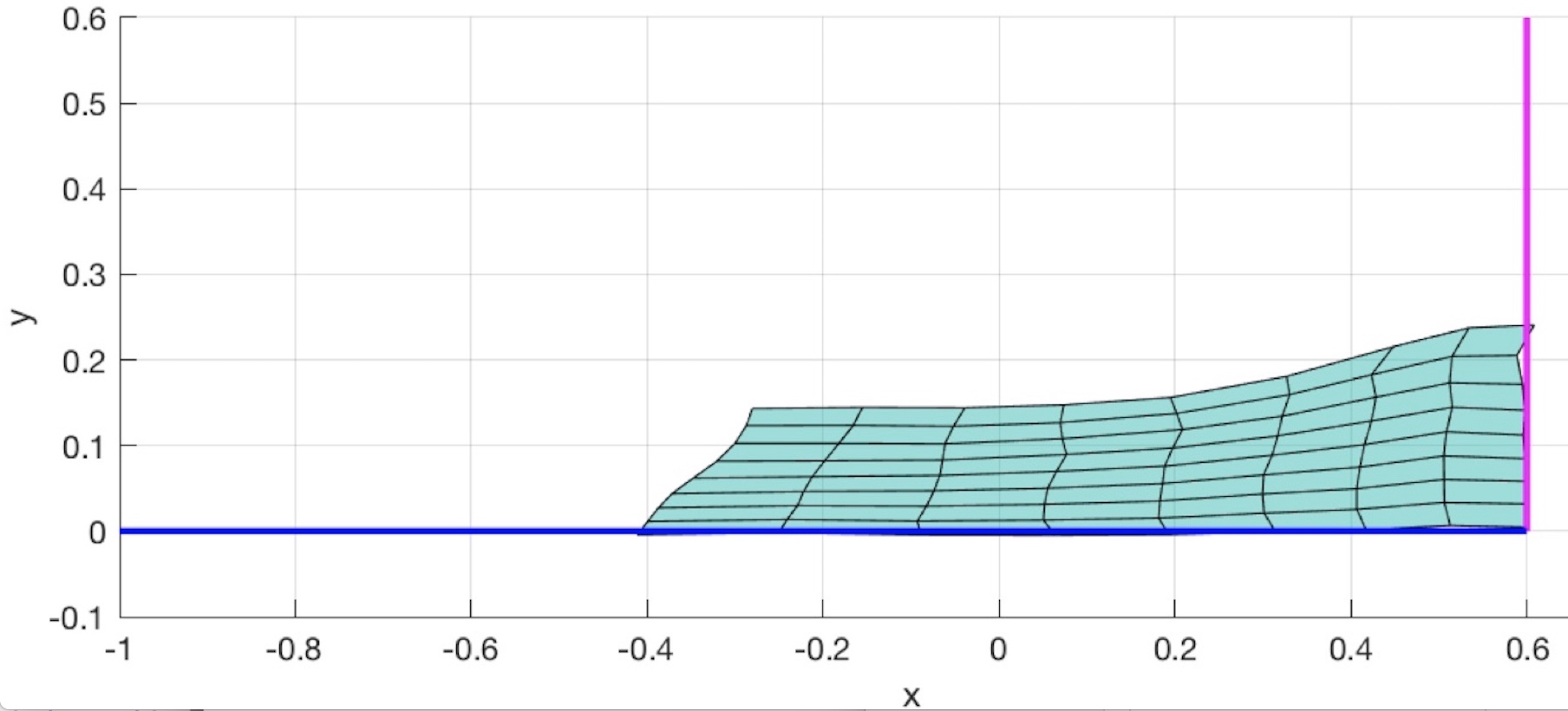} \vspace{-3pt} 
   \includegraphics[width=1.88 in]{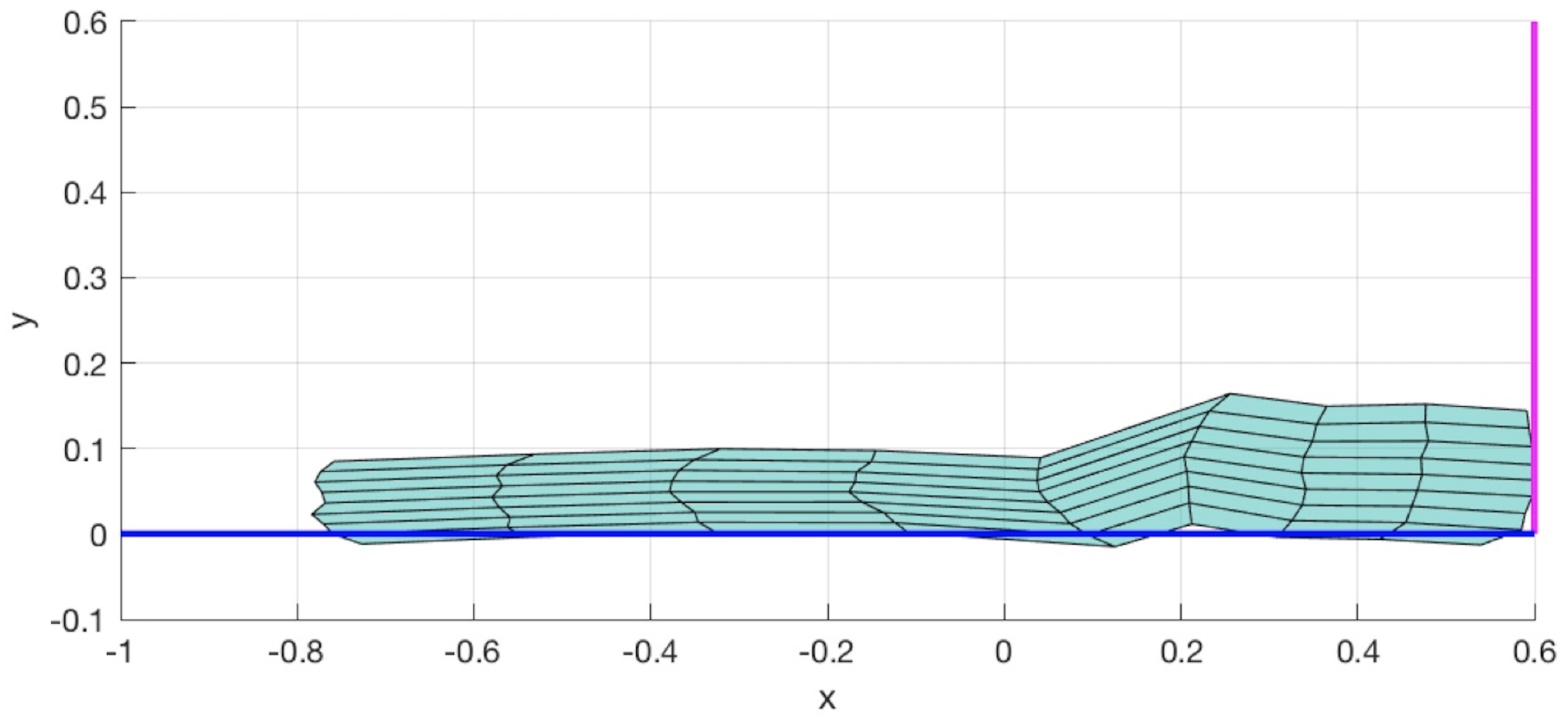} \vspace{-3pt} 
 \caption{\footnotesize Barotropic fluid. \textit{From top to bottom}: motion in vacuum with free boundaries ($\Delta s_1=\Delta s_2 =0.05$, $\Delta t=2 \times 10^{-3}$), and impact against an obstacle of a fluid flowing on a surface ($\Delta s_1=\Delta s_2 =0.05$, $\Delta t=2 \times 10^{-3}$). \textit{From left to right}: after $0.15$s, $0.3$s, and $0.5$s.} \label{convergence_time_contact} 
\end{figure}

\paragraph{Impact against an obstacle of a fluid flowing on a surface.} The values of the impenetrability penalty parameters are $K_1=K_2= 3\times10^7$. 
Given a fixed mesh, with values $\Delta s_1=\Delta s_2= 0.057$m, we repeat the experiment, described in Fig.\,\ref{convergence_time_contact}, with varying values of $\Delta t \in \{5 \times 10^{-3}, \, 2.5 \times 10^{-3}, \, 1.25 \times 10^{-3}, \, 6.25 \times 10^{-4} \}$.
Then, we compute the $L^2$-errors in the position $\varphi_d$ at time $t^N=0.25$s, by comparing $\varphi_d$ with an ``exact solution'' obtained with the time-step $\Delta t_{\rm ref}=3.125 \times 10^{-4}$.  We get the following convergence with respect to $\Delta t$
 
\begin{figure}[H] \centering 
\begin{tabular}{| c | c | c | c | c |}
\hline
$\Delta t$ & $5 \times 10^{-3}$ & $ 2.5 \times 10^{-3}$ & $1.25 \times 10^{-3}$ & $6.25 \times 10^{-4}$ \\
\hline
$\| \varphi_d - \varphi_{\rm ref} \|_{L^2}$ & $7.5\times 10^{-2}$  & $ 3.5 \times 10^{-2}$  &  $1.5\times 10^{-2}$ &  $5\times 10^{-3}$ \\
\hline
$ \text{rate} $  &    & 1.101 &  1.226 &  1.586 \\
\hline
\end{tabular}
\end{figure}

Similarly, given a fixed time-step $\Delta t=2 \times 10^{-3}$ we repeat the same experiment with varying values of $\Delta s_1=\Delta s_2$ $\in \{0.4, \, 0.2,\, 0.1,\, 0.05 \}$. The ``exact solution'' is chosen with $\Delta s_{1;\rm ref}=\Delta s_{2;\rm ref} = 0.025$m. We compute the $L^2$-errors in the position $\varphi_d$ at time $t^N=0.1$s. Therefore we get the following convergence with respect to $\Delta s_1 = \Delta s_2$ 

\begin{figure}[H] \centering 
\begin{tabular}{| c | c | c | c | c | c |c|} 
\hline 
$\Delta s_1=\Delta s_2 $ &  $0.4$  &$0.2$  & $0.1$ & $0.05$   \\
\hline
$\| \varphi_d - \varphi_{\rm ref} \|_{L^2}$ & $6.38\times 10^{-2}$  & $3.84\times 10^{-2}$   &  $1.75\times 10^{-2} $ &  $6.8\times 10^{-3} $   \\
\hline
$ \text{rate} $  &    & 0.733 &  1.132 &  1.366 \\
\hline 
\end{tabular}
\end{figure}

An illustration of the test used for the numerical convergence is given in Fig.\,\ref{convergence_time_contact}.

\begin{figure}[H] \centering 
  \includegraphics[width=1.7 in]{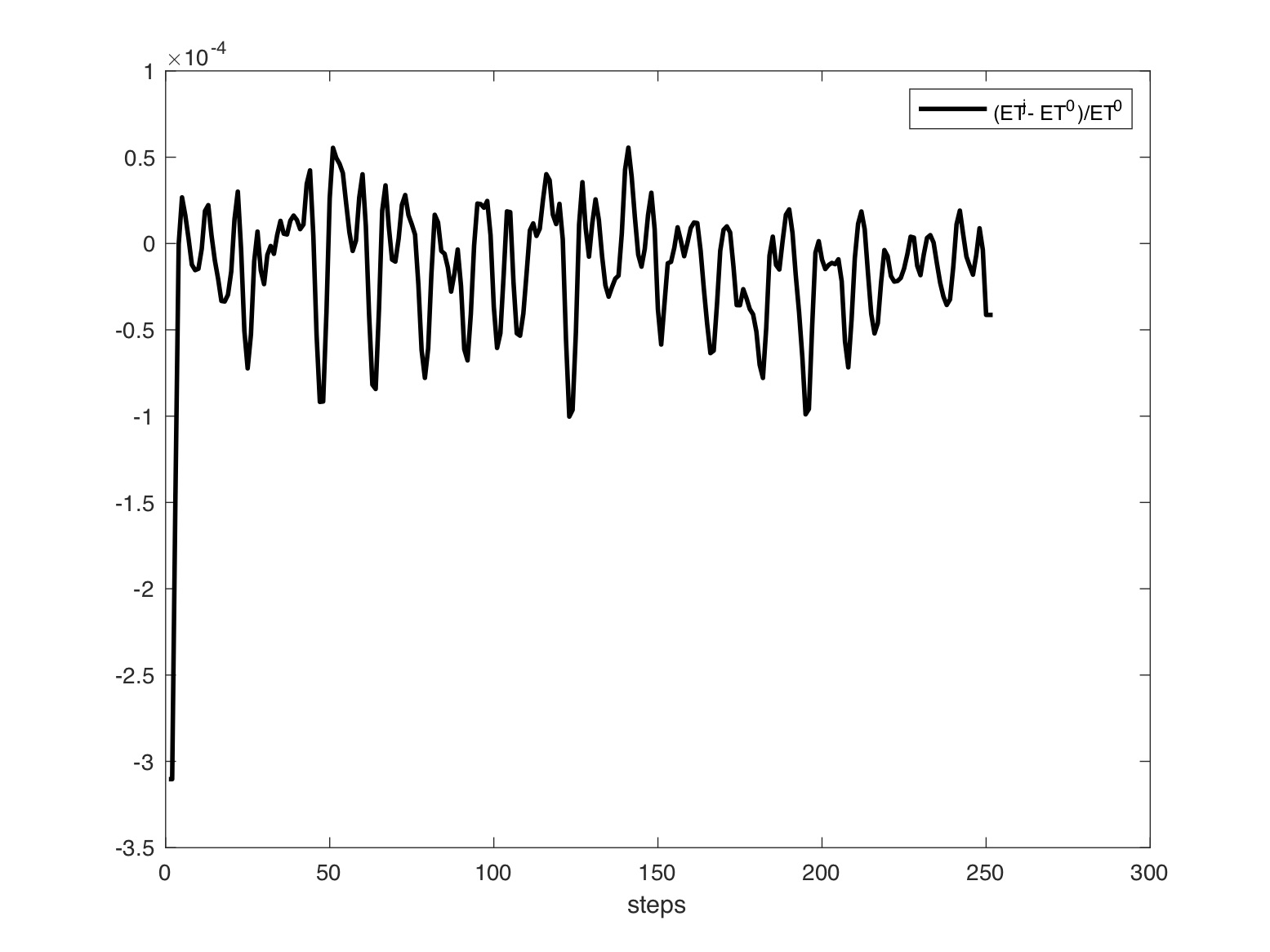} \vspace{-3pt} \qquad
    \includegraphics[width=1.7 in]{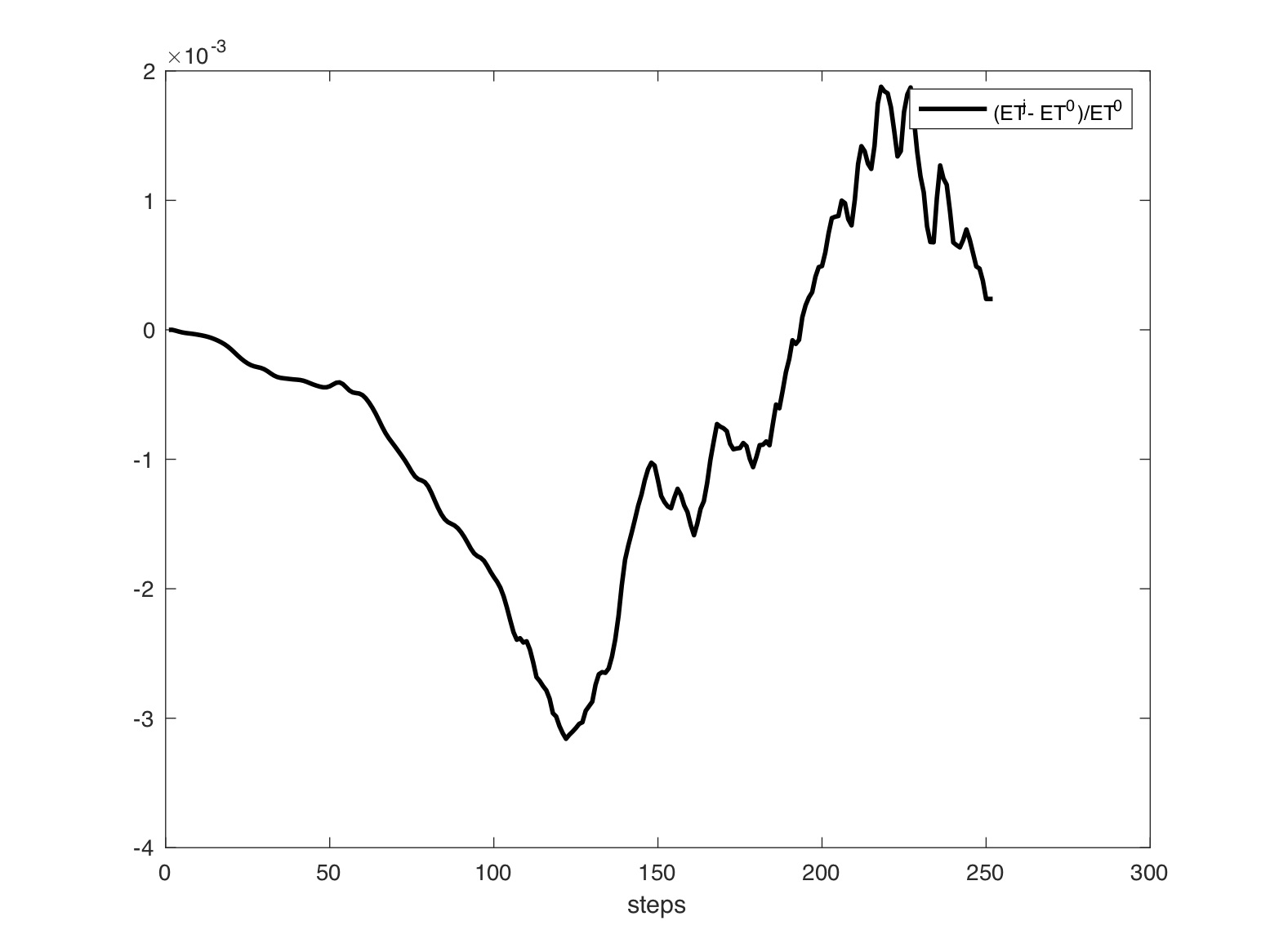} \vspace{-3pt} \qquad
 \caption{\footnotesize Relative error in the energy $(ET^j-ET^0)/ET^0$, with $ET^j$ the total energy at time $t^j$. \textit{From left to right}: motion in vacuum with free boundaries ($\Delta s_1=\Delta s_2 =0.05$, $\Delta t=2 \times 10^{-3}$), and impact against an obstacle of a fluid flowing on a surface ($\Delta s_1=\Delta s_2 =0.05$, $\Delta t=2 \times 10^{-3}$). } \label{relative_error_energy} 
\end{figure}

\section{3D discrete barotropic and incompressible fluid models} \label{3D_wave_prop_fluid}

In this section we indicate how the developments made in \S\ref{2D_wave_prop_fluid} extend to the 3D case. The general discrete multisymplectic framework (discrete configuration bundle, discrete first jet, discrete multisymplectic form, etc...) have been already explained in a general setting in \S\ref{2D_wave_prop_fluid}.
We assume that $ \mathcal{B} $ is a parallelepiped in $ \mathbb{R} ^3 $ and take $ \mathcal{M}= \mathbb{R} ^3$.
  
\subsection{Multisymplectic discretizations}

\subsubsection{Discrete configuration bundle}

\begin{figure}[H] \centering 
  \includegraphics[width=4.9 in]{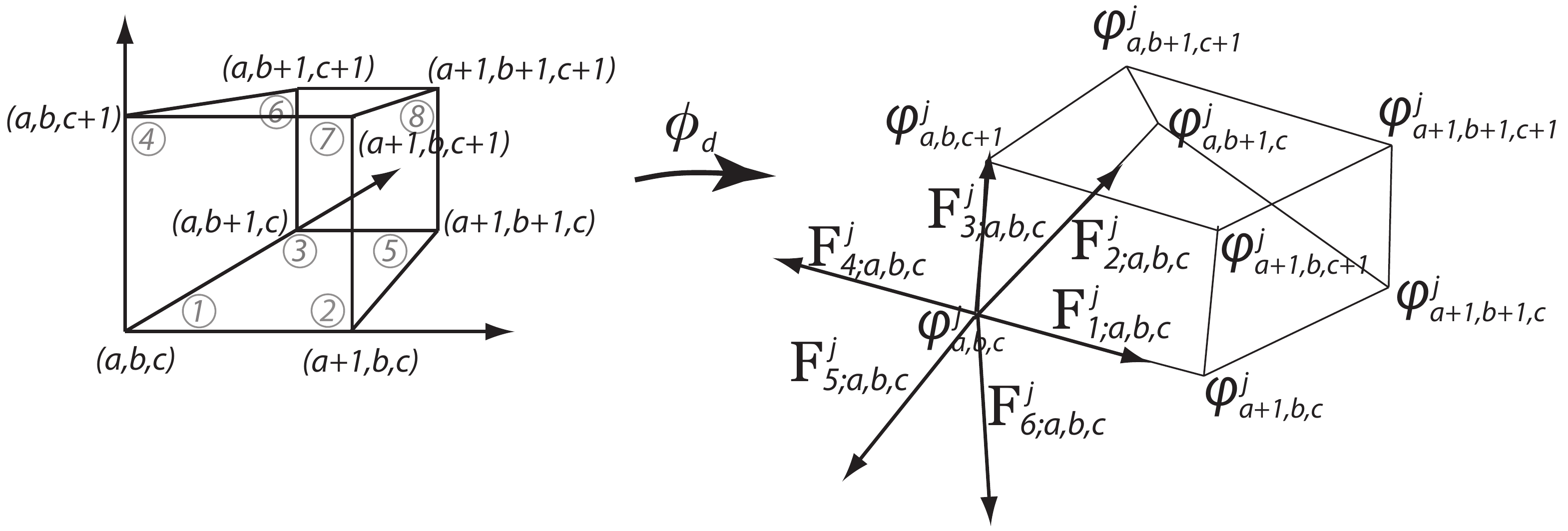} \vspace{-3pt}  
 \caption{\footnotesize Discrete field $\phi_d=\varphi_d\circ  \phi _{ \mathcal{X} _d}$ evaluated on $\mbox{\mancube}_{a,b,c}^j$ at time $t^j$. }\label{phi_evalutation_3D} 
 \end{figure}

The discrete parameter space $ \mathcal{U} _d$ is now decomposed in a set of elements $\,\mbox{\mancube}_{a,b,c}^{\,j} $ defined by 16 pairs of indices, see Fig.\,\ref{phi_evalutation_3D} for the eight pairs of indices in $\,\mbox{\mancube}_{a,b,c}^j$ at time $t^j$. 


As in \eqref{choice_of_DBSC},  we consider discrete base-space configurations of the form
\begin{equation}\label{3D_disc_base_space_conf}
\phi_{ \mathcal{X} _d}: \mathcal{U} _d \ni (j,a,b,c) \mapsto s^j_{a,b,c}=(t^j,z_a^j,z_b^j,z_c^j) \in    \mathcal{X} _d\subset \mathbb{R} \times \mathcal{B}.
\end{equation}
The discrete field $\varphi_d$ evaluated at $s^j_{a,b,c}$ is denoted $\varphi_{a,b,c}^j:= \varphi_d(s^j_{a,b,c}) \in  \mathbb{R} ^3 $.

\subsubsection{Discrete Jacobian}\label{DJ3D}

Given a discrete base space configuration $ \phi _{ \mathcal{X} _d}$ of the form \eqref{3D_disc_base_space_conf} and a discrete field $ \varphi _d$, we define the following six vectors $\mathbf{F}_{\ell;a,b,c}^j \in \mathbb{R} ^3 $, $\ell=1,...,6$ at each node $(j,a,b,c) \in  \mathcal{U} _d$, see Fig.\,\ref{phi_evalutation_3D} on the right:
\begin{align*} 
\mathbf{F}_{1;a,b,c}^j&= \frac{\varphi_{a+1,b,c}^j - \varphi_{a,b,c}^j}{ | s _{a+1,b,c} - s_{a,b,c}| } , &   \mathbf{F}_{2;a,b,c}^j &= \frac{\varphi_{a,b+1,c}^j- \varphi_{a,b,c}^j}{|s _{a,b+1,c} - s _{a,b,c}| }, \\
\mathbf{F}_{3;a,b,c}^j &= \frac{\varphi_{a,b,c+1}^j- \varphi_{a,b,c}^j}{| s_{a,b,c+1} -s _{a,b,c}| }, &  
\mathbf{F}_{4;a,b,c}^j &=  \frac{\varphi_{a-1,b,c}^j - \varphi_{a,b,c}^j}{| s _{a,b,c} - s _{a-1,b,c}| }= - \mathbf{F}^j_{1,a-1,b,c}\\
\mathbf{F}_{5;a,b,c}^j &= \frac{\varphi_{a,b-1,c}^j- \varphi_{a,b,c}^j}{ | s _{a,b,c} - s _{a,b-1,c}|} = - \mathbf{F}
^j_{2,a,b-1,c}, &  \mathbf{F}_{6;a,b,c}^j&= \frac{\varphi_{a,b,c-1}^j - \varphi_{a,b,c}^j}{|s _{a,b,c} -s_{a,b,c-1}| }= - \mathbf{F}
 ^j_{3;a,b,c-1}.
\end{align*} 
Based on these definitions, the discrete gradient is constructed as follows.

\begin{definition}  \label{3D_gradient_definition}
The  \textit{discrete gradient deformations} of a discrete field $ \varphi _d$ at the element $\,\mbox{\mancube}_{a,b,c}^j$ are the $3 \times 3$ matrices $\mathbf{F}^\ell(\mbox{\mancube}_{a,b,c}^j)$, $\ell=1,...,8$, defined at the eight nodes at time $t^j$ of $\;\mbox{\mancube}_{a,b,c}^j$ as follows:
\begin{equation}\label{3D_gradient_def}
{\small
\begin{aligned}
& \mathbf{F}_1( \mbox{\mancube}_{a,b,c}^j) = \left[\mathbf{F}_{1;a,b,c}^j \; \; \mathbf{F}_{2;a,b,c}^j \; \; \mathbf{F}_{3;a,b,c}^j \right],\\
&\mathbf{F}_2( \mbox{\mancube}_{a,b,c}^j) = \left[\mathbf{F}_{2;a+1,b,c}^j \;\; \mathbf{F}_{4;a+1,b,c}^j  \; \;  \mathbf{F}_{3;a+1,b,c}^j \right] ,\\
&  \mathbf{F}_3( \mbox{\mancube}_{a,b,c}^j)= \left[\mathbf{F}_{5;a,b+1,c}^j  \; \; \mathbf{F}_{1;a,b+1,c}^j \; \; \mathbf{F}_{3;a,b+1,c}^j \right], \\
&\mathbf{F}_4( \mbox{\mancube}_{a,b,c}^j)= \left[ \mathbf{F}_{2;a,b,c+1}^j  \; \; \mathbf{F}_{1;a,b,c+1}^j \; \; \mathbf{F}_{6;a,b,c+1}^j \right], \\
&  \mathbf{F}_5( \mbox{\mancube}_{a,b,c}^j) = \left[\mathbf{F}_{4;a+1,b+1,c}^j  \; \; \mathbf{F}_{5;a+1,b+1,c}^j  \; \; \mathbf{F}_{3;a+1,b+1,c}^j \right],\\
& \mathbf{F}_6( \mbox{\mancube}_{a,b,c}^j) = \left[\mathbf{F}_{1;a,b+1,c+1}^j  \; \; \mathbf{F}_{5;a,b+1,c+1}^j \; \; \mathbf{F}_{6;a,b+1,c+1}^j \right],\\
& \mathbf{F}_7( \mbox{\mancube}_{a,b,c}^j)= \left[\mathbf{F}_{4;a+1,b,c+1}^j \; \; \mathbf{F}_{2;a+1,b,c+1}^j  \; \; \mathbf{F}_{6;a+1,b,c+1}^j \right],\\
& \mathbf{F}_8( \mbox{\mancube}_{a,b,c}^j) = \left[\mathbf{F}_{5;a+1,b+1,c+1}^j  \; \; \mathbf{F}_{4;a+1,b+1,c+1}^j \; \; \mathbf{F}_{6;a+1,b,c+1}^j \right].
\end{aligned} }
\end{equation}
The ordering $\ell=1$ to $\ell =8$ is respectively associated to the nodes $(j,a,b,c)$, $(j,a+1,b,c)$, $(j,a,b+1,c)$, $(j,a,b,c+1)$, $(j,a+1,b+1,c)$, $(j,a,b+1,c+1)$, $(j,a+1,b,c+1)$, $(j,a+1,b+1,c+1)$, see Fig.\,\ref{phi_evalutation_3D} on the left.
\end{definition}

Then we define the Jacobian in each node, as follows

\begin{definition}
The discrete Jacobians of a discrete field $ \varphi _d $ at the element $\,\mbox{\mancube}_{a,b,c}^j$ are the numbers $J_\ell(\mbox{\mancube}_{a,b}^j)$, $\ell=1,...,8$, defined at the eight nodes at time $t^j$ of $\;\mbox{\mancube}_{a,b,c}^j$ as follows:
\begin{equation} \label{3D_Jacobian}
J_1 (\mbox{\mancube}_{a,b,c}^j)=  (\mathbf{F}_{1;a,b,c}^j \times \mathbf{F}_{2;a,b,c}^j)\cdot \mathbf{F}_{3;a,b,c}^j =\mathrm{det} \big( \mathbf{F}_1(\mbox{\mancube}_{a,b,c}^j)\big) .
\end{equation}
\end{definition}

See in \S\ref{3D_Jacobian's} for the others Jacobian on $\mbox{\mancube}_{a,b,c}^{\,j}$. We can now establish the link between the discrete Jacobian and the discrete gradient deformation.

\medskip
In terms of the discrete field $ \varphi _d $, the discrete Jacobians are 
\[
J_1 ( \mbox{\mancube}_{a,b,c}^j)=  \frac{ (( \varphi_{a+1,b,c}^j - \varphi_{a,b,c}^j)\times (\varphi_{a,b+1,c}^j - \varphi_{a,b,c}^j) ) \cdot (\varphi_{a,b,c+1}^j - \varphi_{a,b,c}^j)  }{ | s_{a+1,b,c} - s_{a,b,c}|  | s_{a,b+1,c} - s_{a,b,c}| | s_{a,b,c+1} - s_{a,b,c}| }
\]
\[
J_2(\mbox{\mancube}_{a,b,c}^j)= \frac{ (( \varphi_{a+1,b+1,c}^j - \varphi_{a+1,b,c}^j)\times (\varphi_{a,b,c}^j - \varphi_{a+1,b,c}^j) ) \cdot (\varphi_{a+1,b,c+1}^j - \varphi_{a+1,b,c}^j)  }{ | s_{a+1,b+1,c} - s_{a+1,b,c}|  | s_{a,b,c} - s_{a+1,b,c}| | s_{a+1,b,c+1} - s_{a+1,b,c}| }
\]
similarly for the other ones.

\subsubsection{Discrete Lagrangian}

The discrete Lagrangian for 3D barotropic fluid models has the same general form as \eqref{Discrete_Lagrangian_2D_fluid_metric}, with the obvious 3D extension of formulas \eqref{kinetic_energy}--\eqref{grav_potential}.

\subsubsection{Discrete variations and discrete Euler-Lagrange equations}

The discrete action functional takes the form
\begin{equation}\label{3D_discrete_action_sum}
S_d(\varphi_d) = \sum_{j=0}^{N-1} \sum_{a=0}^{A-1} \sum_{b=0}^{B-1} \sum_{b=0}^{C-1} \mathcal{L}\big( j^1 \varphi _d (\mbox{\mancube}_{a,b,c}^{\,j}) \big)
\end{equation}
and yield the discrete Euler-Lagrange equations
\begin{equation} \label{3D_DCEL}
\begin{aligned}
&M v_{a,b,c}^{j} + A_{a,b,c}^j + B_{a-1,b,c}^j + C_{a,b-1,c}^j + D_{a,b,c-1}^j 
\\
& \qquad  + E_{a-1,b-1,c}^j + F_{a,b-1,c-1}^j + G_{a-1,b,c-1}^j + H_{a-1,b-1,c-1}^j =0,
\end{aligned}
\end{equation} 
where we have used notations analogous to \eqref{PD_2} and \eqref{PD_1}  for the partial derivative of $ \mathcal{L} _d$. We refer to Appendix \ref{3D_water} for the expressions of $A^j_{a,b,c},..., H^j_{a,b,c}$.
Boundary conditions are deduced the discrete Hamilton principle in a similar way as it was done in \eqref{2D_boundary_cond_space} and \eqref{2D_boundary_cond_time} for the 2D case.

\subsubsection{Discrete multisymplectic form formula and discrete Noether theorem}

Following the general definition \eqref{DCF}, the discrete Cartan forms evaluated at the first jet extension $j^1 \varphi _d( \mbox{\mancube}_{a,b,c}^j)$ of a discrete field $ \varphi _d $ are
\begin{equation}\label{3DCartanForm_extended}
\footnotesize
 \begin{aligned}
\Theta_{\mathcal{L}_d }^{1}& = A_{a,b,c}^j \, {\rm d} \varphi_{a,b,c}^j, &\quad &  \Theta_{\mathcal{L}_d }^{2} =  \frac{M}{8}  v_{a,b,c}^{j} \, {\rm d} \varphi_{a,b,c}^{j+1},
\\
\Theta_{\mathcal{L}_d }^{3} & = B_{a,b,c}^j \, {\rm d} \varphi_{a+1,b,c}^j,&\quad   &  \Theta_{\mathcal{L}_d }^{4}  =  \frac{M}{8}  v_{a+1,b,c}^{j} \, 
{\rm d} \varphi_{a+1,b,c}^{j+1},
\\
\Theta_{\mathcal{L}_d }^{5} & = C_{a,b,c}^j \, {\rm d} \varphi_{a,b+1,c}^j, &\quad  &  \Theta_{\mathcal{L}_d }^{6}  = \frac{M}{8}  v_{a,b+1,c}^{j+1} \, 
{\rm d} \varphi_{a,b+1,c}^{j+1},
\\
\Theta_{\mathcal{L}_d }^{7}  & = D_{a,b,c}^j \, {\rm d} \varphi_{a,b,c+1}^j, &\quad  &  \Theta_{\mathcal{L}_d }^{8}  = \frac{M}{8}  v_{a,b,c+1}^{j} \, 
{\rm d} \varphi_{a,b,c+1}^{j+1},
\\ 
\Theta_{\mathcal{L}_d }^{9} & = E_{a,b,c}^j \, {\rm d} \varphi_{a+1,b+1,c}^j, &\quad &  \Theta_{\mathcal{L}_d }^{10} =  \frac{M}{8}  v_{a+1,b+1,c}^{j} \, {\rm d} \varphi_{a+1,b+1,c}^{j+1},
\\
\Theta_{\mathcal{L}_d }^{11}  & = F_{a,b,c}^j \, {\rm d} \varphi_{a,b+1,c+1}^j,&\quad   &  \Theta_{\mathcal{L}_d }^{12} =  \frac{M}{8}  v_{a,b+1,c+1}^{j} \, 
{\rm d} \varphi_{a,b+1,c+1}^{j+1},
\\
\Theta_{\mathcal{L}_d }^{13}  & = G_{a,b,c}^j \, {\rm d} \varphi_{a+1,b,c+1}^j, &\quad  &  \Theta_{\mathcal{L}_d }^{14} = \frac{M}{8}  v_{a+1,b,c+1}^{j+1} \, 
{\rm d} \varphi_{a+1,b,c+1}^{j+1},
\\
\Theta_{\mathcal{L}_d }^{15} & = H_{a,b,c}^j \, {\rm d} \varphi_{a+1,b+1,c+1}^j, &\quad  &  \Theta_{\mathcal{L}_d }^{16} = \frac{M}{8}  v_{a+1,b+1,c+1}^{j} \, 
{\rm d} \varphi_{a+1,b+1,c+1}^{j+1}.
\end{aligned}
\end{equation}

With these forms, the discrete multisymplectic form formula and conservation laws in the presence of a symmetry group (discrete Noether theorem) can be derived in a similar way as it was done in 2D in \S \ref{2D_mult_discret}.

\subsubsection{Symmetries for barotropic fluids}

Exactly as in \S\ref{Sym_barotropic}, the discrete Lagrangian is $SE(3)$ invariant and hence the discrete covariant Noether theorem holds with the covariant discrete momentum maps $J_{ \mathcal{L} _d}^\mathtt{p} : J^1 \mathcal{Y} _d \rightarrow \mathfrak{se}(3) ^*$, $\mathtt{p}=1,...,16$. From this, the discrete momentum map

\vspace{-0.3cm}
{\small\begin{equation}\label{total_momentum_map_3D} 
\begin{aligned} 
\mathbf{J} _d^j=\mathbf{J} _d ( \boldsymbol{\varphi }^j, \boldsymbol{\varphi }^{j+1}) &= \sum_{a=0}^{A-1}\sum_{b=0}^{B-1}  \sum_{c=0}^{C-1}\left( J_{ \mathcal{L} _d}^2 + J_{ \mathcal{L} _d}^4+J_{ \mathcal{L} _d}^6+J_{ \mathcal{L} _d}^8+J_{ \mathcal{L} _d}^{10}+ J_{ \mathcal{L} _d}^{12}+ J_{ \mathcal{L} _d}^{14}+J_{ \mathcal{L} _d}^{16} \right) \\
&= - \sum_{a=0}^{A-1}\sum_{b=0}^{B-1}\sum_{c=0}^{C-1} \left( J_{ \mathcal{L} _d}^1 + J_{ \mathcal{L} _d}^3+J_{ \mathcal{L} _d}^5+J_{ \mathcal{L} _d}^7+J_{ \mathcal{L} _d}^9+J_{ \mathcal{L} _d}^{11}+J_{ \mathcal{L} _d}^{13}+J_{ \mathcal{L} _d}^{15}\right),
\end{aligned}
\end{equation} }
is preserved, as explained in \S\ref{Sym_barotropic}. In 3D, the expression \eqref{3D_momentum_map_fluid} extends as
\[
\mathbf{J}_d^j =  \begin{bmatrix}
\vspace{0.2cm}\displaystyle\sum_{a=0}^{A-1} \sum_{b=0}^{B-1}\sum_{a=0}^{C-1} \mathbf{J}_r\big(j^1 \varphi _d (\mbox{\mancube}_{a,b}^j)\big)  
\\ 
\displaystyle\sum_{a=0}^{A-1} \sum_{b=0}^{B-1} \sum_{a=0}^{C-1} \mathbf{J}_l \big(j^1 \varphi _d (\mbox{\mancube}_{a,b}^j)\big)  
\end{bmatrix}
\] 
with
\begin{align*} 
\mathbf{J}_r\big(j^1 \varphi _d (\mbox{\mancube}_{a,b,c}^j)\big)& = \sum_{\alpha =a}^{a+1} \sum_{\beta =b}^{b+1} \sum_{\alpha =c}^{c+1} \varphi_{\alpha,\beta,\gamma}^j \times \left(\frac{M}{8} v_{\alpha,\beta,\gamma}^j \right) \in \mathbb{R} ^3\\ 
\mathbf{J}_l\big(j^1 \varphi _d (\mbox{\mancube}_{a,b,c}^j)\big)&= \sum_{\alpha =a}^{a+1} \sum_{\beta =b}^{b+1}\sum_{\alpha =c}^{c+1} \; \frac{M}{8} v_{\alpha,\beta,\gamma}^j \in \mathbb{R}  ^3.
\end{align*} 
 
\subsubsection{Incompressible ideal hydrodynamics}

As done in 2D (see \S \ref{2D_penalty_method}), associated to the equality constraint $J=1$ we consider a penalty function 
 \begin{equation}\label{3D_augmenting_fonc}
\Phi_{d0}\big(j^1 \varphi _d (\mbox{\mancube})\big) := \frac{1}{8}  \sum_{\ell=1}^8 \frac{r}{2} \big(J_\ell(\mbox{\mancube})- 1\big)^2,
\end{equation} 
where $r$ is the penalty parameter. 

\subsection{Numerical simulations}

In this section we illustrate the performance of an explicit-in-time integrator in 3D \footnote{Note that with the multisymplectic variational integrators we can equally move in time and in space, see \cite{DeGBRa2014}.}, as it was done in \S \ref{2D_examples}.

\subsubsection{Example 3: barotropic fluid motion in vacuum with free boundaries}

Consider a barotropic fluid with properties $\rho_0= 997 \, \mathrm{kg/m}^3$, $\gamma =6 $, $A = \tilde{A} \rho_0^{-\gamma}$ with $\tilde A = 3.041\times 10^4$ Pa, and $B = 3.0397\times 10^4$ Pa. The size of the mesh at time $t^0$ is $2 \, \mathrm{m} \times 2 \, \mathrm{m}\times 2 \, \mathrm{m}$, with $\Delta s_1= \Delta s_2= \Delta s_3= 0.333$ m. We consider both the compressible barotropic fluid and the incompressible case with penalty parameter $r= 10^5$ and $r= 10^7$. The time step is $\Delta t=10^{-3}$. There are no exterior forces.

\medskip

Initial perturbations are applied at time $t^1$ on nodes $(1,0,1)$ and $(1,0,2)$, in a similar way with the test made in dimension 2, see Fig.\,\ref{water_perturb}.

  \begin{figure}[H] \centering 
\includegraphics[width=1.8 in]{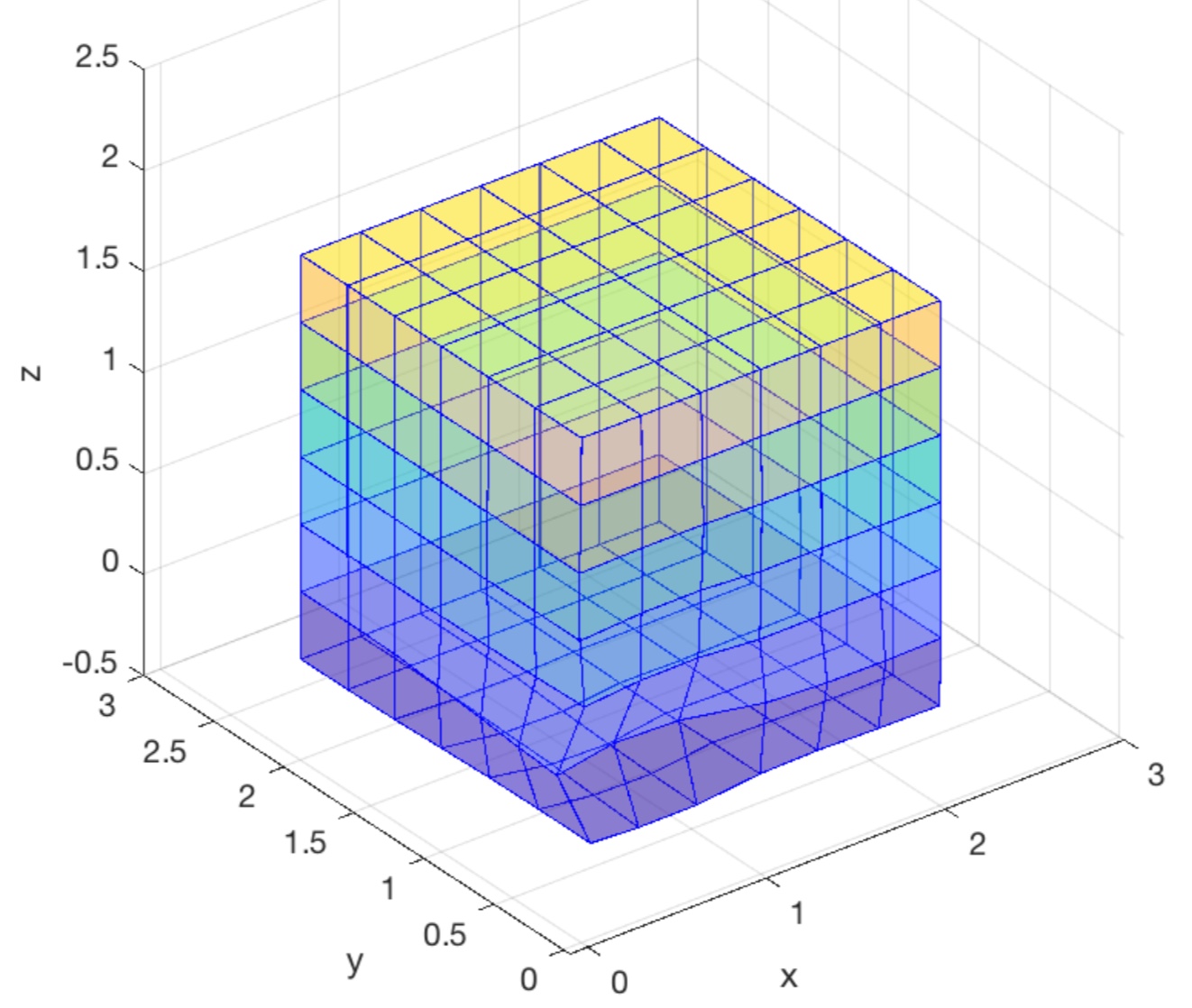} \vspace{-3pt} 
   \includegraphics[width=1.8 in]{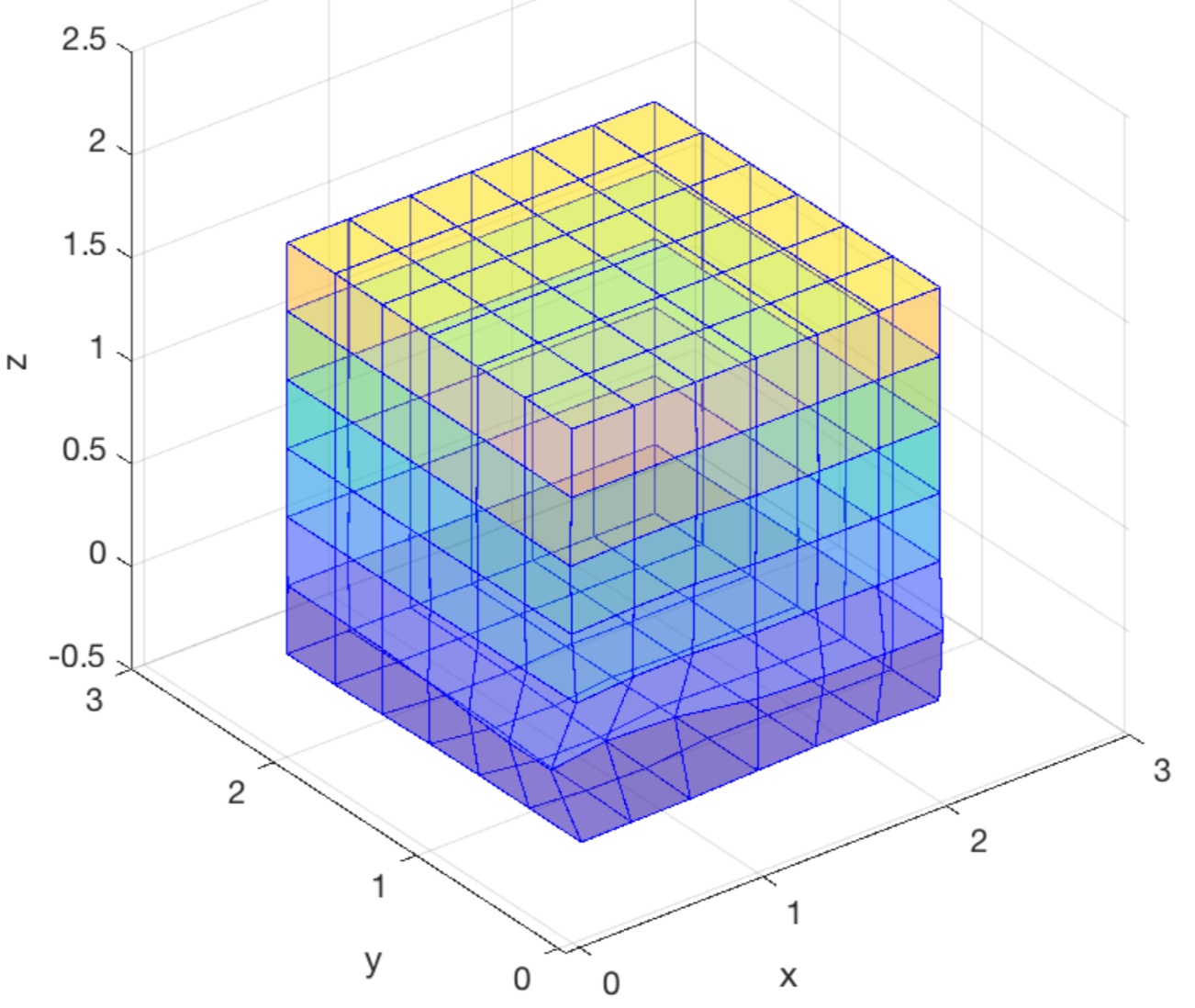} \vspace{-3pt} 
   \includegraphics[width=1.8 in]{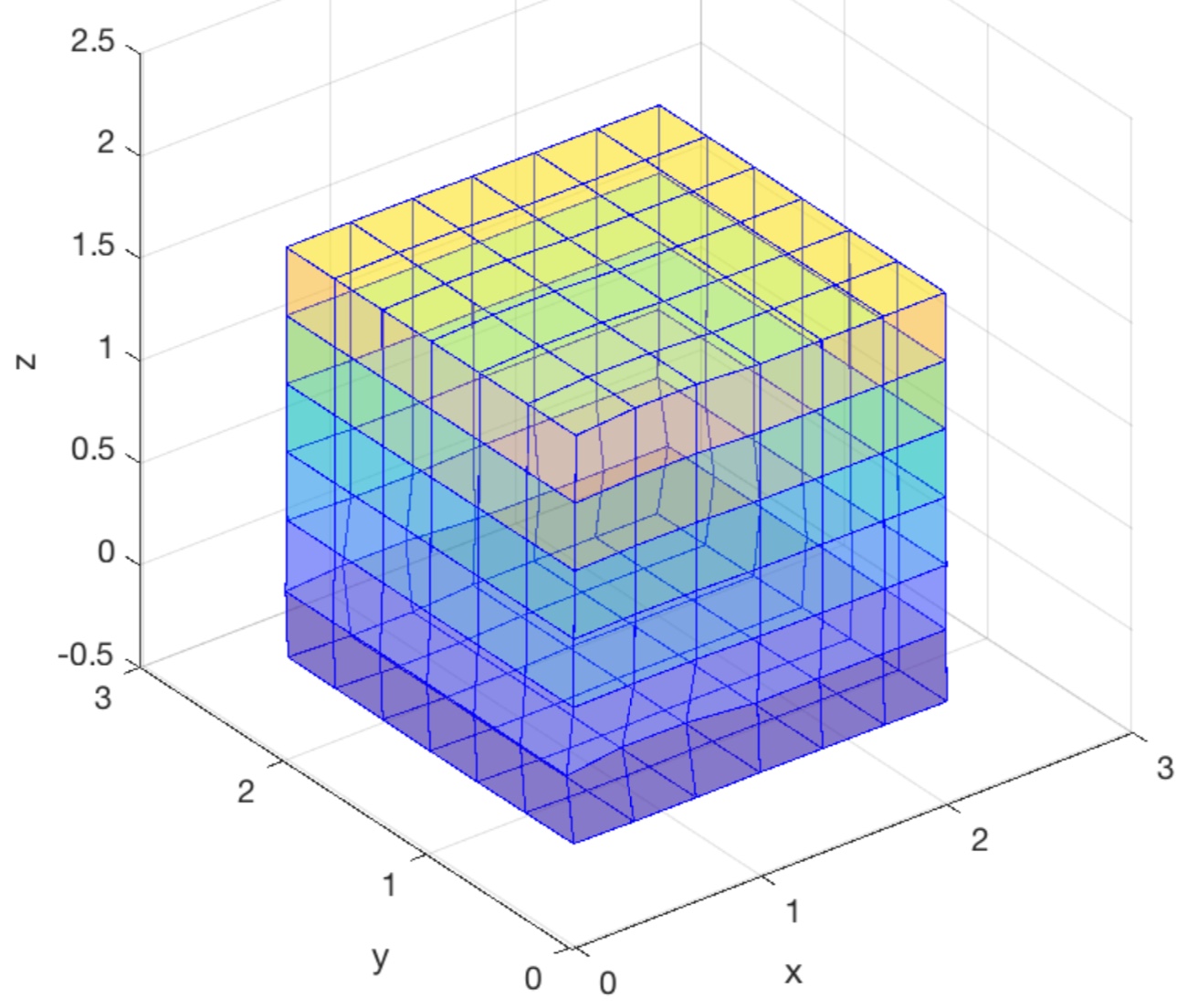} \vspace{-3pt} 
   \\
  \includegraphics[width=1.8 in]{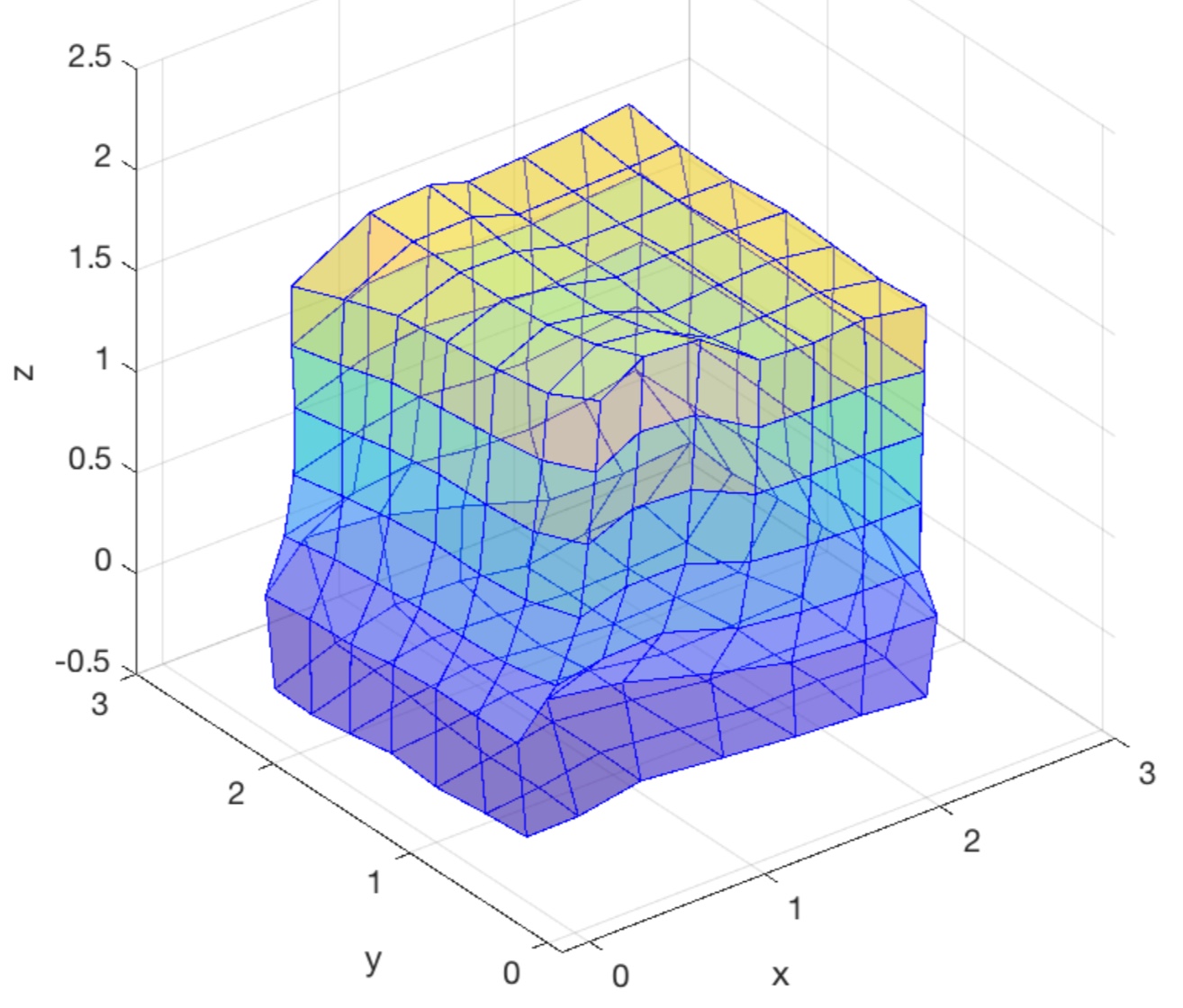} \vspace{-3pt} 
     \includegraphics[width=1.8 in]{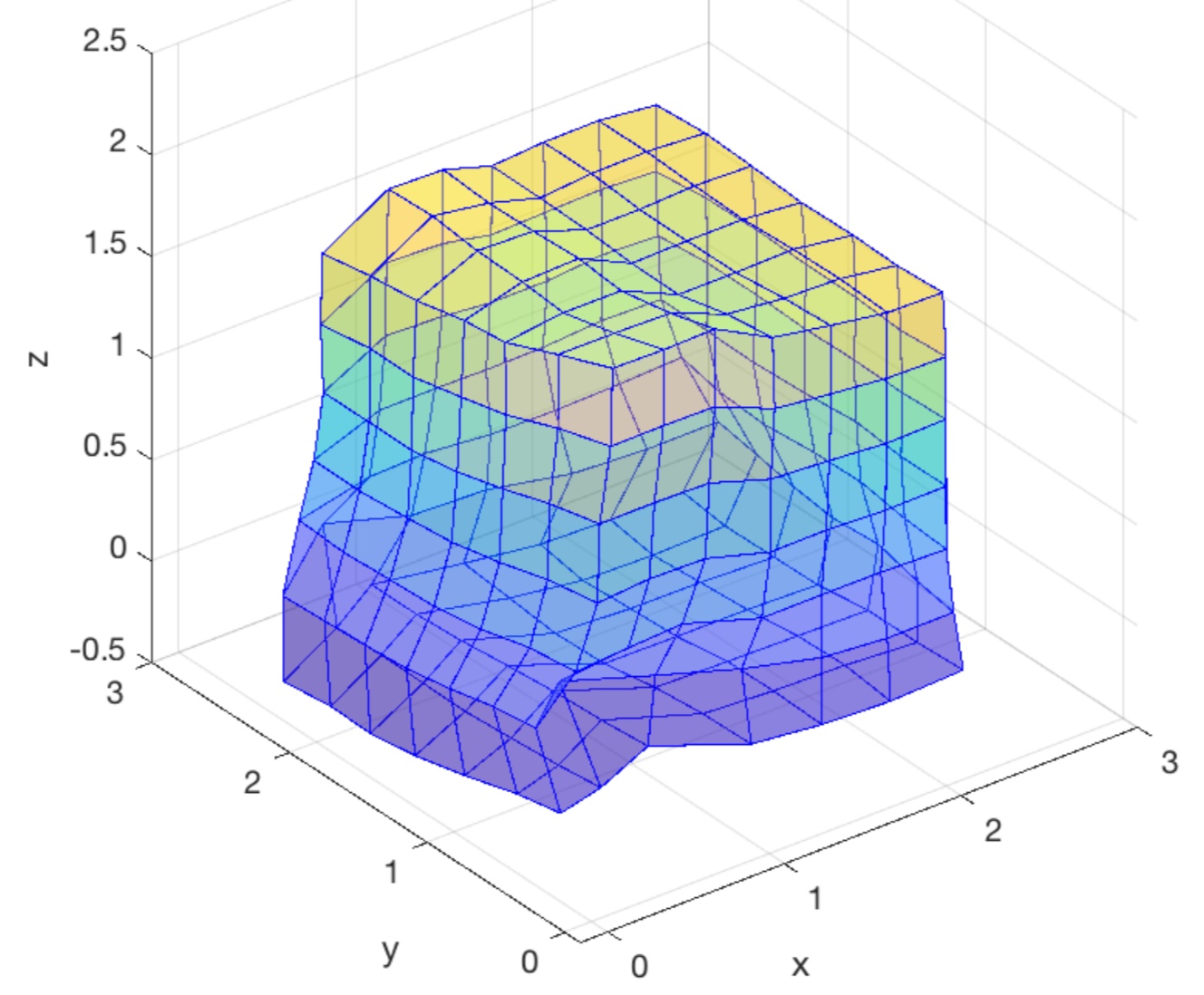} \vspace{-3pt} 
     \includegraphics[width=1.8 in]{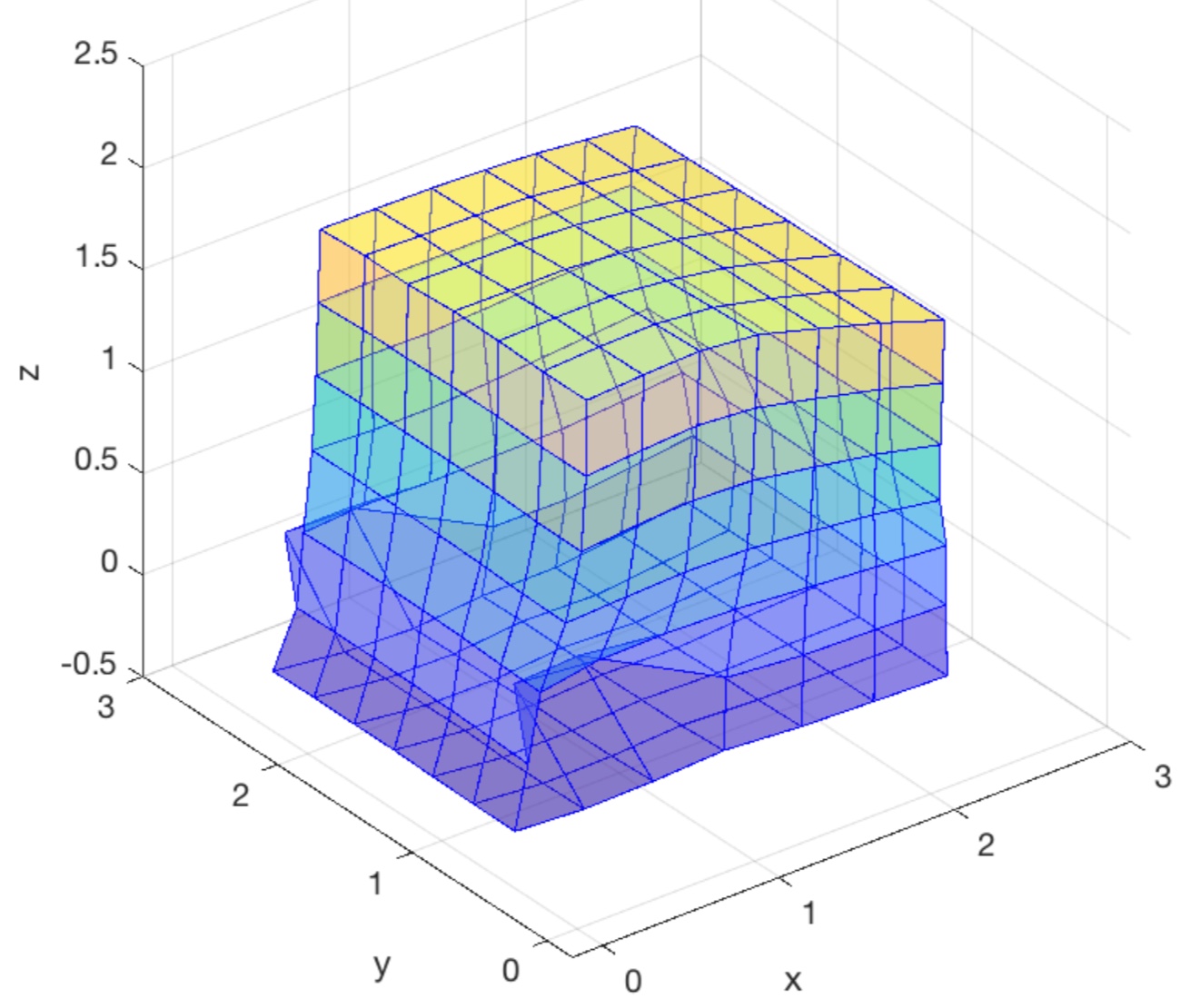} \vspace{-3pt}
   \caption{\footnotesize  \textit{Left to right}: Discrete barotropic and incompressible ideal fluid model ($r= 10^5$ and $r= 10^7$). \textit{Top to bottom}: after  $0.1$s and $2$s.} \label{3D_water_perturb} 
 \end{figure}

 \begin{figure}[H] \centering 
  \includegraphics[width=1.8 in]{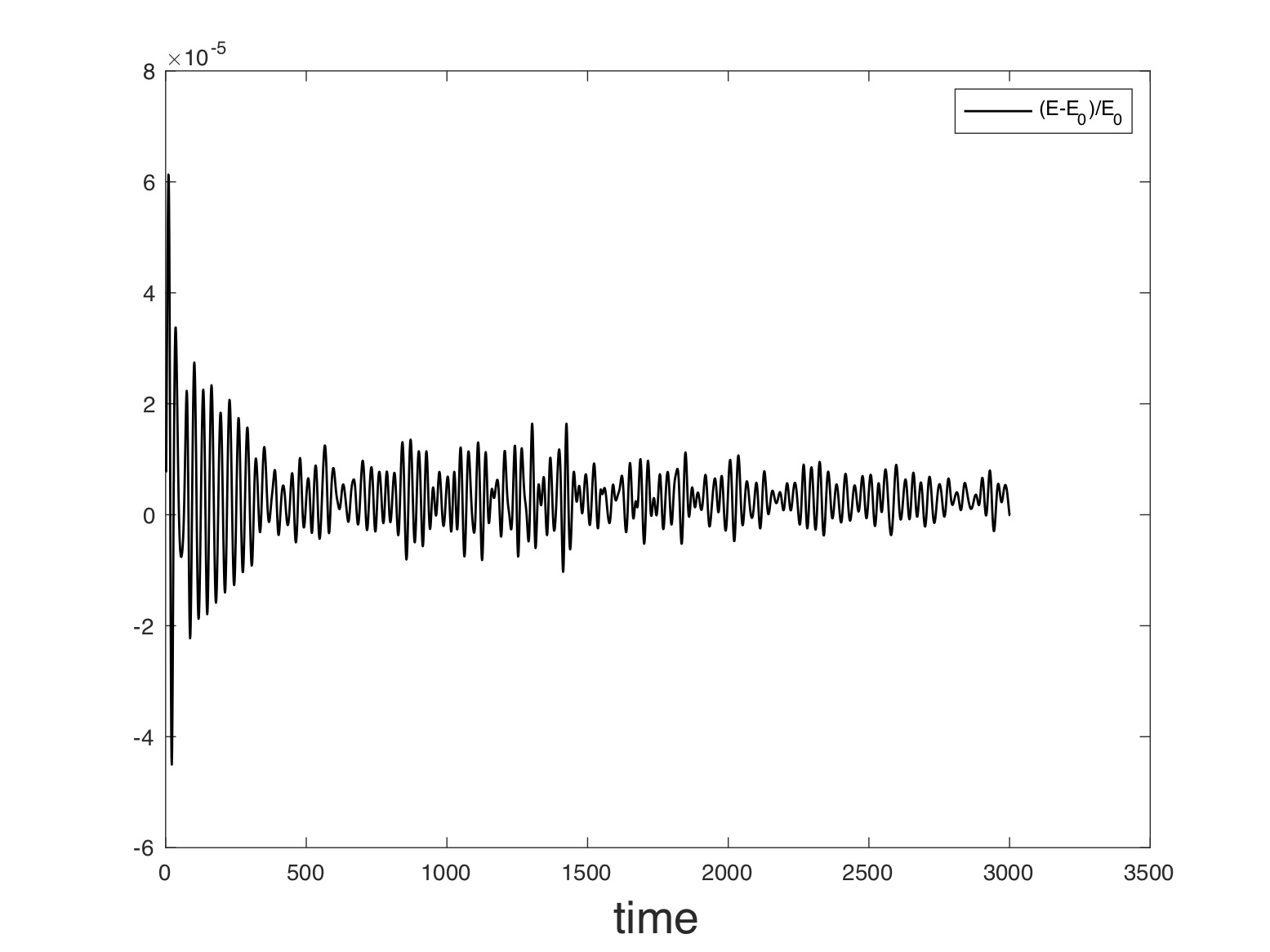} \vspace{-3pt} \;
    \includegraphics[width=1.8 in]{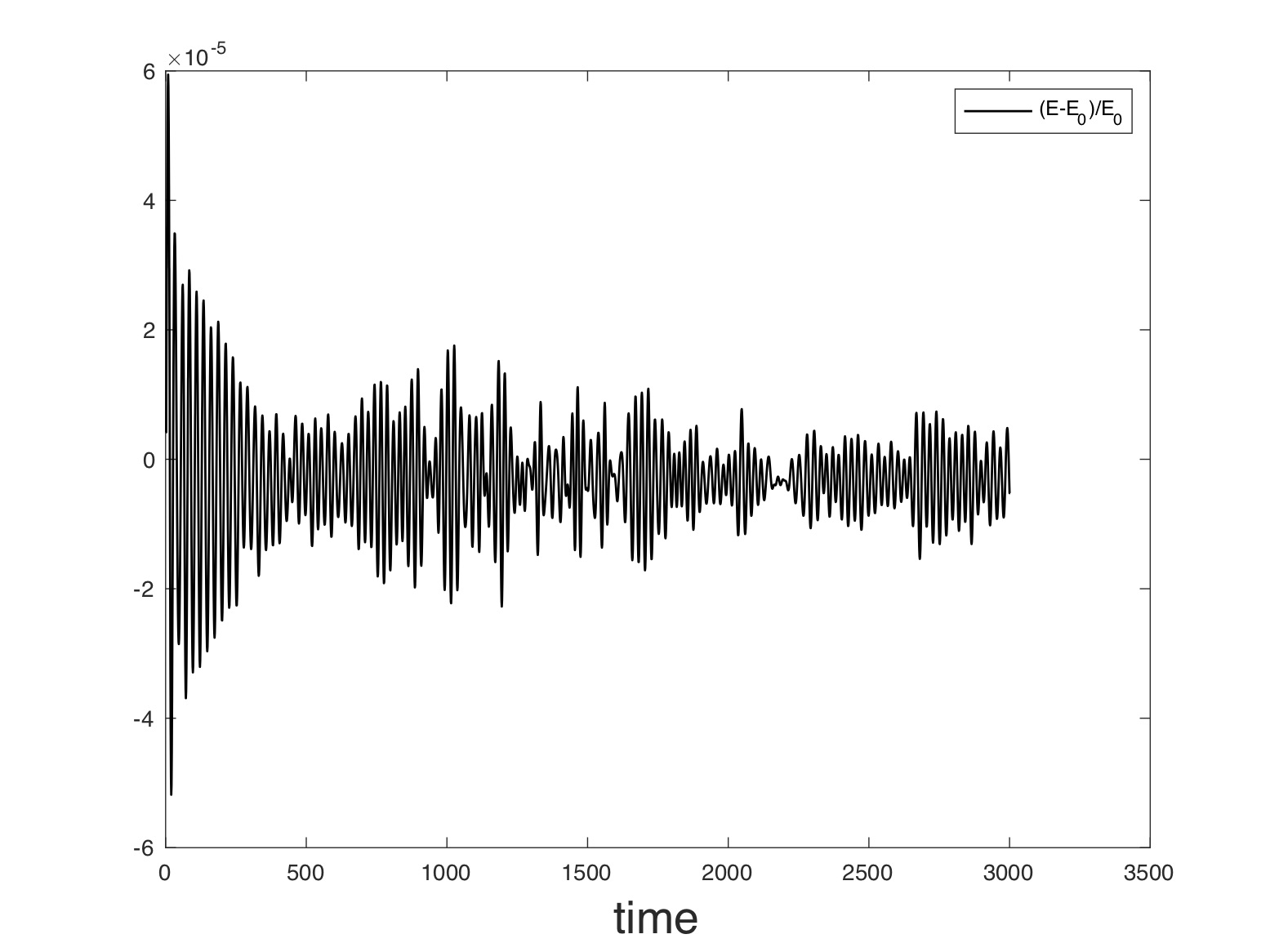} \vspace{-3pt} \;
    \includegraphics[width=1.8 in]{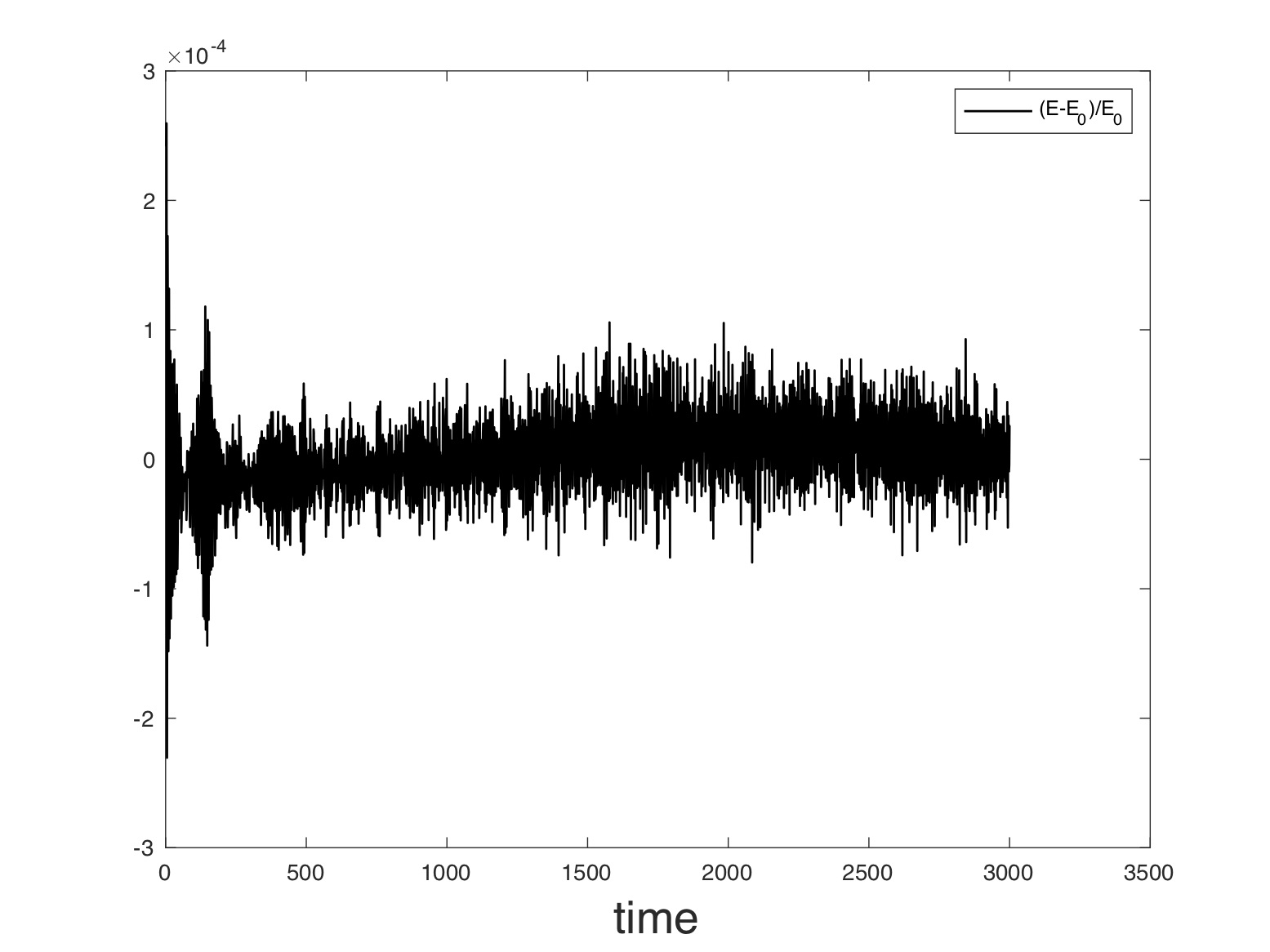} \vspace{-3pt} 
  \\
  \includegraphics[width=1.8 in]{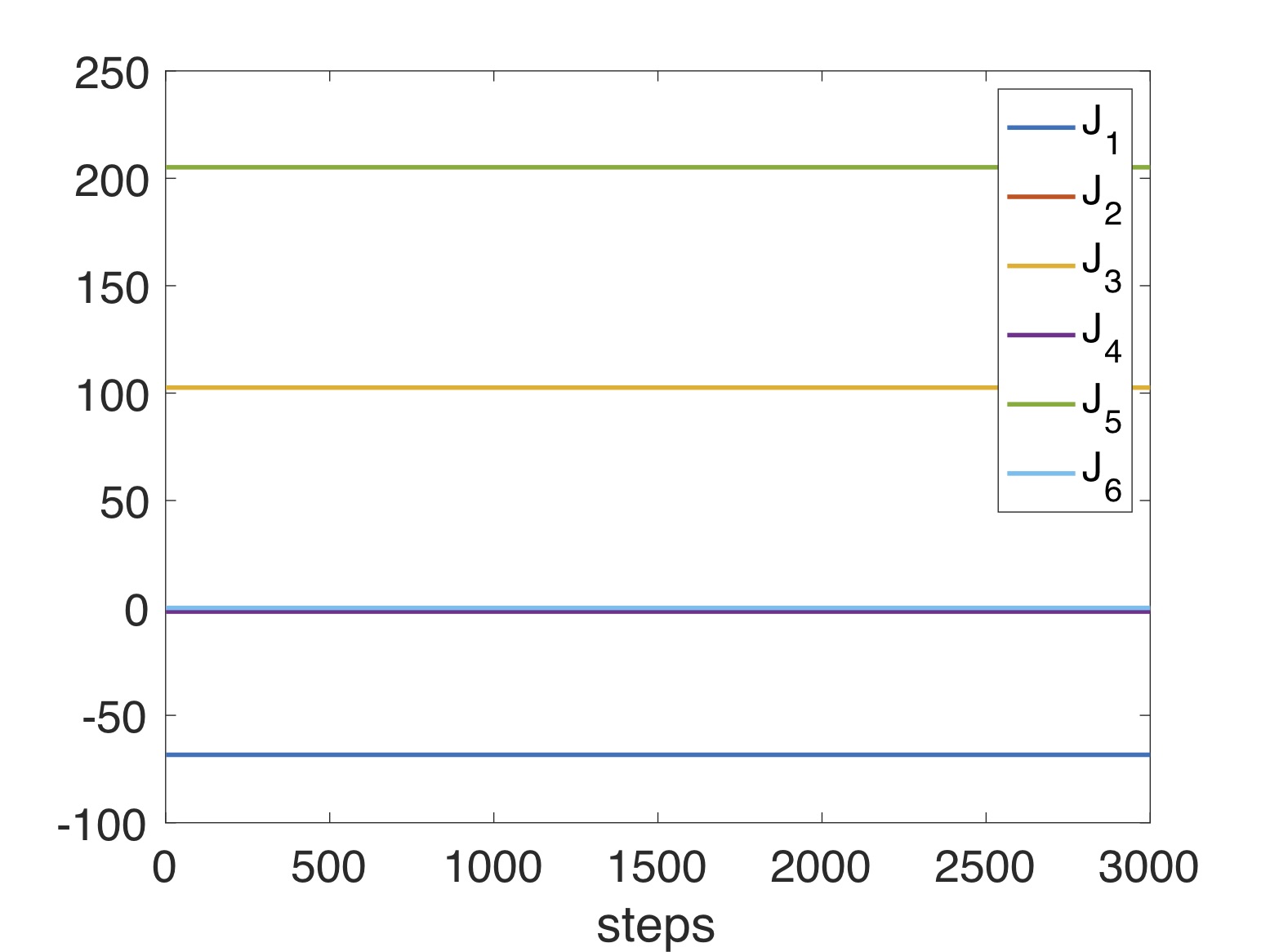} \vspace{-3pt} \;  \includegraphics[width=1.8 in]{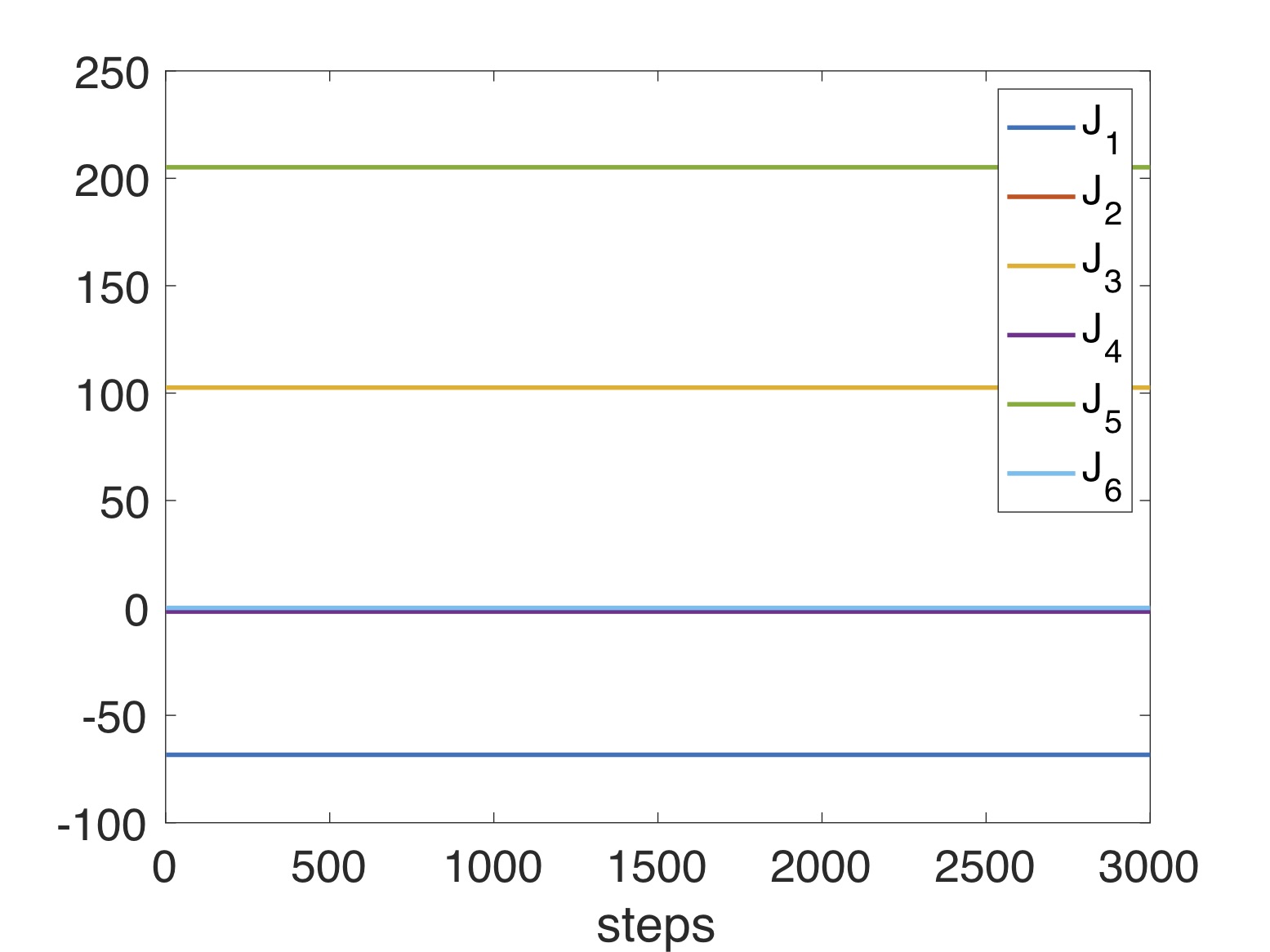} \vspace{-3pt} \;  \includegraphics[width=1.8 in]{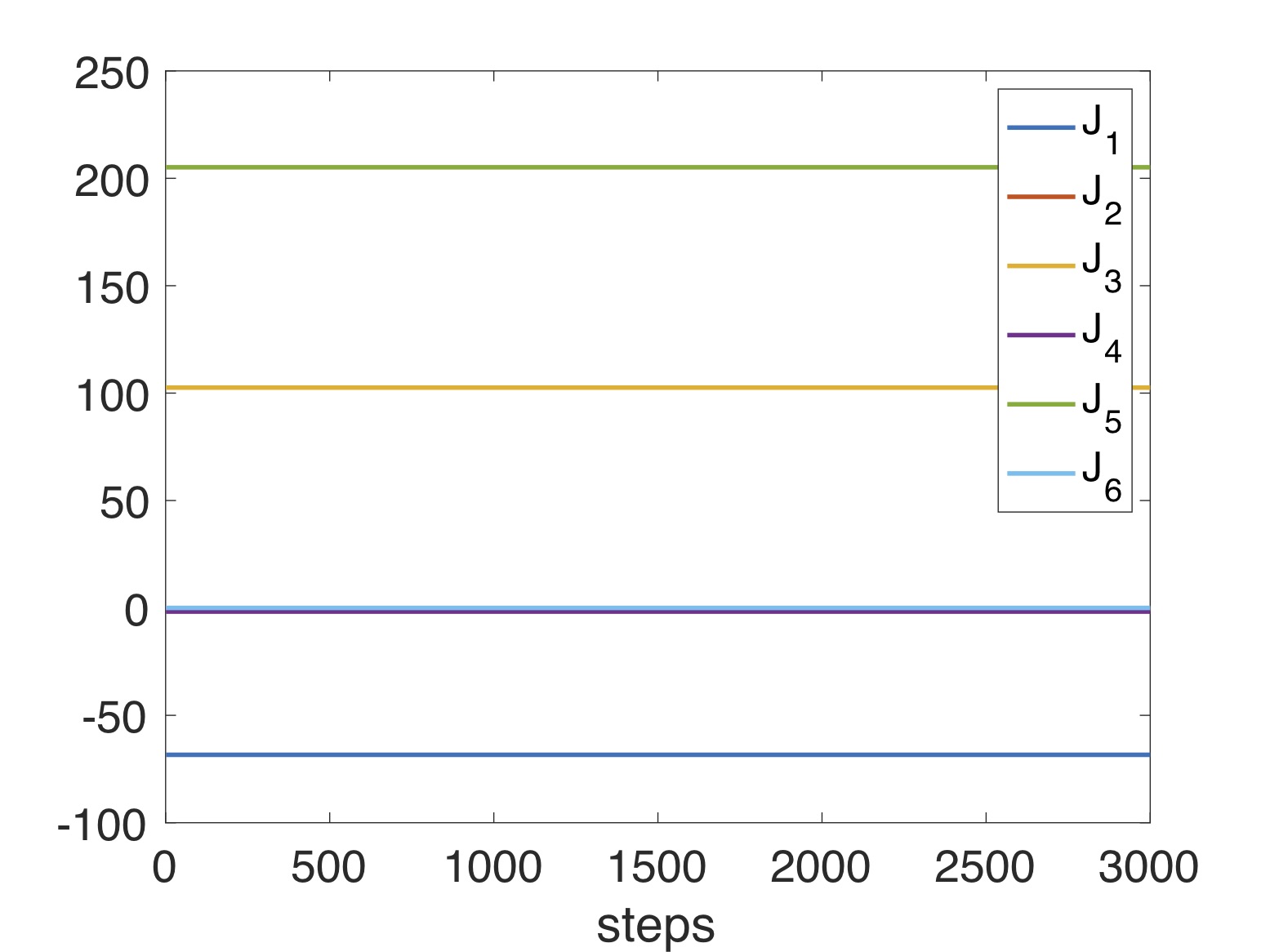} \vspace{-3pt}
   \caption{\footnotesize  \textit{Left to right}: barotropic and incompressible ideal fluid model ($r= 10^5$ and $r=10^7$). \textit{Top to bottom}: Relative energy and momentum map evolution during $3$s.}\label{water_3D_energy_momentum}
  \end{figure}

The main interest of this test in vacuum with free boundaries is to exhibit the perfect preservation of the symmetries, see the figures above.

 \subsubsection{Example 4: Impact against an obstacle of a fluid flowing on a surface} 
 
As explain in \S \ref{2D_water_impact}, in this example the problems are to find the extremum of the action subject to 
an equality constraint associated to incompressibility and inequality constraints imposing the fluid to stay on a surface and outside of an obstacle.

\medskip 

Let $(\mathcal{P}_1)$, resp., $(\mathcal{P}_2)$ denote the problem to solve for barotropic fluid, resp., incompressible ideal fluid. These two problems are already described in \S \ref{2D_water_impact}.

\medskip

Consider a barotropic fluid with properties $\rho_0= 997 \, \mathrm{kg/m}^2$, $\gamma =7 $, $A = \tilde{A} \rho_0^{-\gamma}$ with $\tilde A = 3.041\times 10^4$ Pa, and $B = 3.0397\times 10^4$ Pa. The size of the discrete reference configuration at time $t^0$ is $1.6 \, \mathrm{m} \times 1 \, \mathrm{m}  \times 0.4 \, \mathrm{m}$, with time-step $\Delta t=5\times 10^{-5}$ and space-steps $\Delta s_1=0.1$m, $\Delta s_2= 0.2$m, $\Delta s_3= 0.1$m. The value of the impenetrability penalty coefficients are $K_1=K_2= 5 \times 10^9$. We consider both the compressible barotropic fluid and the incompressible case with penalty given by $r=10^8$ and $r=10^9$.

\medskip 

As in 2D, the initial motion of the fluid is only due to the gravity.

\begin{figure}[H] \centering 
    \includegraphics[width=2.5 in]{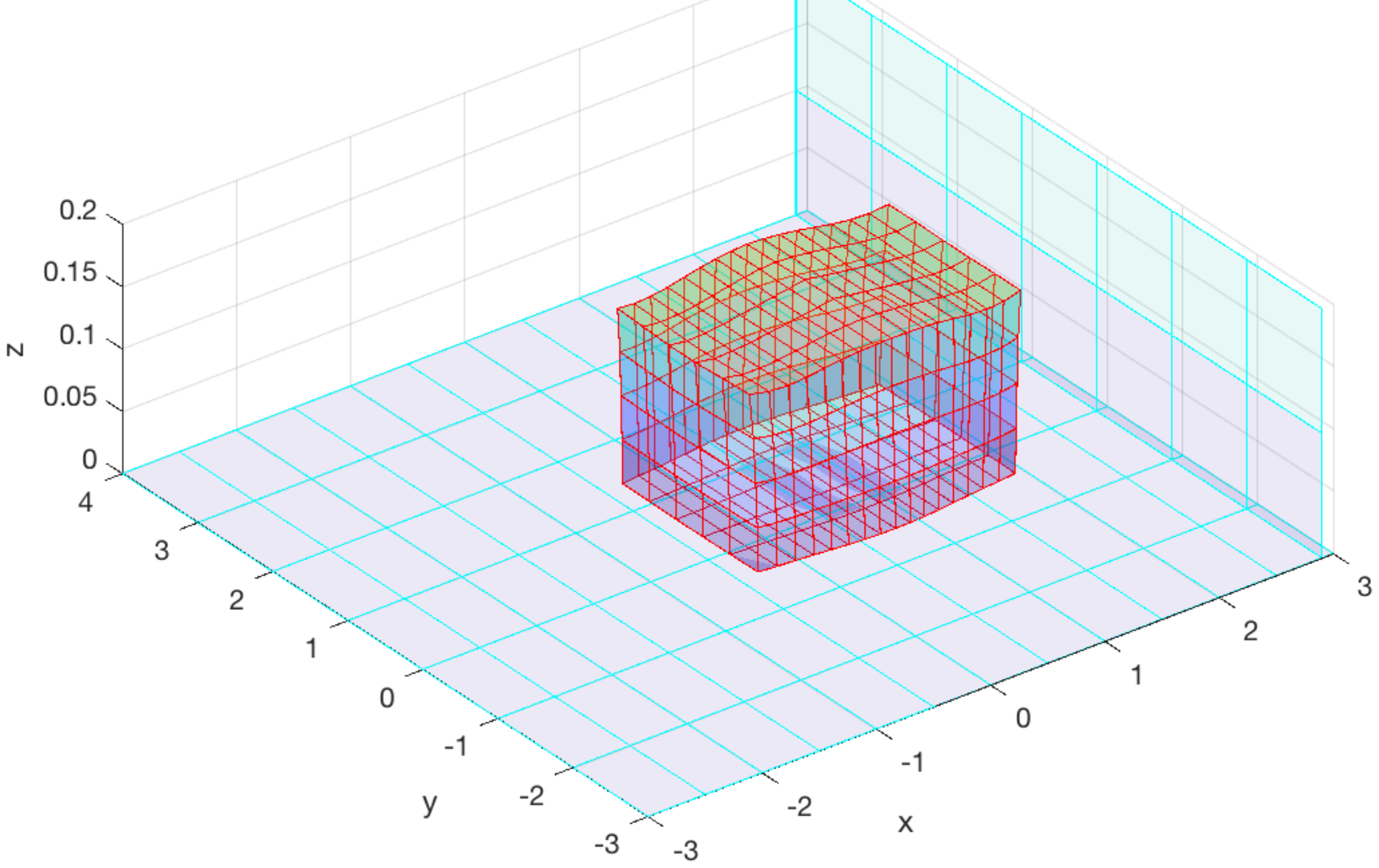} \vspace{2pt}
     \includegraphics[width=2.5 in]{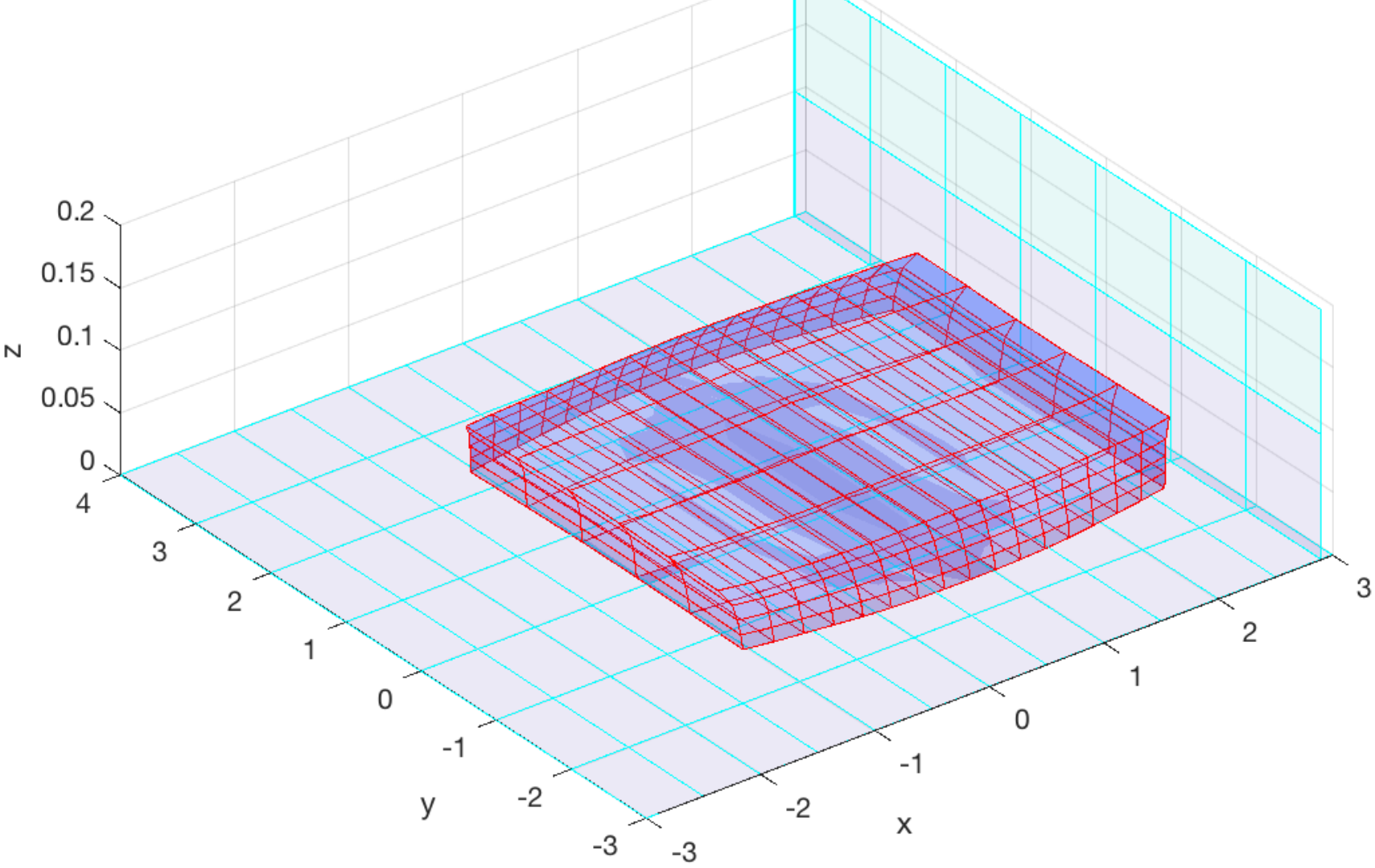} \vspace{2pt} 
    \\
      \includegraphics[width=2.5 in]{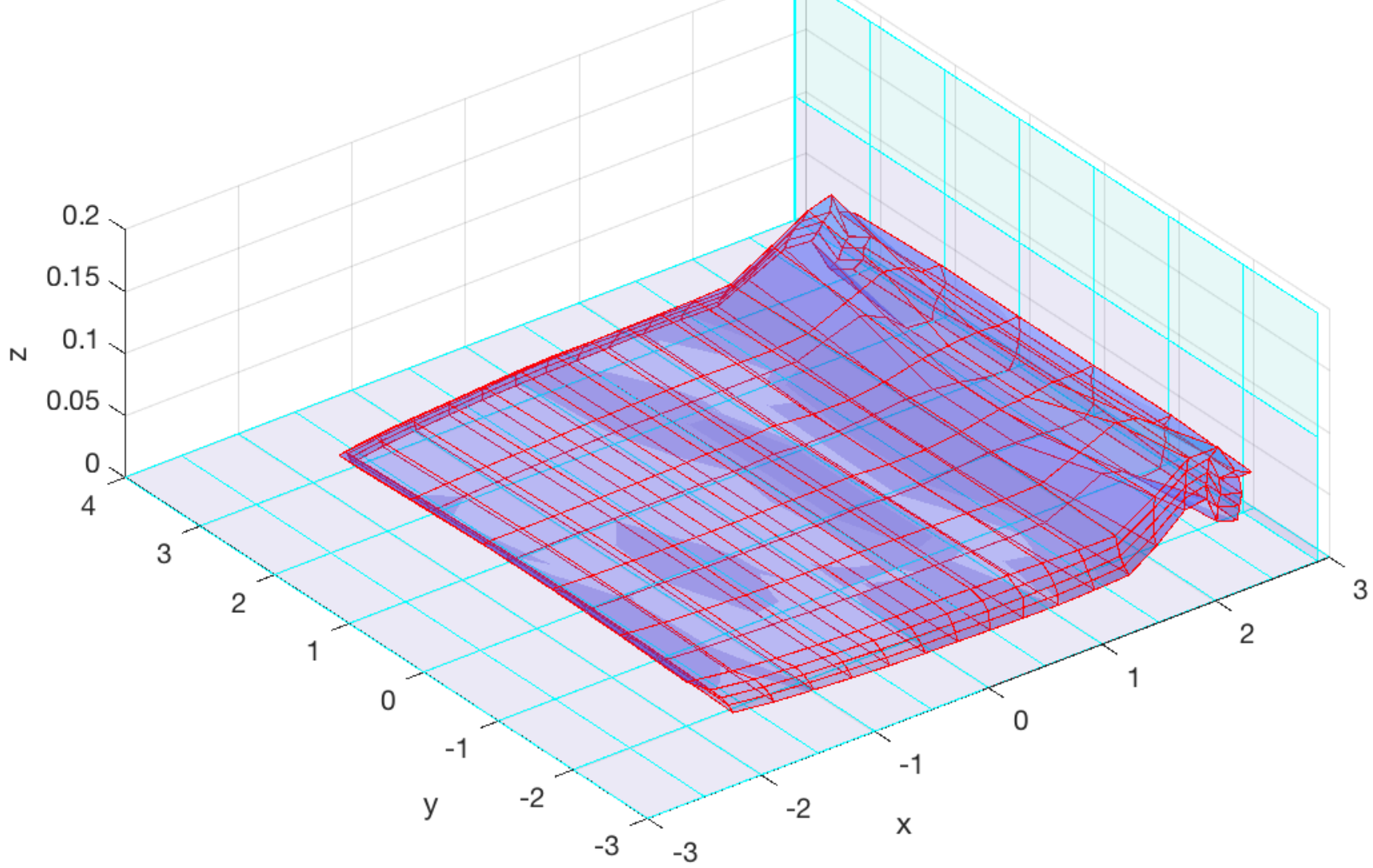}
     \includegraphics[width=2.5 in]{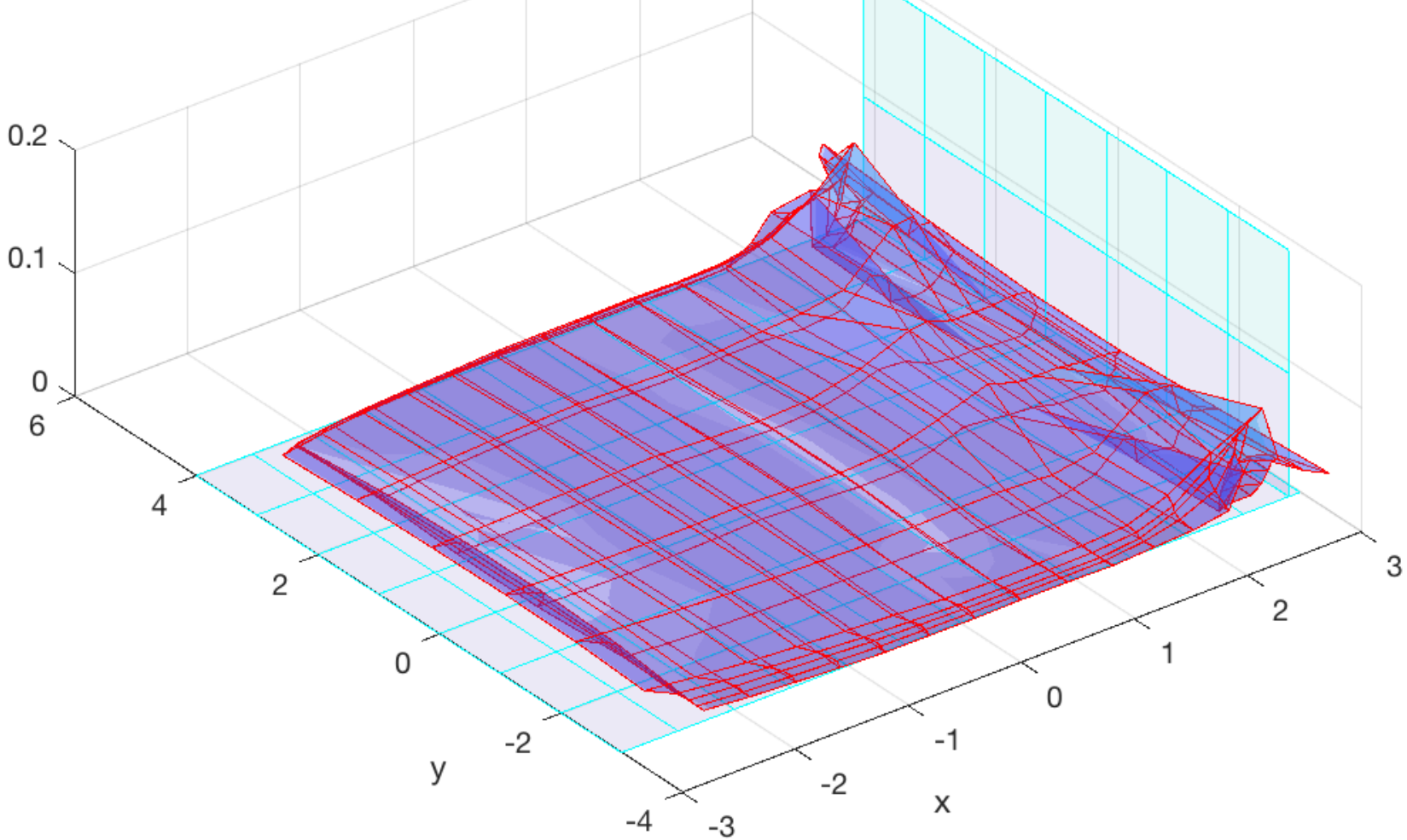} \vspace{-3pt}
   \caption{\footnotesize Barotropic fluid model with impact  after $0.4$s, $0.8$s,  $1.1$s, $1.4$s.} \label{3D_baro_water_contact} 
 \end{figure}
 
 \begin{figure}[H] \centering 
      \includegraphics[width=2.5 in]{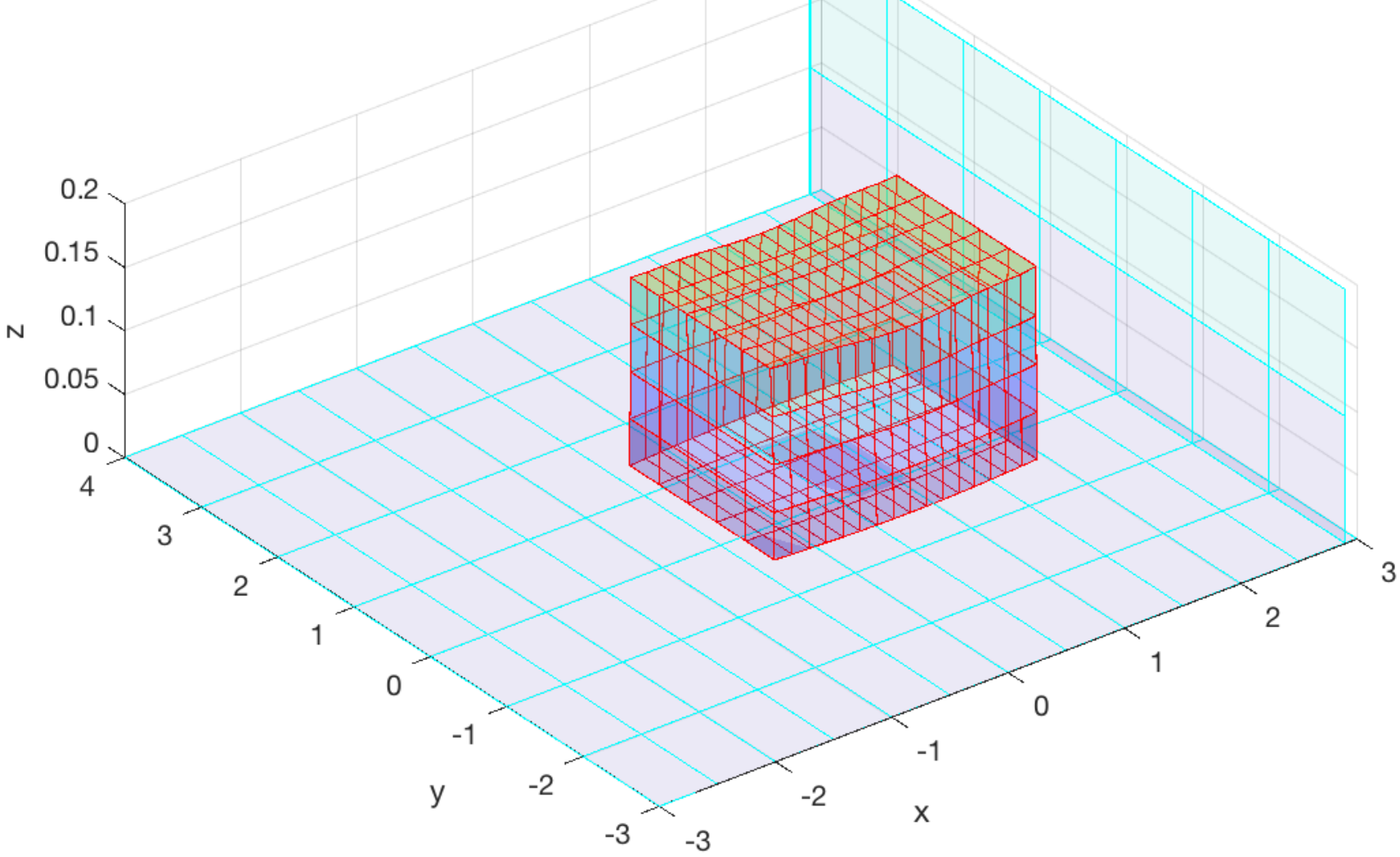} \vspace{2pt} \quad
    \includegraphics[width=2.5 in]{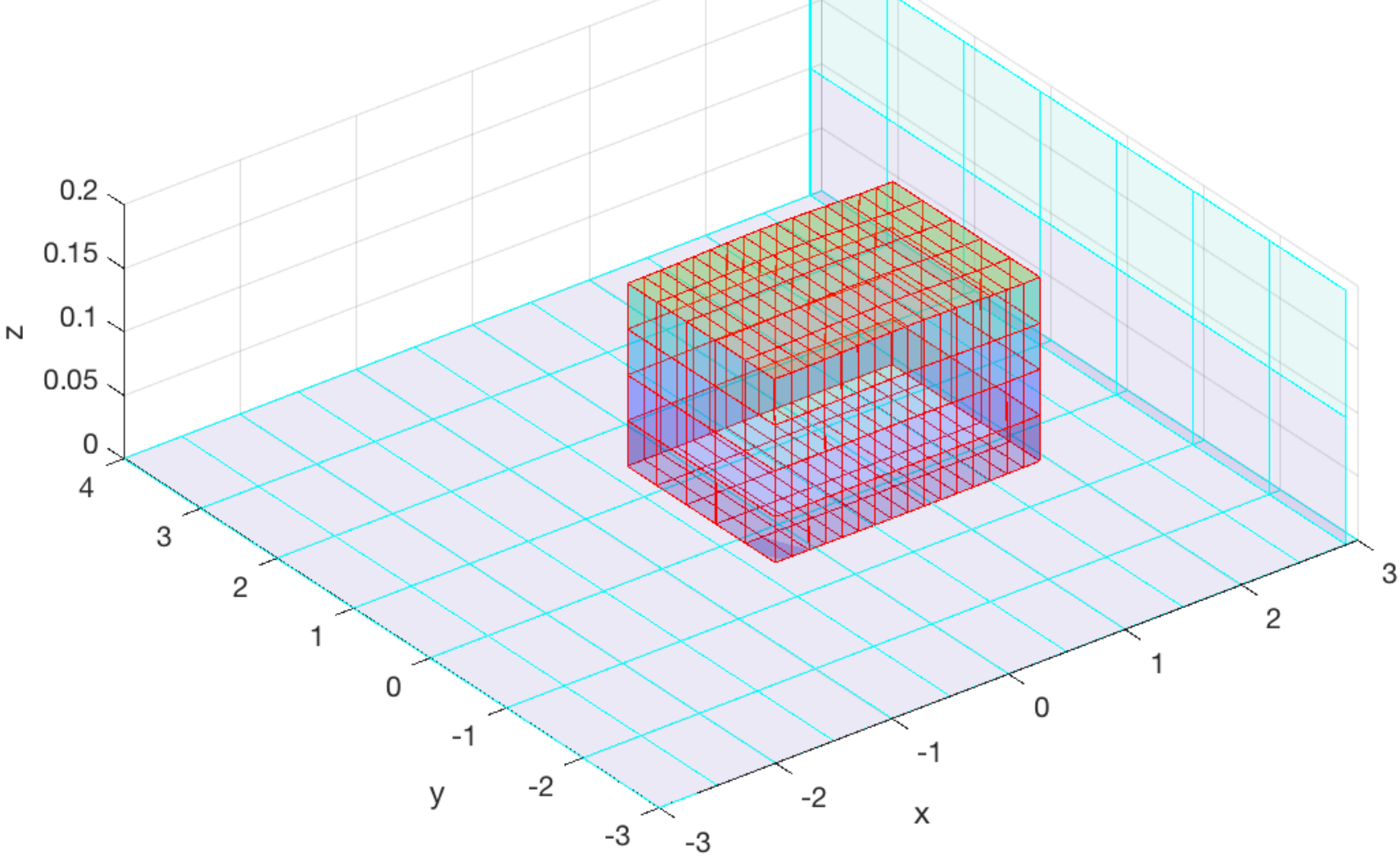} \vspace{2pt}
    \\
     \includegraphics[width=2.5 in]{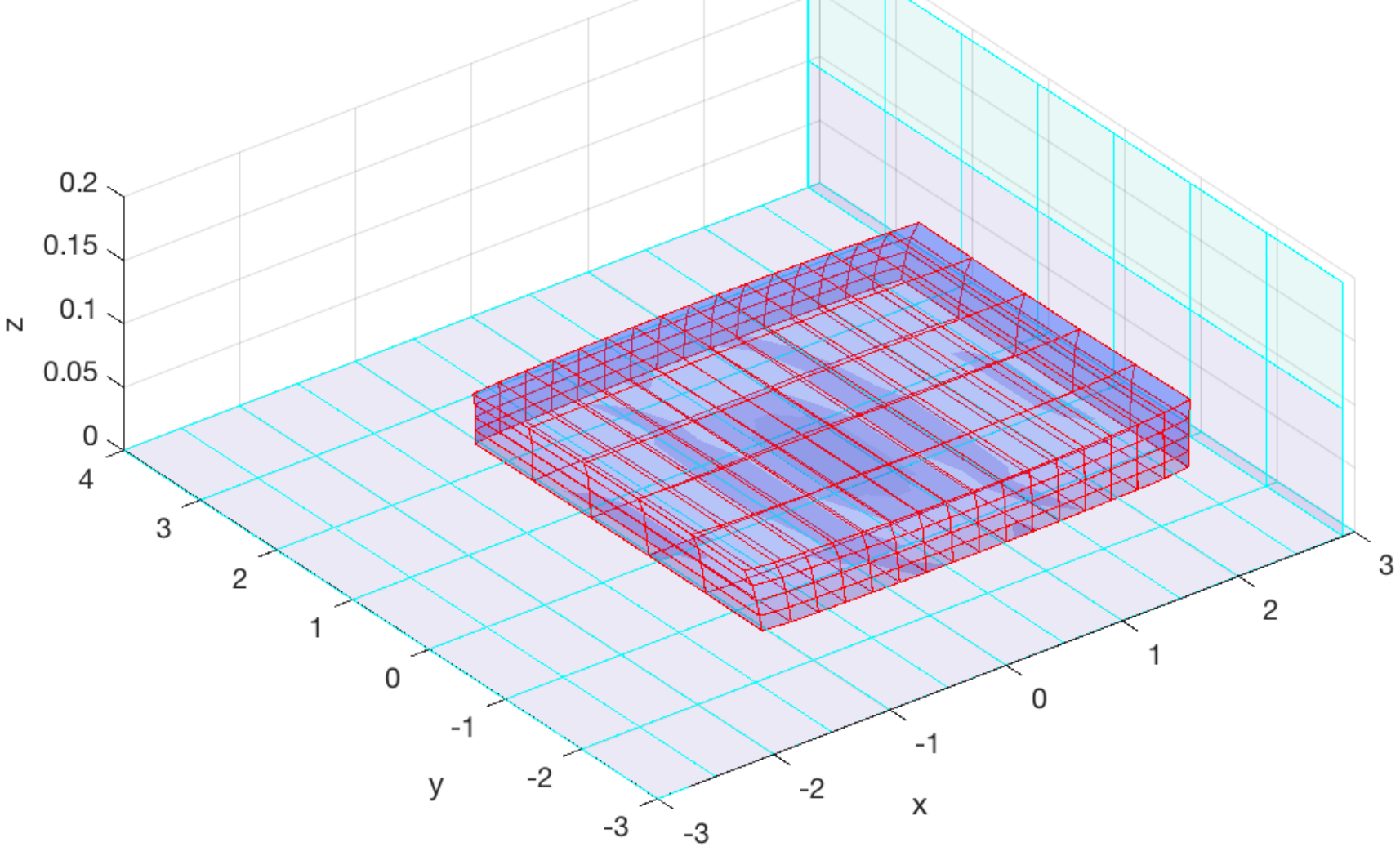} \vspace{2pt} \quad
      \includegraphics[width=2.5 in]{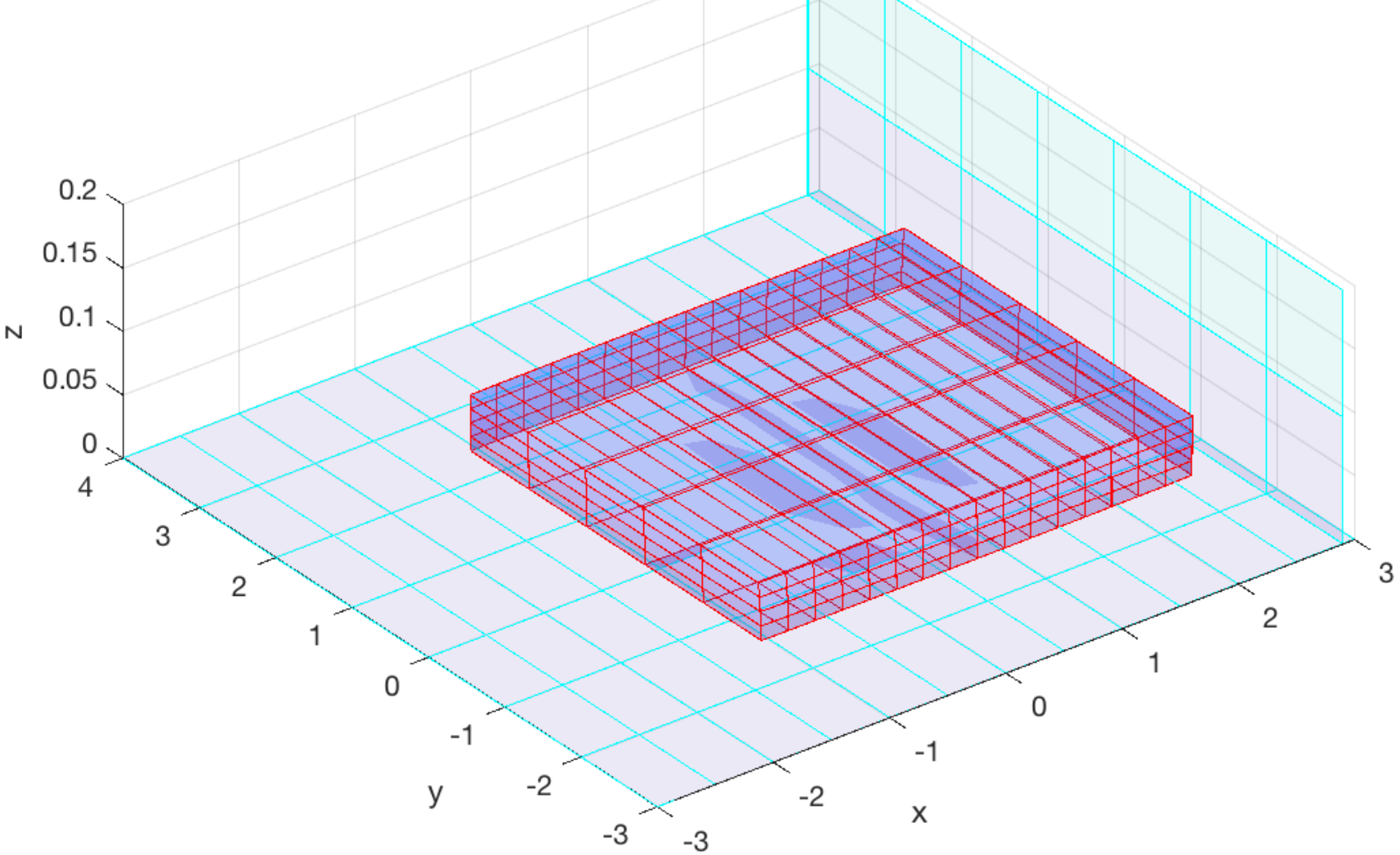} \vspace{2pt} 
    \\
      \includegraphics[width=2.5 in]{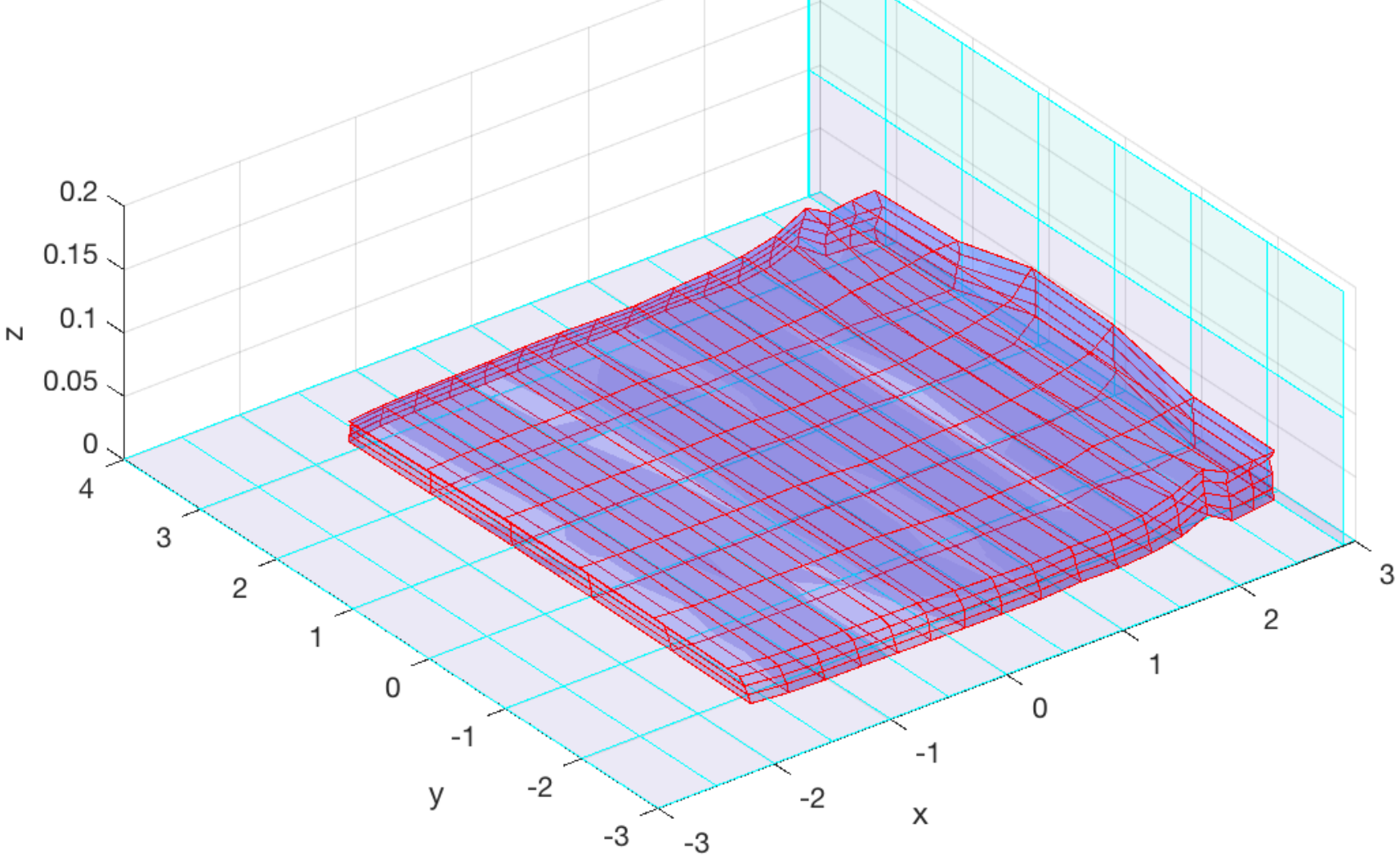} \vspace{-3pt} \quad
       \includegraphics[width=2.5 in]{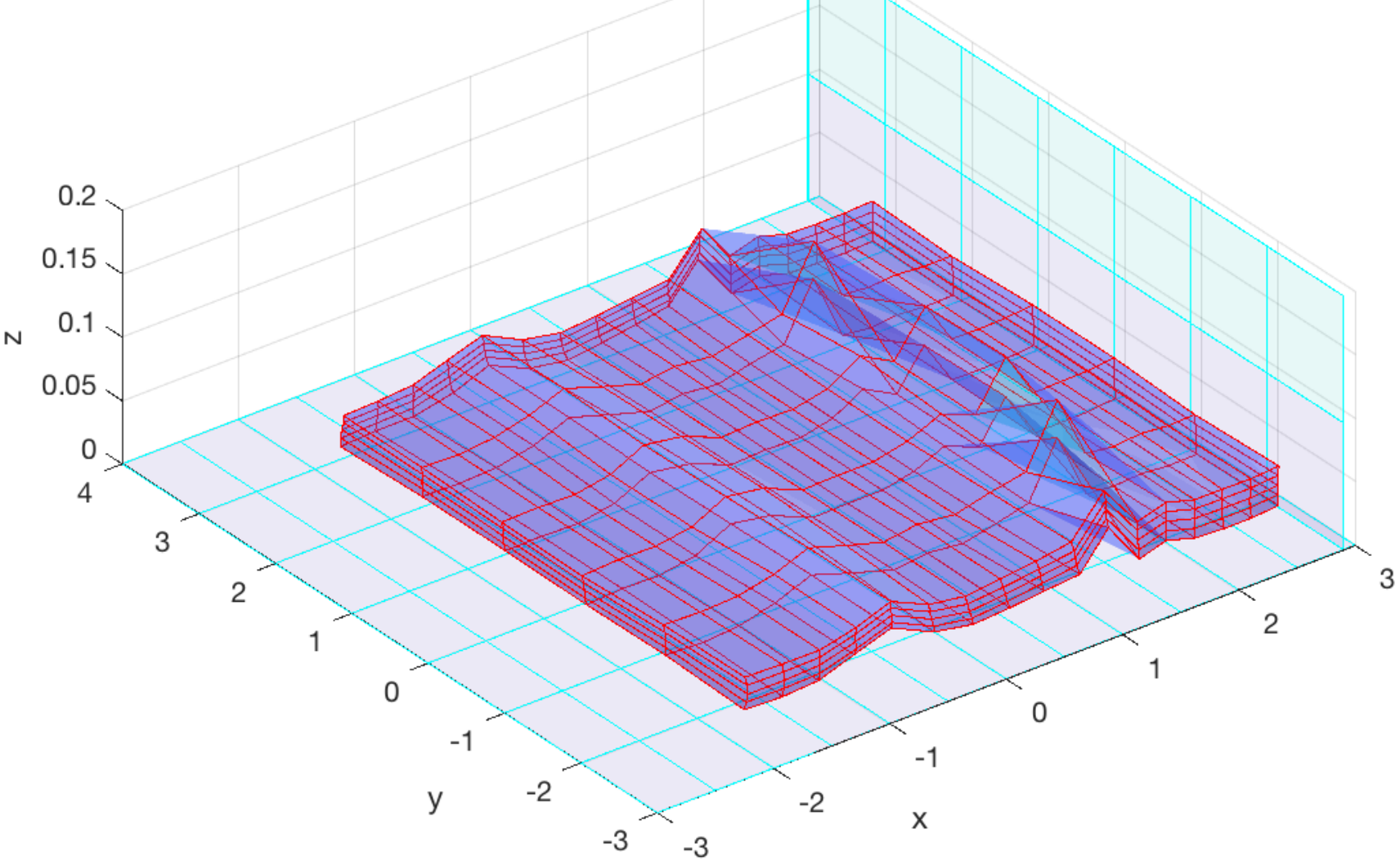} \vspace{-3pt}
      \\
           \includegraphics[width=2.9 in]{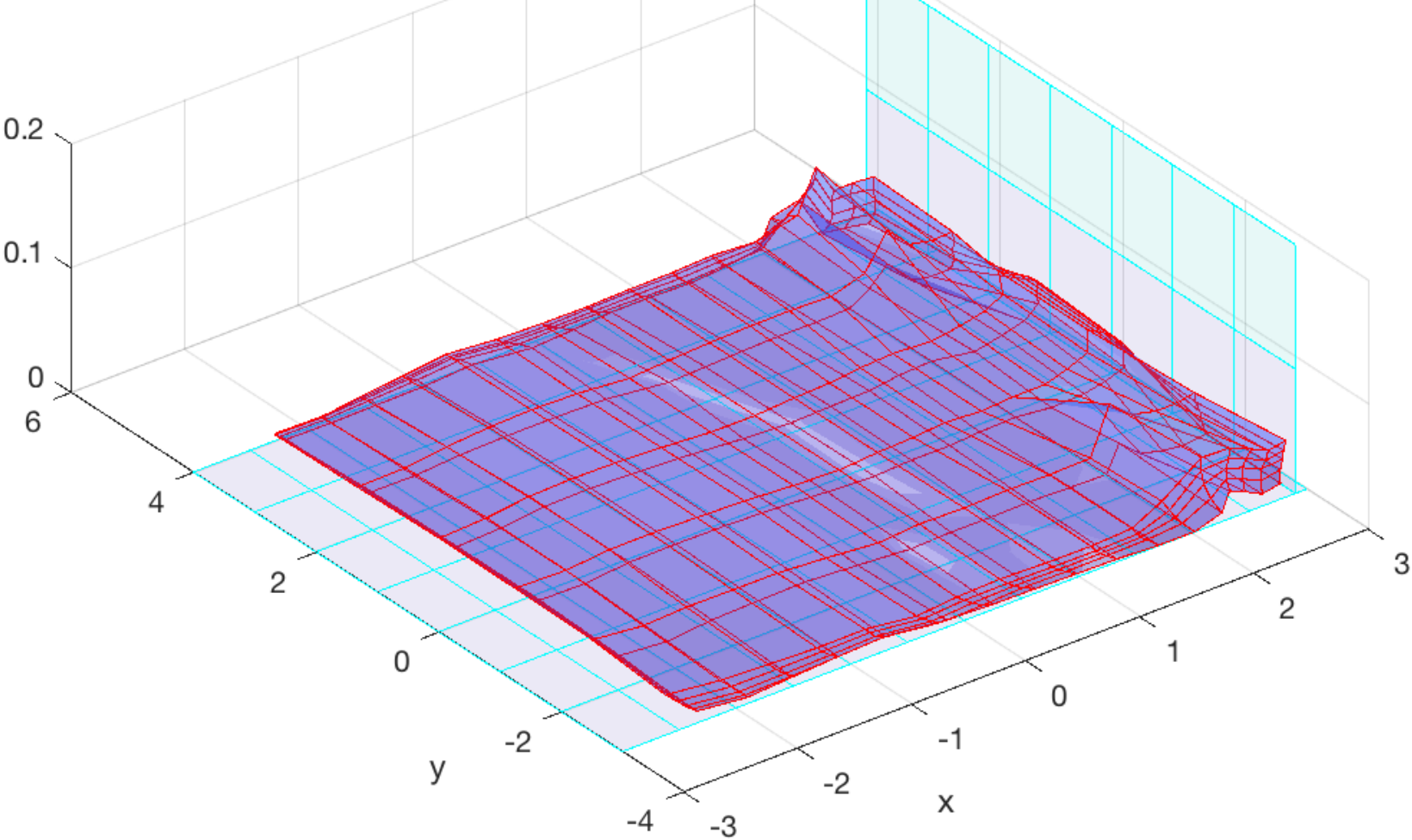} \vspace{-3pt}
     \includegraphics[width=2.75 in]{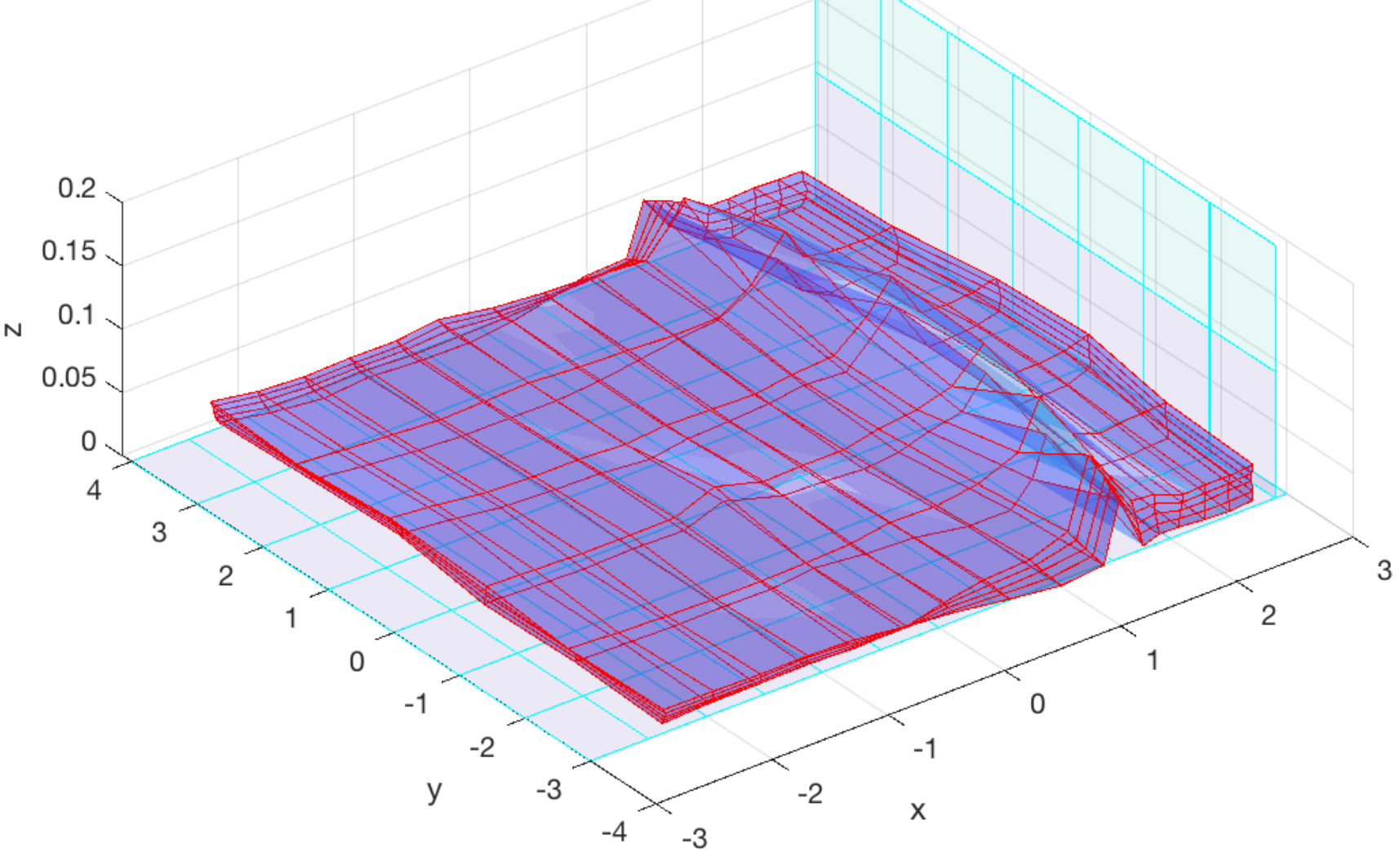} \vspace{-3pt}
   \caption{\footnotesize Fluid impact.  \textit{Left to right}: From less ($r=10^8$) to more ($r= 10^9$) incompressibility.\textit{Top to bottom}:  after $0.4$s, $0.8$s,  $1.1$s, $1.4$ s.} \label{3D_icomp_water_contact} 
 \end{figure}
 
 Note the large differences in behavior between the three tests that correspond to barotropic fluid (Fig.\,\ref{3D_baro_water_contact}) and to fluids which are more or less incompressible (Fig.\,\ref{3D_icomp_water_contact}).

\begin{figure}[H] \centering 
\includegraphics[width=1.8 in]{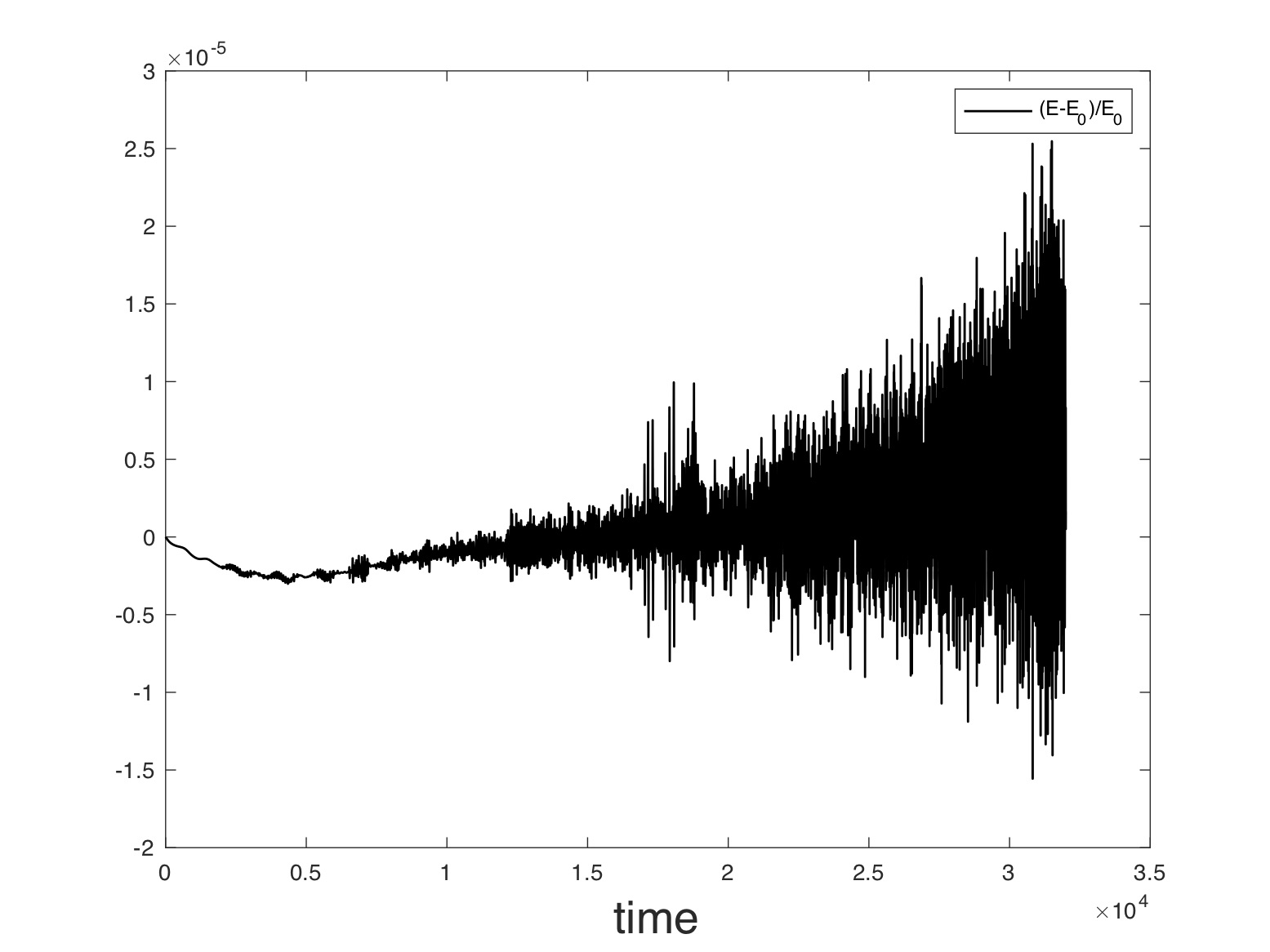} \vspace{-3pt} \quad
\includegraphics[width=1.8 in]{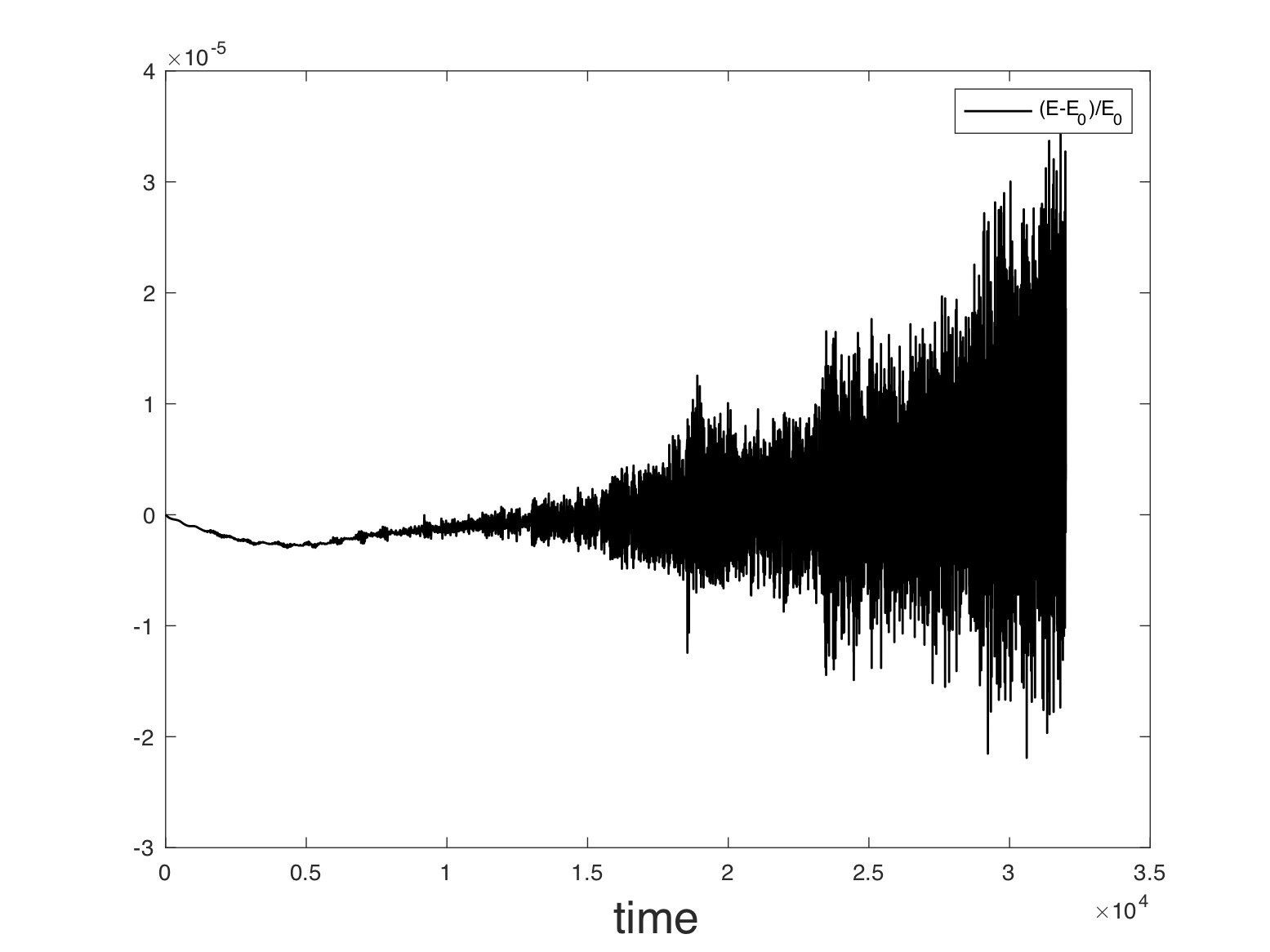} \vspace{-3pt} 
\includegraphics[width=1.8 in]{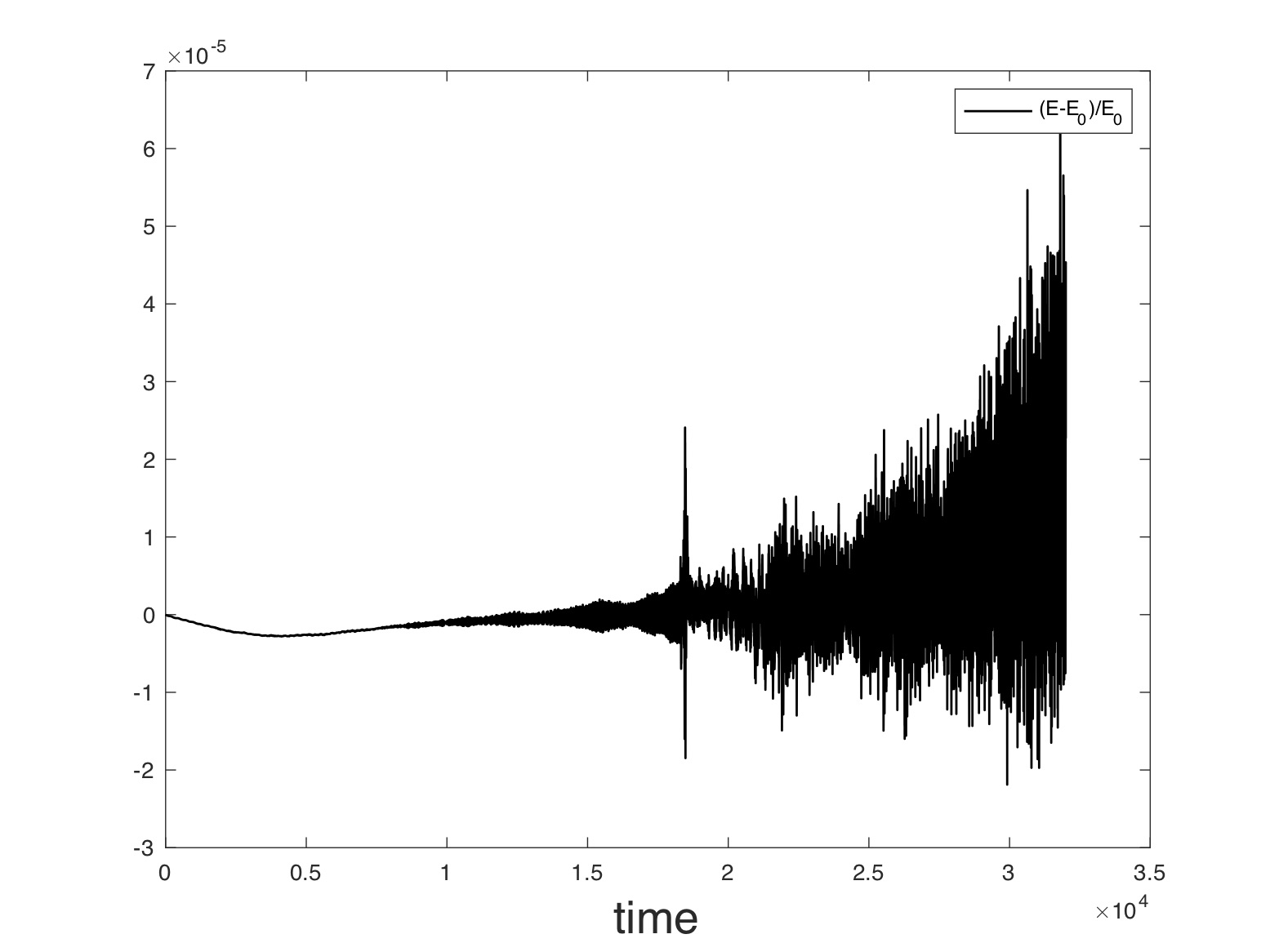} \vspace{-3pt} \\ 
\includegraphics[width=1.8 in]{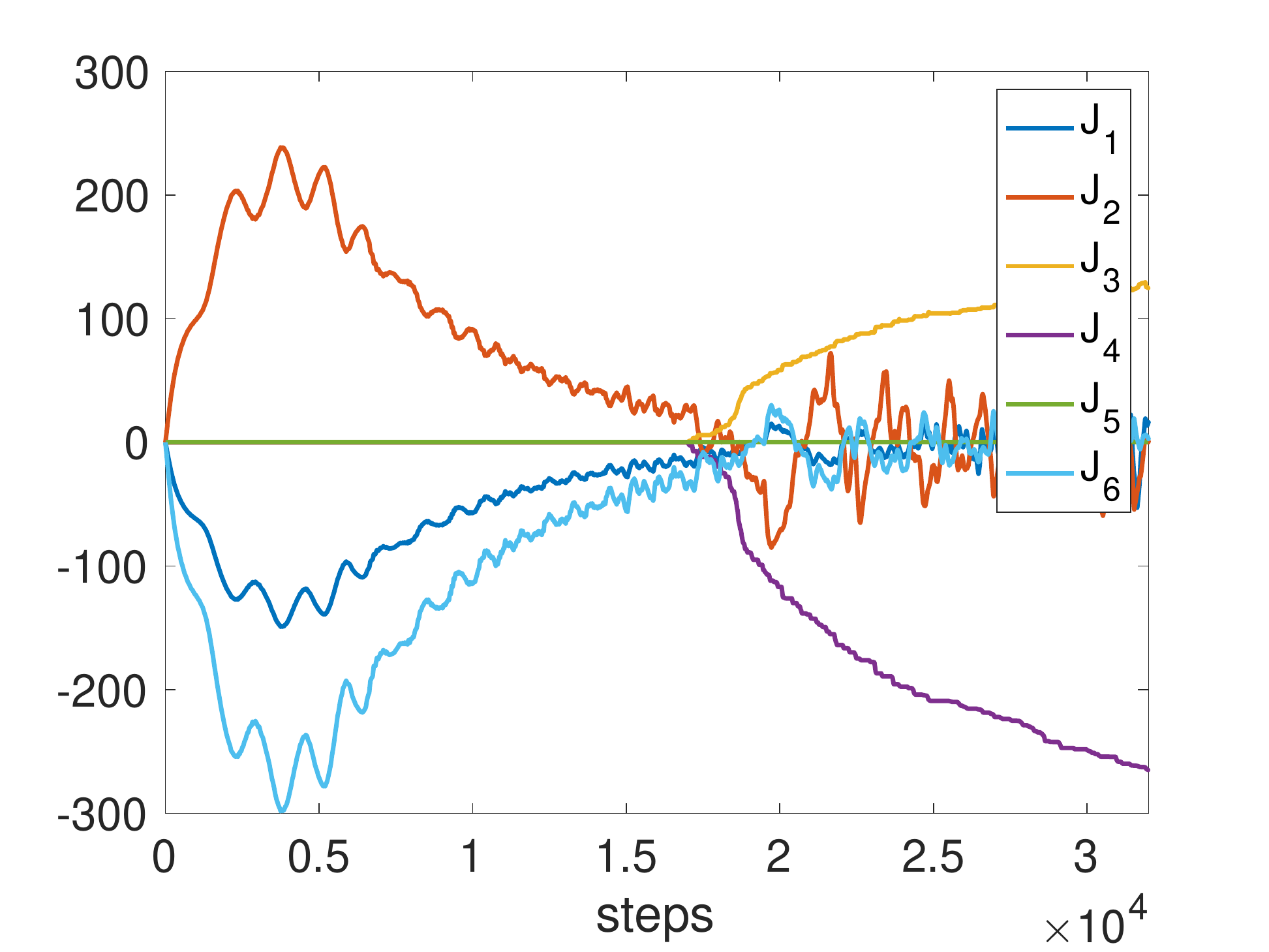} \vspace{-3pt} \quad
\includegraphics[width=1.8 in]{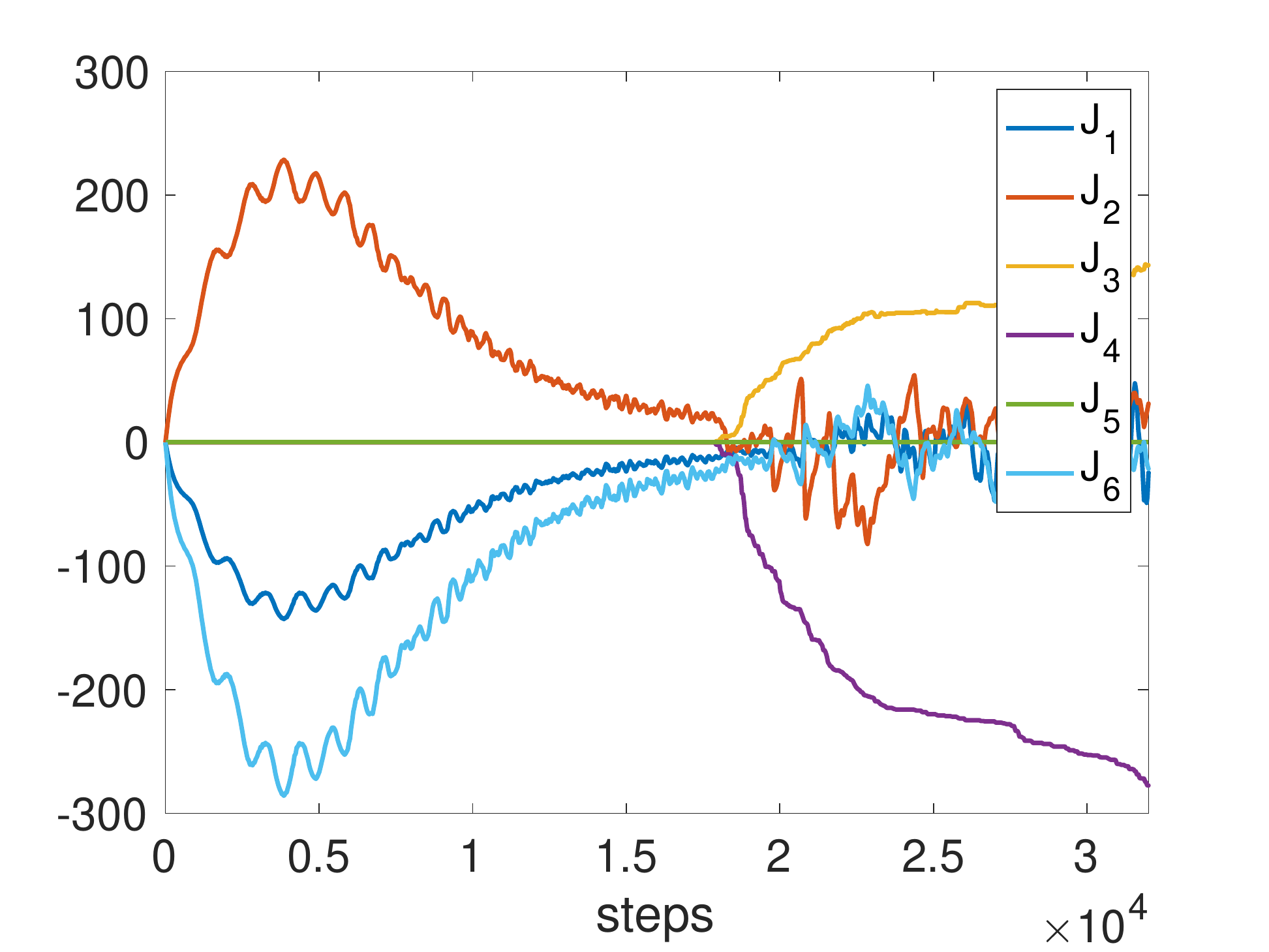} \vspace{-3pt} 
\includegraphics[width=1.8 in]{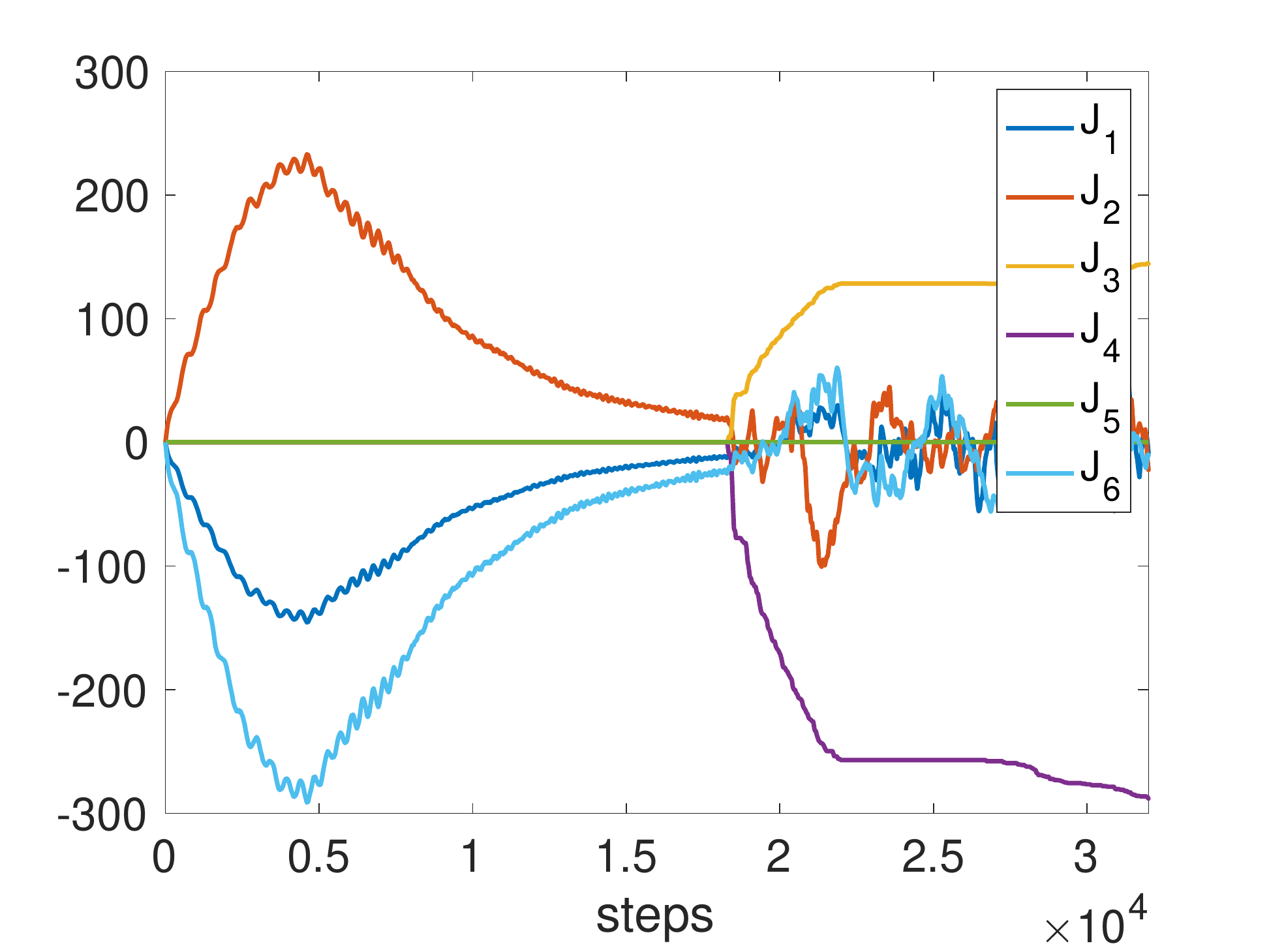} \vspace{-3pt} 
\caption{\footnotesize  \textit{From left to right}: Barotropic, incompressible ideal ($r=10^8$),  incompressible ideal ($r= 10^9$) fluid impact model. \textit{From top to bottom}: Relative energy and momentum map evolution.} \label{water_contact_energy_J_3D}
\end{figure}
We observe that three components of the momentum map are preserved until the contact with the obstacle. These three components correspond to the symmetries of the gravitational term, given by the three dimensional subgroup of $SE(3)$ consisting of rotations around the vertical axis and translations in the horizontal plane.  After contact only one symmetry with respect to one axis of translation is conserved, namely, the translation parallel to the obstacle wall. After impact, the total energy becomes more unstable, as observed in the 2D test, see Fig.\,\ref{water_contact_energy_J}.

\medskip

Results concerning impact between fluid and solid are promising, however further investigations in numerics are necessary. In particular, we must develop an implicit integrator in order to increase the time-step and the performance of the integrator when the stability of the flow begin to be lost (e.g., when the flow is perturbed by the impact). 

 \subsubsection{Convergence tests} \label{3D_convergence}

Consider a barotropic fluid model with properties $\rho_0= 997 \, \mathrm{kg/m}^2$, $\gamma =6 $, $A = \tilde{A} \rho_0^{-\gamma}$ with $\tilde A = 3.041\times 10^4$ Pa, and $B = 3.0397\times 10^4$ Pa. The size of the discrete reference configuration at time $t^0$ is $0.4 \, \mathrm{m} \times 0.4 \, \mathrm{m} \times 0.4 \, \mathrm{m}$. We consider the \textit{explicit} integrator to study the convergence with respect to $\Delta t$ and $\Delta s_i$, $i=1,2,3$.

\paragraph{Barotropic fluid flowing freely over a surface.} Given a fixed mesh, with values $\Delta s_1=\Delta s_2= \Delta s_3= 0.1$m, we impose the gravity and one impenetrability constraint. We vary the time-steps as $\Delta t \in \{ 2.5 \times 10^{-4}, \, 1.25 \times 10^{-4}, \, 6.25 \times 10^{-5}, \, 3.125 \times 10^{-5} \}$.
We compute the $L^2$-errors in the position $\varphi_d$ at time $t^N=0.5$s, by comparing $\varphi_d$ with an ``exact solution'' obtained with the time-step $\Delta t_{\rm ref}=7.8 \times 10^{-6}s$.  That is, for each value of $\Delta t$ we calculate
\begin{equation}\label{L_2norm_3D}
\| \varphi_d - \varphi_{\rm ref} \|_{L^2} = \left( \sum_a \sum_b \sum_c \| \varphi_{a,b,c}^N - \varphi_{{\rm ref};a,b,c}^N \|^2 \right)^{1/2}.
\end{equation}
This yields the following convergence with respect to $\Delta t$ 
\begin{figure}[H] \centering 
\begin{tabular}{| c | c | c | c | c |}
\hline
$\Delta t$ & $2.5 \times 10^{-4}$ & $1.25 \times 10^{-4}$ & $ 6.25 \times 10^{-5}$ & $3.125 \times 10^{-5}$  \\
\hline
$\| \varphi_d - \varphi_{\rm ref} \|_{L^2}$ &  $9 \times 10^{-3}$  & $4.4\times 10^{-3}$  & $2.1\times 10^{-3}$  &  $9.7 \times 10^{-4}$ \\
\hline
$ \text{rate} $  &    & 1.03 & 1.07  & 1.11  \\
\hline
\end{tabular}
\end{figure}

Given a fixed time-step $\Delta t=3.125 \times 10^{-5}$, we vary the space-steps as $\Delta s_1=\Delta s_2=\Delta s_3$ $\in \{0.2, \, 0.1,\, 0.05,\,0.025 \}$. The ``exact solution'' is chosen with $\Delta s_{1;\rm ref}=\Delta s_{2;\rm ref}  =\Delta s_{3;\rm ref} = 0.0125$m. We compute the $L^2$-errors in the position $\varphi_d$ at time $t^N=0.1$s. We get the following convergence with respect to $\Delta s_1=\Delta s_2=\Delta s_3$
\begin{figure}[H] \centering 
\begin{tabular}{| c | c | c | c |c |c|c|}
\hline
$\Delta s_1=\Delta s_2 =\Delta s_3 $ & $0.2$  & $0.1$  & $0.05$    & $0.025$\\
\hline
$\| \varphi_d - \varphi_{\rm ref} \|_{L^2}$ & $0.0753$   &  $0.0532 $ &  $0.0269$ &  $0.0120$ \\
\hline
$ \text{rate} $  &  & 0.5  & 0.98 & 1.16 \\
\hline
\end{tabular}
\end{figure}

\medskip

\begin{figure}[H] \centering 
\includegraphics[width=2.7 in]{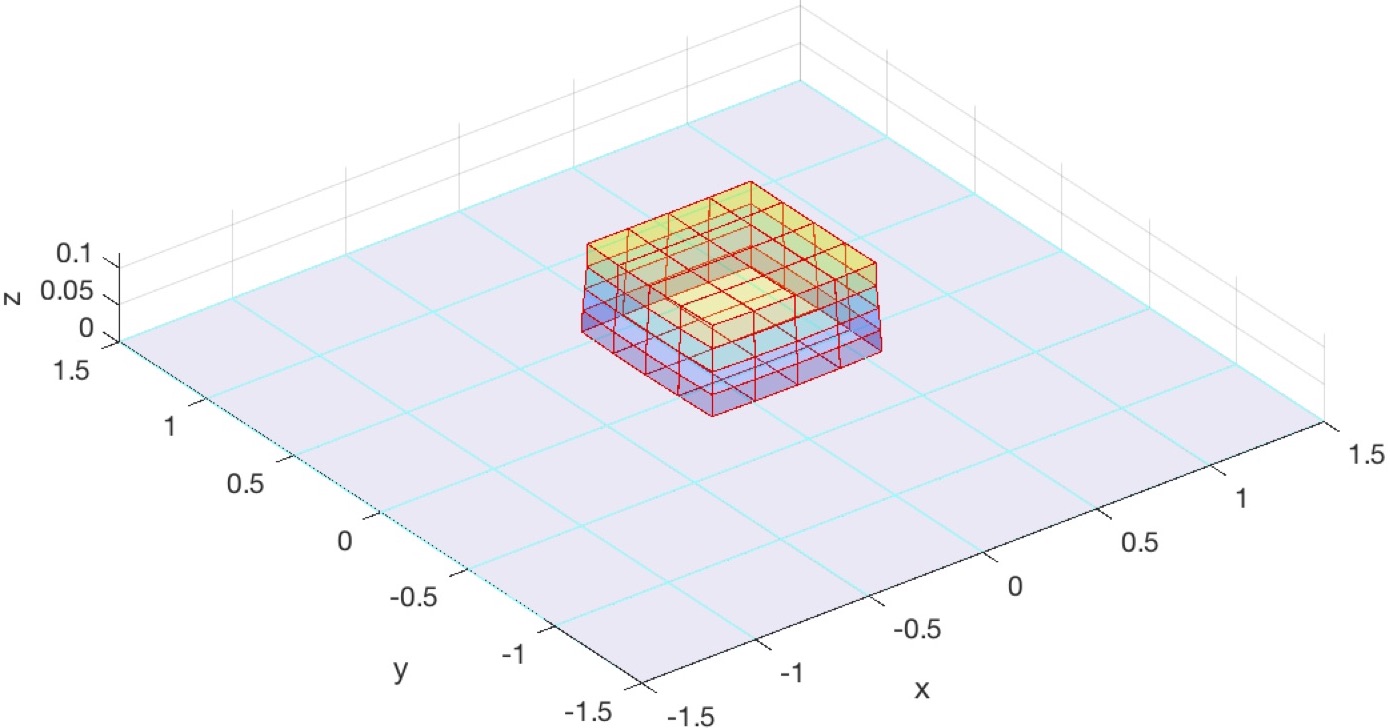} \vspace{-1pt} 
 \includegraphics[width=2.7 in]{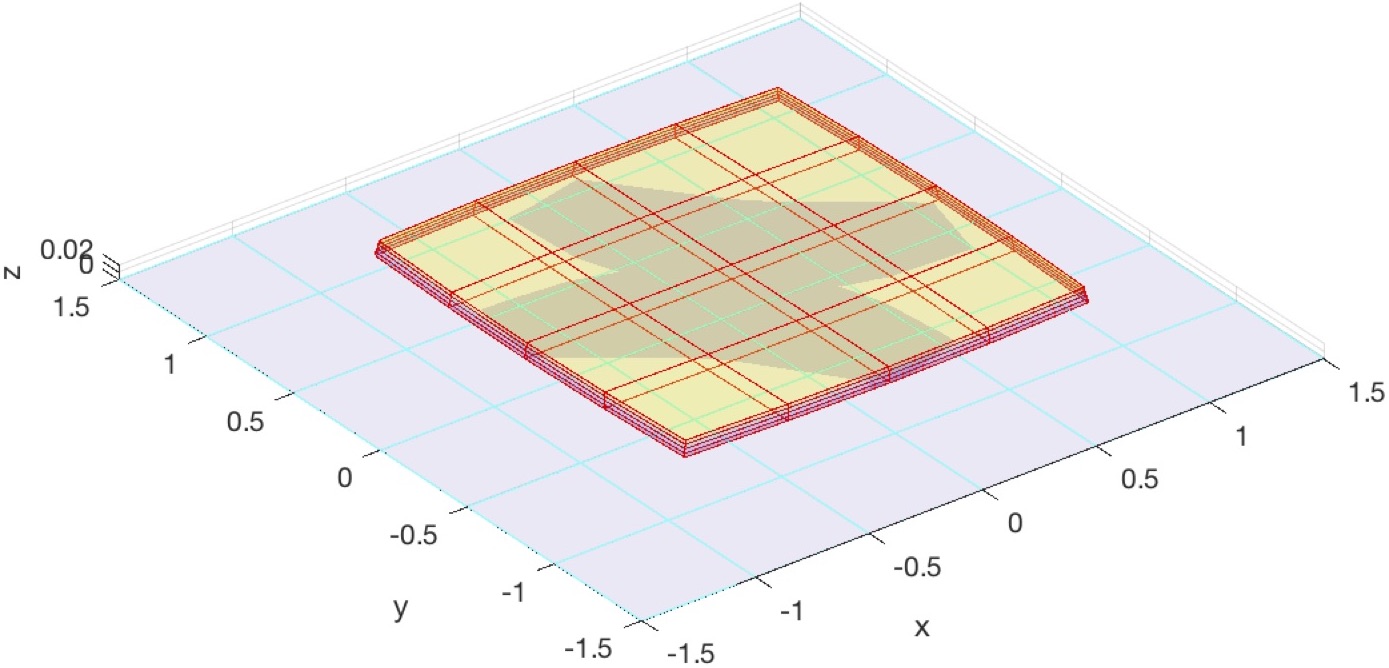} \vspace{-1pt} 
 \\
 \includegraphics[width=2.7 in]{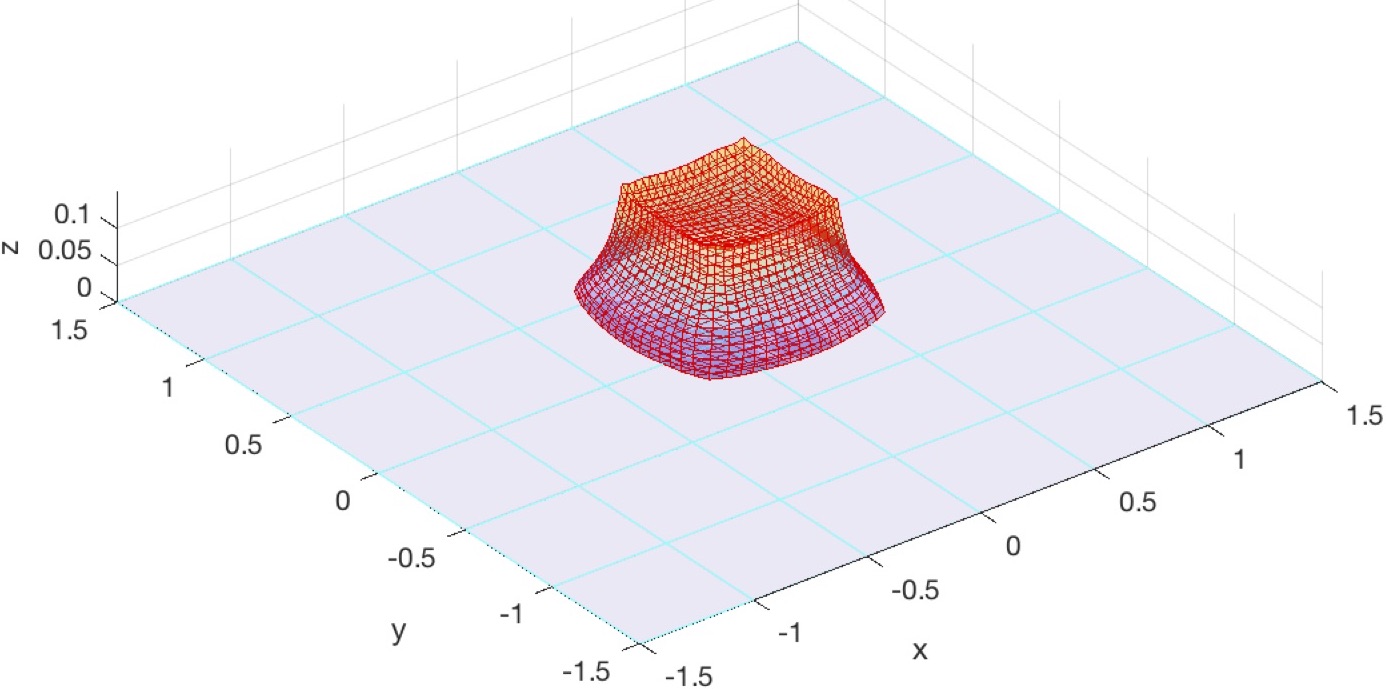} \vspace{-1pt} 
 \includegraphics[width=2.9 in]{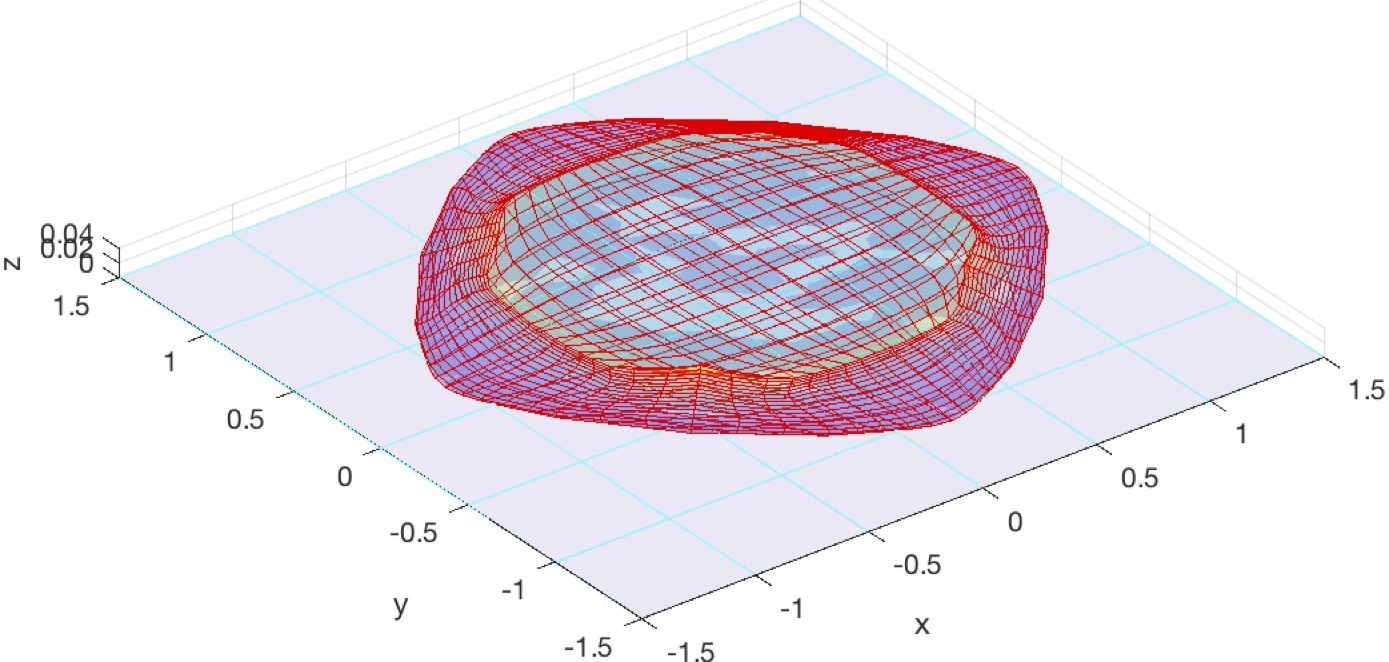} \vspace{-1pt} 
 \caption{\footnotesize Barotropic fluid flowing freely over a surface with $\Delta t=3.125 \times 10^{-5}$. \textit{From left to right}: after $0.25$s and $0.5$s. \textit{From top to bottom}: with $\Delta s_i=0.1$m and $\Delta s_i=0.025$m.} \label{convergence_freely_flowing_3D} 
\end{figure}

An illustration of the test used for the numerical convergence is given in Fig.\,\ref{convergence_freely_flowing_3D}.

\begin{figure}[H] \centering 
  \includegraphics[width=1.4 in]{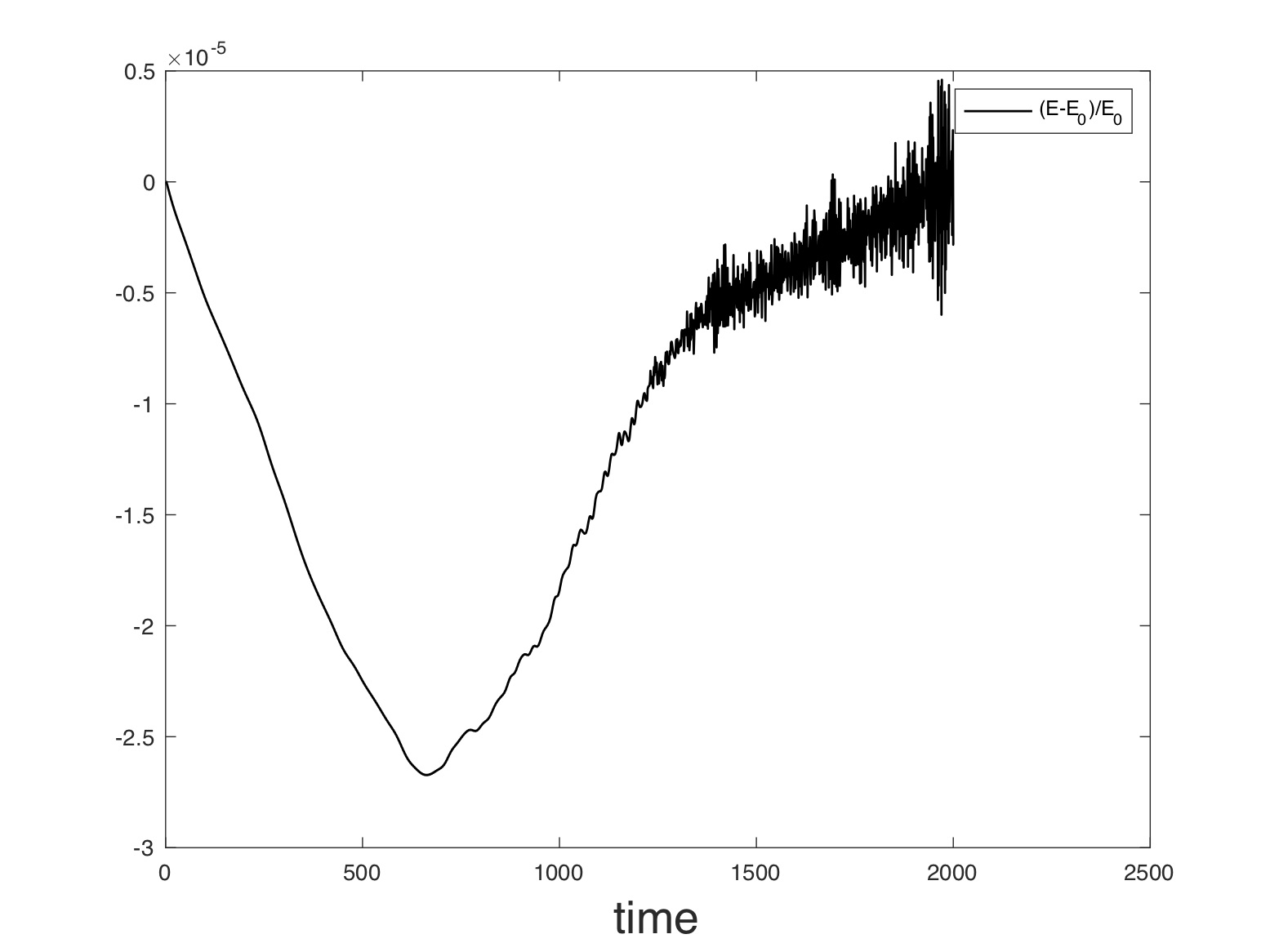} \vspace{-3pt}\qquad 
    \includegraphics[width=1.4 in]{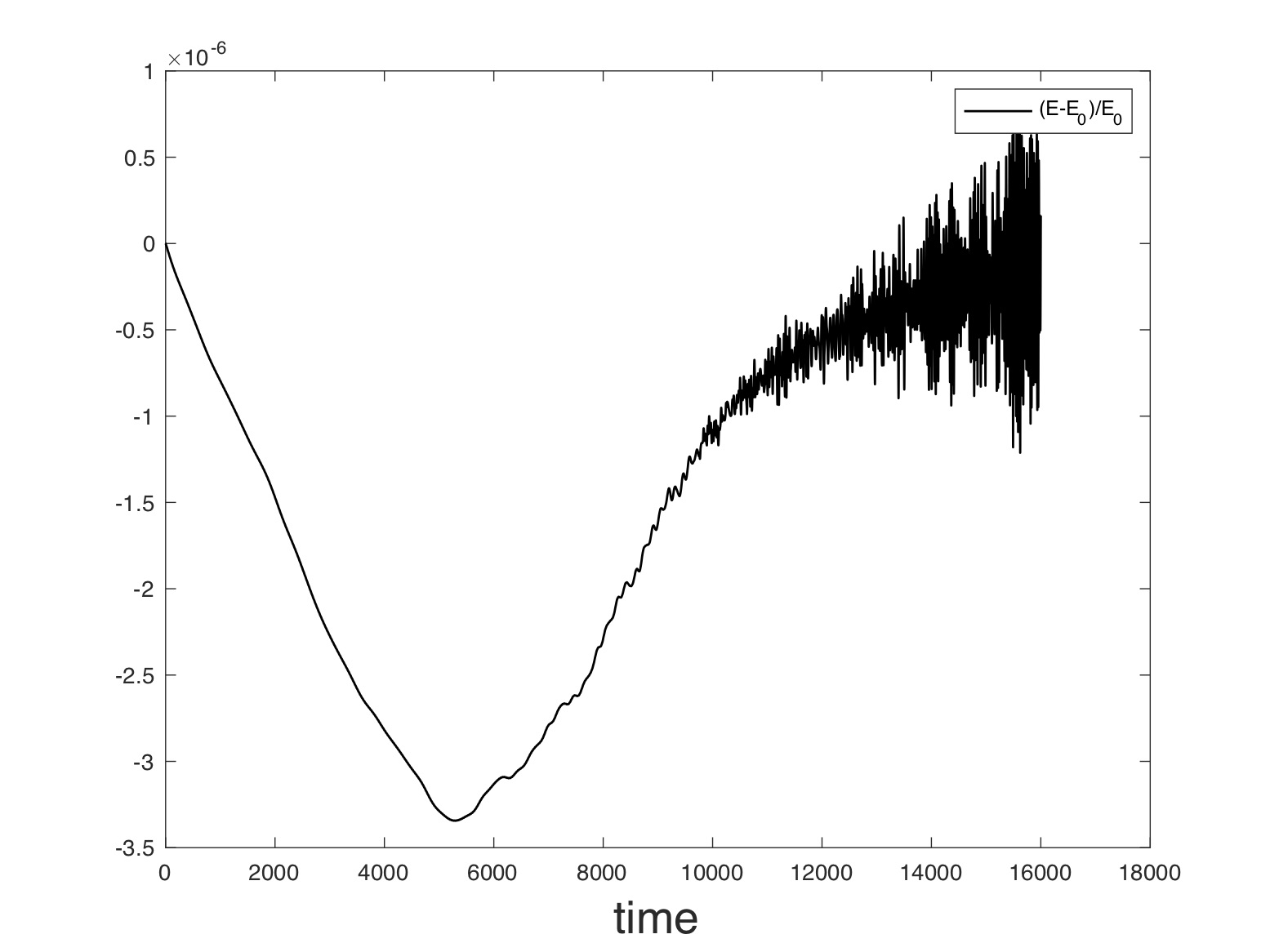} \vspace{-3pt}\qquad   \includegraphics[width=1.4 in]{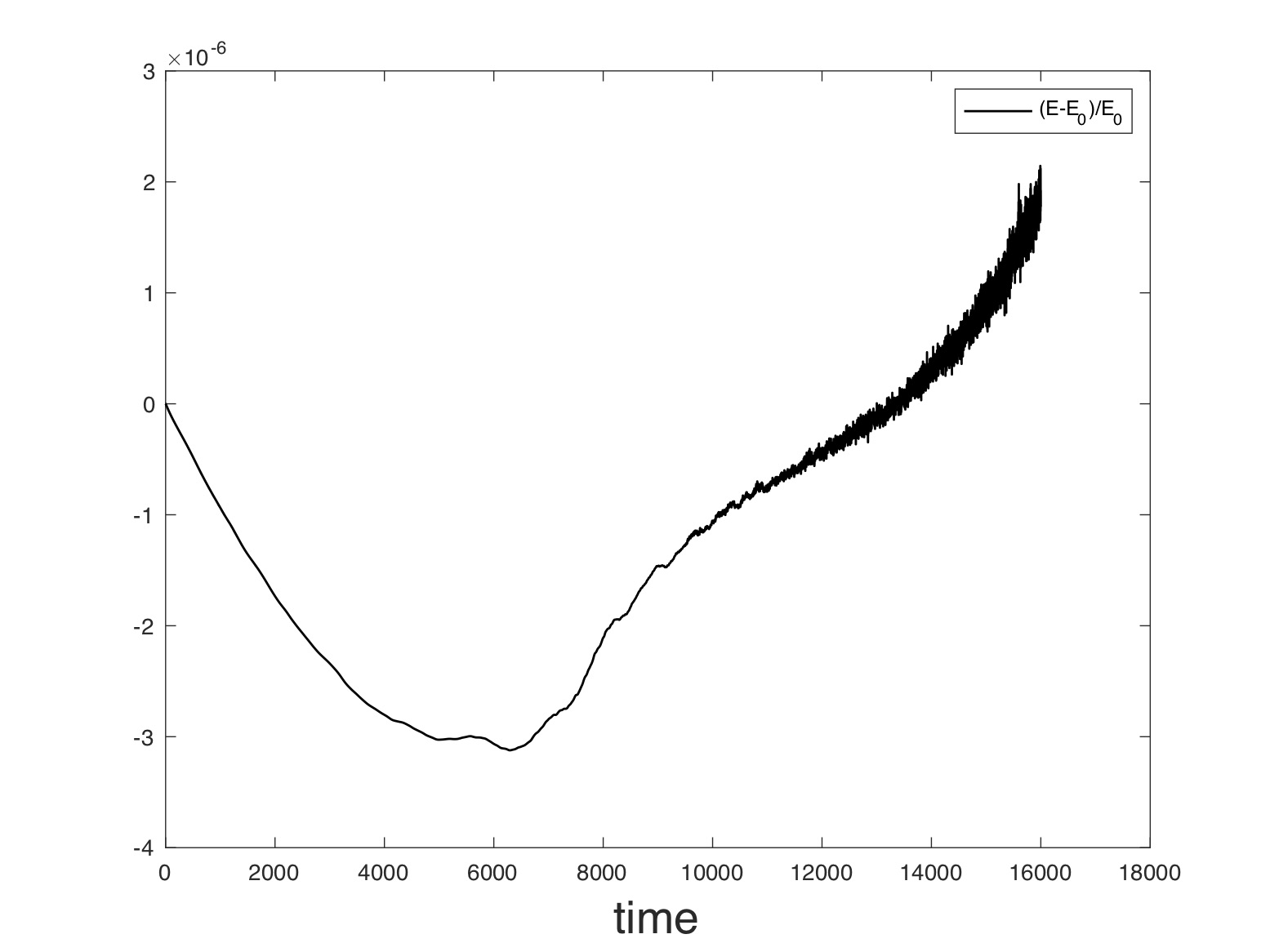} \vspace{-3pt} 
 \caption{\footnotesize Relative error in the energy $(ET^j-ET^0)/ET^0$. \textit{From left to right}: $(\Delta s_i =0.1,\Delta t = 2.5 \times 10^{-4})$, $(\Delta s_i =0.1,\Delta t =3.125 \times 10^{-5})$, $(\Delta t =3.125 \times 10^{-5}, \Delta s_i=0.025)$.} \label{relative_error_energy_3D} 
\end{figure}

\section{Concluding remarks and future directions} \label{conclusion}

This paper has presented new Lagrangian schemes for the regular motion of barotropic and incompressible fluid models, which preserve the momenta associated to symmetries, up to machine precision, and satisfy the nearly constant energy property of symplectic integrators, see Fig.\,\ref{water_perturb_sym_energ} and \ref{water_3D_energy_momentum}. The schemes are derived by discretization of the geometric and variational structures underlying the spacetime formulation of continuum mechanics, seen as a particular instance of field theory. We have illustrated how this approach can naturally accommodate incompressibility and fluid impact against an obstacle by appropriately augmenting the discrete Lagrangian and thanks to our definition of the discrete Jacobian.

\medskip

An important task for the future is the study of contacts, and of friction with heat exchange, between liquid and gas, which are important phenomena encountered for instance in ocean-atmosphere coupling.

\medskip 

Thanks to the clear link between the expressions of the discrete Jacobian and discrete gradient deformation in \S\ref{DJ2D} \& \S\ref{DJ3D}, see also \cite{Demoures2019}, it is possible to treat the coupling of fluid and elasticity dynamics, based on appropriate variational formulations such as, for instance, those developed in \cite{FaGBPu2020,GBPu2020} for fluid-structure interaction.

\medskip 

The results in this paper also make it possible to further develop multisymplectic integrators to model flows interacting mightily with obstacles, which are problems commonly met in engineering and biological applications, for example, when rocks fall into the reservoir of a dam or with blood flow in arteries which impacts heart valves.

\appendix 

\section{Appendix}\label{appendix}

\subsection{Derivatives of the discrete Lagrangian for 2D barotropic fluid}\label{water_DCEL}

\paragraph{Explicit integrator:} The partial derivatives of the discrete Lagrangian for a 2D barotropic fluid with internal energy $W( \rho  _0, J)$ are listed below, where the discrete pressure at time $t^j$ and spatial positon $\ell$, with the ordering $\ell=1$ to $\ell=4$ respectively associated to the nodes $(j,a,b)$, $(j,a+1,b)$, $(j,a,b+1)$, $(j,a+1,b+1)$, is defined as
\[
P_\ell ( \mbox{\mancube}_{a,b}^j)= - \rho_0 \frac{\partial W(  \mbox{\mancube}_{a,b}^j) }{\partial J_\ell(  \mbox{\mancube}_{a,b}^j)} = P_W\big( \rho_0  , J_\ell( \mbox{\mancube}_{a,b}^j)\big).
\] 
For the special case of internal energy given in \eqref{energy_water}, we have
\begin{equation}\label{pressure_2D}
P_\ell ( \mbox{\mancube}_{a,b}^j) =A  \Big(  \frac{ \rho  _0}{J_\ell ( \mbox{\mancube}_{a,b}^j)} \Big) ^\gamma -  B.
\end{equation}
{\footnotesize
\begin{align*}
A_{a,b}^j  = &  - \frac{M}{4} v_{a,b}^j + \frac{ \Delta t}{4 } P_1( \mbox{\mancube}_{a,b}^j) \left(\varphi_{a+1,b}^j - \varphi_{a,b+1}^j \right)\times \frac{\mathbf{n}_1( \mbox{\mancube}_{a,b}^j)}{|\mathbf{n}_1( \mbox{\mancube}_{a,b}^j) |}  
\\
 & +\frac{ \Delta t}{4 }P_2 ( \mbox{\mancube}_{a,b}^j)\left(\varphi_{a+1,b}^j - \varphi_{a+1,b+1}^j \right)\times \frac{\mathbf{n}_2( \mbox{\mancube}_{a,b}^j)}{|\mathbf{n}_2( \mbox{\mancube}_{a,b}^j) |} 
\\
 & +\frac{ \Delta t}{4 } P_3 ( \mbox{\mancube}_{a,b}^j) \left(\varphi_{a+1,b+1}^j - \varphi_{a,b+1}^j \right)\times \frac{\mathbf{n}_3( \mbox{\mancube}_{a,b}^j)}{|\mathbf{n}_3( \mbox{\mancube}_{a,b}^j) |} 
 \end{align*}
\begin{align*}
 B_{a,b}^j = & - \frac{M}{4} v_{a+1,b}^j + \frac{ \Delta t}{4 }P_1 ( \mbox{\mancube}_{a,b}^j) \left(\varphi_{a,b+1}^j - \varphi_{a,b}^j \right)\times \frac{\mathbf{n}_1( \mbox{\mancube}_{a,b}^j)}{|\mathbf{n}_1( \mbox{\mancube}_{a,b}^j) |} 
  \\
 & +\frac{ \Delta t}{4 } P_2 ( \mbox{\mancube}_{a,b}^j)   \left(\varphi_{a+1,b+1}^j - \varphi_{a,b}^j \right)\times \frac{\mathbf{n}_2( \mbox{\mancube}_{a,b}^j)}{|\mathbf{n}_2( \mbox{\mancube}_{a,b}^j) |} 
\\
 & + \frac{ \Delta t}{4 } P_4 ( \mbox{\mancube}_{a,b}^j)   \left(\varphi_{a+1,b+1}^j - \varphi_{a,b+1}^j \right)\times \frac{\mathbf{n}_4( \mbox{\mancube}_{a,b}^j)}{|\mathbf{n}_4( \mbox{\mancube}_{a,b}^j) |} 
\end{align*}
\begin{align*}
 C_{a,b}^j = & - \frac{M}{4} v_{a,b+1}^j +\frac{ \Delta t}{4 } P_1 ( \mbox{\mancube}_{a,b}^j) \left(\varphi_{a,b}^j - \varphi_{a+1,b}^j \right)\times \frac{\mathbf{n}_1( \mbox{\mancube}_{a,b}^j)}{|\mathbf{n}_1( \mbox{\mancube}_{a,b}^j) |} 
\\
 & + \frac{ \Delta t}{4 } P_3 ( \mbox{\mancube}_{a,b}^j)  \left(\varphi_{a,b}^j - \varphi_{a+1,b+1}^j \right)\times \frac{\mathbf{n}_3( \mbox{\mancube}_{a,b}^j)}{|\mathbf{n}_3( \mbox{\mancube}_{a,b}^j) |} 
 \\
 & +\frac{ \Delta t}{4 } P_4 ( \mbox{\mancube}_{a,b}^j) \left(\varphi_{a+1,b}^j - \varphi_{a+1,b+1}^j \right)\times \frac{\mathbf{n}_4( \mbox{\mancube}_{a,b}^j)}{|\mathbf{n}_4( \mbox{\mancube}_{a,b}^j) |}
\end{align*} 
\begin{align*} 
D_{a,b}^j  =& - \frac{M}{4} v_{a+1,b+1}^j  +\frac{ \Delta t}{4 }P_2 ( \mbox{\mancube}_{a,b}^j)   \left(\varphi_{a,b}^j - \varphi_{a+1,b}^j \right)\times \frac{\mathbf{n}_2( \mbox{\mancube}_{a,b}^j)}{|\mathbf{n}_2( \mbox{\mancube}_{a,b}^j) |} 
\\
& +\frac{ \Delta t}{4 } P_3 ( \mbox{\mancube}_{a,b}^j)  \left(\varphi_{a,b+1}^j - \varphi_{a,b}^j \right)\times \frac{\mathbf{n}_3( \mbox{\mancube}_{a,b}^j)}{|\mathbf{n}_3( \mbox{\mancube}_{a,b}^j) |} 
\\
& + \frac{ \Delta t}{4 } P_4 ( \mbox{\mancube}_{a,b}^j)   \left(\varphi_{a,b+1}^j - \varphi_{a+1,b}^j \right)\times \frac{\mathbf{n}_4( \mbox{\mancube}_{a,b}^j)}{|\mathbf{n}_4( \mbox{\mancube}_{a,b}^j) |} 
\end{align*}}
with
{\footnotesize
\begin{align*}
\mathbf{n}_1( \mbox{\mancube}_{a,b}^j)  &= (\varphi_{a+1,b}^j - \varphi_{a,b}^j) \times (\varphi_{a,b+1}^j- \varphi_{a,b}^j); 
\\
\mathbf{n}_2( \mbox{\mancube}_{a,b}^j) &=(\varphi_{a+1,b+1} - \varphi_{a+1,b}^j) \times (\varphi_{a,b}^j- \varphi_{a+1,b}^j). \quad 
\\
\mathbf{n}_3( \mbox{\mancube}_{a,b}^j)  &=(\varphi_{a,b} - \varphi_{a,b+1}^j) \times (\varphi_{a+1,b+1}^j- \varphi_{a,b+1}^j); \quad 
\\
\mathbf{n}_4( \mbox{\mancube}_{a,b}^j)  &= (\varphi_{a,b+1}^j - \varphi_{a+1,b+1}^j) \times (\varphi_{a+1,b}^j- \varphi_{a+1,b+1}^j).  
\end{align*}}

\paragraph{Implicit integrator:} The partial derivatives of the Lagrangian for a 2D barotropic fluid with internal energy $W( \rho  _0, J)$, discretized under the mid-point rules, are given by

{\footnotesize
\begin{align*}
\mathbb{A}_{a,b}^j  = &   \frac{ \Delta t}{4^2 } P_1( \mbox{\mancube}_{a,b}^j) \left((\varphi_{a+1,b}^j+ \varphi_{a+1,b}^{j+1})  - (\varphi_{a,b+1}^j +\varphi_{a,b+1}^{j+1}) \right)\times \frac{\mathbb{n}_1( \mbox{\mancube}_{a,b}^j)}{|\mathbb{n}_1( \mbox{\mancube}_{a,b}^j) |}  
\\
 & +\frac{ \Delta t}{4^2 }P_2 ( \mbox{\mancube}_{a,b}^j)\left( (\varphi_{a+1,b}^j + \varphi_{a+1,b}^{j+1}) - (\varphi_{a+1,b+1}^j + \varphi_{a+1,b+1}^{j+1}) \right)\times \frac{\mathbb{n}_2( \mbox{\mancube}_{a,b}^j)}{|\mathbb{n}_2( \mbox{\mancube}_{a,b}^j) |} 
\\
 & +\frac{ \Delta t}{4^2 } P_3 ( \mbox{\mancube}_{a,b}^j) \left( (\varphi_{a+1,b+1}^j + \varphi_{a+1,b+1}^{j+1}) - (\varphi_{a,b+1}^j + \varphi_{a,b+1}^{j+1}) \right)\times \frac{\mathbb{n}_3( \mbox{\mancube}_{a,b}^j)}{|\mathbb{n}_3( \mbox{\mancube}_{a,b}^j) |} 
 \end{align*}
\begin{align*}
 \mathbb{B}_{a,b}^j = &  \frac{ \Delta t}{4^2 }P_1 ( \mbox{\mancube}_{a,b}^j) \left( (\varphi_{a,b+1}^j + \varphi_{a,b+1}^{j+1}) - (\varphi_{a,b}^j + \varphi_{a,b}^{j+1}) \right)\times \frac{\mathbb{n}_1( \mbox{\mancube}_{a,b}^j)}{|\mathbb{n}_1( \mbox{\mancube}_{a,b}^j) |} 
  \\
 & +\frac{ \Delta t}{4^2 } P_2 ( \mbox{\mancube}_{a,b}^j)   \left( (\varphi_{a+1,b+1}^j + \varphi_{a+1,b+1}^{j+1}) - (\varphi_{a,b}^j + \varphi_{a,b}^{j+1}) \right)\times \frac{\mathbb{n}_2( \mbox{\mancube}_{a,b}^j)}{|\mathbb{n}_2( \mbox{\mancube}_{a,b}^j) |} 
\\
 & + \frac{ \Delta t}{4^2 } P_4 ( \mbox{\mancube}_{a,b}^j)   \left((\varphi_{a+1,b+1}^j + \varphi_{a+1,b+1}^{j+1}) - (\varphi_{a,b+1}^j +\varphi_{a,b+1}^{j+1}) \right)\times \frac{\mathbb{n}_4( \mbox{\mancube}_{a,b}^j)}{|\mathbb{n}_4( \mbox{\mancube}_{a,b}^j) |} 
\end{align*}
\begin{align*}
 \mathbb{C}_{a,b}^j = & \frac{ \Delta t}{4^2 } P_1 ( \mbox{\mancube}_{a,b}^j) \left( (\varphi_{a,b}^j + \varphi_{a,b}^{j+1}) - (\varphi_{a+1,b}^j + \varphi_{a+1,b}^{j+1}) \right)\times \frac{\mathbb{n}_1( \mbox{\mancube}_{a,b}^j)}{|\mathbb{n}_1( \mbox{\mancube}_{a,b}^j) |} 
\\
 & + \frac{ \Delta t}{4^2 } P_3 ( \mbox{\mancube}_{a,b}^j)  \left((\varphi_{a,b}^j + \varphi_{a,b}^{j+1}) - (\varphi_{a+1,b+1}^j + \varphi_{a+1,b+1}^{j+1}) \right)\times \frac{\mathbb{n}_3( \mbox{\mancube}_{a,b}^j)}{|\mathbb{n}_3( \mbox{\mancube}_{a,b}^j) |} 
 \\
 & +\frac{ \Delta t}{4^2 } P_4 ( \mbox{\mancube}_{a,b}^j) \left((\varphi_{a+1,b}^j+\varphi_{a+1,b}^{j+1}) - (\varphi_{a+1,b+1}^j + \varphi_{a+1,b+1}^{j+1}) \right)\times \frac{\mathbb{n}_4( \mbox{\mancube}_{a,b}^j)}{|\mathbb{n}_4( \mbox{\mancube}_{a,b}^j) |}
\end{align*} 
\begin{align*} 
\mathbb{D}_{a,b}^j  =& \frac{ \Delta t}{4^2 }P_2 ( \mbox{\mancube}_{a,b}^j)   \left( (\varphi_{a,b}^j + \varphi_{a,b}^{j+1}) - (\varphi_{a+1,b}^j + \varphi_{a+1,b}^{j+1}) \right)\times \frac{\mathbb{n}_2( \mbox{\mancube}_{a,b}^j)}{|\mathbb{n}_2( \mbox{\mancube}_{a,b}^j) |} 
\\
& +\frac{ \Delta t}{4^2 } P_3 ( \mbox{\mancube}_{a,b}^j)  \left((\varphi_{a,b+1}^j + \varphi_{a,b+1}^{j+1}) - (\varphi_{a,b}^j + \varphi_{a,b}^{j+1}) \right)\times \frac{\mathbb{n}_3( \mbox{\mancube}_{a,b}^j)}{|\mathbb{n}_3( \mbox{\mancube}_{a,b}^j) |} 
\\
& + \frac{ \Delta t}{4^2 } P_4 ( \mbox{\mancube}_{a,b}^j)   \left( (\varphi_{a,b+1}^j + \varphi_{a,b+1}^{j+1}) - (\varphi_{a+1,b}^j+\varphi_{a+1,b}^{j+1}) \right)\times \frac{\mathbb{n}_4( \mbox{\mancube}_{a,b}^j)}{|\mathbb{n}_4( \mbox{\mancube}_{a,b}^j) |} 
\end{align*}}
with
{\footnotesize
\begin{align*}
\mathbb{n}_1( \mbox{\mancube}_{a,b}^j)  &= \left((\varphi_{a+1,b}^j + \varphi_{a+1,b}^{j+1}) - (\varphi_{a,b}^j + \varphi_{a,b}^{j+1}) \right) \times \left((\varphi_{a,b+1}^j + \varphi_{a,b+1}^{j+1})- ( \varphi_{a,b}^j + \varphi_{a,b}^{j+1}) \right); 
\\
\mathbb{n}_2( \mbox{\mancube}_{a,b}^j) &=\left((\varphi_{a+1,b+1}^j+\varphi_{a+1,b+1}^{j+1}) - (\varphi_{a+1,b}^j + \varphi_{a+1,b}^{j+1})\right) \times \left((\varphi_{a,b}^j + \varphi_{a,b}^{j+1})- (\varphi_{a+1,b}^j + \varphi_{a+1,b}^{j+1})\right). 
\\
\mathbb{n}_3( \mbox{\mancube}_{a,b}^j)  &= \left((\varphi_{a,b}^j + \varphi_{a,b}^{j+1}) - (\varphi_{a,b+1}^j + \varphi_{a,b+1}^{j+1})\right) \times \left((\varphi_{a+1,b+1}^j + \varphi_{a+1,b+1}^{j+1})- (\varphi_{a,b+1}^j + \varphi_{a,b+1}^{j+1}) \right); 
\\
\mathbb{n}_4( \mbox{\mancube}_{a,b}^j)  &= \left( (\varphi_{a,b+1}^j + \varphi_{a,b+1}^{j+1}) - (\varphi_{a+1,b+1}^j + \varphi_{a+1,b+1}^{j+1})\right) \times \left((\varphi_{a+1,b}^j + \varphi_{a+1,b}^{j+1})- (\varphi_{a+1,b+1}^j + \varphi_{a+1,b+1}^{j+1})\right). 
\end{align*}}
Where the pressures $P_\ell ( \mbox{\mancube}_{a,b}^j)$, $\ell=1,2,3,4$, were defined in \eqref{pressure_2D}.

\subsection{Discrete Euler-Lagrange equations for 2D barotropic fluid}\label{2D_DCEL_geo}
 
{\footnotesize
\begin{align*}
 \rho_0\left( \frac{ v_{a,b}^j - v_{a,b}^{j-1}}{\Delta t} \right) 
 = \frac{1}{4 \Delta s_1 \Delta s_2} &  \left\{ \color{red} - \left( P_1( \mbox{\mancube}_{a,b}^j) \left(\varphi_{a,b+1}^j - \varphi_{a+1,b}^j \right)\times \frac{\mathbf{n}_1( \mbox{\mancube}_{a,b}^j)}{|\mathbf{n}_1( \mbox{\mancube}_{a,b}^j) |}  
 \right. \right.
 \\
 &  \color{red} \qquad + P_2 ( \mbox{\mancube}_{a-1,b}^j)   \left(\varphi_{a-1,b}^j - \varphi_{a,b+1}^j  \right)\times \frac{\mathbf{n}_2( \mbox{\mancube}_{a-1,b}^j)}{|\mathbf{n}_2( \mbox{\mancube}_{a-1,b}^j) |} 
 \\
 &  \color{red} \qquad + P_3 ( \mbox{\mancube}_{a,b-1}^j)  \left(\varphi_{a+1,b}^j - \varphi_{a,b-1}^j \right)\times \frac{\mathbf{n}_3( \mbox{\mancube}_{a,b-1}^j)}{|\mathbf{n}_3( \mbox{\mancube}_{a,b-1}^j) |} 
 \\
 &  \color{red} \qquad \left. + P_4 ( \mbox{\mancube}_{a-1,b-1}^j)   \left( \varphi_{a,b-1}^j - \varphi_{a-1,b}^j \right)\times \frac{\mathbf{n}_4( \mbox{\mancube}_{a-1,b-1}^j)}{|\mathbf{n}_4( \mbox{\mancube}_{a-1,b-1}^j) |} \right) 
 \\
 &  \color{blue} + \left( P_2 ( \mbox{\mancube}_{a,b}^j)\left(\varphi_{a+1,b}^j - \varphi_{a+1,b+1}^j \right)\times \frac{\mathbf{n}_2( \mbox{\mancube}_{a,b}^j)}{|\mathbf{n}_2( \mbox{\mancube}_{a,b}^j) |} \right.
 \\
 &  \color{blue} \left. \qquad + P_4 ( \mbox{\mancube}_{a,b-1}^j) \left(\varphi_{a+1,b-1}^j - \varphi_{a+1,b}^j \right)\times \frac{\mathbf{n}_4( \mbox{\mancube}_{a,b-1}^j)}{|\mathbf{n}_4( \mbox{\mancube}_{a,b-1}^j) |} \right)
 \\
 &  \color{magenta} + \left( P_3 ( \mbox{\mancube}_{a,b}^j) \left(\varphi_{a+1,b+1}^j - \varphi_{a,b+1}^j \right)\times \frac{\mathbf{n}_3( \mbox{\mancube}_{a,b}^j)}{|\mathbf{n}_3( \mbox{\mancube}_{a,b}^j) |} \right.
 \\
 &  \color{magenta} \qquad +\left.  P_4 ( \mbox{\mancube}_{a-1,b}^j)   \left(\varphi_{a,b+1}^j - \varphi_{a-1,b+1}^j \right)\times \frac{\mathbf{n}_4( \mbox{\mancube}_{a-1,b}^j)}{|\mathbf{n}_4( \mbox{\mancube}_{a-1,b}^j) |} \right)
 \\
 &  \color{cyan} + \left( P_1 ( \mbox{\mancube}_{a-1,b}^j) \left(\varphi_{a-1,b+1}^j - \varphi_{a-1,b}^j \right)\times \frac{\mathbf{n}_1( \mbox{\mancube}_{a-1,b}^j)}{|\mathbf{n}_1( \mbox{\mancube}_{a-1,b}^j) |} \right.
 \\
 &  \color{cyan} \left. \qquad + P_3 ( \mbox{\mancube}_{a-1,b-1}^j)  \left(\varphi_{a-1,b}^j - \varphi_{a-1,b-1}^j \right)\times \frac{\mathbf{n}_3( \mbox{\mancube}_{a-1,b-1}^j)}{|\mathbf{n}_3( \mbox{\mancube}_{a-1,b-1}^j) |} \right)
 \\
 & + \left( P_2 ( \mbox{\mancube}_{a-1,b-1}^j)   \left(\varphi_{a-1,b-1}^j - \varphi_{a,b-1}^j \right)\times \frac{\mathbf{n}_2( \mbox{\mancube}_{a-1,b-1}^j)}{|\mathbf{n}_2( \mbox{\mancube}_{a-1,b-1}^j) |} \right.
 \\
 &\left. \qquad \left.+ P_1 ( \mbox{\mancube}_{a,b-1}^j) \left(\varphi_{a,b-1}^j - \varphi_{a+1,b-1}^j \right)\times \frac{\mathbf{n}_1( \mbox{\mancube}_{a,b-1}^j)}{|\mathbf{n}_1( \mbox{\mancube}_{a,b-1}^j) |} \right) \right\}
 \end{align*}
Interpretation of the discrete Euler-Lagrange equations: \textit{discrete balance of momentum}
\[
\rho_0 (v^j-v^{j-1})/\Delta t = - \left[(P_{\rm ext}  \cdot   L_{\rm ext})\mathbf{n}_{\rm ext} - (P_{\rm int}   \cdot  L_{\rm int})\mathbf{n}_{\rm int}\right]_j / \textrm{area} ,
\]
where $\mathbf{n}_{\rm ext}$, $\mathbf{n}_{\rm int}$ are unit vectors that point outward of the boundaries.}
 
 \begin{figure}[H] \centering 
 \includegraphics[width=3.8 in]{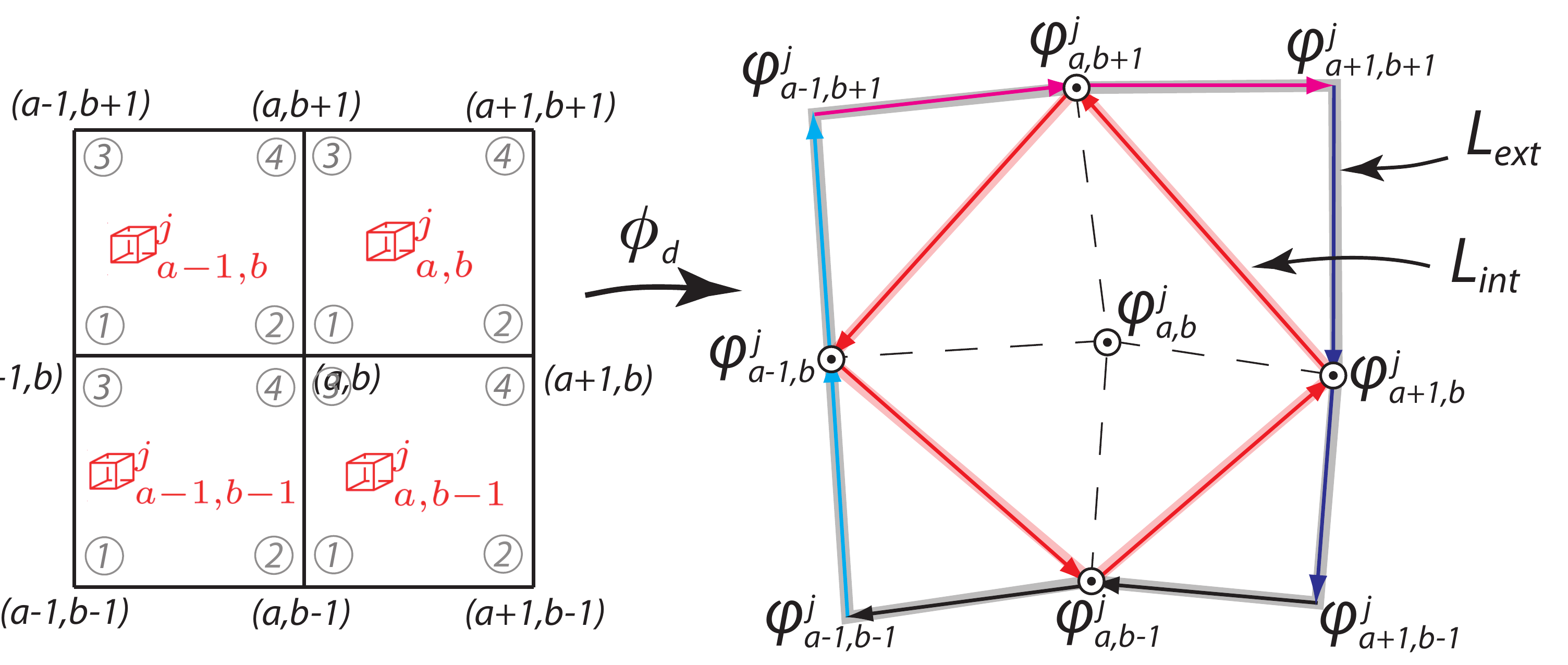} \vspace{-3pt}  
 \caption{\footnotesize Internal and external lengths $L_{\rm int}$, $ L_{\rm ext}$.
} \label{2D_DCEL_fig} 
 \end{figure} 

\begin{figure}[H] \centering 
\includegraphics[width=2 in]{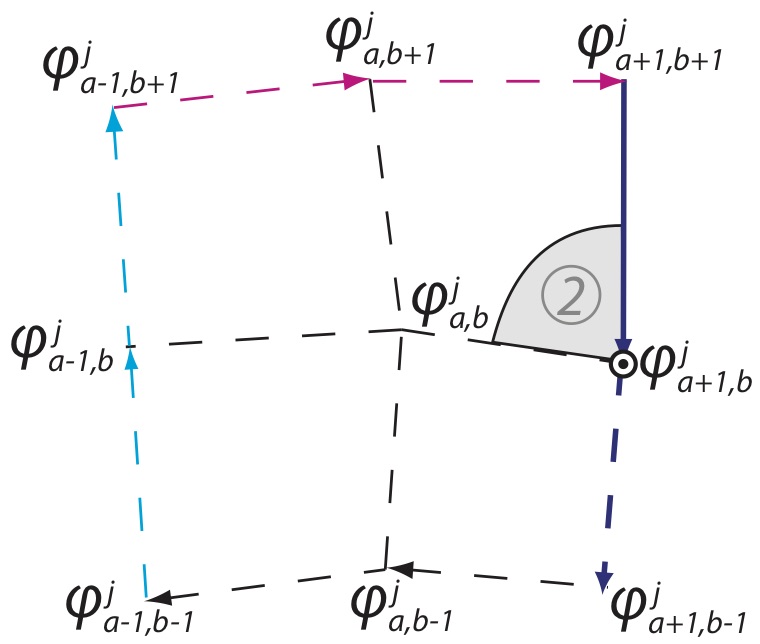} \qquad  \includegraphics[width=2 in]{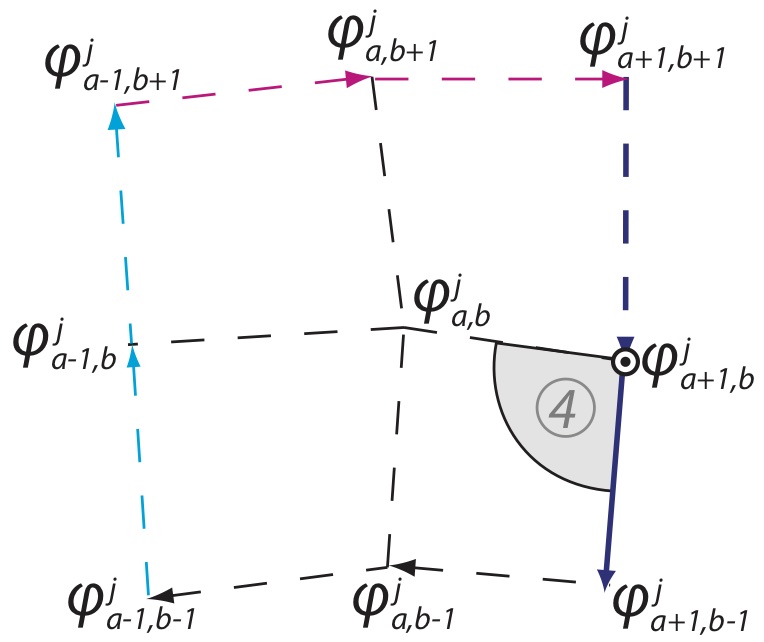} \vspace{-3pt}  
\caption{\footnotesize On the left: the force  \textcolor{blue}{$P_2 ( \mbox{\mancube}_{a,b}^j)\left(\varphi_{a+1,b}^j - \varphi_{a+1,b+1}^j \right)\times \frac{\mathbf{n}_2( \mbox{\mancube}_{a,b}^j)}{|\mathbf{n}_2( \mbox{\mancube}_{a,b}^j) |}$} associated to the pressure $P_2( \mbox{\mancube}_{a,b}^j)$ and the length $| \varphi_{a+1,b}^j - \varphi_{a+1,b+1}^j  |$. 
On the right \textcolor{blue}{$P_4 ( \mbox{\mancube}_{a,b-1}^j) \left(\varphi_{a+1,b-1}^j - \varphi_{a+1,b}^j \right)\times \frac{\mathbf{n}_4( \mbox{\mancube}_{a,b-1}^j)}{|\mathbf{n}_4( \mbox{\mancube}_{a,b-1}^j) |}$} associated to the pressure $P_4 ( \mbox{\mancube}_{a,b-1}^j) $ and the length $| \varphi_{a+1,b-1}^j - \varphi_{a+1,b}^j |$.} \label{DCEL_detail12} 
\end{figure}  

\begin{figure}[H] \centering 
\includegraphics[width=2 in]{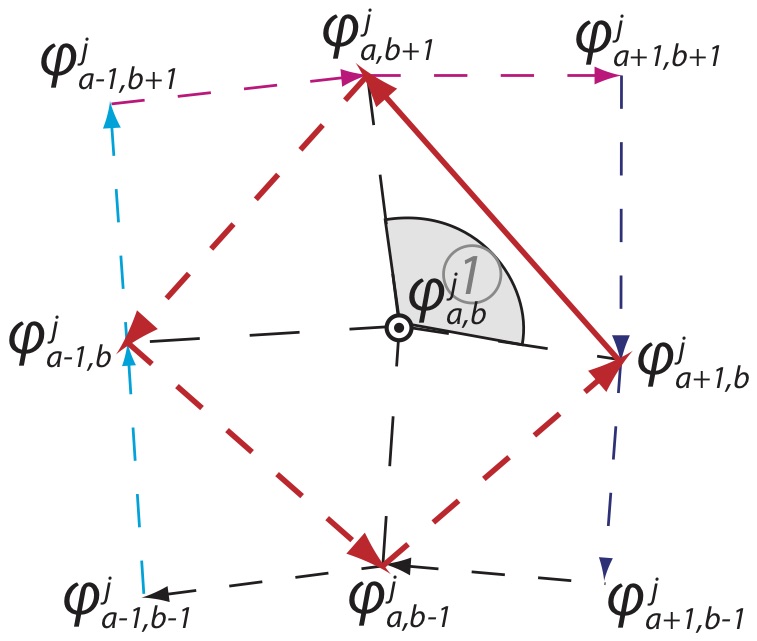} \qquad  \includegraphics[width=2 in]{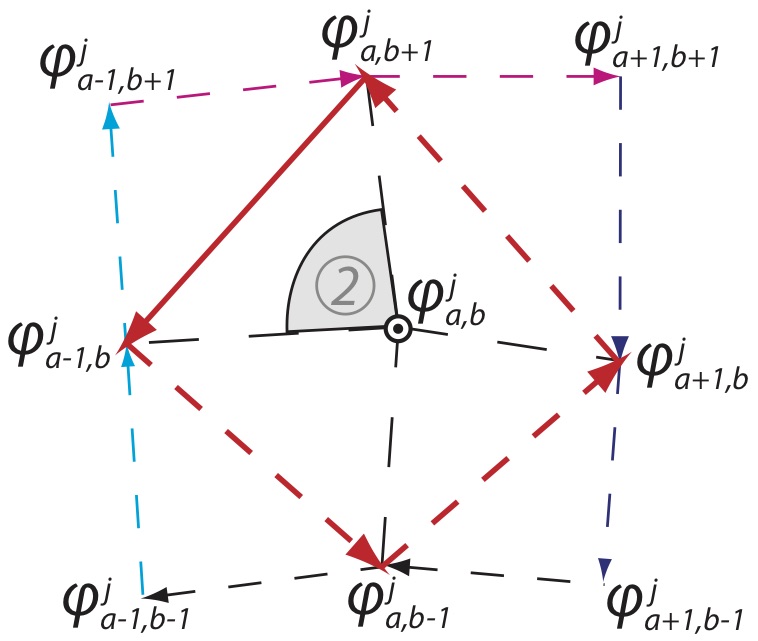} \vspace{-3pt}  
\caption{\footnotesize On the left: the force  \textcolor{red}{$P_1( \mbox{\mancube}_{a,b}^j) \left(\varphi_{a,b+1}^j - \varphi_{a+1,b}^j \right)\times \frac{\mathbf{n}_1( \mbox{\mancube}_{a,b}^j)}{|\mathbf{n}_1( \mbox{\mancube}_{a,b}^j) |} $} associated to the pressure $P_1( \mbox{\mancube}_{a,b}^j)$ and the length $| \varphi_{a,b+1}^j - \varphi_{a+1,b}^j |$. 
 On the right \textcolor{red}{$P_2 ( \mbox{\mancube}_{a-1,b}^j)   \left(\varphi_{a-1,b}^j - \varphi_{a,b+1}^j  \right)\times \frac{\mathbf{n}_2( \mbox{\mancube}_{a-1,b}^j)}{|\mathbf{n}_2( \mbox{\mancube}_{a-1,b}^j) |}$} associated to the pressure $P_2 ( \mbox{\mancube}_{a-1,b}^j) $ and the length $|\varphi_{a-1,b}^j - \varphi_{a,b+1}^j |$.} \label{DCEL_detail12} 
\end{figure}

\subsection{Discrete Jacobian on $\mbox{\mancube}_{a,b,c}^j$ in 3D} \label{3D_Jacobian's}

{\footnotesize
\begin{align*}
 J_2(\mbox{\mancube}_{a,b,c}^j) &=   (\mathbf{F}_{2;a+1,b,c}^j \times \mathbf{F}_{4;a+1,b,c}^j)\cdot \mathbf{F}_{3;a+1,b,c}^j  , \quad J_5(\mbox{\mancube}_{a,b,c}^j) =  (\mathbf{F}_{4;a+1,b+1,c}^j \times \mathbf{F}_{5;a+1,b+1,c}^j)\cdot \mathbf{F}_{3;a+1,b+1,c}^j  
 \\
 J_3(\mbox{\mancube}_{a,b,c}^j) &=   (\mathbf{F}_{5;a,b+1,c}^j \times \mathbf{F}_{1;a,b+1,c}^j)\cdot \mathbf{F}_{3;a,b+1,c}^j  ,  \quad J_6(\mbox{\mancube}_{a,b,c}^j) = (\mathbf{F}_{1;a,b+1,c+1}^j \times \mathbf{F}_{5;a,b+1,c+1}^j)\cdot \mathbf{F}_{6;a,b+1,c+1}^j
 \\
 J_4(\mbox{\mancube}_{a,b,c}^j) &= (\mathbf{F}_{2;a,b,c+1}^j \times \mathbf{F}_{1;a,b,c+1}^j)\cdot \mathbf{F}_{6;a,b,c+1}^j , \quad 
 J_7(\mbox{\mancube}_{a,b,c}^j) = (\mathbf{F}_{4;a+1,b,c+1}^j \times \mathbf{F}_{2;a+1,b,c+1}^j)\cdot \mathbf{F}_{6;a+1,b,c+1}^j  
 \\
  J_8(\mbox{\mancube}_{a,b,c}^j) & = (\mathbf{F}_{5;a+1,b+1,c+1}^j \times \mathbf{F}_{4;a+1,b+1,c+1}^j)\cdot \mathbf{F}_{6;a+1,b+1,c+1}^j .
 \end{align*}}

\subsection{Derivatives of the discrete Lagrangian for 3D barotropic fluid}\label{3D_water}

With the same definition as in \S\ref{water_DCEL}, the partial derivatives of the discrete Lagrangian for a 3D barotropic fluid with internal energy $W( \rho  _0, J)$ are listed below.

{\footnotesize
\begin{align*}
A_{a,b,c}^j = & - \frac{M}{8} v_{a,b,c}^j + \frac{\Delta t }{8} \left\{ -P_1  \left( \mathbf{n}_1^{(12)} + \mathbf{n}_1^{(23)} + \mathbf{n}_1^{(31)} \right) + P_4 \mathbf{n}_4^{(12)} +   P_3 \mathbf{n}_3^{(23)} 
 +  P_{ 2 } \mathbf{n}_2^{(31)} \right\}( \mbox{\mancube}_{a,b,c}^j),
\\
B_{a,b,c}^j = & - \frac{M}{8} v_{a+1,b,c}^j  + \frac{\Delta t  }{8} \left\{ -P_2 \left( \mathbf{n}_2^{(12)} + \mathbf{n}_2^{(23)} + \mathbf{n}_2^{(31)} \right)   + P_7 \mathbf{n}_7^{(12)} +  P_1  \mathbf{n}_1^{(23)} 
+ P_5  \mathbf{n}_5^{(31)} \right\}( \mbox{\mancube}_{a,b,c}^j) ,
\\
C_{a,b,c}^j = & - \frac{M}{8} v_{a,b+1,c}^j  + \frac{\Delta t  }{8} \left\{ -P_3 \left( \mathbf{n}_3^{(12)} + \mathbf{n}_3^{(23)} + \mathbf{n}_3^{(31)} \right)  +  P_6  \mathbf{n}_6^{(12)} + P_5  \mathbf{n}_5^{(23)} +  P_1 \mathbf{n}_1^{(31)} \right\}( \mbox{\mancube}_{a,b,c}^j) ,
\\
D_{a,b,c}^j = & - \frac{M}{8} v_{a,b,c+1}^j  + \frac{\Delta t  }{8} \left\{ -P_4  \left( \mathbf{n}_4^{(12)} + \mathbf{n}_4^{(23)} + \mathbf{n}_4^{(31)} \right)  + P_1  \mathbf{n}_1^{(12)} +   P_7  \mathbf{n}_7^{(23)} +   P_6 \mathbf{n}_6^{(31)} \right\}( \mbox{\mancube}_{a,b,c}^j), 
\\
E_{a,b,c}^j = & - \frac{M}{8} v_{a+1,b+1,c}^j  + \frac{\Delta t  }{8} \left\{ -P_5 \left( \mathbf{n}_5^{(12)} + \mathbf{n}_5^{(23)} + \mathbf{n}_5^{(31)} \right)  + P_8  \mathbf{n}_8^{(12)} + P_2  \mathbf{n}_2^{(23)}  + P_3 \mathbf{n}_3^{(31)} \right\}( \mbox{\mancube}_{a,b,c}^j) , 
\\
F_{a,b,c}^j = & - \frac{M}{8} v_{a,b+1,c+1}^j  +  \frac{\Delta t  }{8} \left\{ - P_6 \left( \mathbf{n}_6^{(12)} + \mathbf{n}_6^{(23)} + \mathbf{n}_6^{(31)} \right)   + P_3 \mathbf{n}_3^{(12)} + P_4 \mathbf{n}_4^{(23)}  +  P^j_8 \mathbf{n}_8^{(31)}\right\} ( \mbox{\mancube}_{a,b,c}^j) ,
\\
G_{a,b,c}^j = & - \frac{M}{8} v_{a+1,b,c+1}^j  + \frac{\Delta t  }{8} \left\{ - P_7 \left( \mathbf{n}_7^{(12)} + \mathbf{n}_7^{(23)} + \mathbf{n}_7^{(31)} \right) + P_2 \mathbf{n}_2^{(12)} + P_8 \mathbf{n}_8^{(23)} + P_4  \mathbf{n}_4^{(31)} \right\} ( \mbox{\mancube}_{a,b,c}^j), 
\\
H_{a,b,c}^j = & - \frac{M}{8} v_{a+1,b+1,c+1}^j  + \frac{\Delta t  }{8} \left\{- P_8 \left( \mathbf{n}_8^{(12)} + \mathbf{n}_8^{(23)} + \mathbf{n}_8^{(31)} \right)  + P_5 \mathbf{n}_5^{(12)} + P_6 \mathbf{n}_6^{(23)} 
 + P_7 \mathbf{n}_7^{(31)} \right\} ( \mbox{\mancube}_{a,b,c}^j),
\end{align*}}
where 
{\footnotesize
\begin{equation*}
\begin{aligned}
\mathbf{n}_1( \mbox{\mancube}_{a,b,c}^j) &= (( \varphi_{a+1,b,c}^j - \varphi_{a,b,c}^j)\times (\varphi_{a,b+1,c}^j - \varphi_{a,b,c}^j) ) \cdot (\varphi_{a,b,c+1}^j - \varphi_{a,b,c}^j),
\\
\mathbf{n}_2( \mbox{\mancube}_{a,b,c}^j) &= (( \varphi_{a+1,b+1,c}^j - \varphi_{a+1,b,c}^j)\times (\varphi_{a,b,c}^j - \varphi_{a+1,b,c}^j) ) \cdot (\varphi_{a+1,b,c+1}^j - \varphi_{a+1,b,c}^j),
\\
\mathbf{n}_3( \mbox{\mancube}_{a,b,c}^j) &= (( \varphi_{a,b,c}^j - \varphi_{a,b+1,c}^j)\times (\varphi_{a+1,b+1,c}^j - \varphi_{a,b+1,c}^j) ) \cdot (\varphi_{a,b+1,c+1}^j - \varphi_{a,b+1,c}^j),
\\
\mathbf{n}_4( \mbox{\mancube}_{a,b,c}^j) &= (( \varphi_{a,b+1,c+1}^j - \varphi_{a,b,c+1}^j)\times (\varphi_{a+1,b,c+1}^j - \varphi_{a,b,c+1}^j) ) \cdot (\varphi_{a,b,c}^j - \varphi_{a,b,c+1}^j),
\\
\mathbf{n}_5( \mbox{\mancube}_{a,b,c}^j) &= (( \varphi_{a,b+1,c}^j - \varphi_{a+1,b+1,c}^j)\times (\varphi_{a+1,b,c}^j - \varphi_{a+1,b+1,c}^j) ) \cdot (\varphi_{a+1,b+1,c+1}^j - \varphi_{a+1,b+1,c}^j),
\\
\mathbf{n}_6( \mbox{\mancube}_{a,b,c}^j) &= (( \varphi_{a+1,b+1,c+1}^j - \varphi_{a,b+1,c+1}^j)\times (\varphi_{a,b,c+1}^j - \varphi_{a,b+1,c+1}^j) ) \cdot (\varphi_{a,b+1,c}^j - \varphi_{a,b+1,c+1}^j),
\\
\mathbf{n}_7( \mbox{\mancube}_{a,b,c}^j) &= (( \varphi_{a,b,c+1}^j - \varphi_{a+1,b,c+1}^j)\times (\varphi_{a+1,b+1,c+1}^j - \varphi_{a+1,b,c+1}^j) ) \cdot (\varphi_{a+1,b,c}^j - \varphi_{a+1,b,c+1}^j),
\\
\mathbf{n}_8( \mbox{\mancube}_{a,b,c}^j) &= (( \varphi_{a+1,b,c+1}^j - \varphi_{a+1,b+1,c+1}^j)\times (\varphi_{a,b+1,c+1}^j - \varphi_{a+1,b+1,c+1}^j) ) \cdot (\varphi_{a+1,b+1,c}^j - \varphi_{a+1,b+1,c+1}^j),
\end{aligned}
\end{equation*}}
and
{\footnotesize
\begin{equation*}
\begin{aligned}
\mathbf{n}_1^{(12)} &=( \varphi_{a+1,b,c}^j - \varphi_{a,b,c}^j)\times (\varphi_{a,b+1,c}^j - \varphi_{a,b,c}^j),\\
\mathbf{n}^{(23)} _1&= ( \varphi_{a,b+1,c}^j - \varphi_{a,b,c}^j)\times (\varphi_{a,b,c+1}^j - \varphi_{a,b,c}^j),\\
\mathbf{n}^{(31)}_1&= ( \varphi_{a,b,c+1}^j - \varphi_{a,b,c}^j)\times (\varphi_{a+1,b,c}^j - \varphi_{a,b,c}^j),\\
J_1 ( \mbox{\mancube}_{a,b,c}^j)&=  \frac{ (( \varphi_{a+1,b,c}^j - \varphi_{a,b,c}^j)\times (\varphi_{a,b+1,c}^j - \varphi_{a,b,c}^j) ) \cdot (\varphi_{a,b,c+1}^j - \varphi_{a,b,c}^j)  }{ | s_{a+1,b,c} - s _{a,b,c}|  | s _{a,b+1,c} -s_{a,b,c}| | s _{a,b,c+1} -s _{a,b,c}| }.
\end{aligned}
\end{equation*}
Note that in $\mathbf{n}_1^{(12)}$, $\mathbf{n}^{(23)} _1$, $\mathbf{n}^{(31)}_1$ we adopt the usual expressions $(12)$, $(23)$, $(31)$ for permutation of the elements of the set $\{( \varphi_{a+1,b,c}^j - \varphi_{a,b,c}^j)$, $(\varphi_{a,b+1,c}^j - \varphi_{a,b,c}^j)$, $(\varphi_{a,b,c+1}^j - \varphi_{a,b,c}^j)\}$ which compose $\mathbf{n}_1$.}

\newpage
 
\subsection{Discrete Euler-Lagrange equations  for 3D barotropic fluid}\label{3D_DCEL_geo}

 {\footnotesize
\begin{align*}
 \rho_0\left( \frac{ v_{a,b,c}^j - v_{a,b,c}^{j-1}}{\Delta t} \right) 
 = & \frac{1}{8 \Delta s_1 \Delta s_2 \Delta s_3}  \left\{ \color{red} -  P_1  \left( \mathbf{n}_1^{(12)} + \mathbf{n}_1^{(23)} + \mathbf{n}_1^{(31)} \right)  ( \mbox{\mancube}_{a,b,c}^j)               \right.
 \\
& \color{red} -  P_2  \left( \mathbf{n}_2^{(12)} + \mathbf{n}_2^{(23)} + \mathbf{n}_2^{(31)} \right)  ( \mbox{\mancube}_{a-1,b,c}^j)  -  P_3  \left( \mathbf{n}_3^{(12)} + \mathbf{n}_3^{(23)} + \mathbf{n}_3^{(31)} \right)  ( \mbox{\mancube}_{a,b-1,c}^j) 
\\
& \color{red} -  P_4  \left( \mathbf{n}_4^{(12)} + \mathbf{n}_4^{(23)} + \mathbf{n}_4^{(31)} \right)  ( \mbox{\mancube}_{a,b,c-1}^j)  -  P_5  \left( \mathbf{n}_5^{(12)} + \mathbf{n}_5^{(23)} + \mathbf{n}_5^{(31)} \right)  ( \mbox{\mancube}_{a-1,b-1,c}^j) 
\\
& \color{red} -  P_6  \left( \mathbf{n}_6^{(12)} + \mathbf{n}_6^{(23)} + \mathbf{n}_6^{(31)} \right)  ( \mbox{\mancube}_{a,b-1,c-1}^j)  -  P_7  \left( \mathbf{n}_7^{(12)} + \mathbf{n}_7^{(23)} + \mathbf{n}_7^{(31)} \right)  ( \mbox{\mancube}_{a-1,b,c-1}^j) 
\\
&  \color{red} -  P_8  \left( \mathbf{n}_8^{(12)} + \mathbf{n}_8^{(23)} + \mathbf{n}_8^{(31)} \right)  ( \mbox{\mancube}_{a-1,b-1,c-1}^j) 
\\
& + \left( P_4 \mathbf{n}_4^{(12)} +   P_3 \mathbf{n}_3^{(23)} 
 +  P_{ 2 } \mathbf{n}_2^{(31)}\right)( \mbox{\mancube}_{a,b,c}^j)
 \\
& + \left( P_7 \mathbf{n}_7^{(12)} +  P_1  \mathbf{n}_1^{(23)} 
+ P_5  \mathbf{n}_5^{(31)} \right)( \mbox{\mancube}_{a-1,b,c}^j)
 \\
& + \left( P_6  \mathbf{n}_6^{(12)} + P_5  \mathbf{n}_5^{(23)} +  P_1 \mathbf{n}_1^{(31)} \right)( \mbox{\mancube}_{a,b-1,c}^j)
 \\
& + \left( P_1  \mathbf{n}_1^{(12)} +   P_7  \mathbf{n}_7^{(23)} +   P_6 \mathbf{n}_6^{(31)} \right)( \mbox{\mancube}_{a,b,c-1}^j)
 \\
& + \left( P_8  \mathbf{n}_8^{(12)} + P_2  \mathbf{n}_2^{(23)}  + P_3 \mathbf{n}_3^{(31)} \right)( \mbox{\mancube}_{a-1,b-1,c}^j)
 \\
& + \left( P_3 \mathbf{n}_3^{(12)} + P_4 \mathbf{n}_4^{(23)}  +  P_8 \mathbf{n}_8^{(31)} \right)( \mbox{\mancube}_{a,b-1,c-1}^j)
 \\
& + \left( P_2 \mathbf{n}_2^{(12)} + P_8 \mathbf{n}_8^{(23)} + P_4  \mathbf{n}_4^{(31)} \right)( \mbox{\mancube}_{a-1,b,c-1}^j)
 \\
& \left.+ \left( P_5 \mathbf{n}_5^{(12)} + P_6 \mathbf{n}_6^{(23)} 
 + P_7 \mathbf{n}_7^{(31)}  \right)( \mbox{\mancube}_{a-1,b-1,c-1}^j) \right\}
\end{align*}
Interpretation of the discrete Euler-Lagrange equations: \textit{discrete balance of momentum}
\[  \rho_0 (v^j-v^{j-1})/\Delta t = - \left[(P_{\rm ext}  \cdot   S_{\rm ext})\mathbf{n}_{\rm ext} - (P_{\rm int}   \cdot  S_{\rm int})\mathbf{n}_{\rm int}\right]_j / \textrm{volume}.
\]
 
\begin{figure}[H] \centering 
\includegraphics[width=4.7 in]{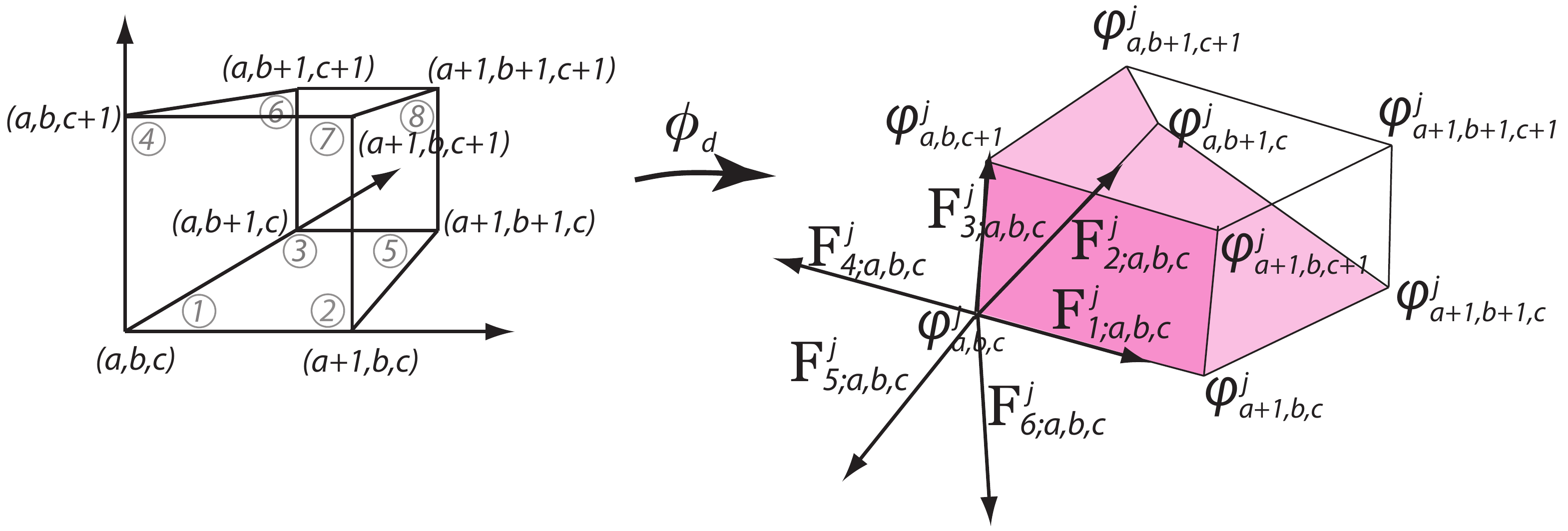} \vspace{-3pt}  
\caption{\footnotesize In red: the internal surfaces $S_{\rm int}$ associated to the cell $ \mbox{\mancube}_{a,b,c}^j$.} \label{3D_DCEL} 
\end{figure} 

{\footnotesize

\bibliographystyle{new}

}

\end{document}